\documentclass[preprint]{elsarticle}
\usepackage{lineno,hyperref}
\modulolinenumbers[5]
\makeatletter
\def\ps@pprintTitle{%
 \let\@oddhead\@empty
 \let\@evenhead\@empty
 \let\@oddfoot\@empty
 \let\@evenfoot\@empty
}
\makeatother

\usepackage[utf8]{inputenc}
\usepackage{amsthm, amsmath, amsfonts, amssymb}
\usepackage[fleqn]{mathtools} 
\usepackage{csquotes}    
\usepackage{mathrsfs}
\usepackage{xcolor}
\usepackage{verbatim}

\newtheorem{lem}{Lemma}[section]
\newtheorem{thm}{Theorem}[section]

\newtheorem{corollary}{Corollary}[section]

\newtheorem{definition}{Definition}[section]
\newcommand{\mycomment}[1]{}

\newcommand{\br}[2]{\left(#1 , #2\right)}
\newcommand{\nr}[1]{\| #1 \|_{L^2}^2}
\newcommand{\nrm}[1]{\| #1 \|_{L^2}}
\newcommand{\norm}[1]{{\left\lVert#1\right\rVert}_{L^{\infty}(\Omega)}}
\newcommand{\innerproduct}[1]{\left\langle #1 \right\rangle}

\usepackage[a4paper, total={6in, 8in}, margin=1.5 cm]{geometry}
\begin{document}
\begin{frontmatter}
\title{ Well-posedness and Fingering Patterns in $A + B \rightarrow C$ Reactive Porous Media Flow }

\author[mymainaddress]{Sahil Kundu\corref{equalcontribution}}
\ead{sahil.20maz0009@iitrpr.ac.in}

\author[mysecondaryaddress]{Surya Narayan Maharana\corref{equalcontribution}}
\ead{surya.maharana@gmail.com}

\author[mymainaddress]{Manoranjan Mishra\corref{mycorrespondingauthor}}
\ead{manoranjan@iitrpr.ac.in}

\address[mymainaddress]{Department of Mathematics, Indian Institute of Technology Ropar, Rupnagar, India}
\address[mysecondaryaddress]{Nonlinear Physical Chemistry Unit, Université libre de Bruxelles, 1050 Brussels, Belgium}

\cortext[equalcontribution]{These authors contributed equally to this work.}

\begin{abstract}
The convection-diffusion-reaction system governing incompressible reactive fluids in porous media is studied, focusing on the \( A + B \to C \) reaction coupled with density-driven flow. The time-dependent Brinkman equation describes the velocity field, incorporating permeability variations modeled as an exponential function of the product concentration. Density variations are accounted for using the Oberbeck-Boussinesq approximation, with density as a function of reactants and product concentration. The existence and uniqueness of weak solutions are established via the Galerkin approach, proving the system's well-posedness. A maximum principle ensures reactant nonnegativity with nonnegative initial conditions, while the product concentration is shown to be bounded, with an explicit upper bound derived in a simplified setting. Numerical simulations are performed using the finite element method to explore reactive fingering instabilities and illustrate the effects of density stratification, differential product mobility, and two- or three-dimensionality. Two cases with initial flat and elliptic interfaces further demonstrate the theoretical result that solutions continuously depend on initial and boundary conditions. These theoretical and numerical findings provide a foundation for understanding chemically induced fingering patterns and their implications in applications such as carbon dioxide sequestration, petroleum migration, and rock dissolution in karst reservoirs.
\end{abstract}

\begin{keyword}
Darcy-Brinkman, Porous media flow, Reactive flow, Existence, Uniqueness, Density fingering instability 
\end{keyword}

\end{frontmatter}


\section{Introduction}
Reactive flow in porous media plays a vital role in a wide range of natural and engineered processes, where fluid transport is intricately linked to chemical reactions and the inherent heterogeneity of the porous structure. A particularly striking feature of such systems is the complex interfacial dynamics that emerge from nonlinear interactions between flow, transport, and reaction processes. These dynamics are often influenced by instabilities such as the Rayleigh–Taylor (RT) instability \cite{rayleigh1882investigation,taylor1950instability}, which arises when fluids of different densities interact under the influence of gravity. These instabilities are primarily driven by the relative acceleration between the fluids, giving rise to finger-like structures at the interface. In miscible, non-reactive systems, such instabilities typically manifest symmetrically: the heavier fluid sinks while the lighter fluid rises, resulting in bidirectional finger propagation \cite{Paoli_2022, Gopalakrishnan_2021}. However, in reactive scenarios, chemical reactions can significantly modify local concentration gradients and, consequently, the density distribution near the interface. This leads to asymmetries in the resulting flow structures and the emergence of more complex interfacial behavior \cite{Almarcha_2010}. Although non-reactive flow in porous media is critically important for various geophysical and industrial applications—such as groundwater contamination \cite{van1988transport,leblanc1984sewage}, ice formation in seawater \cite{feltham2006sea}, petroleum migration \cite{simmons2001variable}, and carbon sequestration aimed at mitigating climate change \cite{emami2015convective, huppert2014fluid, de2021influence, Szulczewski2012}—the coupling between reactive transport and fluid motion in heterogeneous media is equally consequential. It governs key processes in the convective dissolution of \(CO_2 \) \cite{loodts2014control, THOMAS2016230}, reactive transport in geological formations \cite{stone2004}, Earth's mantle convection \cite{davaille1999}, and astrophysical events such as supernova explosions \cite{schmidt2006}, among others.

In this study, we consider two incompressible reacting fluids, \( A \) and \( B \), moving with a common velocity vector \( \boldsymbol{u} \) while undergoing a generic reaction, denoted as \( A + B \rightarrow C \), where \( C \) is the product fluid. The concentrations of the reactants and the product are denoted by \( a \), \( b \), and \( c \), respectively. The interactions among the reactants and the product are described by reaction rates \( R_{\Phi},~ \Phi \in \{a, b, c\} \), and mass diffusion characterized by fluxes \( F_{\Phi},~ \Phi \in \{a, b, c\} \). The reaction terms are expressed as:  
$$ R_{a} = - k ab,~  R_{b} = - k ab,~ R_{c} =  k ab,$$  
where \( k \) is the reaction rate constant.  Following Fick's law, the diffusion fluxes are given by:  
$$ F_{a} = - d_{A} \boldsymbol{\nabla} a,~  F_{b} = - d_{B} \boldsymbol{\nabla} b,~  F_{c} =  d_{C} \boldsymbol{\nabla} c .$$
Here, \( d_A \), \( d_B \), and \( d_C \) are positive diffusion constants for the respective species.
Mostly, this reactive system is analyzed in porous media, where the velocity vector \( \boldsymbol{u} \) satisfies Darcy's law, the conservation of momentum equation, and the conservation of mass equation (see \cite{dewit_2020} and references therein).
However, although Darcy’s law is widely used due to its geophysical significance, it is inadequate for modeling flow in domains with vugs and cracks formed in chemically active or reactive carbonate rocks, where porosity exceeds 0.75 (see \cite{Mccurdy2019,szymczak_ladd_2014} and references therein). In such cases, the Brinkman equation is a suitable alternative, as it extends Darcy’s law by including the Laplacian of the velocity vector, allowing for a more accurate representation of complex flow behaviors \cite{brinkman1949calculation, morales2017darcy}. Therefore, we adopt the Brinkman equation as the momentum conservation equation \cite{KUNDU2024128532,Hou16, John2015, keim2016, Layton2013}.

The coupling between the Brinkman equation and the convection diffusion reaction system occurs through changes in density \cite{Almarcha_2010} or permeability \cite{nagatsu2014hydrodynamic}, as reactant and product concentrations influence these properties. Under the Oberbeck-Boussinesq approximation, density variations are assumed to affect the system solely through the gravity term, which is a reasonable assumption for geophysical applications such as geological carbon dioxide sequestration. Additionally, the density of the mixture is modeled as a linear function of the reactant concentrations \( a \) and \( b \), and the product concentration \( c \), expressed as  
\[
\rho = 1 + R_{A}a + R_{B}b + R_{C}c,
\]  
where \( R_{A} \), \( R_{B} \), and \( R_{C} \) are the respective expansion coefficients \cite{loodts2014control}.
To account for heterogeneity, permeability is modeled as  
\[
K(c) = e^{-\alpha c},
\]  
where \( \alpha \) is a positive constant representing the mobility ratio of reactants to the product \cite{nagatsu2014hydrodynamic}.
The system incorporates no-slip boundary conditions for the velocity and Neumann-type boundary conditions for the reactant and product concentrations. The following mathematical model system, comprising PDEs, describes the fingering phenomena in \( A + B \rightarrow C \) reactive porous media flow.

\noindent \textbf{Mathematical Model:}\label{sec:main}
\begin{subequations}\label{model}
\begin{align}
      \boldsymbol{\nabla} \cdot \boldsymbol{u} = 0 \quad \text{in}\ (0,T) \times \Omega,\label{model1}\\
      \dfrac{\partial \boldsymbol{u}}{\partial t}+\dfrac{\mu}{K(c)}\boldsymbol{u}=-\boldsymbol{\nabla} p+\mu_{e} \Delta \boldsymbol{u} - \rho(a,b,c) \boldsymbol{g}\quad \text{in}\ (0,T) \times \Omega,\label{model2}\\
      \dfrac{\partial a}{\partial t}+\boldsymbol{u}\cdot \boldsymbol{\nabla} a=d_A\Delta a - k  ab \quad \text{in}\ (0,T) \times \Omega,\label{model3}\\
      \dfrac{\partial b}{\partial t}+\boldsymbol{u}\cdot \boldsymbol{\nabla} b = d_B\Delta b - k  ab \quad \text{in}\ (0,T) \times \Omega,\label{model4}\\
      \dfrac{\partial c}{\partial t}+\boldsymbol{u}\cdot \boldsymbol{\nabla} c=d_C\Delta c + k  ab \quad \text{in}\ (0,T) \times \Omega.\label{model5}
    \end{align}  
\end{subequations}
Here, $\mu$, $K$, $p$, $\mu_{e}$, and $g(\boldsymbol{x})$ denote the viscosity, permeability, pressure, effective viscosity, and gravitational field, respectively. The system of PDEs given by equations \eqref{model1} through \eqref{model5} is supplemented with the following boundary and initial data:
\begin{subequations}\label{boundry and initial data}
    \begin{align}
     \boldsymbol{u}(t,\boldsymbol{x}) = 0 , ~\dfrac{\partial a}{\partial \boldsymbol{\eta} }(t,\boldsymbol{x}) = \dfrac{\partial b}{\partial \boldsymbol{\eta} }(t,\boldsymbol{x}) = \dfrac{\partial c}{\partial \boldsymbol{\eta} }(t,\boldsymbol{x}) = 0  
~ \text{ for } \hspace{5pt} (t,\boldsymbol{x}) \in (0,T) \times \partial\Omega.   \label{bc's}\\
        \boldsymbol{u}(0,\boldsymbol{x}) = \boldsymbol{u}_0(\boldsymbol{x}),  a\left( \boldsymbol{x},0\right) = a_{0}\left(\boldsymbol{x}\right), b\left( \boldsymbol{x},0\right) = b_{0}\left(\boldsymbol{x}\right), c\left( \boldsymbol{x},0\right) = c_{0}\left(\boldsymbol{x}\right) \hspace{5pt} \text{ for  } \boldsymbol{x} \in \Omega.\label{ic's}
    \end{align}
\end{subequations}
Where `$\boldsymbol{\eta}$' denotes the unit normal vector. We further assume that there exist positive constants $M_A$, $M_B$, and $M_C$ such that
   \begin{align}\label{intial data assumptation}
    0 \leq a_{0}(\boldsymbol{x}) \leq M_A,~ 0 \leq b_{0}(\boldsymbol{x}) \leq M_B,~ 0 \leq c_{0}(\boldsymbol{x}) \leq M_C ~ ~ a. e. \text{ on } \Omega.
\end{align}

The key questions for understanding the well-posedness of the above model are whether the solution exists, is unique, depends continuously on initial and boundary conditions, and accurately represents the physical phenomena involving fingering instability. To our knowledge, prior research has not addressed the well-posedness of the \( A + B \rightarrow C \) reaction-diffusion system. Gálfi and Rácz \cite{galfi1988} provided an analytical solution under specific conditions but focused on the propagation and width of the reaction front, not its well-posedness. This raises two inquiries: (1) Has prior research established the well-posedness of the \( A + B \rightarrow C \) reaction-diffusion system, alone or coupled with porous media convection? (2) What is needed to bridge the gap between theoretical understanding and experimental observations? The existence of global weak solutions for porous media convection governed by Darcy or Stokes equations has been studied under specific conditions. For example, Amirat et al. \cite{Amirat_2007} analyzed nonreactive compressible flow with steady Darcy velocity, constant permeability, and pressure-dependent density, while Chen and Ewing \cite{Chen_1999} established weak solutions for incompressible flow with space-dependent permeability undergoing a first-order reaction. Nazer et al. \cite{Nazer_2024} demonstrated well-posedness for weakly compressible Darcy flow with pressure-dependent permeability, and Mucha et al. \cite{MUCHA2024104139} investigated compressible flow governed by Stokes equations with multiple irreversible reactions. These studies focus on mathematical analysis but overlook physical interpretations or real-world fluid interactions. Parallel research has also addressed non-Newtonian fingering \cite{Bansal_2023}, diffusive finger stability \cite{Dag_2009}, magnetic Rayleigh-Taylor finger stabilization \cite{Jiang_2018}, elastic-viscous fingering \cite{He_2012}, and mixing length bounds \cite{Kalinin_2024}, but without addressing the well-posedness of underlying models. Efforts to date neither establish well-posedness nor fully explain reactive density fingering patterns in heterogeneous porous media governed by Brinkman's or Darcy's equations. Bridging this gap requires investigating mathematical well-posedness while connecting theory with experimentally observed fingering patterns in heterogeneous porous media flows.

The rest of this paper is organized to justify and explain the answers to these key questions. Section \ref{fspr} provides an overview of the functional spaces and preliminary results used in our analysis. Section \ref{Exis} establishes the well-posedness of the \( A + B \rightarrow C \) reaction-diffusion-convection system for a broader class of initial conditions by addressing the existence of solutions, while Section \ref{continuous} focuses on the uniqueness of solutions and their continuous dependence on initial data. Sections \ref{semi} and \ref{num} present the semi-discrete formulation and numerical methods, including validations that confirm Gálfi and Rácz’s theory \cite{galfi1988}. Section \ref{pattern} demonstrates the continuous dependence of solutions on initial data by analyzing scenarios inspired by prior experiments and simulations, such as reactants separated by a flat horizontal interface \cite{Almarcha_2010} or an elliptical blob of reactant \( A \) in a pool of \( B \) \cite{Jha_2023,SaloSalgado2024}, in two- and three-dimensional domains. Comparisons of pre-existing density stratification with reaction-induced instability highlight how initial density gradients amplify or suppress fingering. Section \ref{pattern} also explores permeability variations on fingering patterns, showing that differential product mobility can significantly reduce finger growth. These findings, consistent with experiments \cite{Binda_2017} demonstrating the role of precipitation in density fingering, are framed within a mathematical model emphasizing permeability heterogeneity \cite{nagatsu2014hydrodynamic}. Finally, the implications of our findings, particularly the impact of three-dimensionality, density stratification, and permeability heterogeneity on reactive fingering, are summarized.
\section{Functional Spaces and Preliminary Results} \label{fspr}
In this section, we introduce standard notations and essential results that will be frequently referenced throughout the paper. The inner product in a given inner product space $X$ is denoted by $(\cdot, \cdot)_{X}$. For a Banach space $X$ and its dual space $X^{*}$, the duality pairing is denoted by $\langle\cdot, \cdot\rangle_{X}$. Here, $W^{m,p}(\Omega)$ denotes the standard Banach space with the corresponding norm. For $p=2$, $W^{m,2}(\Omega)$ is denoted by $H^{m}(\Omega)$. The space $H_{0}^{m}(\Omega)$ is the closure of $\mathcal{C}_{c}^{\infty}(\Omega)$ in $H^{m}(\Omega)$. Let $\mathcal{V} = \{\boldsymbol{u} \in \mathcal{C}_{c}^{\infty}(\Omega), \boldsymbol{\nabla} \cdot \boldsymbol{u} = 0   \}$, then we define $S$, and $V$ to be the closure of $\mathcal{V}$ in $(L^{2}(\Omega))^{d}$ and $(H_{0}^{1}(\Omega))^d$ spaces,  respectively. Then, ${L^{p}(0, T; X)}$ and ${L^{\infty}(0, T; X)}$ are time-dependent Banach spaces, where $X$ is a given Banach space equipped with the norm $\| \cdot \|_{X}$. For convenience, we denote the inner product in $L^2(\Omega)$ and $(L^2(\Omega))^d$ by $(\cdot, \cdot)$. Similarly, we denote the $L^2$ norm by $\nrm{\cdot}$, both in $L^2(\Omega)$ and $(L^2(\Omega))^{d}$.

\begin{thm}[Gagliardo–Nirenberg, cf. { \cite[Lemma~1]{migorski2019nonmonotone}}, { \cite[Lemma~1.1]{garcke2019}}]\label{gagliardo} 
If $\Omega \subset \mathbb{R}^d \,(d=2,3)$ is a domain with $\mathcal{C}^{1}$ boundary, then there exists a constant $M>0$ depending only on $\Omega$ such that, in the case $n=2$:
$$
\|\boldsymbol{\phi}\|_{L^{4}} \leq M\|\boldsymbol{\phi}\|_{L^{2}}^{1 / 2}\|\boldsymbol{\nabla} \boldsymbol{\phi}\|_{L^{2}}^{1 / 2} \hspace{5pt};\hspace{5pt} \forall \hspace{5pt} \boldsymbol{\phi} \in H_{0}^{1}(\Omega)
$$
and in the case $n=3$:
$$
\|\boldsymbol{\phi}\|_{L^{4}} \leq M\|\boldsymbol{\phi}\|_{L^{2}}^{1 / 4}\|\boldsymbol{\nabla} \boldsymbol{\phi}\|_{L^{2}}^{3 / 4}\hspace{5pt};\hspace{5pt} \forall \hspace{5pt} \boldsymbol{\phi} \in H_{0}^{1}(\Omega).
$$  
\end{thm}

\begin{lem}[Theorem 7, \cite{bresch2007effect}]\label{lemma0}
Let $\boldsymbol{u} \in L^2(0,T;V)$, and suppose that $\zeta(t,\boldsymbol{x})$ is a solution to the following equation:
 \begin{align}\label{eq3.1}
     \innerproduct{\frac{\partial \zeta}{\partial t},\phi} + \br{\boldsymbol{u}\cdot \boldsymbol{\nabla}\zeta}{\phi} + d\br{\boldsymbol{\nabla}\zeta}{\boldsymbol{\nabla}\phi} = 0 ~ ~ a.  e. ~ on ~ (0, T), ~ ~ \forall \phi \in H^1(\Omega).
 \end{align}
  with initial data $\zeta_{0}(\boldsymbol{x})$ satisfying $0 < m \leq \zeta_{0}(\boldsymbol{x}) \leq M$. Then $\zeta(t,\boldsymbol{x})$ satisfies $0 < m \leq \zeta(t,\boldsymbol{x}) \leq M $, for almost every $(t,\boldsymbol{x}) \in (0,T) \times \Omega$.
\end{lem}

\begin{definition}\label{def super solution}
    We call $ \overline{\zeta} \in H^1(\Omega)$ a weak super-solution of equation \eqref{eq3.1}, if 
    \begin{align*}
         \innerproduct{\frac{\partial \overline{\zeta}}{\partial t},\phi} + \br{\boldsymbol{u}\cdot \boldsymbol{\nabla}\overline{\zeta}}{\phi} + \br{\boldsymbol{\nabla}\overline{\zeta}}{\boldsymbol{\nabla}\phi} \geq 0 ~ ~a.  e. ~ on ~ (0, T),
    \end{align*}
for all $\phi \in H^1(\Omega),~ \text{and}~  \phi \geq 0 ~ a.e. \text{ on } \Omega.$
\end{definition}

\begin{definition}\label{def sub solution}
    We call $\underline{\zeta} \in H^1(\Omega)$ a weak sub-solution of equation \eqref{eq3.1}, if
    \begin{align*}
         \innerproduct{\frac{\partial\underline{\zeta}}{\partial t},\phi} + \br{\boldsymbol{u}\cdot \boldsymbol{\nabla}\underline{\zeta}}{\phi} + \br{\boldsymbol{\nabla}\underline{\zeta}}{\boldsymbol{\nabla}\phi} \leq 0 ~ ~a.  e. ~ on ~ (0, T),
    \end{align*}
for all $\phi \in H^1(\Omega),~\text{and}~  \phi \geq 0 ~ a.e. \text{ on } \Omega.$
\end{definition}

\begin{definition}[Weak solution to the problem \eqref{model}--\eqref{boundry and initial data}]\label{def11}
We say that a quadruplet \((a, b, c, \boldsymbol{u})\) is a weak solution to the problem \eqref{model}--\eqref{boundry and initial data} if
\[
a, b, c \in L^2(0,T; H^1(\Omega)), \quad 
\boldsymbol{u} \in L^2(0,T; V),
\]
and the initial condition \eqref{ic's} is satisfied, along with the following variational formulations: for almost every \( t \in (0,T) \), and for all test functions \(\phi_1, \phi_2, \phi_3 \in H^1(\Omega)\) and \( \boldsymbol{v} \in V \),
\begin{subequations}
\begin{align}\label{weak1}
\left\langle \frac{\partial \boldsymbol{u}(t)}{\partial t}, \boldsymbol{v} \right\rangle 
&+ \mu \big(e^{\alpha c(t)} \boldsymbol{u}(t), \boldsymbol{v}\big) 
+ \mu_e \big(\nabla \boldsymbol{u}(t), \nabla \boldsymbol{v}\big) \nonumber \\
&+ \big((1 + R_A a(t) + R_B b(t) + R_C c(t)) \boldsymbol{g}, \boldsymbol{v} \big) = 0,
\end{align}
\begin{align}\label{weak2}
\left\langle \frac{\partial a(t)}{\partial t}, \phi_1 \right\rangle 
&+ \big( \boldsymbol{u}(t) \cdot \nabla a(t), \phi_1 \big) 
+ d \big( \nabla a(t), \nabla \phi_1 \big) 
+ k \big( a(t) b(t), \phi_1 \big) = 0,
\end{align}
\begin{align}\label{weak3}
\left\langle \frac{\partial b(t)}{\partial t}, \phi_2 \right\rangle 
&+ \big( \boldsymbol{u}(t) \cdot \nabla b(t), \phi_2 \big) 
+ d \big( \nabla b(t), \nabla \phi_2 \big) 
+ k \big( a(t) b(t), \phi_2 \big) = 0,
\end{align}
\begin{align}\label{weak4}
\left\langle \frac{\partial c(t)}{\partial t}, \phi_3 \right\rangle 
&+ \big( \boldsymbol{u}(t) \cdot \nabla c(t), \phi_3 \big) 
+ d \big( \nabla c(t), \nabla \phi_3 \big) 
- k \big( a(t) b(t), \phi_3 \big) = 0.
\end{align}
\end{subequations}

\noindent In addition, we obtain the further regularity:
\[
a, b, c \in L^{\infty}(0,T; L^2(\Omega)), \quad \boldsymbol{u} \in L^{\infty}(0,T; S),
\]
as a consequence of the energy estimates derived in the analysis.
\end{definition}

\section{Existence of a weak solution}\label{Exis}

\begin{thm}\label{th21}
For any initial condition $ (a_{0}, b_0, c_{0}, \boldsymbol{u}_{\boldsymbol{0}}) \in  (L^2(\Omega))^3  \times S$, there exists a solution  $(a, b, c, \boldsymbol{u}) \in (L^2(0,T; H^1(\Omega))^3 \times  L^2(0,T; V)$ of the problem \eqref{model}-\eqref{boundry and initial data}, in the sense of Definition \ref{def11}. 
\end{thm}

To establish the existence of a solution, we employ the Galerkin method. We use a special basis of $H^{1}(\Omega)$ composed with eigenvectors of the negative Laplace operator associated with the Neumann boundary condition corresponding to eigen values  $\{\lambda_j\}_{j=1}^{\infty}$ and an arbitrary basis of $V$. Let the basis for  $V$ be given by ($\boldsymbol{z}_1,\boldsymbol{z}_2,.....$) and let the basis for $H^{1}(\Omega)$ be given by ($w_1,w_2,.....$). Let $V_n = span(\boldsymbol{z}_1,\boldsymbol{z}_2,.....,\boldsymbol{z}_n) $
and $W_n = span(w_1,w_2,.....w_n) $ are finite dimensional subspaces of $V$ and $H^1(\Omega)$, respectively. Now we look for the functions $\boldsymbol{u}_{n}:[0,T] \rightarrow V_n$, $a_{n}:[0,T] \rightarrow W_n$, $b_{n}:[0,T] \rightarrow W_n$ and $c_{n}:[0,T] \rightarrow W_n$ of the form
$$\boldsymbol{u}_n = \sum_{j=1}^{n} \lambda_{j}^{n}(t)\boldsymbol{z}_j~,~ a_n = \sum_{j=1}^{n} \alpha_{j}^{n}(t) w_{j}~,~  b_{n} = \sum_{j=1}^{n} \beta_{j}^{n}(t) w_{j}~,~ \text{ and }~ c_n = \sum_{j=1}^{n} \gamma_{j}^{n}(t) w_{j}, $$
so that these functions satisfy the following equations:
\begin{align}\label{finite weak 1}
\innerproduct{\dfrac{\partial \boldsymbol{u}_{n}(t)}{\partial t},\boldsymbol{z}_{j}}  +\mu \left( e^{\alpha c_n(t)}\boldsymbol{u}_n(t),\boldsymbol{z}_j\right) + \mu_{e} \br{\boldsymbol{\nabla}\boldsymbol{u}_n(t)}{\boldsymbol{\nabla}\boldsymbol{z}_j} + \br{(1 + R_{A} a_n(t) + R_{B} b_n(t) + R_{C} 
 c_n(t) )\boldsymbol{g}}{\boldsymbol{z}_j}=0,    
\end{align}
\begin{align}\label{finite weak 2}
 \innerproduct{\dfrac{\partial a_{n}(t)}{\partial t},w_{j}} +d \br{\boldsymbol{\nabla} a_{n}(t)} {\boldsymbol{\nabla} w_{j}} + \br{\boldsymbol{u}_n(t)\cdot \boldsymbol{\nabla} a_{n}(t)}{w_{j}} + k \br{a_n(t) b_n(t)}{w_{j}} =0,    
\end{align}
\begin{align}\label{finite weak 3}
\innerproduct{\dfrac{\partial b_{n}(t)}{\partial t},w_{j}} +d \br{\boldsymbol{\nabla} b_{n}(t)} {\boldsymbol{\nabla} w_{j}} + \br{\boldsymbol{u}_n(t)\cdot \boldsymbol{\nabla} b_{n}(t)}{w_{j}} + k \br{a_n(t) b_n(t)}{w_{j}} =0,    
\end{align}
\begin{align}\label{finite weak 4}
 \innerproduct{\dfrac{\partial c_{n}(t)}{\partial t},w_{j}} +d \br{\boldsymbol{\nabla} c_{n}(t)} {\boldsymbol{\nabla} w_{j}} + \br{\boldsymbol{u}_n(t)\cdot \boldsymbol{\nabla} c_{n}(t)}{w_{j}} - k \br{a_n(t) b_n(t)}{w_{j}} =0,   
\end{align}
\noindent for $j=1,2....n$ and for almost every $t\in (0,T)$. Since $V_n$ and $W_n$ are finite-dimensional spaces, the system of equations \eqref{finite weak 1}, \eqref{finite weak 2}, \eqref{finite weak 3}, and \eqref{finite weak 4} reduces to a nonlinear system of first-order ordinary differential equations. Consequently, the local existence of $a_n$, $b_n$, $c_n$, and $\boldsymbol{u}_n$ follows from Carathéodory's existence theorem. 

Next, we establish a priori estimates for $a_n$, $b_n$, $c_n$, $\boldsymbol{u}_n$, and their derivatives, then utilize these to apply a continuation argument, to extend the solutions for each $t \in (0, T]$, for a given $T > 0$.


\begin{lem}\label{lemma sub sup solution}
    A weak sub-solution of equation \eqref{eq3.1} satisfies the maximum principle, while a weak super-solution satisfies the minimum principle. Specifically, if the initial conditions satisfy $\overline{\zeta}(0, \boldsymbol{x}) \geq m$ and $\underline{\zeta}(0, \boldsymbol{x}) \leq M$ for almost every $\boldsymbol{x} \in \Omega$, then the super-solution $\overline{\zeta}(t, \boldsymbol{x})$ remains bounded below by $m$, and the sub-solution $\underline{\zeta}(t, \boldsymbol{x})$ remains bounded above by $M$ for almost every $(t, \boldsymbol{x}) \in (0, T) \times \Omega$.
\end{lem}
\begin{proof}
As $\underline{\zeta}$ is a weak sub-solution of equation \eqref{eq3.1}, it follows by definition \eqref{def sub solution}
\begin{align*}
         \innerproduct{\frac{\partial\underline{\zeta}}{\partial t},\phi} + \br{\boldsymbol{u}\cdot \boldsymbol{\nabla}\underline{\zeta}}{\phi} + d\br{\boldsymbol{\nabla}\underline{\zeta}}{\boldsymbol{\nabla}\phi} \leq 0, ~ ~ \forall \phi \in H^1(\Omega),~  \phi \geq 0 ~ a.e. \text{ on } \Omega.
    \end{align*}
Define $\underline{\zeta}^{+} = \max(0,\underline{\zeta} - M)$. From above inequality, we can write 
\begin{align*}
         \innerproduct{\frac{\partial(\underline{\zeta}-M)}{\partial t},\phi} + \br{\boldsymbol{u}\cdot \boldsymbol{\nabla}(\underline{\zeta}-M)}{\phi} + d\br{\boldsymbol{\nabla}(\underline{\zeta}-M)}{\boldsymbol{\nabla}\phi} \leq 0, 
    \end{align*}
for all $\phi \in H^1(\Omega),~  \phi \geq 0 ~ a.e. \text{ on } \Omega$. Now using $\phi = \underline{\zeta}^{+}$, we have 
\begin{align*}
    \innerproduct{\frac{\partial\underline{\zeta}^{+}}{\partial t},\underline{\zeta}^{+}} + \br{\boldsymbol{u}\cdot \boldsymbol{\nabla}\underline{\zeta}^{+}}{\underline{\zeta}^{+}} + d\br{\boldsymbol{\nabla}\underline{\zeta}^{+}}{\boldsymbol{\nabla}\underline{\zeta}^{+}} \leq 0,  ~ a.  e. ~ on ~ (0, T).
\end{align*}
Utilizing integration by parts with boundary conditions \eqref{bc's} and the continuity equation \eqref{model1}, the second term of the above inequality vanishes.
    \begin{align*}
    \frac{1}{2} \frac{d}{dt} \nr{\underline{\zeta}^{+}} + d \nr{\boldsymbol{\nabla}\underline{\zeta}^{+}} \leq 0,  ~ a.  e. ~ on ~ (0, T).
 \end{align*}
     Integrating with respect to $t$, we have 
 \begin{align*}
     \nr{\underline{\zeta}^{+}(t)} + 2d\int_{0}^{t} \nr{\boldsymbol{\nabla}\underline{\zeta}^{+}(\tau)}d\tau \leq  \nr{\underline{\zeta}^{+}(t=0)}, \text{ for } t \in (0,T). 
 \end{align*}
 Utilizing $\underline{\zeta}^{+}(t=0) = 0$, we deduce that $\underline{\zeta}^{+}(t) = 0$, which implies $\underline{\zeta}(t,\boldsymbol{x}) \leq M$ for almost every $(t,\boldsymbol{x}) \in (0,T) \times \Omega$. Similarly, we can prove $\overline{\zeta}(t,\boldsymbol{x}) \geq m$ for almost every $(t,\boldsymbol{x}) \in (0,T) \times \Omega$.
\end{proof}

\begin{lem}\label{max min principle} Let $a_n(0, \boldsymbol{x})$, $b_n(0, \boldsymbol{x})$, and $c_n(0, \boldsymbol{x}) \in L^2(\Omega)$, and assume they satisfy the initial data assumptions \eqref{intial data assumptation}. Then, for $\boldsymbol{u}_n \in L^2(0, T; V)$, the functions $a_n(t, \boldsymbol{x})$, $b_n(t, \boldsymbol{x})$, and $c_n(t, \boldsymbol{x})$ remain nonnegative for all times and satisfy the following bounds: 
    \begin{align*}
      0 \leq a_n(t,\boldsymbol{x}) \leq M_A, ~ 0 \leq b_n(t,\boldsymbol{x}) \leq M_B,~ 0 \leq c_n(t,\boldsymbol{x}) \leq \tilde{M} ~ ~ a.  e. ~ on ~ (0, T) \times \Omega,  
    \end{align*}  
where $\tilde{M} = M_C + k M_A M_B T$ is a positive constant.
\end{lem}
\begin{proof}  For a small positive real number $\epsilon$, let $a_{n}^{\epsilon}$, $b_{n}^{\epsilon}$, and $c_{n}^{\epsilon}$ be the solutions of the system \eqref{finite weak 2}-\eqref{finite weak 4}, with initial data $a_{n}^{\epsilon}(0) = a_{n}(0) + \epsilon$, $b_{n}^{\epsilon}(0) = b_{n}(0) + \epsilon$, and $c_{n}^{\epsilon}(0) = c_{n}(0) + \epsilon$. As $\epsilon \to 0$, we have $a_{n}^{\epsilon} \to a_{n}$, $b_{n}^{\epsilon} \to b_{n}$, and $c_{n}^{\epsilon} \to c_{n}$. Let $t = t_0 \in (0, T)$ denote the earliest time at which one of the concentrations, $a_n^{\epsilon}$, $b_n^{\epsilon}$, or $c_n^{\epsilon}$, becomes negative. Without loss of generality, assume that $a_n^{\epsilon}$ is the first concentration to become negative. Now, as $b_{n}^{\epsilon}$ is positive at $t = t_0$, there exists a $t_1 \in (t_0, T)$ such that $b_{n}^{\epsilon}$ remains non-negative for all $t \in (t_0, t_1)$.
Then from equation \eqref{finite weak 2}, we have 
\begin{align*}
 \innerproduct{\dfrac{\partial a_{n}^{\epsilon}(t)}{\partial t},\phi} + d_A \br{\boldsymbol{\nabla} a_{n}^{\epsilon}(t)} {\boldsymbol{\nabla} \phi} + \br{\boldsymbol{u}_n(t)\cdot \boldsymbol{\nabla} a_{n}^{\epsilon}(t)}{\phi} \geq 0,   
\end{align*}
for all $ \phi \in (H^1)_n, \phi \geq 0$, a.e. $\boldsymbol{x}\in \Omega$.
As $a_n^{\epsilon}$ satisfies the above inequality, it is a weak super-solution of the convection-diffusion equation. By Lemma \eqref{lemma sub sup solution}, we have $a_n^{\epsilon} \geq \epsilon$, which implies that $a_n \geq 0$ for all $t \in (t_0, t_1)$. This leads to a contradiction. 
Thus, we conclude that $a_n$, $b_n$, and $c_n$ are non-negative functions for almost every $(t, \boldsymbol{x}) \in (0,T) \times \Omega$. To obtain the upper bounds for $a_n$ and $b_n$, utilizing the non-negativity of $a_n(t, \boldsymbol{x})$ and $b_n(t, \boldsymbol{x})$ in equations \eqref{finite weak 2} and \eqref{finite weak 3}, we have:
\begin{align*} \innerproduct{\dfrac{\partial a_{n}(t)}{\partial t},\phi} + d_A \br{\boldsymbol{\nabla} a_{n}(t)} {\boldsymbol{\nabla} \phi} + \br{\boldsymbol{u}_n(t)\cdot \boldsymbol{\nabla} a_{n}(t)}{\phi} \leq 0,\\
 \innerproduct{\dfrac{\partial b_{n}(t)}{\partial t},\phi} + d_B\br{\boldsymbol{\nabla} b_{n}(t)} {\boldsymbol{\nabla} \phi} + \br{\boldsymbol{u}_n(t)\cdot \boldsymbol{\nabla} b_{n}(t)}{\phi} \leq 0,   \end{align*}
for all $ \phi \in (H^1)_n, \phi \geq 0$. From Lemma \eqref{lemma sub sup solution}, we conclude that $a_n(t, \boldsymbol{x}) \leq M_A$ and $b_n(t, \boldsymbol{x}) \leq M_B$ for almost every $(t, \boldsymbol{x}) \in (0, T) \times \Omega$. Using these upper bounds in equation \eqref{finite weak 4}, we obtain:
\begin{align*} \innerproduct{\dfrac{\partial c_{n}(t)}{\partial t},\phi} + d_C\br{\boldsymbol{\nabla} c_{n}(t)} {\boldsymbol{\nabla} \phi} + \br{\boldsymbol{u}_n(t)\cdot \boldsymbol{\nabla} c_{n}(t)}{\phi} \leq k\br{M_A M_B}{\phi},   \end{align*}
for all $ \phi \in (H^1)_n, \phi \geq 0$. Substituting $c_n = \tilde{c}_{n} + k M_A M_B t$ into the above inequality, we have:
\begin{align*} \innerproduct{\dfrac{\partial \tilde{c}_{n}(t)}{\partial t},\phi} + d_C\br{\boldsymbol{\nabla} \tilde{c}_{n}(t)} {\boldsymbol{\nabla} \phi} + \br{\boldsymbol{u}_n(t)\cdot \boldsymbol{\nabla} \tilde{c}_{n}(t)}{\phi} \leq 0,   \end{align*}
for all $ \phi \in (H^1)_n, \phi \geq 0$. Using Lemma \eqref{lemma sub sup solution}, we find: \begin{align*} \tilde{c}_{n}(t,\boldsymbol{x}) \leq \tilde{c}_{n}(0,\boldsymbol{x}) \leq c_n(0,\boldsymbol{x}) \leq M_C, \quad \text{for } (t, \boldsymbol{x}) \in (0, T) \times \Omega. \end{align*} This implies that $c_n(t, \boldsymbol{x}) \leq M_C + k M_A M_B T$ for $(t, \boldsymbol{x}) \in (0, T) \times \Omega$. We denote this upper bound by $\tilde{M}$. In the case of equidiffusivity, a simpler upper bound can be derived.
    \end{proof}
    \begin{corollary}
       If we take the diffusion constants of all the chemicals to be the same in the problem given by \eqref{model}-\eqref{boundry and initial data}, then we have $$ c_n(t,\boldsymbol{x}) \leq \frac{M_a + M_b + 2M_c}{2} ~  ~ a.  e. ~ on ~ (0, T) \times \Omega.$$  
    \end{corollary}
\begin{proof}
      Multiplying equation \eqref{finite weak 4} by $2$ and summing it with equations \eqref{finite weak 2} and \eqref{finite weak 3}, we obtain:
    \begin{align*}
        &\innerproduct{\dfrac{\partial(a_n(t)+b_n(t)+2 c_n(t))}{\partial t},\phi}+\br{\boldsymbol{u}_n(t)\cdot \boldsymbol{\nabla}(a_n(t)+b_n(t)+2 c_n(t))}{\phi} \\ &+ d\br{\boldsymbol{\nabla} (a_n(t)+b_n(t)+2c_n(t))}{\boldsymbol{\nabla}\phi} = 0, ~ \forall \phi \in H^1_n(\Omega), \quad \text{a.e. in } (0,T).
    \end{align*}
    Now, $a_n + b_n + 2c_n$ satisfies a convection-diffusion equation. Consequently, by Lemma \eqref{lemma0}, we have:
    $$ a_n(t,\boldsymbol{x}) + b_n(t,\boldsymbol{x}) + 2c_n(t,\boldsymbol{x}) \leq M_a + M_b + 2M_c ~  ~ a.  e. ~ on ~ (0, T) \times \Omega.$$
  Using the fact that $a_n$ and $b_n$ are non-negative functions, we conclude
    $$ c_n(t,\boldsymbol{x}) \leq \frac{M_a + M_b + 2M_c}{2} ~  ~ a.  e. ~ on ~ (0, T) \times \Omega.$$
\end{proof}

\begin{lem}\label{lemma2}
 The concentration sequences $a_n$, $b_n$, and $c_n$ are uniformly bounded in the spaces $L^{\infty}(0, T; L^2(\Omega))$ and $L^2(0, T; H^1(\Omega))$, and the velocity sequence $\boldsymbol{u}_n$ is uniformly bounded in the spaces $L^{\infty}(0, T; S)$ and $L^2(0, T; V)$.
 \end{lem}
\begin{proof}
Multiplying equation \eqref{finite weak 2} by $\alpha_j^n$ and summing over $j = 1, \dots, n$, we obtain:
\begin{align*}
\innerproduct{\frac{\partial a_n(t)}{\partial t},a_n(t)} + \big(\boldsymbol{u}_n\cdot \boldsymbol{\nabla} a_n(t), a_n(t) \big) + d_A \big(\boldsymbol{\nabla}a_n(t),\boldsymbol{\nabla}a_n(t)\big) = - k \big(a_n(t) b_n(t) , a_n(t)\big), ~ ~  a.  e. ~ on ~ (0, T).
\end{align*}
Utilizing the boundary condition and equation \eqref{model1}, the second term in the above equation vanishes after applying integration by parts. Also, using Lemma \eqref{max min principle} in the above equation, we have:
\begin{align}\label{4.5}
    \frac{1}{2}\frac{d}{d t}\|a_n(t)\|_{L^2}^2 + d_A\nr{\boldsymbol{\nabla} a_n(t)} \leq  k \norm{b_n(t)} \nr{a_n(t)} ~ ~a.  e. ~ on ~ (0, T).
\end{align}
Similarly we obtain,
\begin{align}\label{4.6}
    \frac{1}{2}\frac{d}{d t}\|b_n(t)\|_{L^2}^2 + d_B\nr{\boldsymbol{\nabla} b_n(t)} \leq  k \norm{a_n(t)} \nr{b_n(t)} ~ ~a.  e. ~ on ~ (0, T).
    \end{align}
    \begin{align}\label{4.7}
    \frac{1}{2}\frac{d}{d t}\|c_n(t)\|_{L^2}^2 + d_C\nr{\boldsymbol{\nabla} c_n(t)} \leq  \frac{k}{2}\norm{a_n(t)} \big(\nr{b_n(t)}+\nr{c_n(t)}\big), 
\end{align}
$a.  e. ~ on ~ (0, T)$. Multiplying equation \eqref{finite weak 1} by $\lambda_j^n(t)$ and summing over $j = 1, \dots, n$, gives:
\begin{align*}
    &\frac{1}{2}\frac{d}{d t}\nr{\boldsymbol{u}_n(t)} + \mu \br{e^{\alpha c_n(t)}\boldsymbol{u}_n(t)}{\boldsymbol{u}_n(t)} + \mu_e \nr{\boldsymbol{\nabla}\boldsymbol{u}_n(t)}\nonumber\\ &= - \br{(1 + R_{A} a_n(t) + R_{B} b_n(t) + R_{C} 
 c_n(t) )\boldsymbol{g}}{\boldsymbol{u}_n(t)}, ~ a.  e. ~ on ~ (0, T).
\end{align*}
Using Lemma \eqref{max min principle} that $c_n$ is non-negative and applying Hölder's and Young's inequalities, we obtain
\begin{align}\label{4.8}
  \frac{d}{d t}\nr{\boldsymbol{u}_n(t)} + 2\mu \nr{\boldsymbol{u}_n} + 2\mu_e \nr{\boldsymbol{\nabla}\boldsymbol{u}_n(t)}  &\leq \nr{\boldsymbol{g}} + R\norm{\boldsymbol{g}} \big(\nr{a_n(t)} + \nr{b_n(t)} + \nr{c_n(t)} \big)\nonumber \\ &+ \big(1+ 3 R \norm{g}  \big)\nr{\boldsymbol{u}_n(t)}, \quad a.  e. ~ on ~ (0, T).  
\end{align}
Here $R = \max(R_A, R_B, R_C)$. Summing the inequalities from \eqref{4.5} to \eqref{4.8} and using the uniform upper bounds for $\norm{a_n}$ and $\norm{b_n}$, we deduce that
\begin{align*}
    &\frac{d}{dt}\big(\nr{a_n(t)}+\nr{b_n(t)}+\nr{c_n(t)}+\nr{\boldsymbol{u}_n(t)}\big) + 2\mu \nr{\boldsymbol{u}_n}   + 2\mu_e \nr{\boldsymbol{\nabla}\boldsymbol{u}_n(t)}\\ & + 2d_A\nr{\boldsymbol{\nabla} a_n(t)} + 2d_B\nr{\boldsymbol{\nabla} b_n(t)}  + 2 d_C\nr{\boldsymbol{\nabla} c_n(t)} \leq \nr{\boldsymbol{g}} \\ & + \big(1+ 3R \norm{g}  \big)\nr{\boldsymbol{u}_n(t)}   +  R\norm{\boldsymbol{g}} \big(\nr{a_n(t)} + \nr{b_n(t)} + \nr{c_n(t)} \big) \\ &+ k\big( 3 M_A + 2 M_B\big) \nr{a_n(t)}  +  k\big( 3 M_A + 2 M_B\big) \big( \nr{b_n(t)} + \nr{c_n(t)} \big),~ a. e. ~ on ~ (0, T).
\end{align*}
Integrating above inequality  from $0$ to $\tau \in (0,T]$, gives
\begin{align*}
&\nr{a_n(\tau)}+\nr{b_n(\tau)}+\nr{c_n(\tau)}+ \nr{\boldsymbol{u}_n(\tau)} +  2 \int_{0}^{\tau} \mu \nr{\boldsymbol{u}_n(t)}  \\ & +  2 \int_{0}^{\tau}  \mu_e \nr{\boldsymbol{\nabla}\boldsymbol{u}_n(t)}+ 2 \int_{0}^{\tau}\left(d_A \nr{\boldsymbol{\nabla} a_n(t)} + d_B \nr{\boldsymbol{\nabla} b_n(t)}+  d_C \nr{\boldsymbol{\nabla} c_n(t)} \right) \\ & \leq \nr{a_n(0)}+\nr{b_n(0)} +\nr{c_n(0)}  + \nr{\boldsymbol{g}}  +  M  \int_{0}^{\tau}\big(\nr{a_n(t)} + \nr{b_n(t)} + \nr{c_n(t)} + \nr{\boldsymbol{u}_n(t)} \big).
\end{align*}
Where $$M = \max\big(k(3M_A+2M_B), R \norm{\boldsymbol{g}}, 1+ 3 R \norm{\boldsymbol{g}} \big).$$ Applying Grönwall's inequality to the above expression, we conclude that $a_n$, $b_n$, and $c_n$ are uniformly bounded in $L^{\infty}(0, T; L^{2}(\Omega))$ and $L^2(0, T; H^1(\Omega))$, while $\boldsymbol{u}_n$ is uniformly bounded in $L^{\infty}(0, T; S)$ and $L^2(0, T; V)$ spaces.
\end{proof}
\begin{lem}\label{lemma3}
The time derivatives $\frac{\partial a_n}{\partial t}$, $\frac{\partial b_n}{\partial t}$, and $\frac{\partial c_n}{\partial t}$ are uniformly bounded in $L^2(0, T; (H^{1})^*)$, and $\frac{\partial \boldsymbol{u}_n}{\partial t}$ is uniformly bounded in $L^2(0, T; V^{\ast})$.
\end{lem}
\begin{proof}
Using Lemma \eqref{max min principle} and Hölder's inequality in equations \eqref{finite weak 1} and \eqref{finite weak 2}, we derive
\begin{align*}
    \left|\innerproduct{\frac{\partial \boldsymbol{u}_n(t)}{\partial t},\boldsymbol{v}}\right| &\leq \mu e^{\alpha \norm{c_n(t)}} \nrm{\boldsymbol{u}_n(t)} \nrm{\boldsymbol{v}} + \mu_e \nrm{\boldsymbol{\nabla}\boldsymbol{u}_n(t)} \nrm{\boldsymbol{\nabla}\boldsymbol{v}} \\ &+ \norm{\boldsymbol{g}}\nrm{1 + R_{A}a_n(t) + R_{B} b_n(t) + R_{C}c_n(t)} \nrm{\boldsymbol{v}}, \text{ and }\\
  \left|\innerproduct{\frac{\partial a_n(t)}{\partial t},\phi(t)}\right|  &\leq \Big(\norm{a_n(t)}\nrm{\boldsymbol{u}_n} + d_A\nrm{\boldsymbol{\nabla}a_n(t)}\Big) \nrm{\boldsymbol{\nabla}\phi}+ k \nrm{a_n(t) b_n(t)} \nrm{\phi(t)},
\end{align*}
$a. e. ~ on ~ (0, T)$. Taking the supremum over all $\phi \in H^1(\Omega)$ and $\boldsymbol{v} \in V$ such that $\|\phi\|_{H^1(\Omega)} \leq 1$ and $\|\boldsymbol{v}\|_V \leq 1$ in the above equations, and then summing the results, we obtain 
\begin{align*}
 &\left\|\frac{\partial a_n(t)}{\partial t}\right\|_{(H^{1})^*} +  \left\| \frac{\partial \boldsymbol{u}_n(t)}{\partial t}\right\|_{V^{\ast}} \leq  \norm{a_n(t)}\nrm{\boldsymbol{u}_n}  + d_A \nrm{\boldsymbol{\nabla}a_n(t)}  + k \nrm{a_n(t) b_n(t)}\\&+ \mu e^{\alpha \norm{c_n(t)}} \nrm{\boldsymbol{u}_n(t)} + \mu_e \nrm{\boldsymbol{\nabla}\boldsymbol{u}_n(t)}+ \norm{\boldsymbol{g}}\nrm{1 + R_{A}a_n(t) + R_{B} b_n(t) + R_{C}c_n(t)},
\end{align*}
$a. e. ~ on ~ (0, T)$. All terms on the right-hand side of the above inequalities are uniformly bounded in $L^2(0, T)$. Therefore, we conclude that the sequences $ \frac{\partial a_n}{\partial t} $ and $ \frac{\partial \boldsymbol{u}_n}{\partial t} $ are uniformly bounded in $L^2(0, T; (H^{1})^{\ast})$ and $L^2(0, T; V^{\ast})$, respectively. Following the same line of reasoning, we can show that the sequences $ \frac{\partial b_n}{\partial t} $ and $ \frac{\partial c_n}{\partial t} $ are also uniformly bounded in $L^2(0, T; (H^{1})^{\ast})$.
\end{proof}
\noindent \textbf{Proof of the Theorem \ref{th21}} 
\begin{proof}
From Lemmas (\ref{lemma2}) and (\ref{lemma3}), we establish that $a_{n}, b_{n}, c_{n}$ are uniformly bounded in $L^{\infty}\left(0,T;L^2\right)$, and $\boldsymbol{u}_n$ is uniformly bounded in $L^{\infty}\left(0,T;S\right)$ for all $T > 0$. Consequently, employing a continuation argument, we conclude that $a_{n}(t), b_n(t), c_{n}(t)$ and $\boldsymbol{u}_{n}(t)$ exist for all $t \in (0, T)$, for any $T > 0$. Additionally, based on lemmas (2.1)-(2.3), we can deduce the existence of a constant  $M > 0$ such that:
\begin{align*}
    &\left\| a_n\right\|_{L^{\infty}\left(0,T;L^2\right)},   &&\left\| b_n   \right\|_{L^{\infty}\left(0,T;L^2\right)},  
    &&\left\| c_n   \right\|_{L^{\infty}\left(0,T;L^2\right)},  \\
    &\left\|a_n\right\|_{L^{2}\left(0,T;H^1\right)},  
    &&\left\|b_n\right\|_{L^{2}\left(0,T;H^1\right)},  
    &&\left\| c_n\right\|_{L^{2}\left(0,T;H^1\right)},  \\
    &\left\|\boldsymbol{u}_n\right\|_{L^{\infty}(0,T;S)},  
    &&\left\|\boldsymbol{u}_n\right\|_{L^2(0,T;V)},  
    &&\left\|\frac{\partial a_n}{\partial t}\right\|_{L^{2}\left(0,T;(H^{1})^*\right)},  \\
    &\left\|\frac{\partial b_n}{\partial t}\right\|_{L^{2}\left(0,T;(H^{1})^*\right)},  
    &&\left\|\frac{\partial c_n}{\partial t}\right\|_{L^{2}\left(0,T;(H^{1})^*\right)},  
    &&\left\|\frac{\partial \boldsymbol{u}_n}{\partial t}\right\|_{L^2(0,T;V^{*})} \leq M.
\end{align*}
\noindent for all $T>0$. By applying the diagonalization argument with the above estimates, we obtain sub-sequences, which are again labeled by $a_n, b_n$, $c_n$, and $\boldsymbol{u}_n$ that satisfy the following convergence result:
\begin{center}
$\begin{aligned}
    a_n &\rightarrow a  &&\text{weakly in } L^{2}\left(0,T;H^1(\Omega)\right), \\
    b_n &\rightarrow b  &&\text{weakly in } L^{2}\left(0,T;H^1(\Omega)\right), \\
    c_n &\rightarrow c  &&\text{weakly in } L^{2}\left(0,T;H^1(\Omega)\right), \\
    \boldsymbol{u}_n &\rightarrow \boldsymbol{u}  &&\text{weakly in } L^{2}\left(0,T;V\right), \\
    \frac{\partial a_n}{\partial t} &\rightarrow \frac{\partial a}{\partial t}  
    &&\text{weakly in } L^{2}\left(0,T;(H^{1}(\Omega))^*\right), \\
    \frac{\partial b_n}{\partial t} &\rightarrow \frac{\partial b}{\partial t}  
    &&\text{weakly in } L^{2}\left(0,T;(H^{1}(\Omega))^*\right), \\ 
    \frac{\partial c_n}{\partial t} &\rightarrow \frac{\partial c}{\partial t}  
    &&\text{weakly in } L^{2}\left(0,T;(H^{1}(\Omega))^*\right), \\
    \frac{\partial \boldsymbol{u}_n}{\partial t} &\rightarrow \frac{\partial \boldsymbol{u}}{\partial t}  
    &&\text{weakly in } L^{2}\left(0,T;V^{*}\right), 
\end{aligned}$
\end{center}
\noindent as $ n \rightarrow \infty.$
By applying the Aubin-Lions Lemma, we obtain strongly convergent subsequences of the Galerkin solutions: $a_n$, $b_n$, and $c_n$ in $L^2(0, T; L^2(\Omega))$, and $\boldsymbol{u}_n$ in $L^2(0, T; S)$. This strong convergence enables us to pass to the limit in the finite-dimensional equations \eqref{finite weak 1}–\eqref{finite weak 4}.

\noindent \textbf{Passing to the limit:} 
 For passing the limit into the equations \eqref{finite weak 1}-\eqref{finite weak 4}, we fix a $m \in \mathbb{N}$ such that $m \leq n$. Then for all $\phi \in W_m$, from equation \eqref{finite weak 2}, we have 
\begin{align}\label{3.8}
    \innerproduct{\frac{\partial a_n(t)}{\partial t}, \phi} + \br{\boldsymbol{u}_n(t)\cdot \boldsymbol{\nabla}a_n(t)}{\phi} + d \br{\boldsymbol{\nabla}a_n(t)}{\boldsymbol{\nabla}\phi}+k \br{a_n(t) b_n(t)}{\phi} = 0,
\end{align}
since $\phi$ is a linear combination of $w_j$ for $j \in \{1,2,3....m  \}$.
Passing the limit in the linear terms of the equation is straightforward. From Lemmas \eqref{lemma2} and \eqref{lemma3}, we observe that $a_n \boldsymbol{u}_n$ is uniformly bounded in $(L^2(\Omega))^d$. By weak compactness results, there exists a weakly convergent subsequence (still denoted by $a_n \boldsymbol{u}_n$) in $(L^2(\Omega))^d$. Utilizing the results $\nrm{a_n - a} \rightarrow 0$ and $\nrm{\boldsymbol{u}_n - \boldsymbol{u}} \rightarrow 0$ as $n \rightarrow \infty$, we can deduce that $a_n \boldsymbol{u}_n$ converges strongly to $a \boldsymbol{u}$ in $(L^1(\Omega))^d$. Thus, we assume that the weak limit of $a_n \boldsymbol{u}_n$ is indeed $a \boldsymbol{u}$.
Now, using the result that $a_n(t)\boldsymbol{u}_n(t)$ converges weakly to $a(t)\boldsymbol{u}(t)$ in $(L^2(\Omega))^d$ in the above equation, we conclude 
\begin{align*}
    \lim_{n\to\infty} \big(a_n(t) \boldsymbol{u}_n(t)  , \boldsymbol{\nabla}\phi\big) = \big( a(t) \boldsymbol{u}(t) , \boldsymbol{\nabla} \phi\big) ~ ~  a.  e. ~ on ~ (0, T),~ \forall \phi \in W_m. \end{align*}
Using integration by parts, equation \eqref{model1}, and the boundary condition \eqref{bc's}, we obtain: 
\begin{align}\label{3.9}
    \lim_{n\to\infty} \big(\boldsymbol{u}_n(t)\cdot \boldsymbol{\nabla}a_n(t) , \phi\big) = \big( \boldsymbol{u}(t)\cdot \boldsymbol{\nabla}a(t), \phi\big) ~ ~  a.  e. ~ on ~ (0, T),~ \forall \phi \in W_m. \end{align}
Employing Holder's inequality, we write  
\begin{align*}
   \br{a_n(t) b_n(t)- a(t) b(t)}{\phi}  &=   \br{a_n(t) b_n(t) - a_n(t) b(t) + a_n(t) b(t) - a(t) b(t)}{\phi}  \nonumber \\ & \leq \nrm{b_n(t) - b(t)} \nrm{a_n(t) \phi} + \nrm{a_n(t) - a(t)} \nrm{b(t) \phi},
\end{align*}
 $a. e. ~ on ~ (0, T)$. From the strong convergence of $a_n$ and $b_n$ in $L^2(0,T;L^2)$ we have $\left\| a_n(t)- a(t)  \right\|_{L^2} \rightarrow 0 $ and $\left\| b_n(t)- b(t)  \right\|_{L^2} \rightarrow 0 $ as $n \rightarrow \infty$ $~a.  e. ~ on ~ (0, T)$. Having these convergence results and the observation that $\left\|a_n(t) \phi  \right\|_{L^2}$ and $\left\| b(t)  \phi  \right\|_{L^2}$ are uniformly bounded, we get
 \begin{align}\label{3.10}
    \lim_{n \rightarrow \infty} \left\langle a_n(t) b_n(t), \phi \right\rangle = \left\langle a(t) b(t), \phi \right\rangle \quad \text{a.e. on } (0, T),\quad \forall \, \phi \in W_m.
\end{align}
Since $m \in \mathbb{N}$ is chosen arbitrarily, the convergence results in equations \eqref{3.9} and \eqref{3.10} hold for all $\phi \in \bigcup_{m \geq 1} W_m$. Noting that $\bigcup_{m \geq 1} W_m$ is dense in $H^1(\Omega)$, we extend these results for all $\phi \in H^1(\Omega)$. Consequently, using these convergence results in equation \eqref{3.8}, we conclude
\begin{align}
    \innerproduct{\frac{\partial a(t)}{\partial t}, \phi} + \br{\boldsymbol{u}\cdot \boldsymbol{\nabla} a(t)}{\phi} + d \br{\boldsymbol{\nabla}a(t)}{\boldsymbol{\nabla}\phi}+k \br{a(t) b(t)}{\phi} = 0,  
\end{align}
 $ a.  e. ~ on ~ (0,T) \text{ and } \forall ~ \phi \in H^{1}(\Omega).$
 Similarly, passing to the limit in \eqref{finite weak 3} and \eqref{finite weak 4} yields:
 \begin{align}   \innerproduct{\frac{\partial b(t)}{\partial t}, \phi} + \br{\boldsymbol{u}\cdot \boldsymbol{\nabla} b(t)}{\phi} + d \br{\boldsymbol{\nabla}b(t)}{\boldsymbol{\nabla}\phi}+k \br{a(t) b(t)}{\phi} = 0, \\\innerproduct{\frac{\partial c(t)}{\partial t}, \phi} + \br{\boldsymbol{u}\cdot \boldsymbol{\nabla} c(t)}{\phi} + d \br{\boldsymbol{\nabla}c(t)}{\boldsymbol{\nabla}\phi}-k \br{a(t) b(t)}{\phi} = 0, \end{align}$ a.  e. ~ on ~ (0, T) \text{ and } \forall ~ \phi \in H^1(\Omega)$. 
Utilizing Hölder's inequality and the Lipschitz continuity of the exponential function, we arrive at
\begin{align*}
   \left|\big(e^{\alpha c_n}\boldsymbol{u}_n-e^{\alpha c}\boldsymbol{u}, \boldsymbol{v}\big)\right| &= \left|\big(e^{\alpha c_n}\boldsymbol{u}_n - e^{\alpha c_n}\boldsymbol{u} + e^{\alpha c_n}\boldsymbol{u} - e^{\alpha c}\boldsymbol{u}, \boldsymbol{v}  \big)\right|\\ & \leq \left|\big( e^{\alpha c_n}(\boldsymbol{u}_n - \boldsymbol{u}), \boldsymbol{v}  \big)\right|  + \left|(e^{\alpha c_n} - e^{\alpha c})\boldsymbol{u},\boldsymbol{v}) \right|  \\ & \leq e^{\alpha \norm{c_n}} \nrm{\boldsymbol{u}_n - \boldsymbol{u}} \nrm{\boldsymbol{v}} + \nrm{e^{\alpha c_n}-e^{\alpha c}} \nrm{\boldsymbol{u} \cdot \boldsymbol{v}} \\ & \leq e^{\alpha \norm{c_n}} \nrm{\boldsymbol{u}_n - \boldsymbol{u}} \nrm{\boldsymbol{v}} + M \nrm{ c_n- c} \nrm{\boldsymbol{u} \cdot \boldsymbol{v}}. 
\end{align*}
Now, using $\nrm{\boldsymbol{u}_n - \boldsymbol{u}} \to 0$ and $\nrm{c_n - c} \to 0$ as $n \to \infty$, along with the fact that the other terms on the right-hand side of the above inequality are uniformly bounded, we conclude
\begin{align*}
    \lim_{n \to \infty} \big(e^{\alpha c_n(t)}\boldsymbol{u}_n(t),\boldsymbol{v}(t)\big) = \big(e^{\alpha c(t)}\boldsymbol{u}(t), \boldsymbol{v}(t)\big) \quad \text{a.e. on } (0, T), \quad \forall \, \boldsymbol{v} \in V_m.
\end{align*}
Again we see
\begin{align*}
    &\left|\Big((1+ R_{A} a_n(t) +  R_{B} b_n(t) +  R_{C}c_n(t))\boldsymbol{g} - (1+ R_{A} a(t) +  R_{B} b(t) +  R_{C}c(t))\boldsymbol{g},\boldsymbol{v}\Big) \right| \\ &\leq \norm{\boldsymbol{g}}\Big( R_{A}\nrm{a_n(t) - a(t)} + R_{B}\nrm{b_n(t) - b(t)} + R_{B}\nrm{c_n(t) - c(t)}\Big)\nrm{\boldsymbol{v}}
\end{align*}
$\text{a.e. on } (0, T).$
Utilizing the strong convergence of $a_n$, $b_n$, and $c_n$ in $L^2(\Omega)$ space, we conclude:
\begin{align*}
    \lim_{n \to \infty} &\Big(\big(1 + R_{A}a_n(t) + R_{B}b_n(t) + R_{C}c_n(t)\big)\boldsymbol{g}, \boldsymbol{v}\Big) = \Big ( \big(1 + R_{A}a(t) + R_{B}b(t) + R_{C}c(t)\big) \boldsymbol{g},\boldsymbol{v}\Big)\quad \text{a.e. on } (0, T).
\end{align*}
Thus, we have proven the existence of at least one solution to the problem described by equations \eqref{model} and \eqref{boundry and initial data}.
\end{proof}
\section{ Continuous Dependence on the Initial Data and Uniqueness of Solution}\label{continuous}
\begin{thm}\label{continuous_dependence}
 Let \(( a_{i}, b_{i}, c_{i}, \boldsymbol{u}_{i}), i=1,2\) be two solutions to system \eqref{model}-\eqref{boundry and initial data}, corresponding to the initial data \((a_{i0}, b_{i0}, c_{i0}, \boldsymbol{u}_{i0})\) for \(i=1,2\), respectively. Then, the following continuous dependence estimate holds:
 \begin{align*}
  &\nr{a_1(\tau)-a_2(\tau)} + \nr{b_1(\tau)-b_2(\tau)} + \nr{c_1(\tau) - c_2(\tau)} + \nr{\boldsymbol{u}_1(\tau)-\boldsymbol{u}_2(\tau)}  \\  
  &\leq \left(\nr{a_{10} - a_{20}} + \nr{b_{10}- b_{20}} + \nr{c_{10} - c_{20}} + \nr{\boldsymbol{u}_{10} - \boldsymbol{u}_{20}} \right) e^{\int_{0}^{\tau} \phi(s) \, ds},
\end{align*}   
for all \(\tau \in (0,T]\). Here, \(\phi(s)\) is a non-negative function in \(L^1(0, T)\).
\end{thm}
\begin{proof}
Let $( a_1, b_1, c_1, \boldsymbol{u}_1)$ and $( a_2, b_2, c_2, \boldsymbol{u}_2)$ be two solution for problem given by \eqref{model}-\eqref{boundry and initial data}. Then putting these solutions in equations \eqref{weak1}-\eqref{weak4}, and subtracting resulting equations, we have  
\begin{align}\label{26}
    \innerproduct{\frac{\partial \boldsymbol{u}(t)}{\partial t}, \boldsymbol{v}} +  \mu \big( e^{\alpha c_1(t)} \boldsymbol{u}_1(t)- e^{\alpha c_2(t)} \boldsymbol{u}_2(t) , \boldsymbol{v} \big) + \mu_{e}\big( \boldsymbol{\nabla}\boldsymbol{u}(t), \boldsymbol{\nabla}\boldsymbol{v}\big) = - \big((R_{A}a(t) + R_{B}b(t) + R_{C}c(t) )\boldsymbol{g}, \boldsymbol{v} \big),\end{align}
   \begin{align}\label{27}
    \innerproduct{\frac{\partial a(t)}{\partial t}, \phi_1} + \br{\boldsymbol{u}_1(t)\cdot \boldsymbol{\nabla}a_1(t) - \boldsymbol{u}_2(t)\cdot \boldsymbol{\nabla}a_2(t)}{\phi_1} + d_A\br{\boldsymbol{\nabla}a(t)}{\boldsymbol{\nabla}\phi_1}= - k \br{a_1(t) b_1(t) - a_2(t) b_2(t)}{\phi_1}, \end{align}
    \begin{align}\label{28}
   \innerproduct{\frac{\partial b(t)}{\partial t}, \phi_2}  +  \br{\boldsymbol{u}_1(t)\cdot \boldsymbol{\nabla}b_1(t) - \boldsymbol{u}_2(t)\cdot \boldsymbol{\nabla}b_2(t)}{\phi_2} + d_B\br{\boldsymbol{\nabla}b(t)}{\boldsymbol{\nabla}\phi_2}= - k \br{a_1(t) b_1(t) - a_2(t) b_2(t)}{\phi_2},  \end{align}
   \begin{align}\label{29}
     \innerproduct{\frac{\partial c(t)}{\partial t}, \phi_3} + \br{\boldsymbol{u}_1\cdot \boldsymbol{\nabla}c_1(t) - \boldsymbol{u}_2\cdot \boldsymbol{\nabla}c_2(t)}{\phi_3} + d_C \br{\boldsymbol{\nabla}c_1(t)}{\boldsymbol{\nabla}\phi_3} = k \br{a_1(t) b_1(t) - a_2(t) b_2(t)}{\phi_3},
\end{align}
for all $\phi_1, \phi_2, \phi_3 \in H^1(\Omega)$ and for all $\boldsymbol{v} \in V$.
Where $\boldsymbol{u} = \boldsymbol{u}_1 - \boldsymbol{u}_2, a = a_1 - a_2, b = b_1 - b_2, c = c_1 - c_2 ~ a.  e. ~ on ~ (0, T).$ Now, using $\boldsymbol{v} = \boldsymbol{u}(t)$ in equation \eqref{26}, and utilizing the fact that $c_1$ and $c_2$ are positive functions, as well as applying H\"older's and Young's inequalities, we obtain:
\begin{align}\label{33}
  &\frac{1}{2}\frac{d}{dt}\nr{\boldsymbol{u}(t)}  + \mu \nr{\boldsymbol{u}(t)} + \mu_e \nr{\boldsymbol{\nabla} \boldsymbol{u}(t)}
  \leq - \mu \br{e^{\alpha c_1(t)}- e^{\alpha c_2(t)}\boldsymbol{u}_2(t)}{\boldsymbol{u}(t)}\nonumber \\ &  + \frac{\|g\|_{L^{\infty}(\Omega)}}{2} \Big(R_{A}^{2} \nr{a(t)} + R_{B}^{2} \nr{b(t)} + R_{C}^{2} \nr{c(t)} + 3 \nr{\boldsymbol{u}(t)} \Big), ~ a.  e. ~ on ~ (0, T)
\end{align}
Putting $\phi_1 = a(t)$, $\phi_2 = b(t)$, and $\phi_3 = c(t)$ in equation's \eqref{27}, \eqref{28} and \eqref{29} respectively and  rearranging the convection terms, we have  
\begin{align}
  \frac{1}{2}\frac{d}{dt}\nr{a(t)} + d_A \nr{\boldsymbol{\nabla} a(t)} = - \left(\boldsymbol{u}(t)\cdot\boldsymbol{\nabla}a_1(t),a(t)\right) 
       - \left(\boldsymbol{u}_2(t) \cdot \boldsymbol{\nabla}a(t),a(t)\right) 
       - k \br{a_1(t) b_1(t) - a_2(t) b_2(t)}{a(t)},\label{34} \\
    \frac{1}{2}\frac{d}{dt}\nr{b(t)} + d_B \nr{\boldsymbol{\nabla} b(t)} = - \left(\boldsymbol{u}(t)\cdot\boldsymbol{\nabla}b_1(t),b(t)\right) 
       - \left(\boldsymbol{u}_2(t) \cdot \boldsymbol{\nabla}b(t),b(t)\right) 
       - k \br{a_1(t) b_1(t) - a_2(t) b_2(t)}{b(t)},\label{35} \\
    \frac{1}{2}\frac{d}{dt}\nr{c(t)} + d_C \nr{\boldsymbol{\nabla} c(t)} = - \left(\boldsymbol{u}(t)\cdot\boldsymbol{\nabla}c_1(t),c(t)\right) 
       - \left(\boldsymbol{u}_2(t) \cdot \boldsymbol{\nabla}c(t),c(t)\right) 
       + k \br{a_1(t) b_1(t) - a_2(t) b_2(t)}{c(t)},\label{36} 
\end{align}
$ a.  e. ~ on ~ (0, T)$. Utilizing integration by parts with the divergence-free condition on velocity and the boundary conditions, the second term on the right-hand side of equations \eqref{34}, \eqref{35}, and \eqref{36} will vanish. Now, using Holder's and Young's inequalities, we estimate the following terms. 
\begin{align}
  \Big| \big(\boldsymbol{u}(t)\cdot\boldsymbol{\nabla}a_{1}(t),a(t)\big) \Big| 
  &\leq \frac{\epsilon}{2} \|\boldsymbol{\nabla}a(t)\|_{L^2}^2 + \frac{1}{2 \epsilon}\|a_1(t)\|_{L^{\infty}(\Omega)}^2 \|\boldsymbol{u}(t)\|_{L^2}^2, \label{37} \\
  \Big| \big(\boldsymbol{u}(t)\cdot\boldsymbol{\nabla}b_{1}(t),b(t)\big) \Big| 
  &\leq \frac{\epsilon}{2} \|\boldsymbol{\nabla}b(t)\|_{L^2}^2 + \frac{1}{2 \epsilon}\|b_1(t)\|_{L^{\infty}(\Omega)}^2 \|\boldsymbol{u}(t)\|_{L^2}^2, \label{38} \\
  \Big| \big(\boldsymbol{u}(t)\cdot\boldsymbol{\nabla}c_{1}(t),c(t)\big) \Big| 
  &\leq \frac{\epsilon}{2} \|\boldsymbol{\nabla}c(t)\|_{L^2}^2 + \frac{1}{2 \epsilon}\|c_1(t)\|_{L^{\infty}(\Omega)}^2 \|\boldsymbol{u}(t)\|_{L^2}^2, \label{39}\end{align}
  \begin{align}\label{40}
    \big(a_1(t) b_1(t) - a_2(t) b_2(t), a(t) \big) &= \big(a_1(t) b_1(t) - a_2(t) b_1(t) + a_2(t) b_1(t) - a_2(t) b_2(t), a(t) \big)  \nonumber \\
    &\leq \|b_1(t)\|_{L^{\infty}(\Omega)} \|a(t)\|_{L^2}^2 + \frac{\|a_2(t)\|_{L^{\infty}(\Omega)}}{2} \big(\|a(t)\|_{L^2}^2 + \|b(t)\|_{L^2}^2\big),
  \end{align}
   \begin{align}\label{41}
    \big(a_1(t) b_1(t) - a_2(t) b_2(t), b(t) \big)
    &= \big(a_1(t) b_1(t) - a_1(t) b_2(t) + a_1(t) b_2(t) - a_2(t) b_2(t), b(t) \big) \nonumber \\
    &\leq \|a_1(t)\|_{L^{\infty}(\Omega)} \|b(t)\|_{L^2}^2 + \frac{\|b_2(t)\|_{L^{\infty}(\Omega)}}{2} \big(\|a(t)\|_{L^2}^2 + \|b(t)\|_{L^2}^2\big),    
  \end{align}
   \begin{align}\label{42}  
    \big(a_1(t) b_1(t) - a_2(t) b_2(t), c(t) \big)
    &= \big(a_1(t) b_1(t) - a_1(t) b_2(t) + a_1(t) b_2(t) - a_2(t) b_2(t), c(t) \big) \nonumber \\
    &\leq \frac{\|a_1(t)\|_{L^{\infty}(\Omega)}}{2} \big(\|b(t)\|_{L^2}^2 + \|c(t)\|_{L^2}^2\big)  + \frac{\|b_2(t)\|_{L^{\infty}(\Omega)}}{2} \big(\|a(t)\|_{L^2}^2 + \|c(t)\|_{L^2}^2\big), 
  \end{align}
$ a. e. ~ on ~ (0, T)$. Utilizing Lipschitz's continuity of exponential function and  Holder's inequality, we arrive at
\begin{align*}
   \mu\left|\Big(\left(e^{\alpha c_{1}(t)} - e^{\alpha c_{2}(t)}\right) \boldsymbol{u}_2(t), \boldsymbol{u}(t)  \Big) \right| \leq \mu M \alpha \nrm{c(t)} \|\boldsymbol{u}_2(t)\|_{L^4}\|\boldsymbol{u}(t)\|_{L^4}.
   \end{align*}
Applying Gagliardo–Nirenberg's and Young's inequalities yields
\begin{align}\label{42a}
     \mu\left|\Big(\left(e^{\alpha c_{1}(t)} - e^{\alpha c_{2}(t)}\right) \boldsymbol{u}_2(t), \boldsymbol{u}(t)  \Big) \right|  &\leq \mu M \alpha \nrm{c(t)} \nrm{\boldsymbol{u}_2(t)}^{1/4} \nrm{\boldsymbol{\nabla}\boldsymbol{u}_2(t)}^{3/4}  \nrm{\boldsymbol{u}(t)} ^{1/4} \nrm{\boldsymbol{\nabla}\boldsymbol{u}(t)}^{3/4} \nonumber\\ &\leq \frac{1}{2} \nr{c(t)}\nrm{\boldsymbol{\nabla}\boldsymbol{u}_2(t)}^{3/2} + \frac{3\epsilon}{4}\nr{\boldsymbol{\nabla}\boldsymbol{u}(t)} + \frac{\mu^{8} M^{8} \alpha^{8}}{64 \epsilon^{3}} \nr{\boldsymbol{u}_2(t)} 
    \nr{\boldsymbol{u}(t)},  \end{align}
$a. e. ~ on ~ (0, T)$. Using the estimates from inequalities \eqref{37}-\eqref{42a} in equations \eqref{34}, \eqref{35}, and \eqref{36}, and then summing those equations with equation \eqref{33}, we have:
\begin{align*}
&\frac{d}{dt}\Big( \nr{a} + \nr{b} + \nr{c} +  \nr{\boldsymbol{u}} \Big) + 2 \mu \nr{\boldsymbol{u}} +\left(2\mu_e - \frac{3\epsilon}{2}\right) \nr{\boldsymbol{\nabla}\boldsymbol{u}}\\ & + \big( 2 d_A - \epsilon \big) \nr{\boldsymbol{\nabla}a} + \big(  2d_B - \epsilon \big) \nr{\boldsymbol{\nabla}b}  +\big( 2 d_C - \epsilon \big) \nr{\boldsymbol{\nabla}c}\\ &\leq  \frac{1}{\epsilon}\big(\norm{a_1}^{2} + \norm{b_1}^{2} + \norm{c_1}^{2}  \big) \nr{\boldsymbol{u}} +2\|b_1(t)\|_{L^{\infty}(\Omega)} \|a(t)\|_{L^2}^2 \\ &+ \|a_2(t)\|_{L^{\infty}(\Omega)}\big(\|a(t)\|_{L^2}^2 + \|b(t)\|_{L^2}^2\big) + 2\|a_1(t)\|_{L^{\infty}(\Omega)} \|b(t)\|_{L^2}^2 \\&+ \|b_2(t)\|_{L^{\infty}(\Omega)}\big(\|b(t)\|_{L^2}^2 + \|c(t)\|_{L^2}^2\big) + \|a_1(t)\|_{L^{\infty}(\Omega)}\big(\|b(t)\|_{L^2}^2 + \|c(t)\|_{L^2}^2\big) \nonumber\\ & + 2\|b_2(t)\|_{L^{\infty}(\Omega)}\|a(t)\|_{L^2}^2  +  \frac{\mu^{8} M^{8} \alpha^{8}}{32 \epsilon^{3}} \nr{\boldsymbol{u}_2(t)} 
    \nr{\boldsymbol{u}(t)} +  \nr{c(t)}\nrm{\boldsymbol{\nabla}\boldsymbol{u}_2(t)}^{3/2} \\  & +\|g\|_{L^{\infty}(\Omega)} \Big(R_{A}^{2} \nr{a(t)} + R_{B}^{2} \nr{b(t)} + R_{C}^{2} \nr{c(t)} + 3 \nr{\boldsymbol{u}(t)} \Big), ~ a.  e. ~ on ~ (0, T).
\end{align*}
Choosing $\epsilon > 0$ such that $\epsilon = \min\left(d_A, d_B, d_C, \mu_{e}\right)$, and then ignoring the non-negative terms on the left-hand side of the above inequality, we have:
\begin{align}\label{44}
 \frac{d}{dt}\Big(\nr{a(t)} + \nr{b(t)} + \nr{c(t)} + \nr{\boldsymbol{u}(t)}   \Big) \leq \phi(t)   \Big(  \nr{a(t)} + \nr{b(t)} + \nr{c(t)} + \nr{\boldsymbol{u}(t)}\Big)\end{align} 
$ a.  e. ~ on ~ (0, T)$. Where,
\begin{align*}
\phi(t) &=   \left(3 +\frac{1}{\epsilon} \right) \norm{a_1(t)} +  \norm{a_2(t)} + \left(2 +\frac{1}{\epsilon} \right) \norm{b_1(t)} +  \norm{b_2(t)}  + \frac{1}{\epsilon} \norm{c_1(t)}  \\ &  +\left(R_{A}^{2} + R_{B}^{2} + R_{C}^{2} + 3\right) \norm{\boldsymbol{g}}   +  \frac{\mu^{8} M^{8} \alpha^{8}}{32 \epsilon^{3}} \nr{\boldsymbol{u}_2(t)}  
    + \nrm{\boldsymbol{\nabla}\boldsymbol{u}_2(t)}^{3/2} \quad a.  e. ~ on ~ (0, T).    
\end{align*}
By applying Grönwall's inequality in inequality \eqref{44}, we obtain
\begin{align*}
  & \nr{a_1(\tau)-a_2(\tau)} + \nr{b_1(\tau)-b_2(\tau)} + \nr{c_1(\tau) - c_2(\tau)} + \nr{\boldsymbol{u}_1(\tau)-\boldsymbol{u}_2(\tau)}   \\  &\leq\left(\nr{a_{10} - a_{20}} + \nr{b_{10}- b_{20}} + \nr{c_{10} - c_{20}} + \nr{\boldsymbol{u}_{10} - \boldsymbol{u}_{20}} \right) e^{\int_{0}^{\tau} \phi(s)ds},
\end{align*}   
for all $ \tau \in (0,T].$
\end{proof}
\noindent \textbf{Uniqueness of the Weak Solution:} 
\begin{thm}
    Let \( ( a_{i}, b_{i}, c_{i}, \boldsymbol{u}_{i}), \, i=1,2 \), be two solutions to the system \eqref{model}-\eqref{boundry and initial data}, with initial data satisfying assumption \eqref{intial data assumptation}. Then, for almost all \( (t, \boldsymbol{x}) \in (0, T) \times \Omega \), we have    
\[
    a_1(t, \boldsymbol{x}) = a_2(t, \boldsymbol{x}), ~ b_1(t, \boldsymbol{x}) = b_2(t, \boldsymbol{x}), ~ c_1(t, \boldsymbol{x}) = c_2(t, \boldsymbol{x}), ~ \text{and} ~ \boldsymbol{u}_1(t, \boldsymbol{x}) = \boldsymbol{u}_2(t, \boldsymbol{x}).
    \]
\end{thm}
\begin{proof}
    Define \( a = a_1 - a_2 \), \( b = b_1 - b_2 \), \( c = c_1 - c_2 \), and \( \boldsymbol{u} = \boldsymbol{u}_1 - \boldsymbol{u}_2 \). Then, by Theorem \eqref{continuous_dependence},
    \begin{align*}
         \nr{a(t)} + \nr{b(t)} + \nr{c(t)} + \nr{\boldsymbol{u}(t)} 
        \leq  \Big( \nr{a_0} + \nr{b_0} + \nr{c_0} + \nr{\boldsymbol{u}_0} \Big) e^{\int_{0}^{t} \phi(s) \, ds},
    \end{align*}
    for all \( t \in (0, T) \). Since \(\phi(t)\) is non-negative, we have \( e^{\int_{0}^{t} \phi(s) \, ds} \leq e^{\int_{0}^{T} \phi(s) \, ds} \) for all \( t \in (0, T) \). Given the initial conditions \( a_0 = 0 \), \( b_0 = 0 \), \( c_0 = 0 \), and \( \boldsymbol{u}_0 = \boldsymbol{0} \), along with \(\phi \in L^1(0, T)\), we conclude that \( a(t) = b(t) = c(t) = \boldsymbol{u}(t) = 0 \) for all \( t \in (0, T) \), thereby establishing the uniqueness of the solution to the problem given by \eqref{model}-\eqref{boundry and initial data}.
\end{proof}
\section{Semi-discrete formulation}\label{semi}
Let $\mathcal{T}_h$ be a shape-regular family of partitions of the domain $\Omega$ into simplices (triangles or rectangles in 2D and tetrahedrons in 3D).  The mesh size $h$ is defined as $h = \max_{\mathbb{T} \in \mathcal{T}_h} \text{diam}(\mathbb{T})$. The set of all faces (edges in 2D and faces in 3D) of the partition $\mathcal{T}_h$ denoted by `$\Gamma$', with $\Gamma^{\mathrm{o}}$ and $\Gamma^{\partial}$ representing the sets of interior and boundary edges, respectively, such that $\Gamma = \Gamma^{\mathrm{o}} \cup \Gamma^{\partial}$. We define the finite element spaces as follows:
\begin{align*}
    X_{h} = \{ z \in \mathcal{C}(\Omega): z|_{\mathbb{T}} \in \mathcal{P}_{k}(\mathbb{T}), \forall ~ \mathbb{T} \in \mathcal{T}_h  \} \\
    V_h = \{ \boldsymbol{w} \in \mathcal{C}(\Omega): \boldsymbol{w}|_{\mathbb{T}} \in  (X_{h})^d, \forall ~ \mathbb{T} \in \mathcal{T}_h     \}
\end{align*}
where $\mathcal{P}_{k}(\mathbb{T})$ is the space of scalar polynomials of degree at most $k$ in each variable on the element $\mathbb{T}$. We will approximate the velocity by $\boldsymbol{u}_{h} \in V_h$, pressure by $p_h \in X_{h}$, and concentrations by $a_h, b_h, c_h \in X_{h}$, respectively. The velocity \(\boldsymbol{u}_h\) and concentrations \(a_h\), \(b_h\), and \(c_h\) satisfy the following discrete boundary conditions:  
\[
\boldsymbol{n}_e \cdot \boldsymbol{u}_h = 0, \quad \boldsymbol{n}_e \cdot \nabla a_h = 0, \quad \boldsymbol{n}_e \cdot \nabla b_h = 0, \quad \boldsymbol{n}_e \cdot \nabla c_h = 0 \quad \text{on each } e \in \Gamma^{\partial},
\]
where \(\boldsymbol{n}_e\) is the unit normal vector to the edge \(e\). Note that here we adopted the slip boundary condition as described in the experimental work of \cite{Paoli_2022}, rather than the no-slip boundary condition considered in the theoretical framework presented in the previous sections of this paper. Nevertheless, for this slip boundary condition, the well-posedness of the problem can also be established by considering suitable admissible solution spaces, as shown in \cite{migorski2019nonmonotone}.

We seek to find the semi-discrete solution $(\boldsymbol{u}_h, p_h, a_h, b_h, c_h) \in L^2(0,T;V_h) \times (L^2(0,T;X_h))^4$, such that for all test functions $\boldsymbol{w}_h \in V_h$ and $z_h \in X_h$, the following system of equations is satisfied:

\begin{align}\label{eq6}
\br{\boldsymbol{\nabla}\cdot \boldsymbol{u}_h}{z_h} = 0\\
    \innerproduct{\frac{\partial \boldsymbol{u}_h}{\partial t},\boldsymbol{w}_h} + \mu \br{e^{\alpha c_h}\boldsymbol{u}_h}{\boldsymbol{w}_h} + \mu_{e}\br{\boldsymbol{\nabla}\boldsymbol{u}_h}{\boldsymbol{\nabla}\boldsymbol{w}_h} + \br{p_{h}}{\boldsymbol{\nabla}\cdot \boldsymbol{w}_h} + \br{\rho_{h} \boldsymbol{g}}{\boldsymbol{w}_h}  = 0\\
     \innerproduct{\frac{\partial a_h}{\partial t}, z_h} + \br{\boldsymbol{u}_h \cdot \boldsymbol{\nabla}a_h}{z_h} + d \br{\boldsymbol{\nabla}a_h}{\boldsymbol{\nabla}z_h} + k \br{a_h b_h}{z_h} = 0\\
     \innerproduct{\frac{\partial b_h}{\partial t}, z_h} + \br{\boldsymbol{u}_h \cdot \boldsymbol{\nabla}b_h}{z_h} + d \br{\boldsymbol{\nabla}b_h}{\boldsymbol{\nabla}z_h} + k \br{a_h b_h}{z_h} = 0\\
     \innerproduct{\frac{\partial c_h}{\partial t}, z_h} + \br{\boldsymbol{u}_h \cdot \boldsymbol{\nabla}c_h}{z_h} + d \br{\boldsymbol{\nabla}c_h}{\boldsymbol{\nabla}z_h} - k \br{a_h b_h}{z_h} = 0
\end{align}

The initial data $\boldsymbol{u}_h(t=0)$, $a_h(t=0)$, $b_h(t=0)$, and $c_h(t=0)$ are suitable approximations of $\boldsymbol{u}(t=0)$, $a(t=0)$, $b(t=0)$, and $c(t=0)$, respectively.
\section{Numerical method and its validations for 2D simulations} \label{num}
\begin{figure}[h!]
    \centering
    \hspace{0.25 in} (a) \hspace{2 in} (b)\\
    \includegraphics[width=0.32\linewidth,clip]{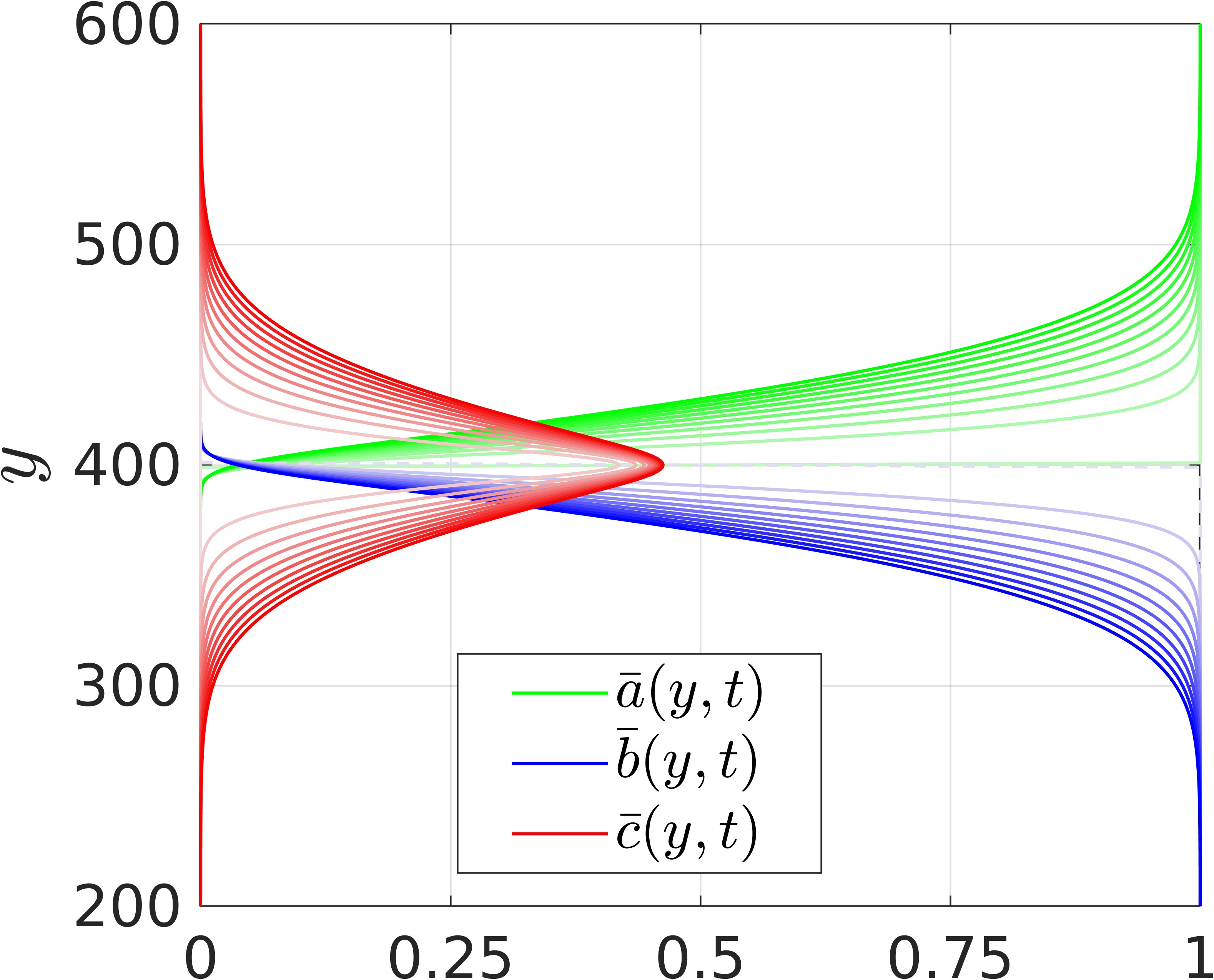}
    \includegraphics[width=0.31\linewidth,clip]{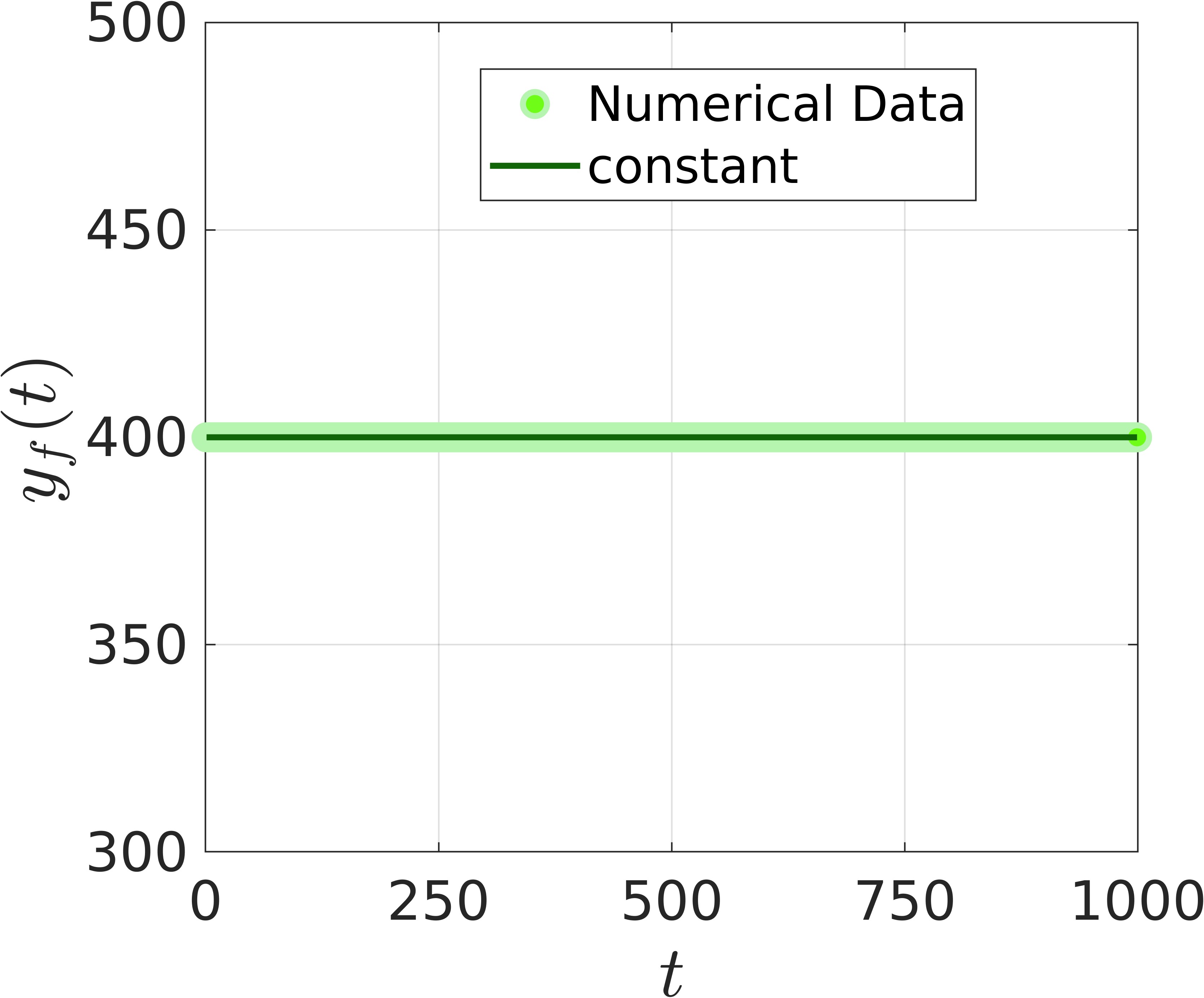}
    
    \hspace{0.13 in} (c) \hspace{2 in} (d) \hspace{2 in} (e)\\
    \includegraphics[width=0.32\linewidth]{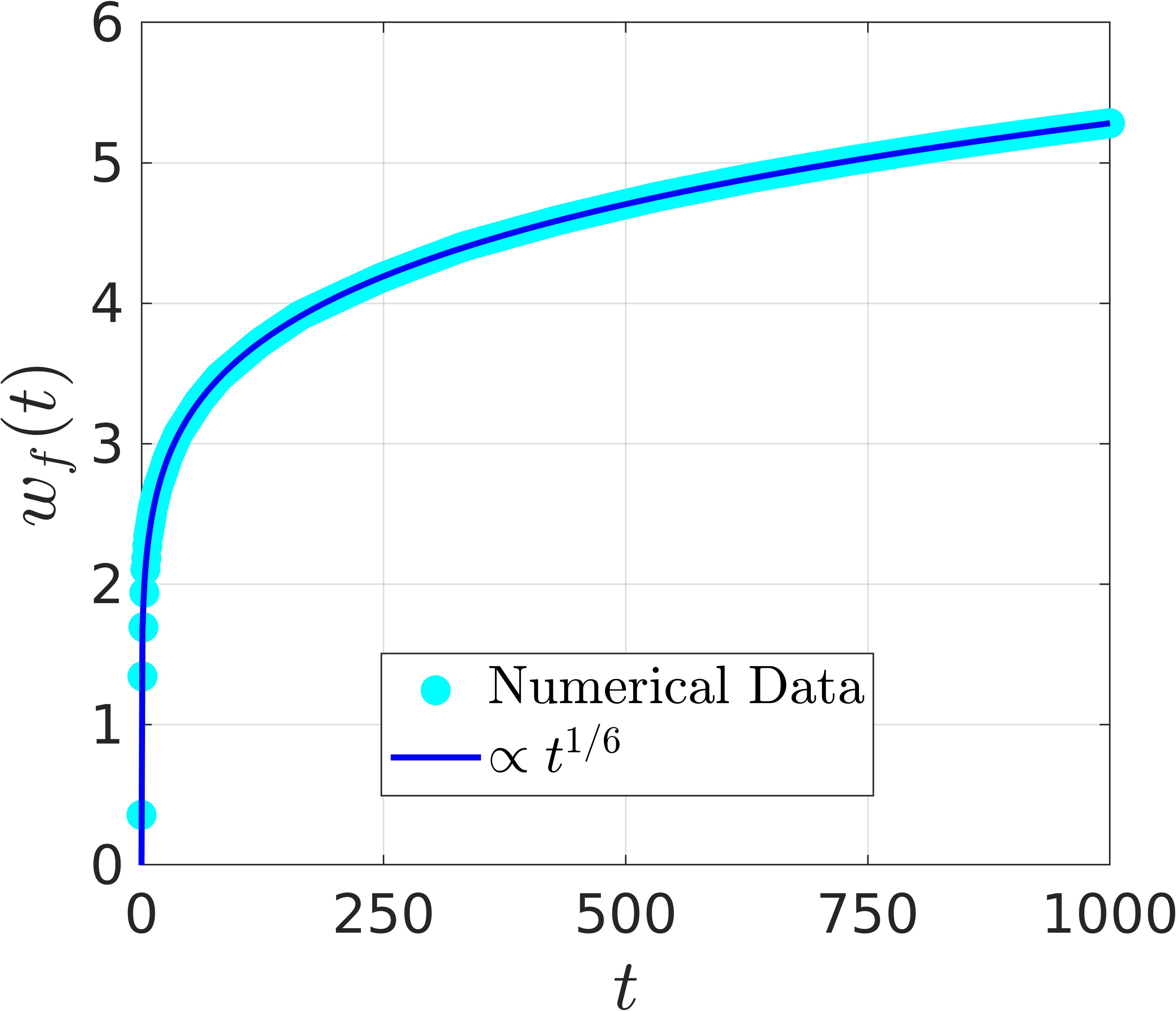}
    \includegraphics[width=0.32\linewidth]{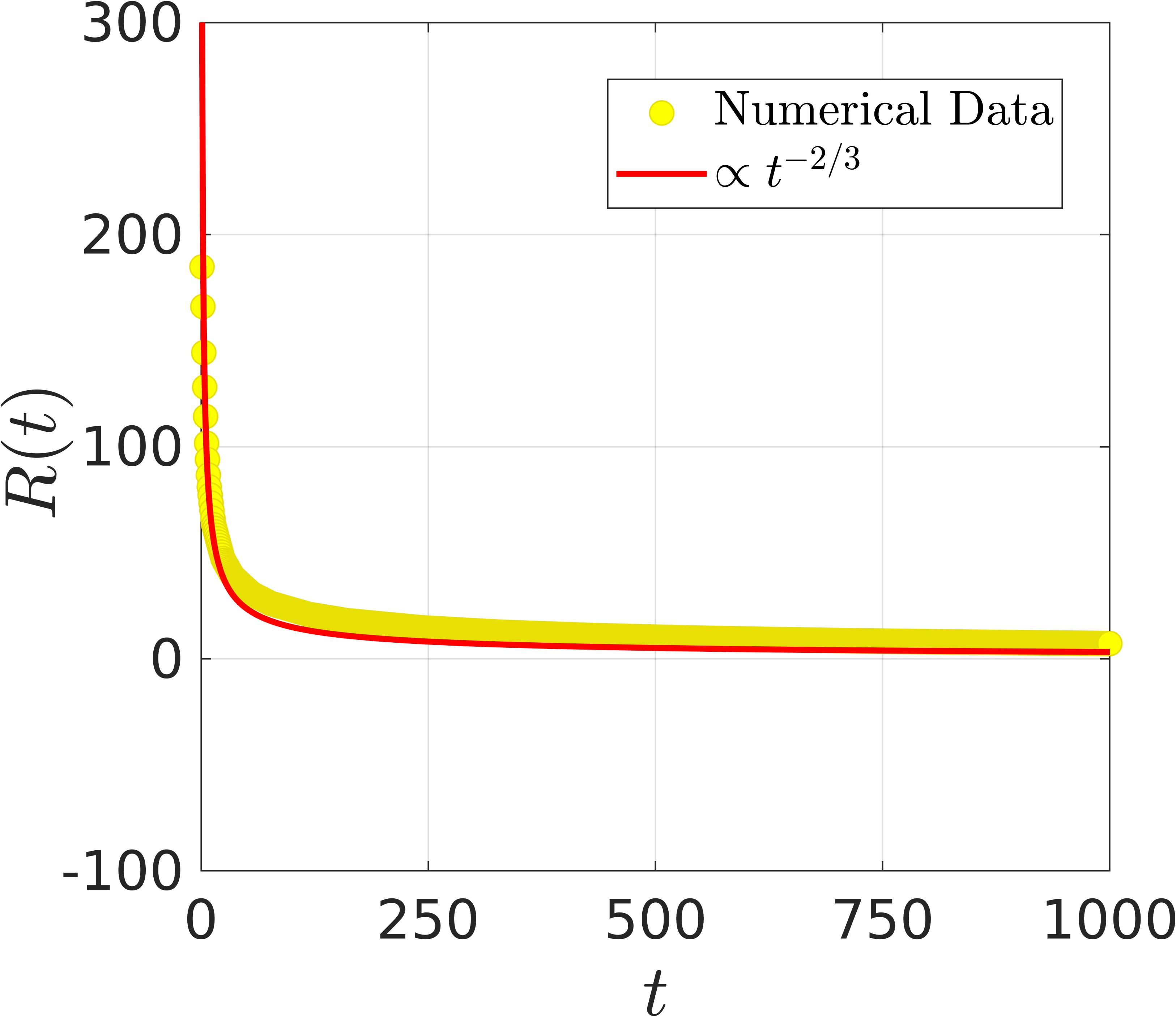}
    \includegraphics[width=0.33\linewidth]{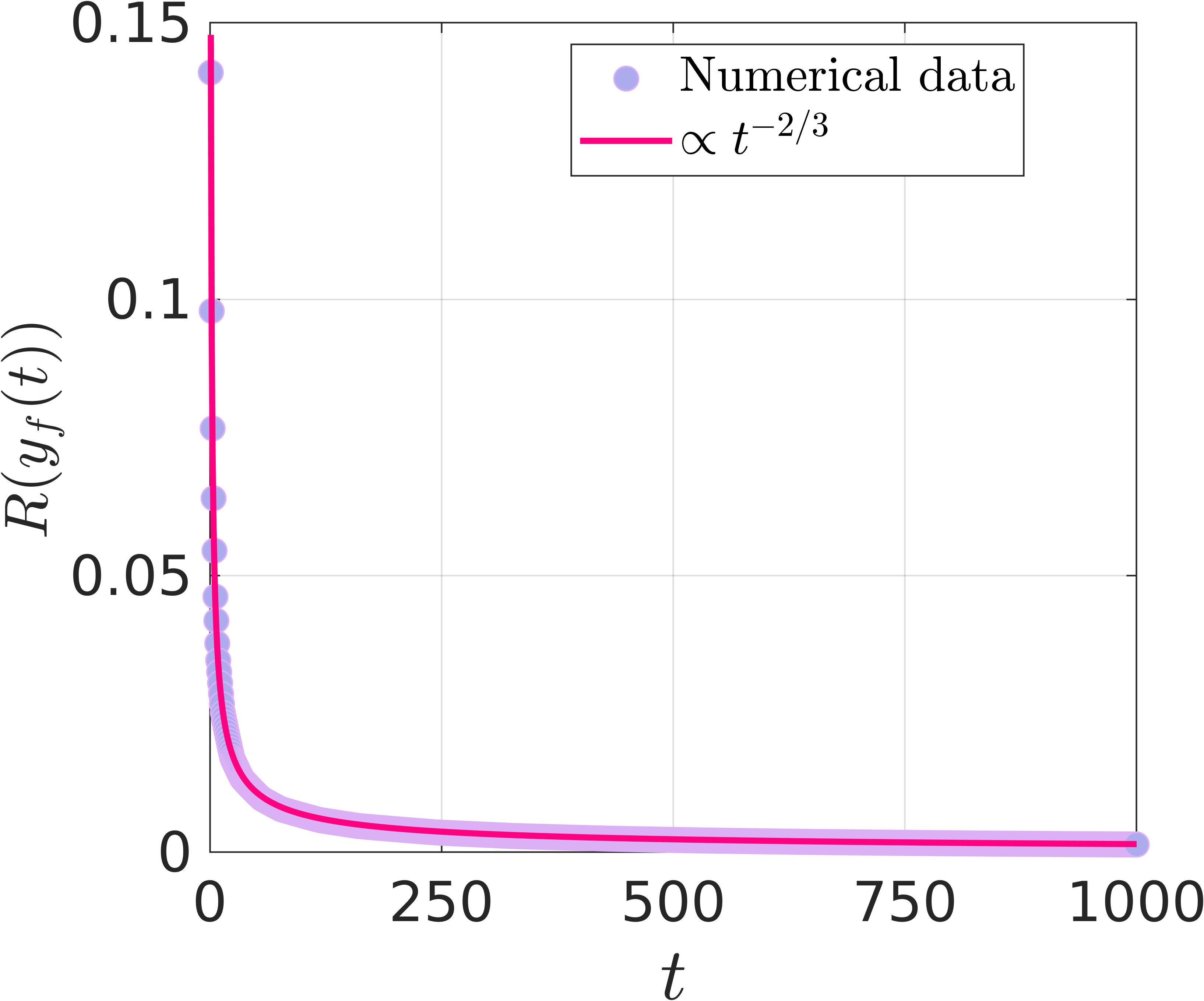}
    
   \vspace{-0.3 cm}\caption{(a) Spatio-temporal evolution of \( x \)-averaged reaction-diffusion profiles \(\bar{a}(y,t)\), \(\bar{b}(y,t)\), and \(\bar{c}(y,t)\). Darker gradients indicate increasing time from \( t = 0 \) to \( t = 1000 \), with intervals of 100 units. The diffusion ratio is fixed at \( d = 1 \). (b) Temporal evolution of the reaction front position \( y_f(t) \), fitted with a constant of 400. Temporal evolution of (c) reaction front width \( w_f(t) \), (d) total reaction rate \( R(t) \), and (e) reaction rate at the front \( R(y_f(t)) \), each with corresponding fitted power laws and proportionality constants of 2.73, 322, and 0.147, respectively.}
    \label{fig:galfi}
\end{figure}

The semi-discrete problem (5.1--5.5) is numerically solved using the Finite Element Method (FEM) solver of \textsc{COMSOL Multiphysics} \textsuperscript{\textregistered}  \cite{COMSOL2024} within a 2D rectangular domain 
$
\Omega = \{ (x, y) \mid x \in [-400, 400], \, y \in [0, 800] \}
$.
The default time stepping in COMSOL is chosen for the simulation, where it employs an adaptive temporal solver utilizing first and second-order backward difference formulae with an initial step size $\Delta t =10^{-6}$. For spatial discretization, we use Lagrange \(\mathcal{P}_2\) elements for the velocity \(\boldsymbol{u}_h\) and \(\mathcal{P}_1\) elements for pressure \(p_h\) and concentrations \(a_h\), \(b_h\), and \(c_h\), defined on a mapped mesh. The fully discretized system is solved using the PARDISO \cite{schenk2002pardiso} direct sparse solver with Bunch-Kaufman pivoting \cite{bunch1977pivoting} and row preordering. The relative tolerance is fixed as 0.005, ensuring robust and efficient computation of fluid flow variables. For the flat interface 2D simulations (discussed in the following sections), the maximum element size is \( h = 1 \), resulting in a total of 10,267,231 degrees of freedom.

\subsection{Validation with classical reaction-diffusion results}
We validate our 2D simulation results against the scalings proposed by G\'alfi and R\'acz \cite{galfi1988}. Their reaction-diffusion theory is established for 1D results, so we average 2D stable solutions to enable direct comparison. Here, with an initial zero velocity field, \(\boldsymbol{u}_0(\boldsymbol{x})=\boldsymbol{0}\), no gravitational field, \(\boldsymbol{g}=\boldsymbol{0}\), and initial concentrations defined as
 \begin{align}\label{e71}
a_0(\boldsymbol{x}) = \frac{\alpha_0}{2} + \frac{\alpha_0}{2} \cdot \operatorname{erf}\left(\frac{y - y_0}{\delta}\right),\quad
b_0(\boldsymbol{x}) = \frac{\beta_0}{2} - \frac{\beta_0}{2} \cdot \operatorname{erf}\left(\frac{y - y_0}{\delta}\right),\quad 
c_0(\boldsymbol{x}) = 0,
\end{align}   
where \(\alpha_0 = \beta_0 = 1\), \(y_0 = 400\), and \(\delta = 10^{-5}\), the problem \eqref{model1}-\eqref{model5} reduces to an uncoupled 2D reaction-diffusion system. We average the 2D concentrations \(a\), \(b\), and \(c\) along the \(x\)-direction to obtain 1D reaction-diffusion profiles \(\bar{a}(y,t)\), \(\bar{b}(y,t)\) and \(\bar{c}(y,t)\), respectively.  For instance,  \(\bar{a}(y,t)\) is defined as $$\bar{a}(y,t) = \frac{1}{L_x} \int_{x_{min}}^{x_{max}} a(x,y,t) \, dx,$$ where $L_x$ is the domain length and $x_{min}$ and $x_{max}$ are the locations of left and right boundary of the domain. The spatio-temporal evolution of 1D reaction-diffusion profiles \(\bar{a}(y,t)\), \(\bar{b}(y,t)\) and \(\bar{c}(y,t)\), are shown in Fig. \ref{fig:galfi}(a).
If the diffusion coefficients are the same for all species, \(A\), \(B\), and \(C\) (i.e., \(d_A = d_B = d_C = d\)), then according to G\'alfi and R\'acz's theory, the reaction front position,
$$
y_f(t) = \frac{\int_{y_{min}}^{y_{max}} y \, \bar{R}(y,t) \, dy}{\int_{y_{min}}^{y_{max}} \bar{R}(y,t) \, dy},
$$
where \(\bar{R}(y,t) = \bar{a}(y,t) \cdot \bar{b}(y,t)\), remains fixed at the initial front position. Our numerical results confirm this, showing \(y_f(t) = \text{constant} = 400\) (Fig. \ref{fig:galfi}(b)).
The width of the reaction front, \(w_f(t)\), where
$$
w_f^2(t) = \frac{\int_{y_{min}}^{y_{max}} (y - y_f)^2 \, \bar{R}(y,t) \, dy}{\int_{y_{min}}^{y_{max}} \bar{R}(y,t) \, dy},
$$
is predicted to evolve as \(w_f(t) \propto t^{1/6}\) according to \cite{galfi1988}, which our numerical data, shown in Fig. \ref{fig:galfi}(c), also support. Finally, the total reaction rate $$ R(t) = \int_{y_{min}}^{y_{max}} \bar{R}(y,t) \, dy $$ and the reaction rate at the front, \( R(y_f(t)) \), are expected to scale as \( \propto t^{-2/3} \) per G\'alfi and R\'acz's theory. Our numerical results are in excellent agreement with these predictions, as shown in Fig. \ref{fig:galfi}(d) and Fig. \ref{fig:galfi}(e), respectively.

\begin{table}[h!]
\centering
\begin{tabular}{ccccccc}
\hline
\hline
$h$ & 0.8 & 1 & 2 & 4 & 8  \\
DoF &16034031&10267231&2573631&646831&163431\\
\hline
\end{tabular}
\caption{Various $h$ and corresponding Degrees of Freedoms (DoF) in the mesh independence test}
\label{tab2}
\end{table}
\subsection{Scheme convergence test for unstable case}
We also test our numerical scheme on an unstable case where the reaction modifies the permeability and density profiles, leading to the formation of falling fingers. For this test, we set the parameters as solute expansion coefficients \( R_A = R_B = 0 \), \( R_C = 2 \), mobility ratio \( \alpha = 2 \), diffusion coefficient \( d = 0.005 \), initial conditions for concentration given by Eq.~\eqref{e71}, and a zero initial velocity field for the model \eqref{model1}-\eqref{model5}.
To assess the mesh independence of the numerical solutions, the mesh is refined by varying the maximum element size $h$, from the set $\{8,4,2,1,0.8\}$ (corresponding degrees of freedom are given in Table \ref{tab2}). For this purpose, we consider the solution obtained from the finest grid, $h = 0.8$, as our reference for the most accurate results and the computed quantities with this element size is denoted by a subscript $e$.
Next, we calculate the relative error in the interfacial length of the reactant $A$, denoted as $E_{A,h}(t)$, the interfacial length of the reactant $B$, denoted as $E_{B,h}(t)$, and the mixing length from the product $C$ concentration, denoted as $E_{ml,h}(t)$, with respect to other element sizes. These errors are defined as follows:
\begin{equation}
\begin{aligned}
    E_{A,h}(t) &= \frac{\left| I_{A,h}(t) - I_{A,e}(t) \right|}{I_{A,e}(t)}, \\
    E_{B,h}(t) &= \frac{\left| I_{B,h}(t) - I_{B,e}(t) \right|}{I_{B,e}(t)}, \\
    E_{ml,h}(t) &= \frac{\left| ml_{h}(t) - ml_{e}(t) \right|}{ml_{e}(t)}.
\end{aligned}
\end{equation}
where the subscripts $h$ denote that the values are computed by varying $h$.
Here, the interfacial lengths, $I_{A}(t)$ and $I_{B}(t)$ are defined as:
\begin{equation}
        I_A(t) = \int_{\Omega} \sqrt{\left(\frac{\partial a}{\partial x}\right)^2 + \left(\frac{\partial a}{\partial y}\right)^2} \, d\Omega, \quad
        I_B(t) = \int_{\Omega} \sqrt{\left(\frac{\partial b}{\partial x}\right)^2 + \left(\frac{\partial b}{\partial y}\right)^2} \, d\Omega.
\end{equation}
The mixing length \(ml(t)\) is computed as the distance between the \(y\)-locations \(y_{\text{upper}}\) and \(y_{\text{lower}}\), where the $y$-averaged product's concentration \(\bar{c}(y, t) > 0.01\) from above and below, respectively:
\begin{equation}
    ml(t) = y_{\text{upper}} - y_{\text{lower}},
\end{equation}
where
\[
y_{\text{upper}} = \max \{ y \mid \bar{c}(y, t) > 0.01 \}, \quad y_{\text{lower}} = \min \{ y \mid \bar{c}(y, t) > 0.01 \}.
\]
The interfacial lengths and mixing lengths are critical measures of fingering instability growth \cite{Gopalakrishnan_2021, Sharma__2019}. Therefore, we incorporate these measures into our error analysis to perform a mesh independence study.

\begin{figure}
    \centering
    \includegraphics[width=0.4\linewidth]{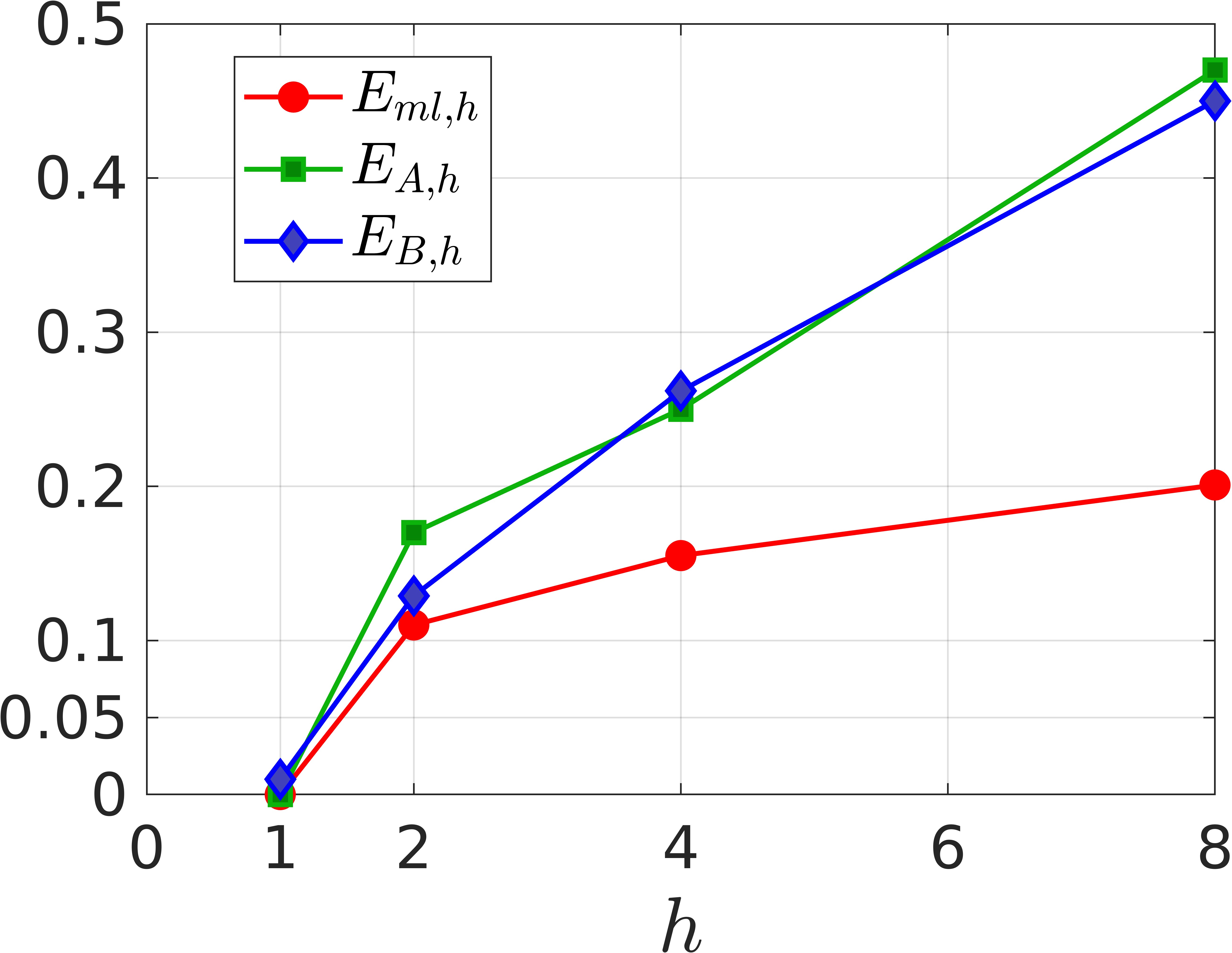}
    \vspace{-0.3 cm}
    \caption{The maximum relative errors $E_{A,h}$, $E_{B,h}$, and $E_{ml,h}$ as functions of the maximum element size $h$.}
    \label{fig:mesh_ind}
\end{figure}

Figure \ref{fig:mesh_ind} presents the $L^\infty$ norm of the relative errors $E_{A,h}(t)$, $E_{B,h}(t)$, and $E_{ml,h}(t)$ over the numerical time interval $[0,T]$ as functions of the maximum element size $h$. These errors are defined as:
\[
E_{A,h} = \| E_{A,h}(t) \|_{L^\infty(0,T)}, \quad
E_{B,h} = \| E_{B,h}(t) \|_{L^\infty(0,T)}, \quad
E_{ml,h} = \| E_{ml,h}(t) \|_{L^\infty(0,T)}.
\]
As $h$ decreases, the errors $E_{A,h}$, $E_{B,h}$, and $E_{ml,h}$ decrease, demonstrating convergence of the numerical solutions.
The smallest error is observed for \( h = 1 \), where the maximum relative errors \( E_{A,h} \), \( E_{B,h} \), and \( E_{ml,h} \) are reduced below \( 2\% \), which is well within acceptable limits, given the highly unstable nature of the flow. Consequently, \( h = 1 \) is used for subsequent analysis.
\section{Towards Pattern Formation} \label{pattern}
From this point onward, we delve into the influence of initial conditions on the evolution of solutions as fingering patterns. While the theory is established for a broad class of initial conditions, we focus our analysis on two specific configurations: (1) a flat interface and (2) an elliptical interface, both motivated by previous studies \cite{Almarcha_2010, Jha_2023}. Additionally, we investigate the effects of various parameters on fingering instabilities to bridge the gap between experimental observations and simulations. In subsequent subsections, we extend this discussion to the three-dimensional case, demonstrating the well-posedness and relevance of the system. This exploration strengthens the connection between theoretical analysis and real-world phenomena.

\subsection{2D initial flat interface}
\begin{figure}[h!]
    \centering
    \includegraphics[width=0.5\linewidth]{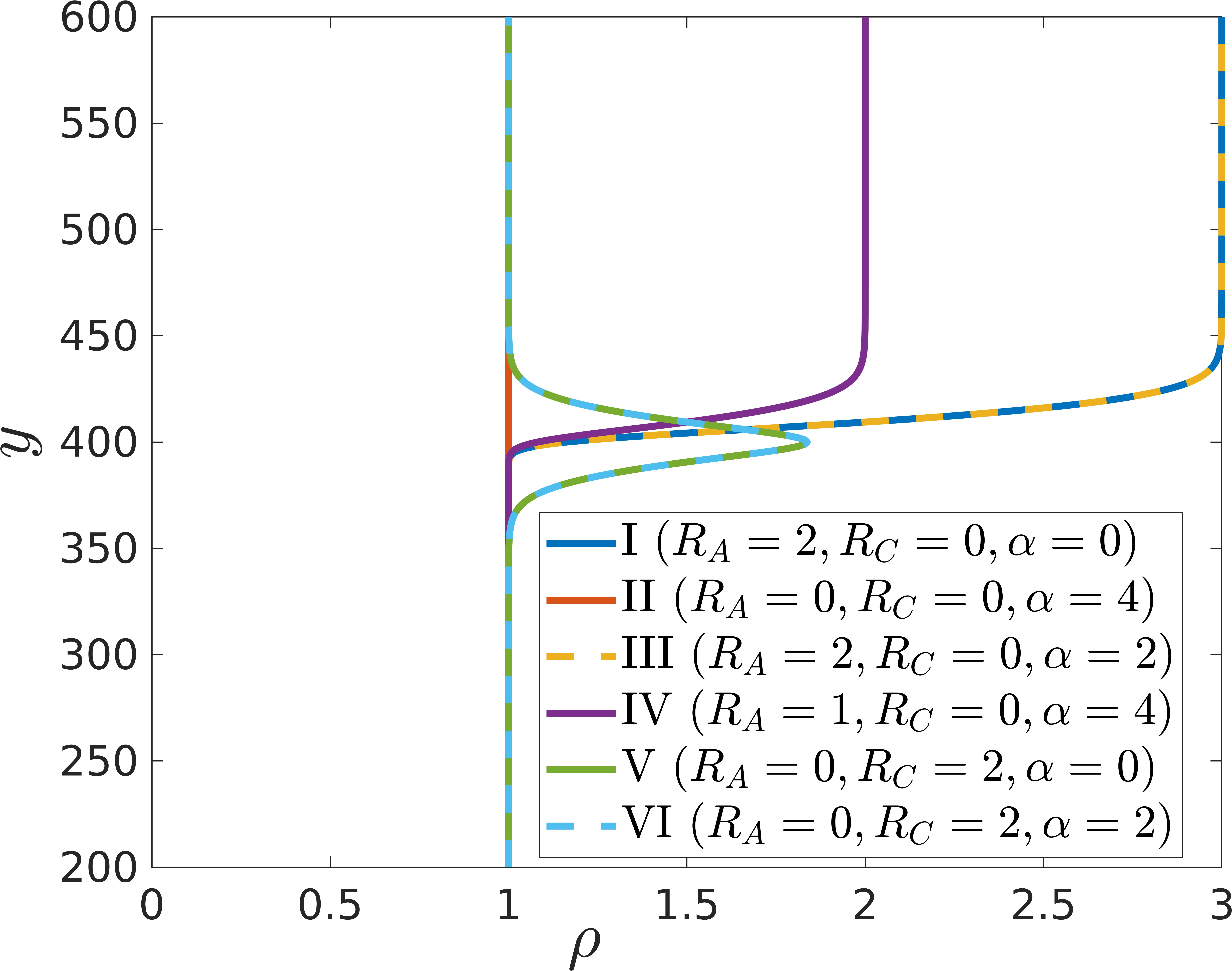}
     \vspace{-0.3 cm} \caption{Typical density profiles \(\rho\)  for each case I to VI as a function of the vertical coordinate \(y\), given by \(\rho = 1 + R_{A} \bar{a} + R_{B} \bar{b} + R_{C} \bar{c}\), where \(\bar{a}(y)\), \(\bar{b}(y)\), and \(\bar{c}(y)\) are reaction-diffusion averaged concentrations at \(t = 100\), shown in Fig. \ref{fig:galfi}(a).}
    \label{fig:rho}
\end{figure}
In non-reactive fluids with two solutes of differing diffusivities, the vertical density profile and resulting instability patterns are influenced by the density ratios $R_A$ and $R_B$. These patterns can be classified into Rayleigh-Taylor (RT), Double-Diffusive Rayleigh-Taylor (DD-RT), and Diffusive Layer Convection Rayleigh-Taylor (DLC-RT) regimes, as discussed in \cite{Gopalakrishnan_2021} and references therein. For simplicity, however, we consider all species to have the same diffusivity, focusing on how classical buoyancy-driven fingering or RT instability is influenced or induced by chemical reactions through local density changes or permeability variations. To this end, we set $R_B = 0$ and examine six cases, for which typical density profiles are shown in Fig. \ref{fig:rho}:
\begin{itemize}
    \item \textbf{Case I ($R_A = 2$, $R_C = 0$, $\alpha = 0$):} The reaction does not induce RT instability but modifies the pre-existing buoyancy-driven RT instability ($R_A = 2$) through reaction effects at the interface.

    \item \textbf{Case II ($R_A = 0$, $R_C = 0$, $\alpha = 4$):} Here, there is no density change due to either reaction or initial stratification, resulting in a constant density profile (Fig. \ref{fig:rho}). However, the reaction modifies the permeability ($\alpha = 4$).

    \item \textbf{Case III ($R_A = 2$, $R_C = 0$, $\alpha = 2$):} The classical RT instability is influenced by reaction-induced permeability changes, while the density profile remains the same as in Case I (Fig. \ref{fig:rho}).

    \item \textbf{Case IV ($R_A = 1$, $R_C = 0$, $\alpha = 4$):} The reactant-induced density is smaller in magnitude compared to Cases I and III. However, the permeability change due to the reaction is more pronounced than in Case III.

    \item \textbf{Case V ($R_A = 0$, $R_C = 2$, $\alpha = 0$):} In this case, the reaction induces RT instability solely by altering the product density ($R_C = 2$), with no permeability change.

    \item \textbf{Case VI ($R_A = 0$, $R_C = 2$, $\alpha = 2$):} Here, the reaction-induced RT instability is affected by permeability changes ($\alpha = 2$).
\end{itemize}

These cases highlight how chemical reactions can induce or modify pattern formation by altering density or permeability. This understanding is central to analyzing reaction-driven instabilities in the current problem.
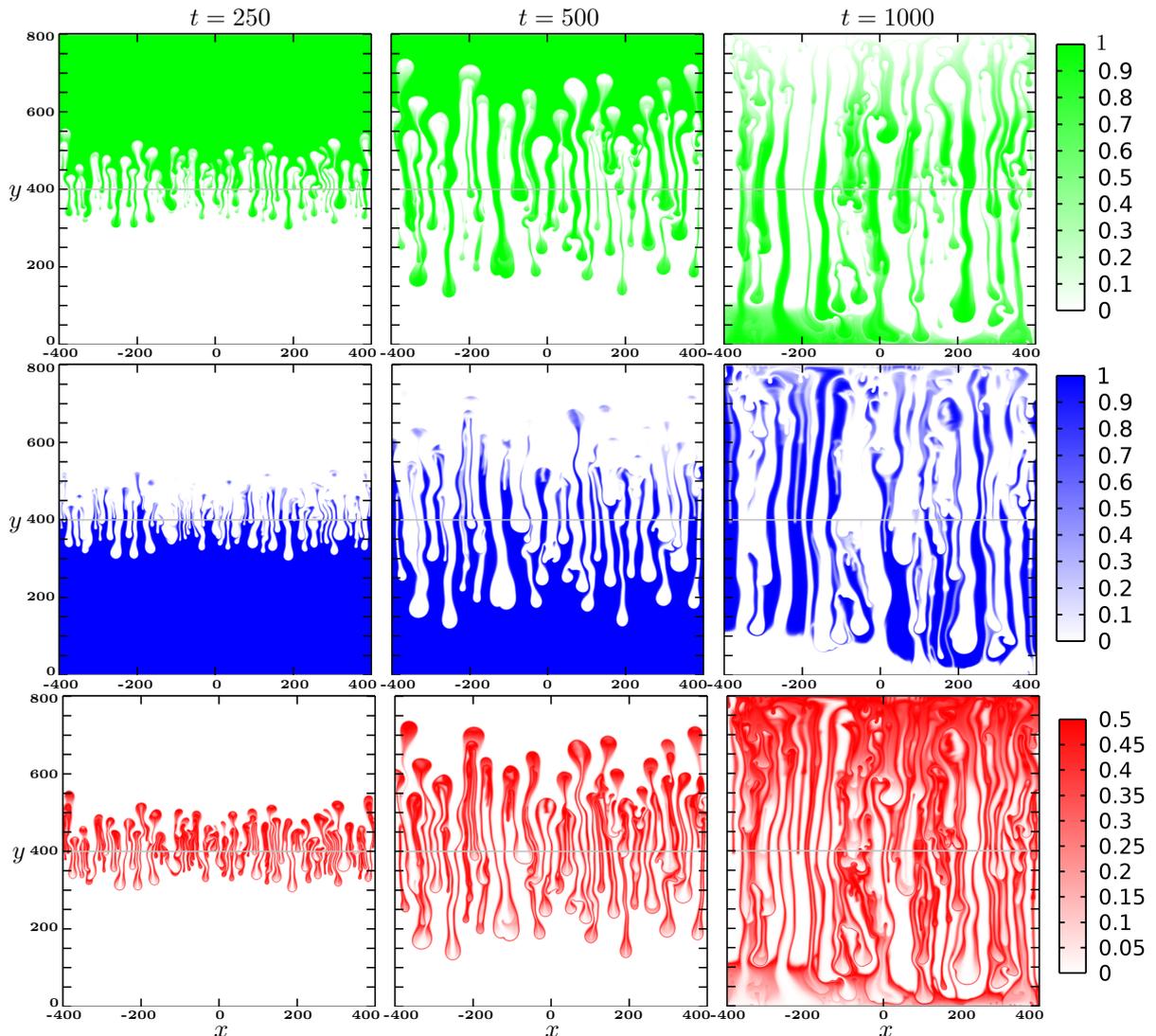
\begin{figure}[h!]
    \centering
    \hspace{0.46 cm}\includegraphics[trim=31 15 0 2,clip,width=0.25\textwidth]{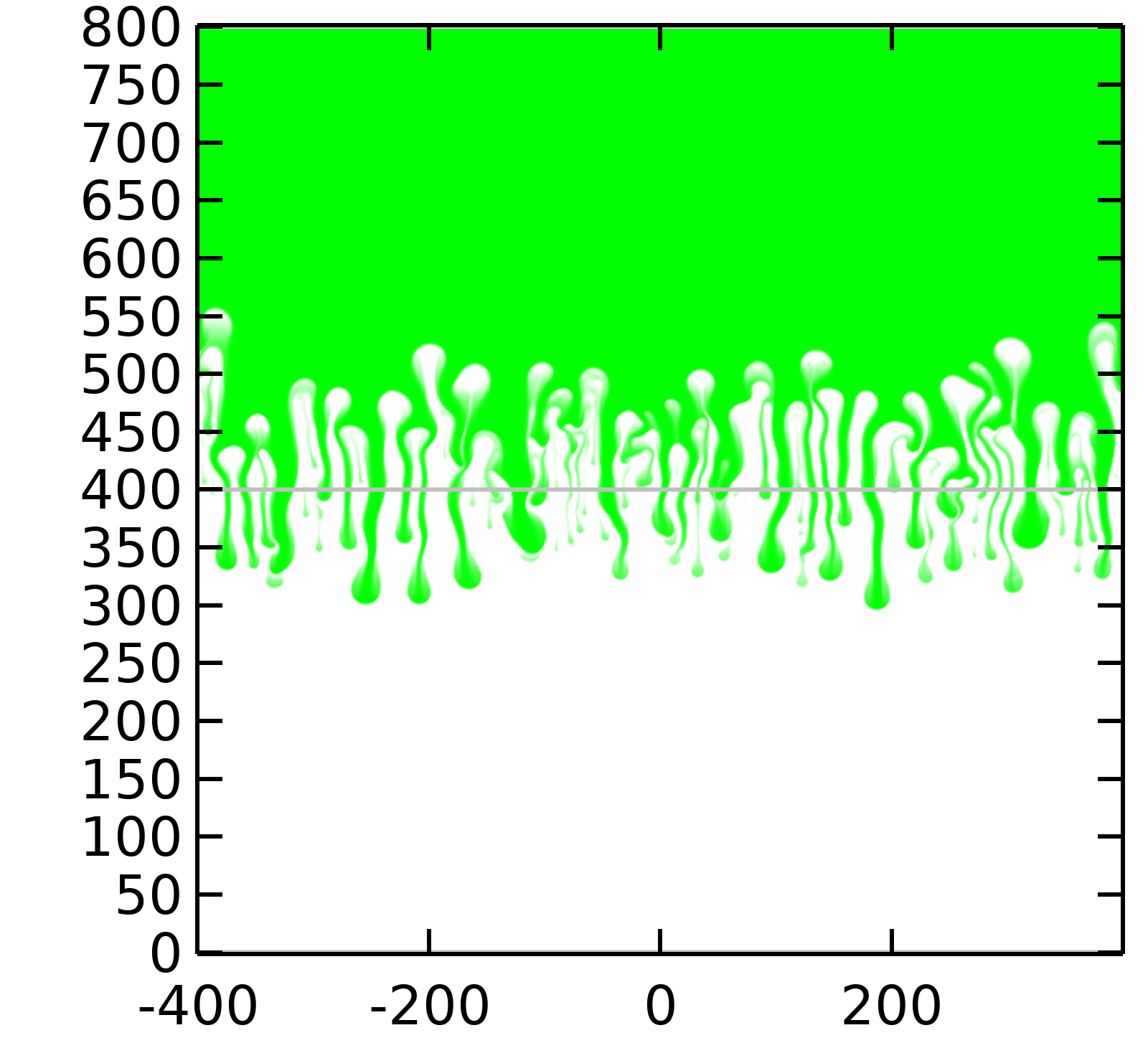}
    \includegraphics[trim=31 15 0 2,clip,width=0.25\textwidth]{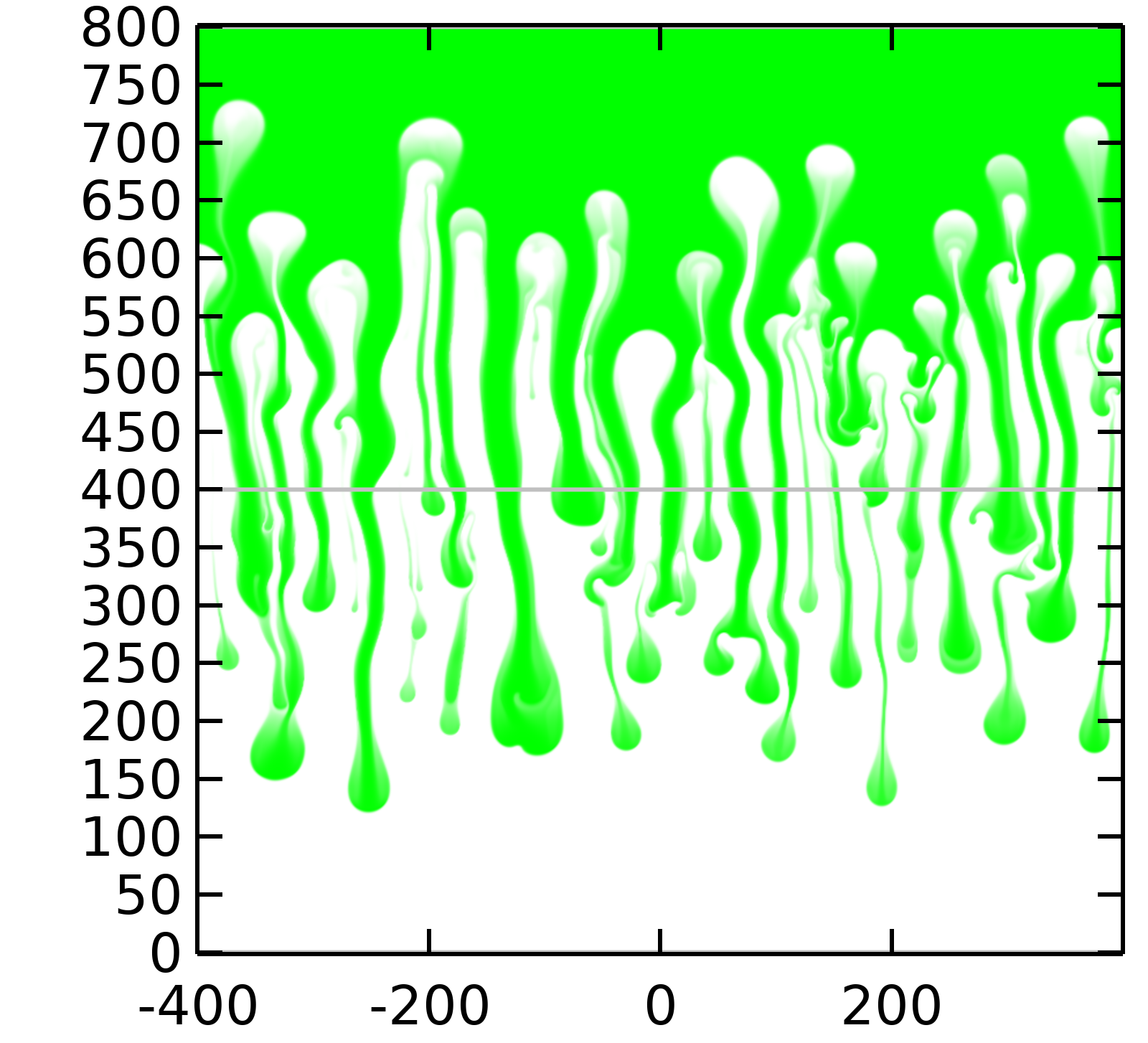}
    \includegraphics[trim=31 15 0 2,clip,width=0.25\textwidth]{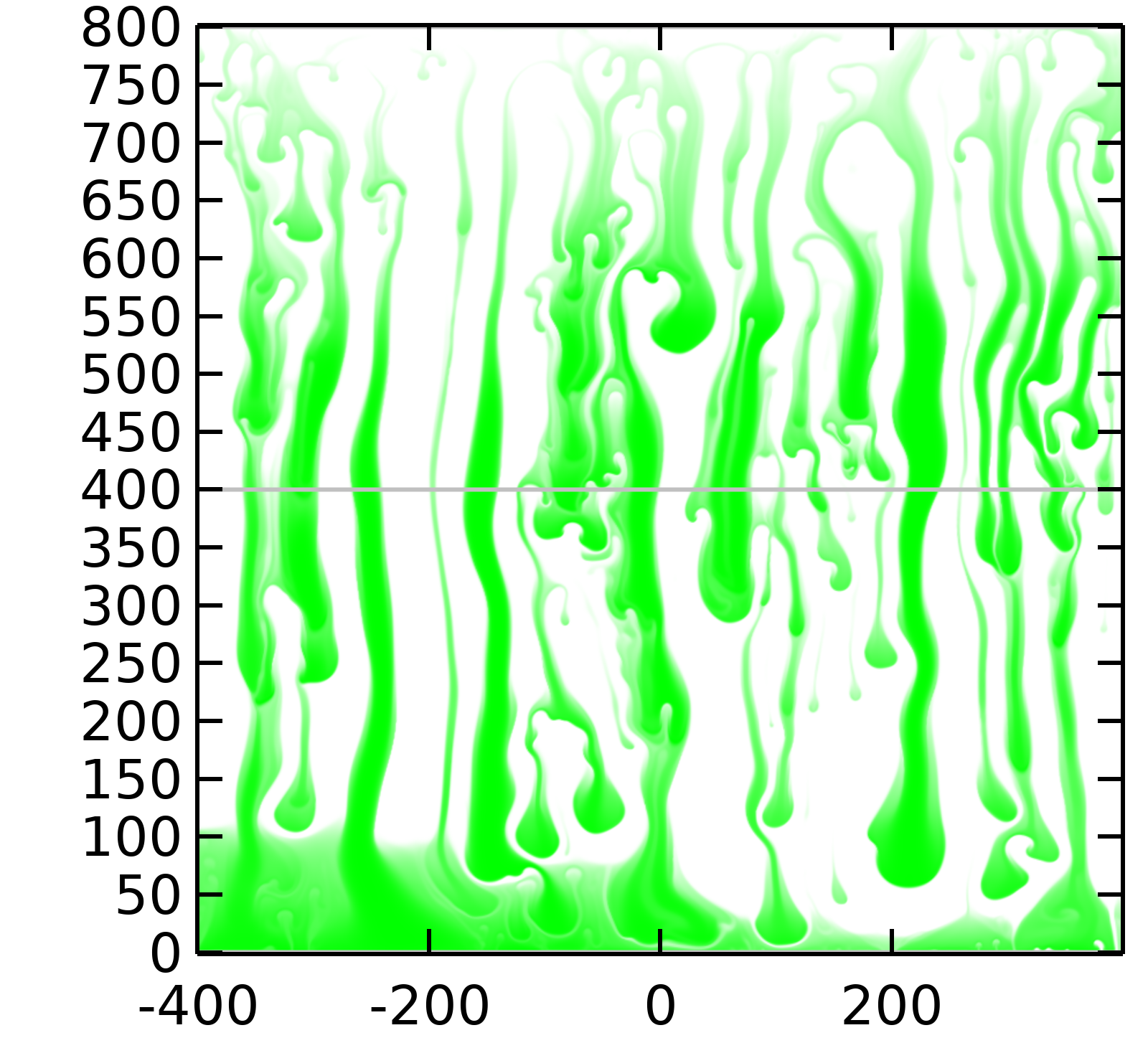}
    \begin{picture}(0,0)
    \put(-410,60){\makebox(0,0)[]{$y$}}
    \put(-325,130){\makebox(0,0)[]{$t=250$}}
    \put(-195,130){\makebox(0,0)[]{$t=500$}}
    \put(-65,130){\makebox(0,0)[]{$t=1000$}}
    \put(-396,1.5){\makebox(0,0)[]{\textbf{\tiny 0}}}
    \put(-400,31.5){\makebox(0,0)[]{\textbf{\tiny 200}}}
    \put(-400,62){\makebox(0,0)[]{\textbf{\tiny 400}}}
    \put(-400,92.5){\makebox(0,0)[]{\textbf{\tiny 600}}}
    \put(-400,122.5){\makebox(0,0)[]{\textbf{\tiny 800}}}
    \put(-333,-3.5){\makebox(0,0)[]{\textbf{\tiny-400~~~~~~-200~~~~~~~~~0~~~~~~~~~200~~~~~~400}}}
    \put(-201,-3.5){\makebox(0,0)[]{\textbf{\tiny-400~~~~~~-200~~~~~~~~0~~~~~~~~~200~~~~~~400}}}
    \put(-70.5,-3.5){\makebox(0,0)[]{\textbf{\tiny-400~~~~~~-200~~~~~~~~~0~~~~~~~~~200~~~~~~400}}}
    \put(-5,10){\includegraphics[trim=120 26 30 35, clip, width=0.08\textwidth]{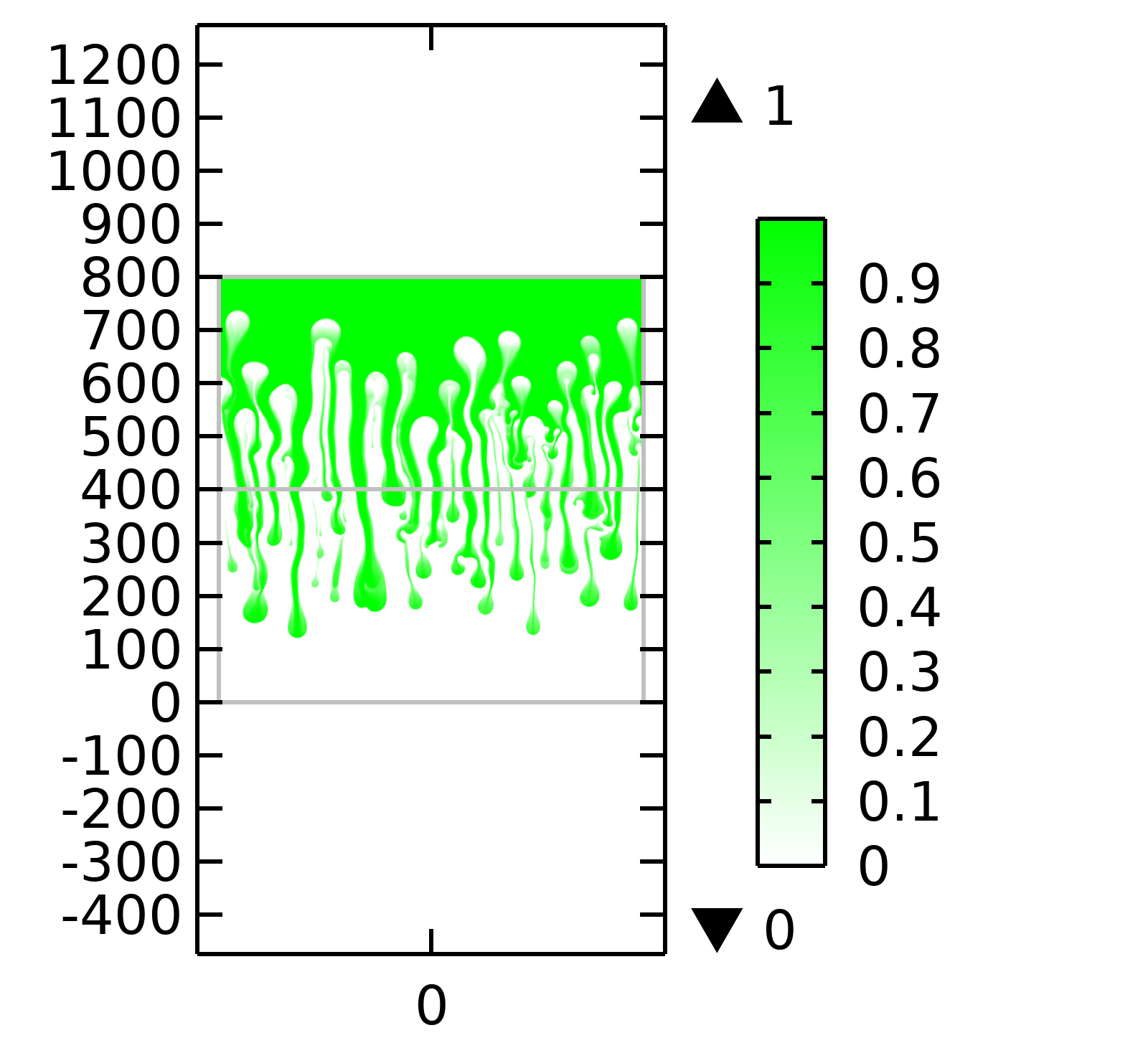}}
    \put(19.5,120){\makebox(0,0)[]{\small 1}}
     \end{picture} \vspace{0.2 cm}

    \hspace{0.35 cm}\includegraphics[trim=31 15 0 2,clip,width=0.25\textwidth]{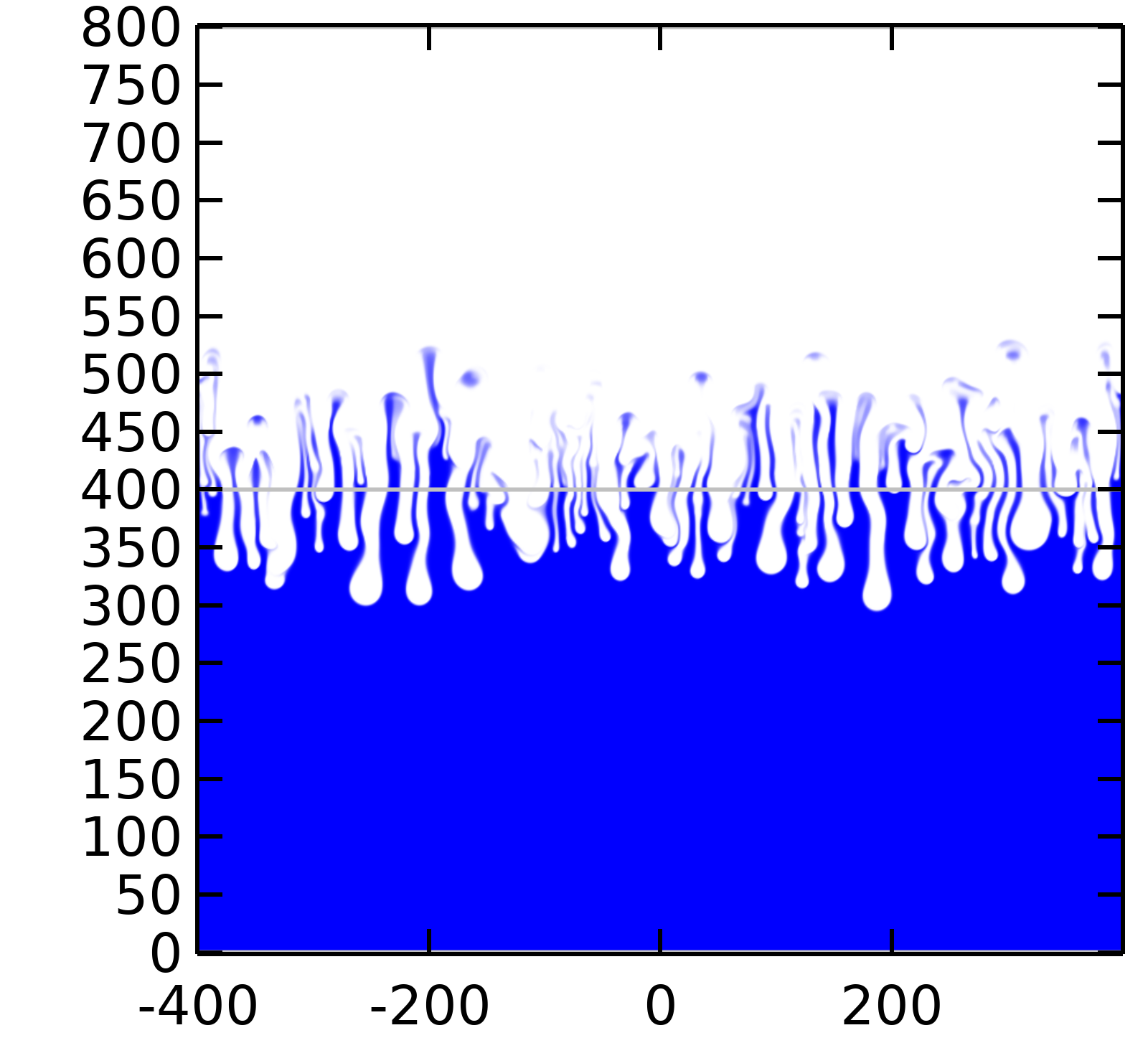}
    \includegraphics[trim=31 15 0 2,clip,width=0.25\textwidth]{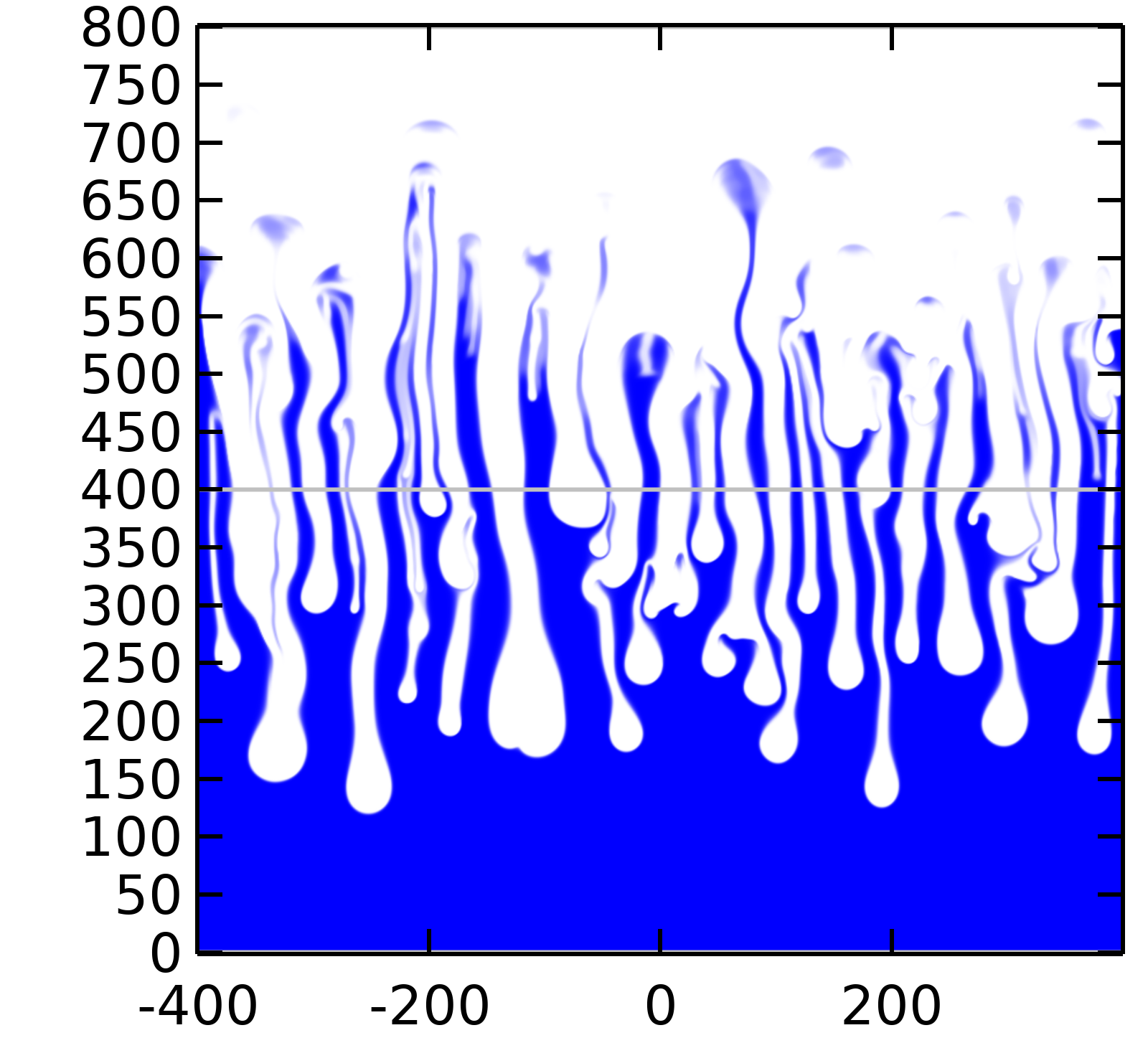}
    \includegraphics[trim=31 15 0 2,clip,width=0.25\textwidth]{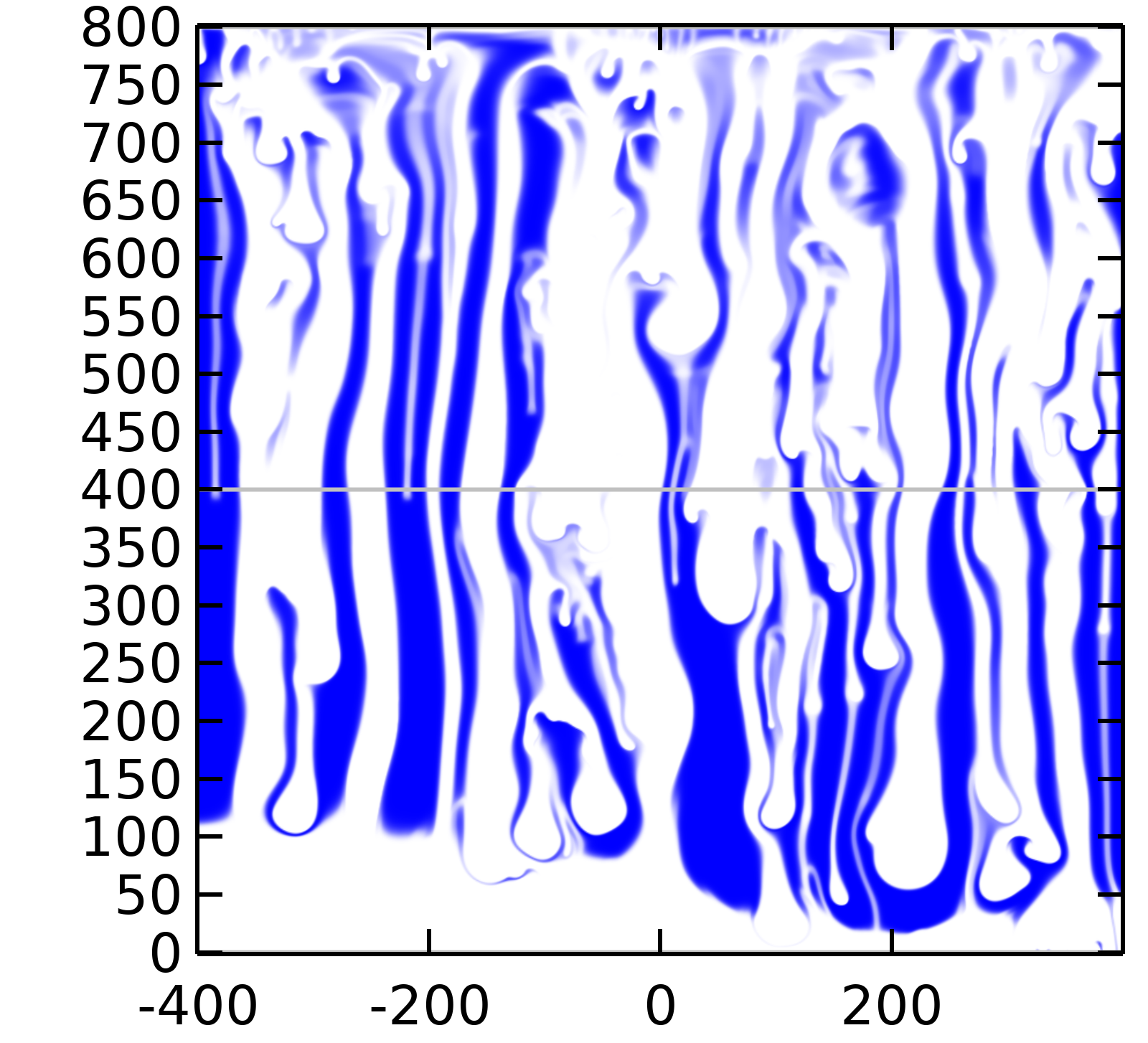}
    \begin{picture}(0,0)
    \put(-410,60){\makebox(0,0)[]{$y$}}
    \put(-396,1.5){\makebox(0,0)[]{\textbf{\tiny 0}}}
    \put(-400,31.5){\makebox(0,0)[]{\textbf{\tiny 200}}}
    \put(-400,62){\makebox(0,0)[]{\textbf{\tiny 400}}}
    \put(-400,92.5){\makebox(0,0)[]{\textbf{\tiny 600}}}
    \put(-400,122.5){\makebox(0,0)[]{\textbf{\tiny 800}}}
    \put(-333,-3.5){\makebox(0,0)[]{\textbf{\tiny-400~~~~~~-200~~~~~~~~~0~~~~~~~~~200~~~~~~400}}}
    \put(-201,-3.5){\makebox(0,0)[]{\textbf{\tiny-400~~~~~~-200~~~~~~~~0~~~~~~~~~200~~~~~~400}}}
    \put(-70.5,-3.5){\makebox(0,0)[]{\textbf{\tiny-400~~~~~~-200~~~~~~~~~0~~~~~~~~~200~~~~~~400}}}
    \put(-5,10){\includegraphics[trim=120 26 30 33, clip, width=0.08\textwidth]{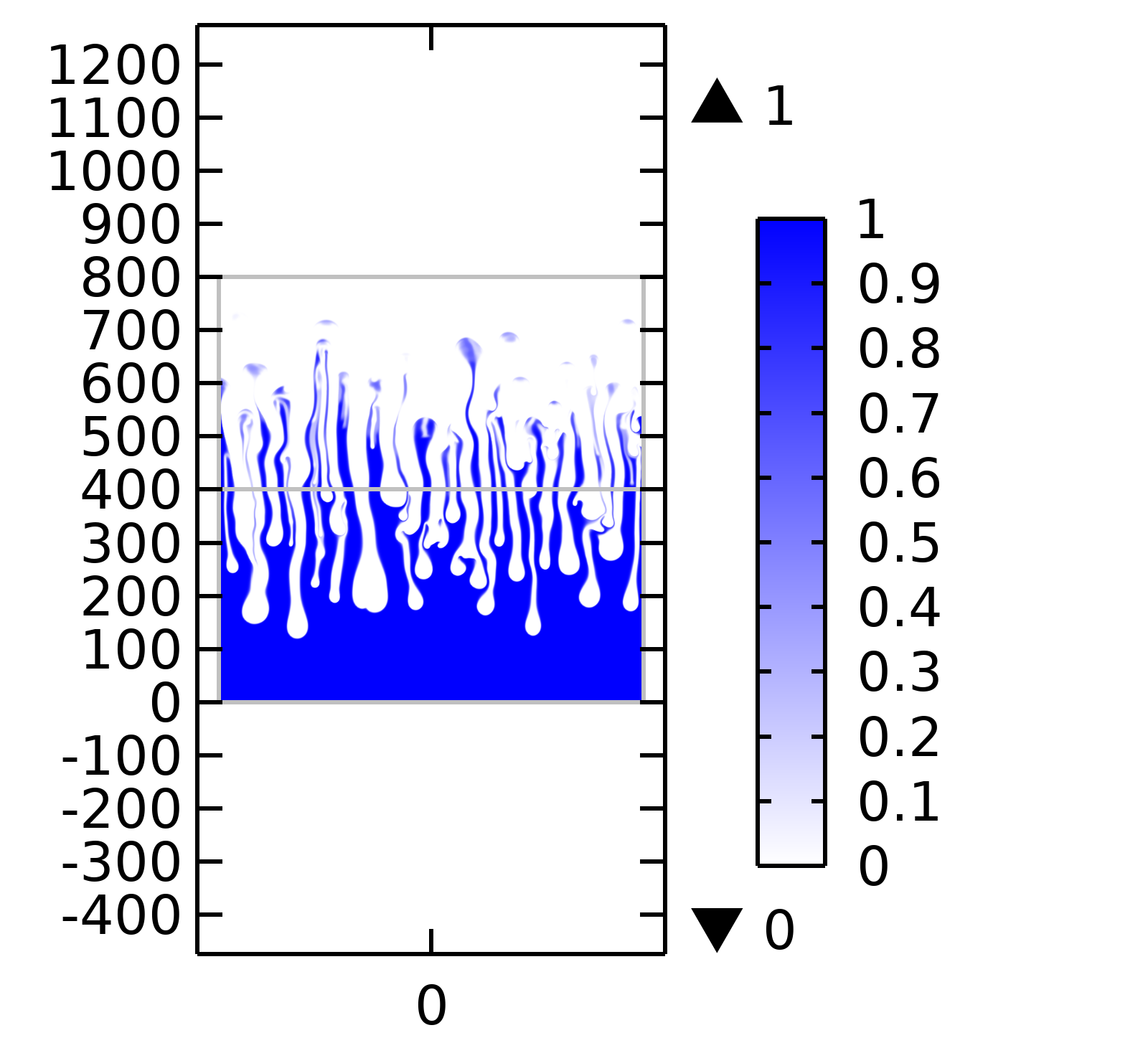}}
    \end{picture}\vspace{0.2 cm}

    \hspace{0.44 cm}\includegraphics[trim=31 15 0 2,clip,width=0.25\textwidth]{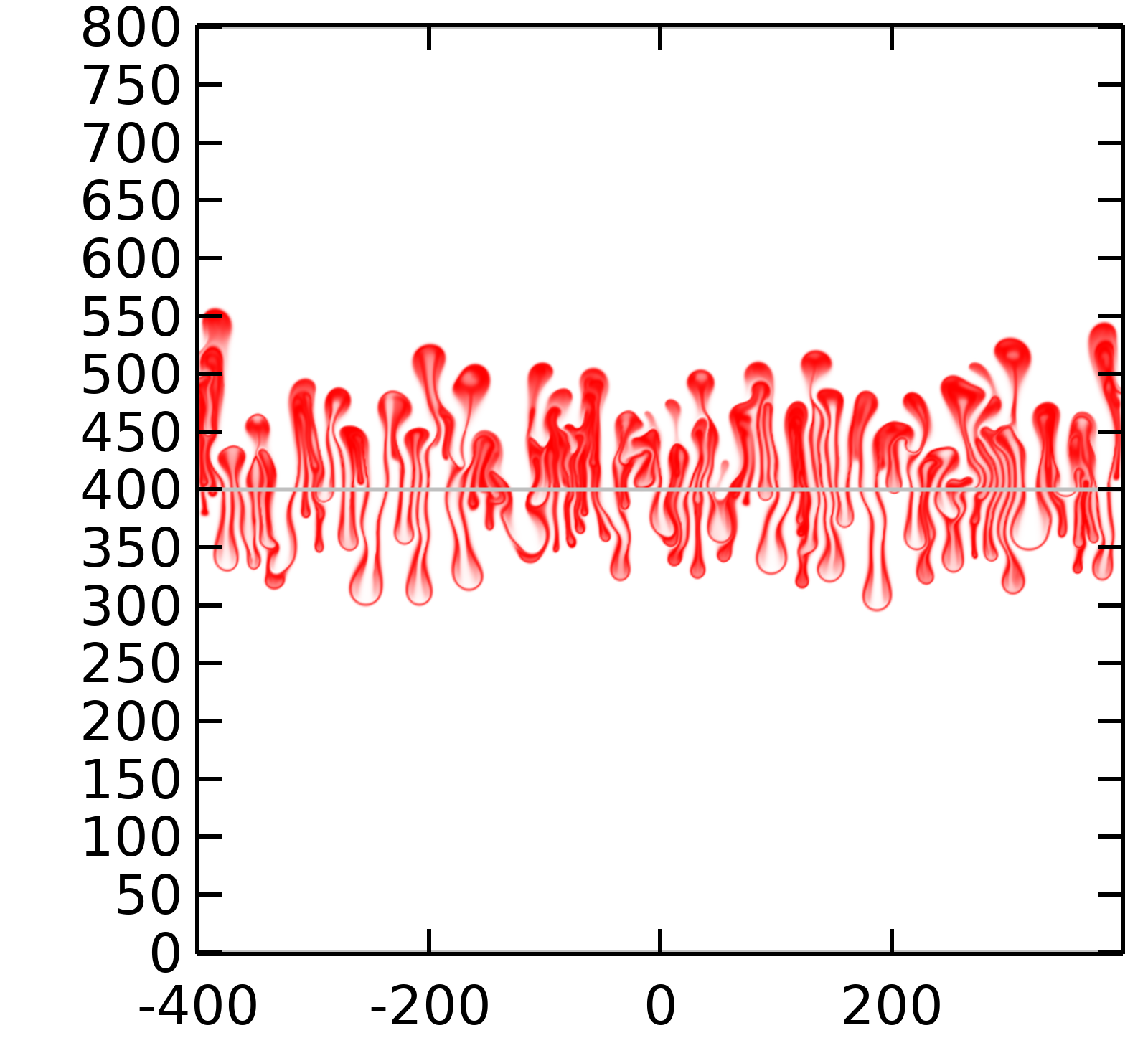}
    \includegraphics[trim=31 15 0 2,clip,width=0.25\textwidth]{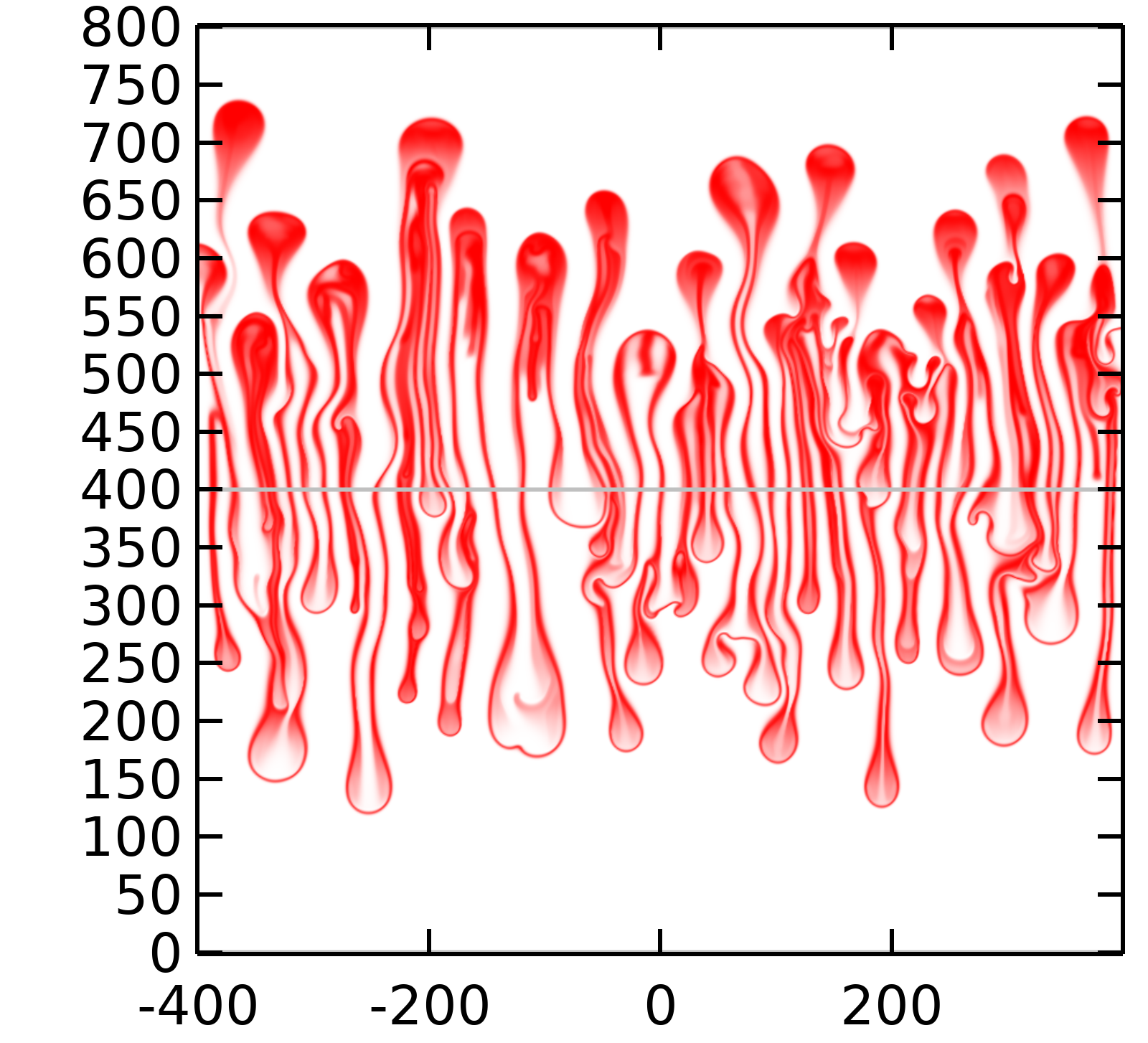}
    \includegraphics[trim=31 15 0 2,clip,width=0.25\textwidth]{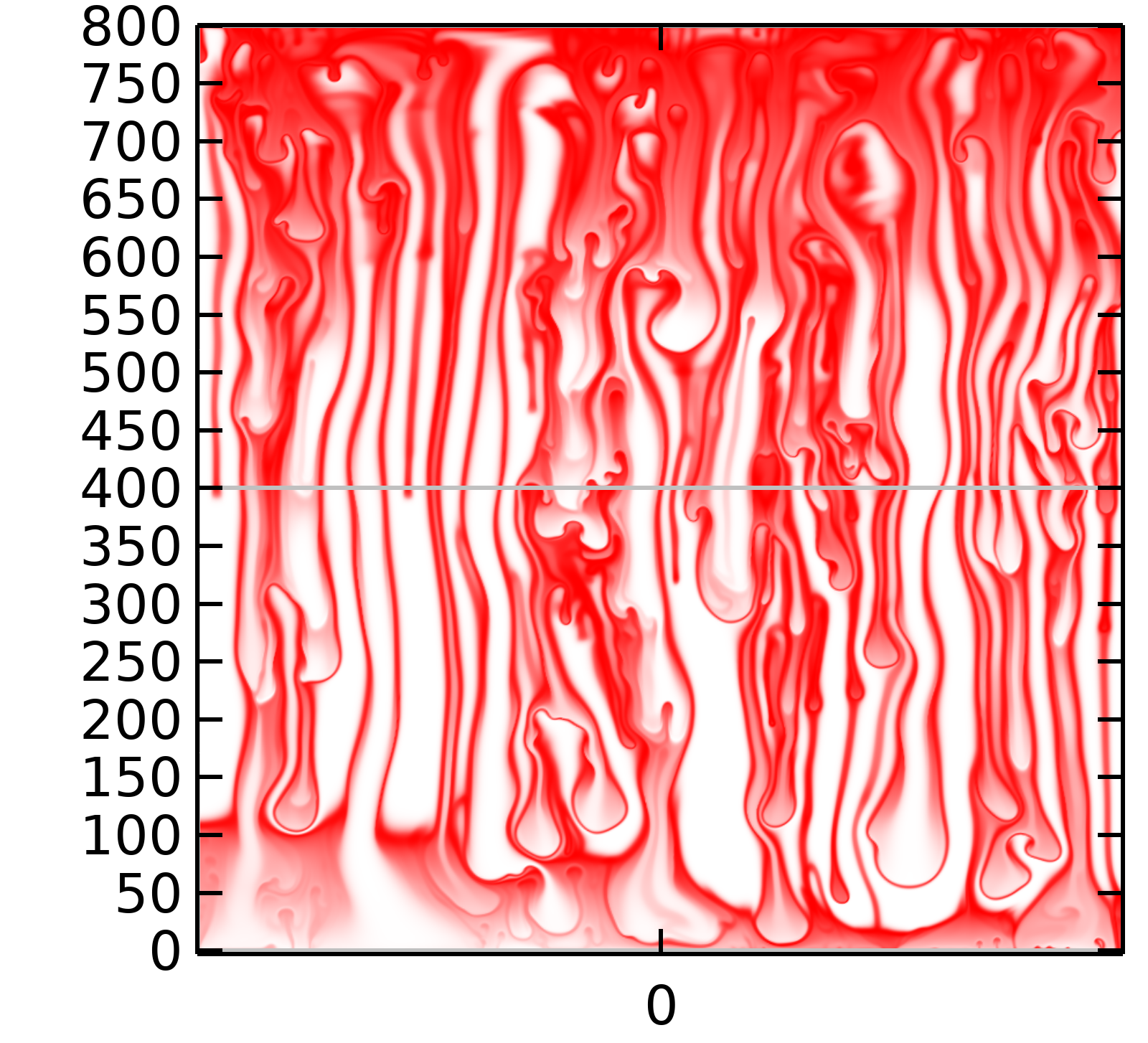}
    \begin{picture}(0,0)
    \put(-410,60){\makebox(0,0)[]{$y$}}
    \put(-330,-10){\makebox(0,0)[]{$x$}}
    \put(-199,-10){\makebox(0,0)[]{$x$}}
    \put(-67,-10){\makebox(0,0)[]{$x$}}
    
    \put(-396,1.5){\makebox(0,0)[]{\textbf{\tiny 0}}}
    \put(-400,31.5){\makebox(0,0)[]{\textbf{\tiny 200}}}
    \put(-400,62){\makebox(0,0)[]{\textbf{\tiny 400}}}
    \put(-400,92.5){\makebox(0,0)[]{\textbf{\tiny 600}}}
    \put(-400,122.5){\makebox(0,0)[]{\textbf{\tiny 800}}}
    \put(-333,-3.5){\makebox(0,0)[]{\textbf{\tiny-400~~~~~~-200~~~~~~~~~0~~~~~~~~~200~~~~~~400}}}
    \put(-201,-3.5){\makebox(0,0)[]{\textbf{\tiny-400~~~~~~-200~~~~~~~~0~~~~~~~~~200~~~~~~400}}}
    \put(-70.5,-3.5){\makebox(0,0)[]{\textbf{\tiny-400~~~~~~-200~~~~~~~~~0~~~~~~~~~200~~~~~~400}}}
    \put(-5,10){\includegraphics[trim=120 26 28 33, clip, width=0.08\textwidth]{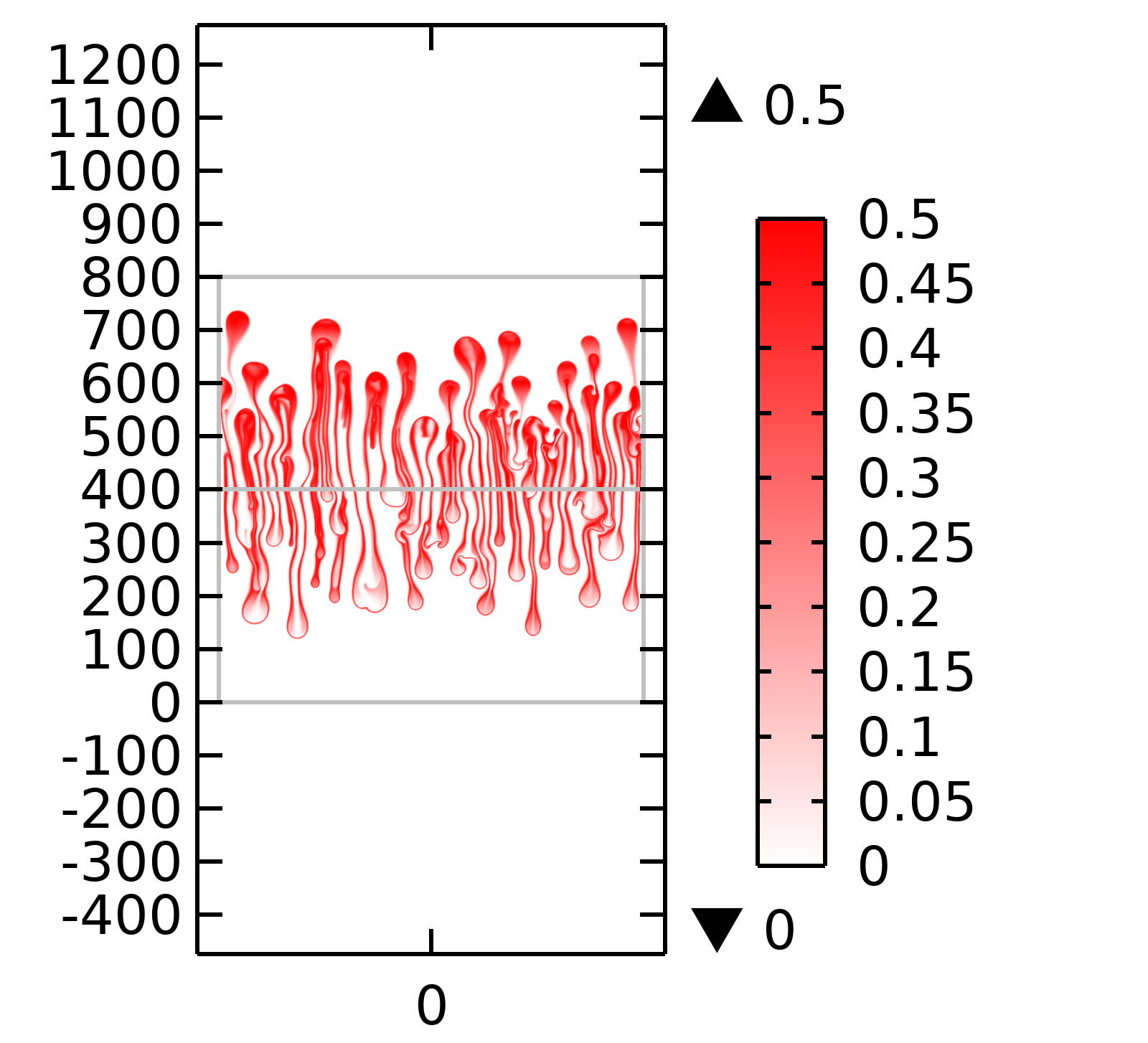}}
    \end{picture}
     \caption{Spatio-temporal evolution of concentration profiles for case-I, with \( R_A = 2 \), \( R_C = 0 \), and \( \alpha = 0 \). The top panel shows concentration \(a\), the middle panel shows \(b\), and the bottom panel shows \(c\). In this configuration, only the more dense fluid, reactant \(A\) (the top fluid), is present, leading to Rayleigh-Taylor fingering as it falls. The reaction occurs at the interface between reactants \(A\) and \(B\).}
    \label{fig:abc}
\end{figure}
\begin{figure}[h!]
    \centering
    \includegraphics[trim=31 14 0 0,clip,width=0.242\textwidth]{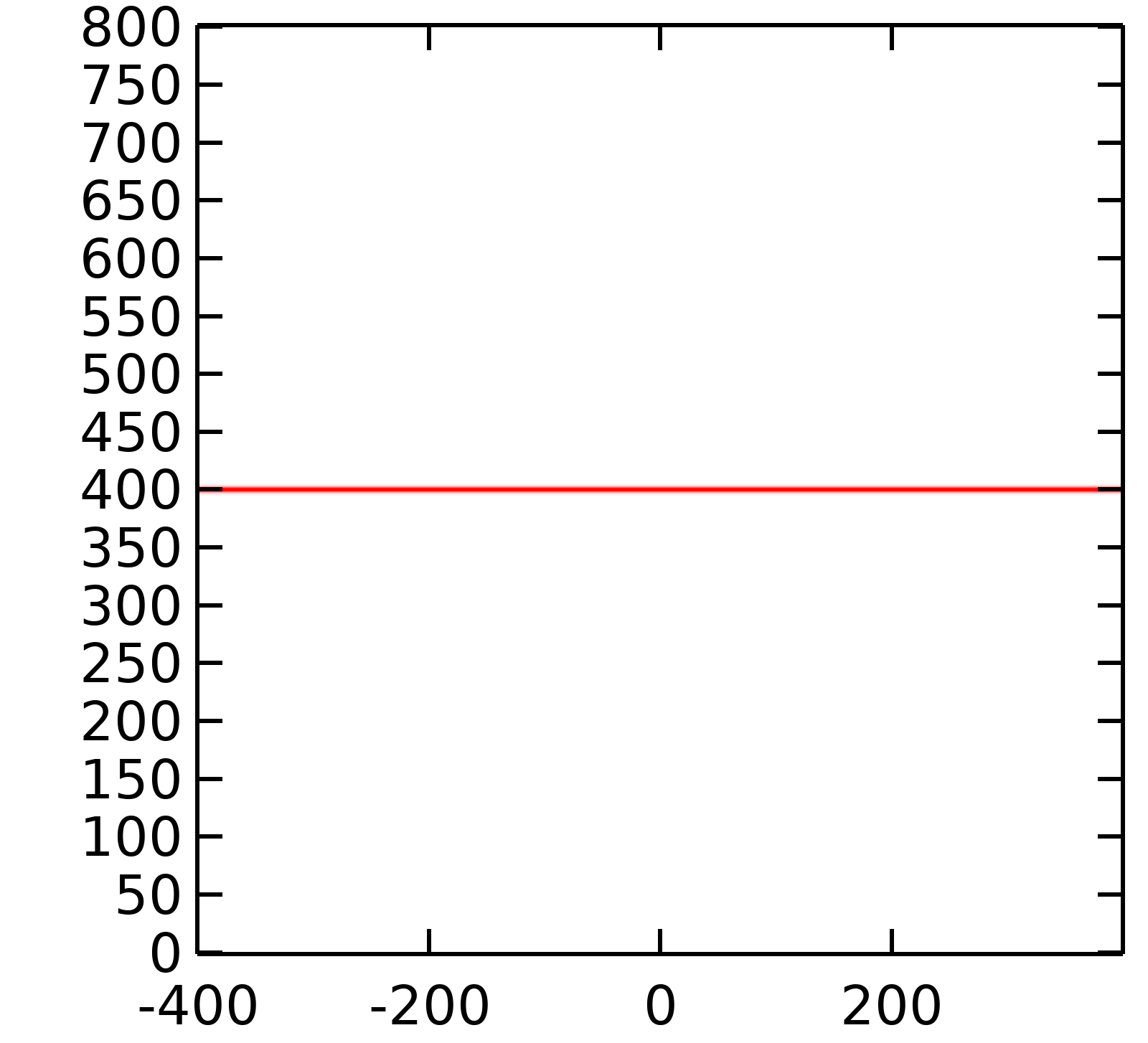}
    \includegraphics[trim=31 14 0 0,clip,width=0.242\textwidth]{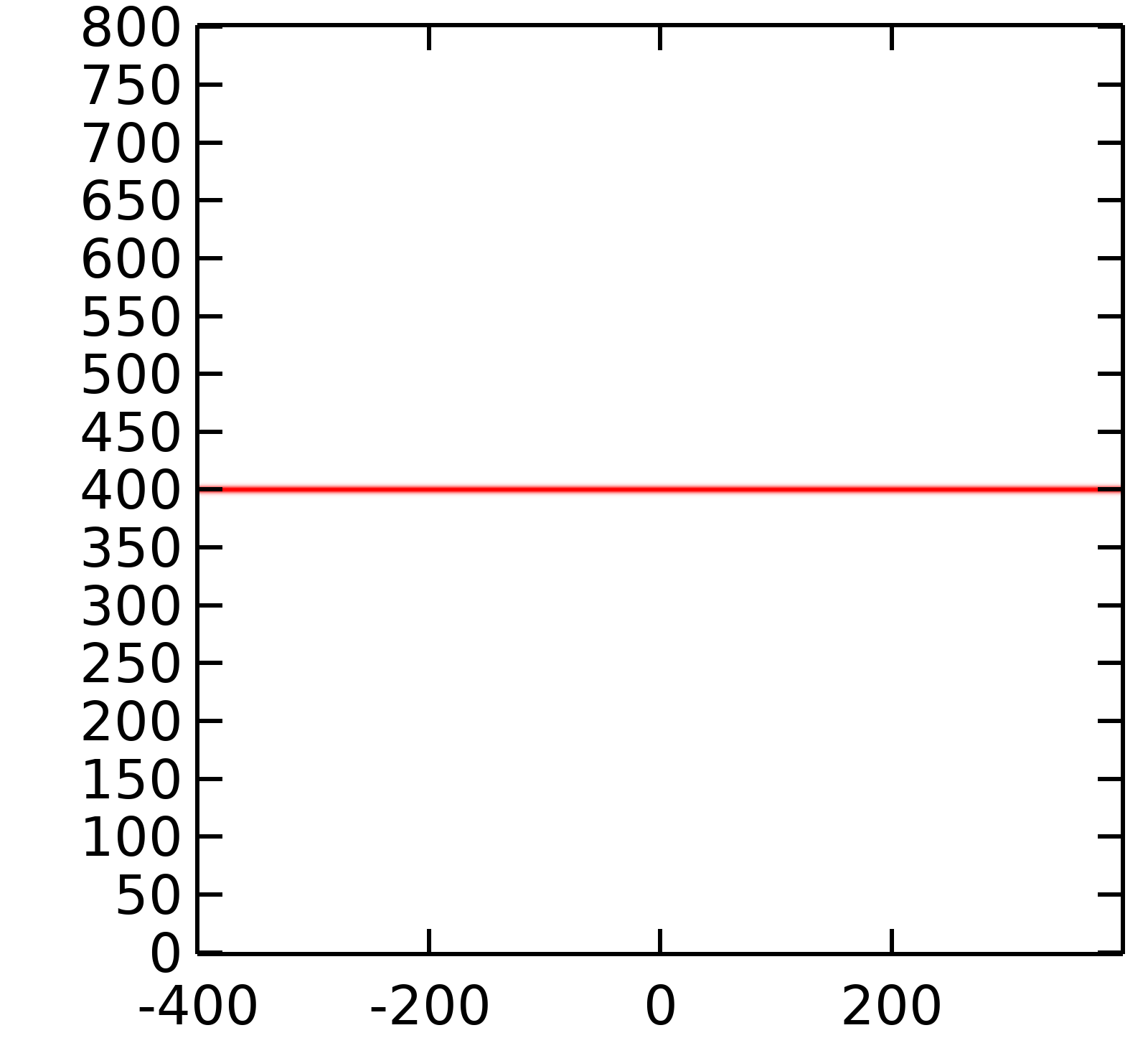}
    \includegraphics[trim=31 14 0 0,clip,width=0.242\textwidth]{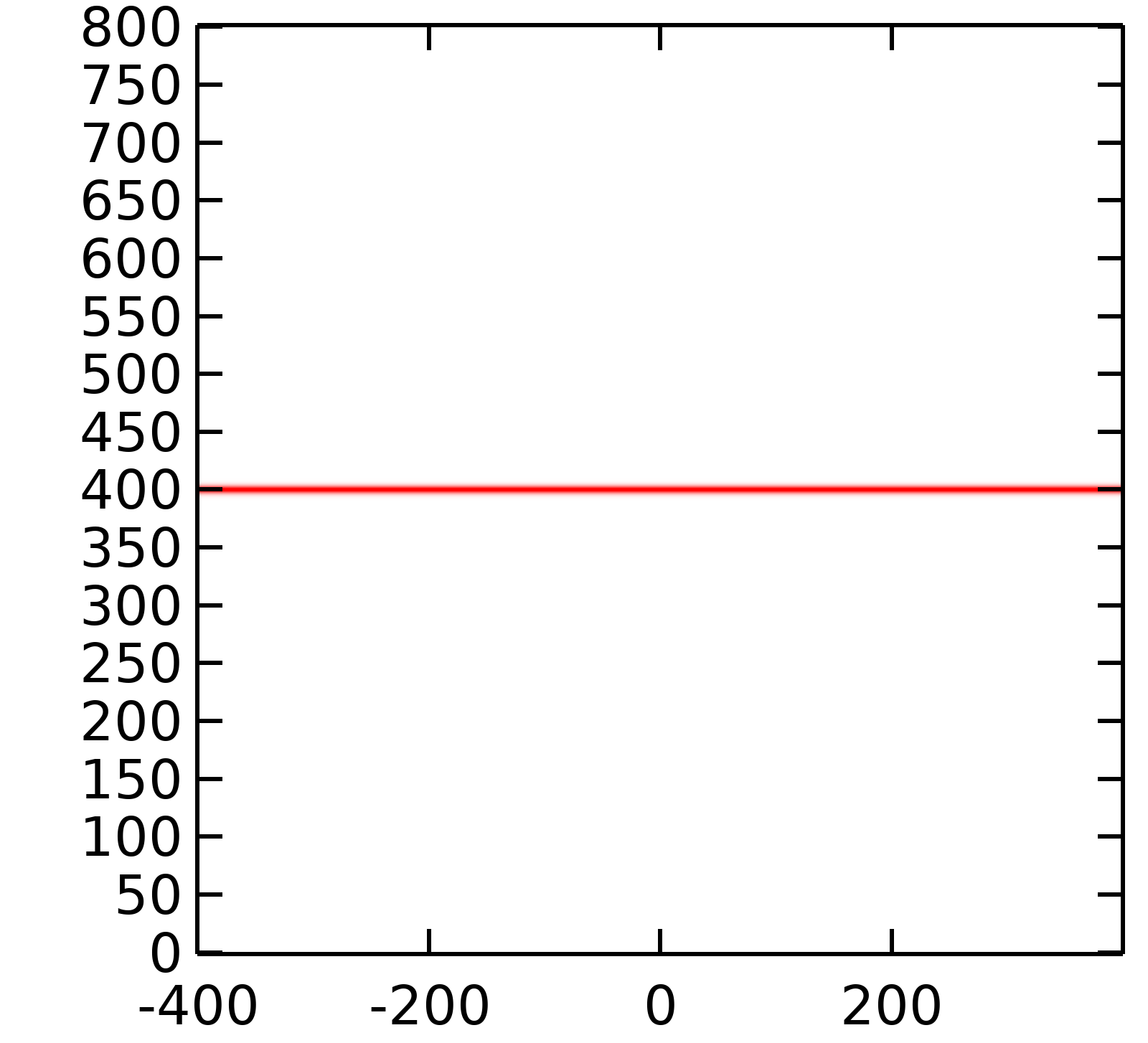}
    \includegraphics[trim=31 14 0 0,clip,width=0.242\textwidth]{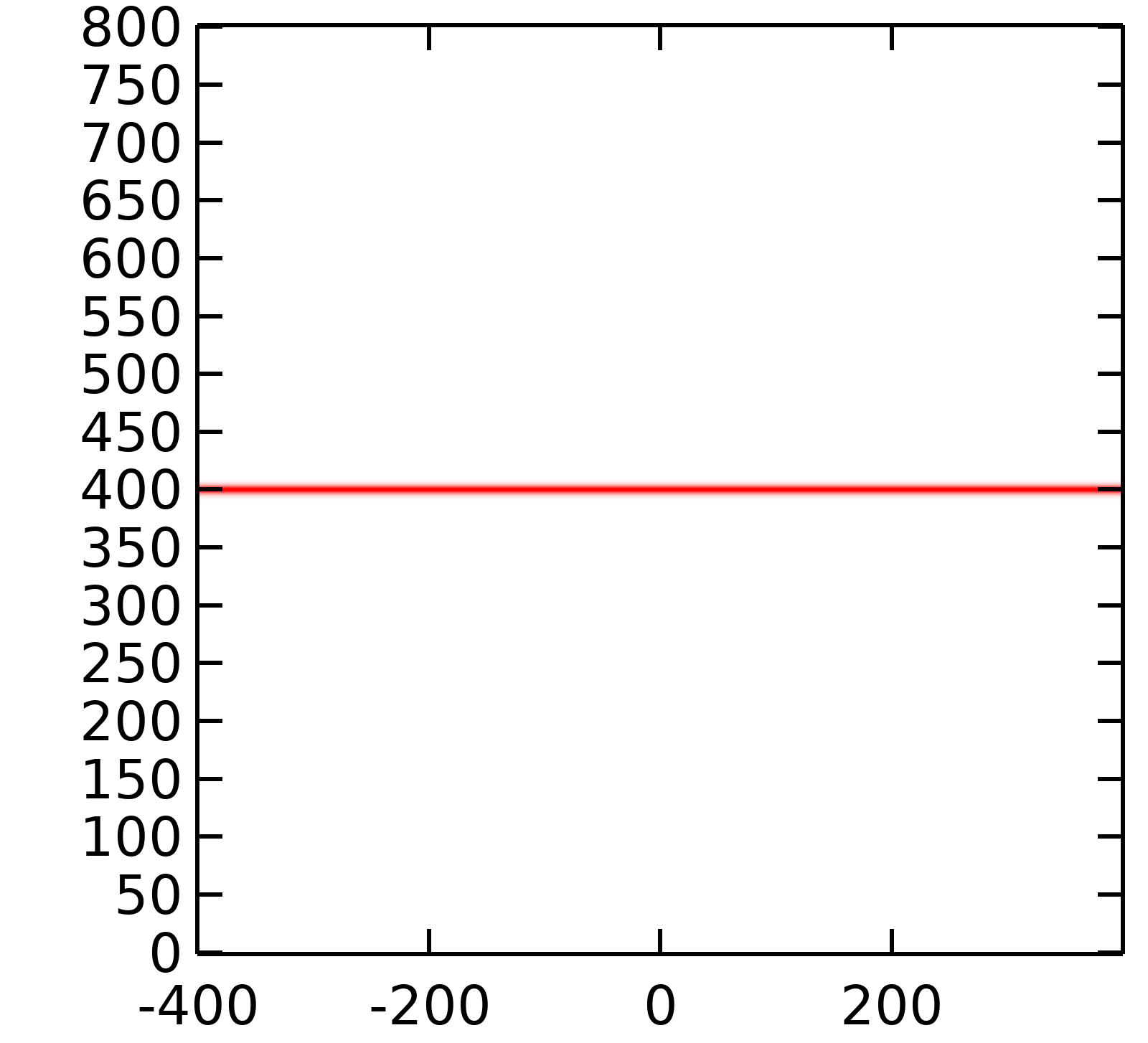}
    \begin{picture}(0,0)
    \put(-450,20){\makebox(0,0)[]{case-II}}
    \put(-450,10){\makebox(0,0)[]{\scriptsize ($R_{A}=0,R_{C}=0,\alpha=4$)}}
   \end{picture}
     
    \includegraphics[trim=31 14 0 0,clip,width=0.242\textwidth]{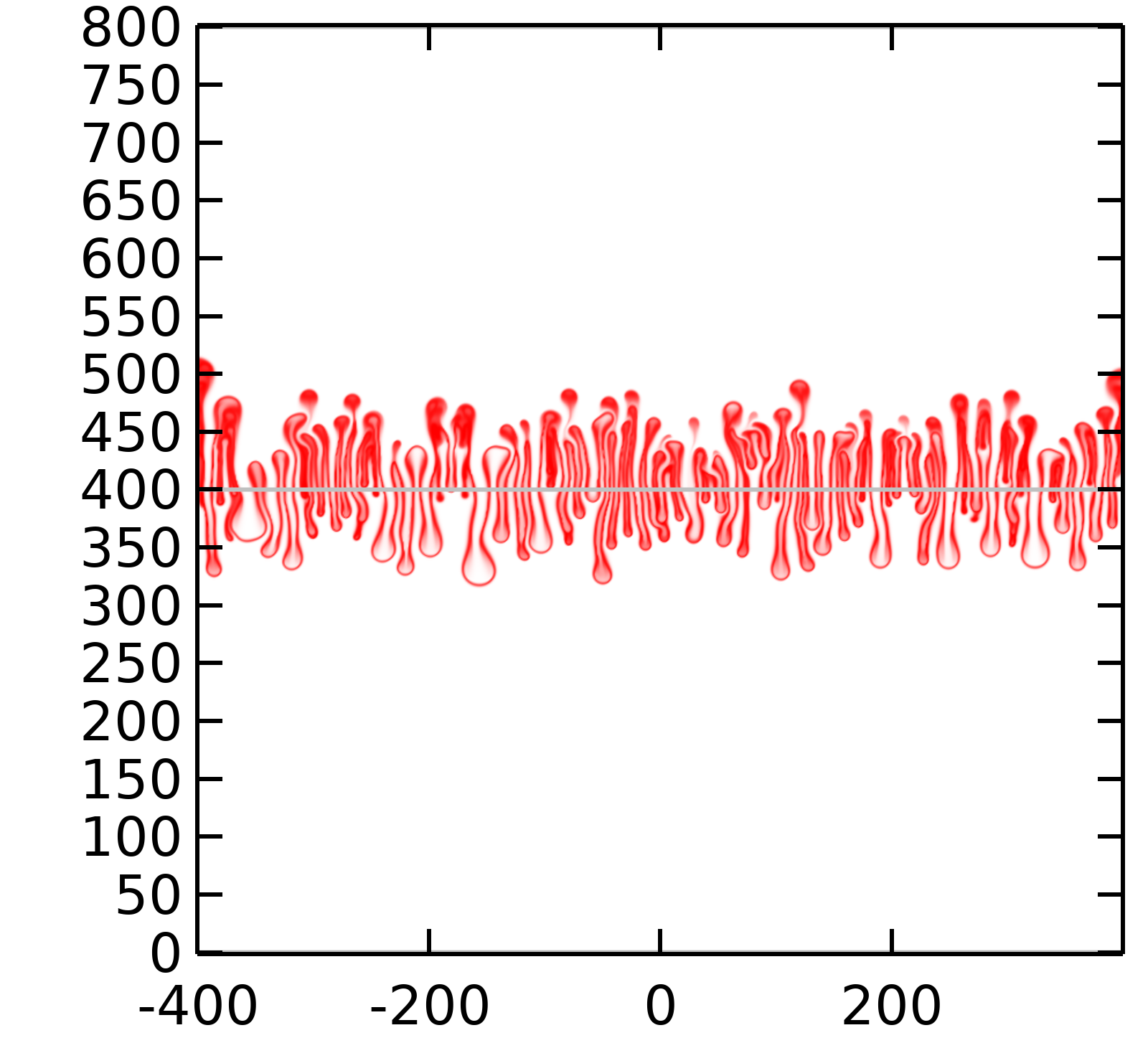}
    \includegraphics[trim=31 14 0 0,clip,width=0.242\textwidth]{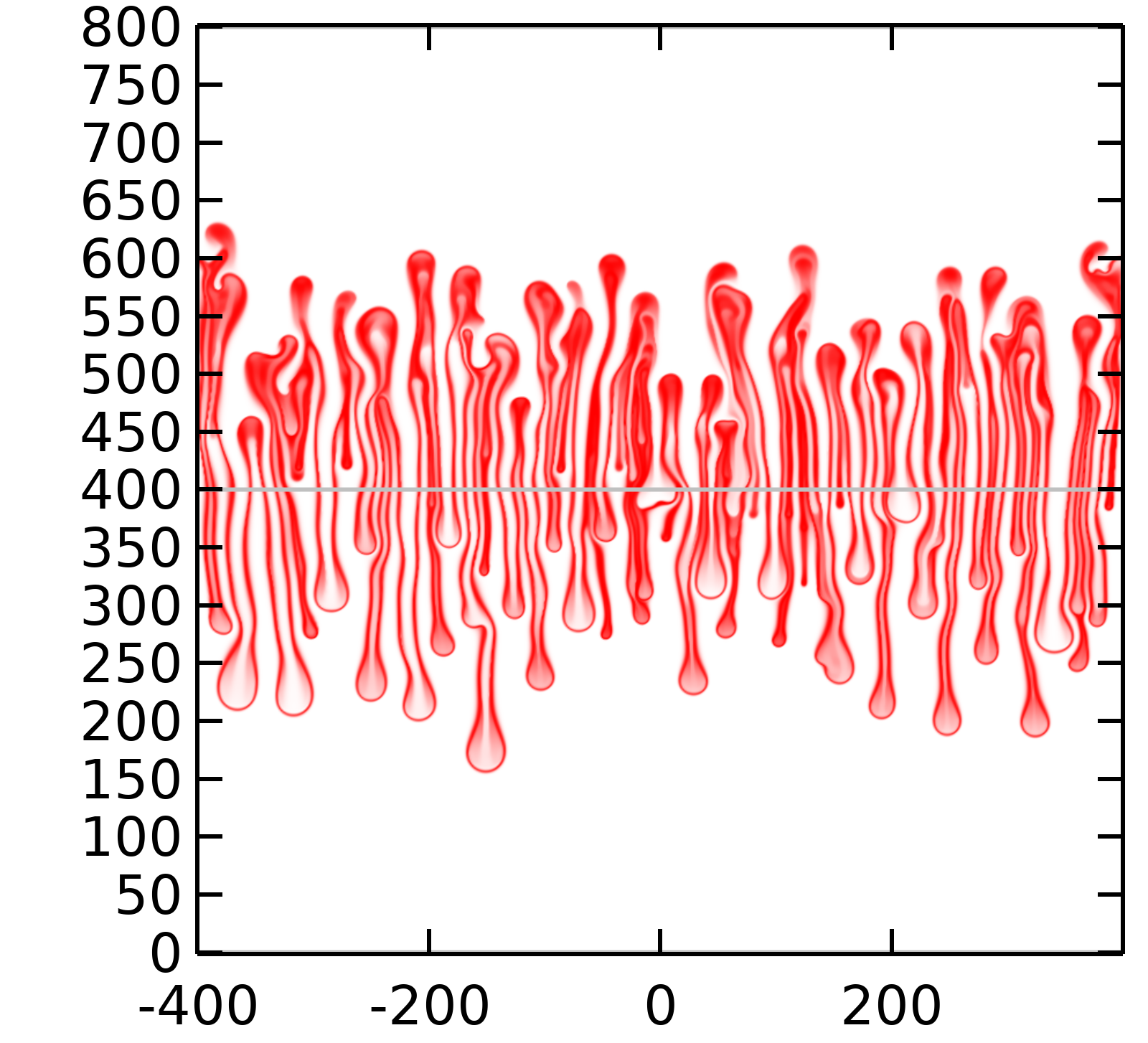}
    \includegraphics[trim=31 14 0 0,clip,width=0.242\textwidth]{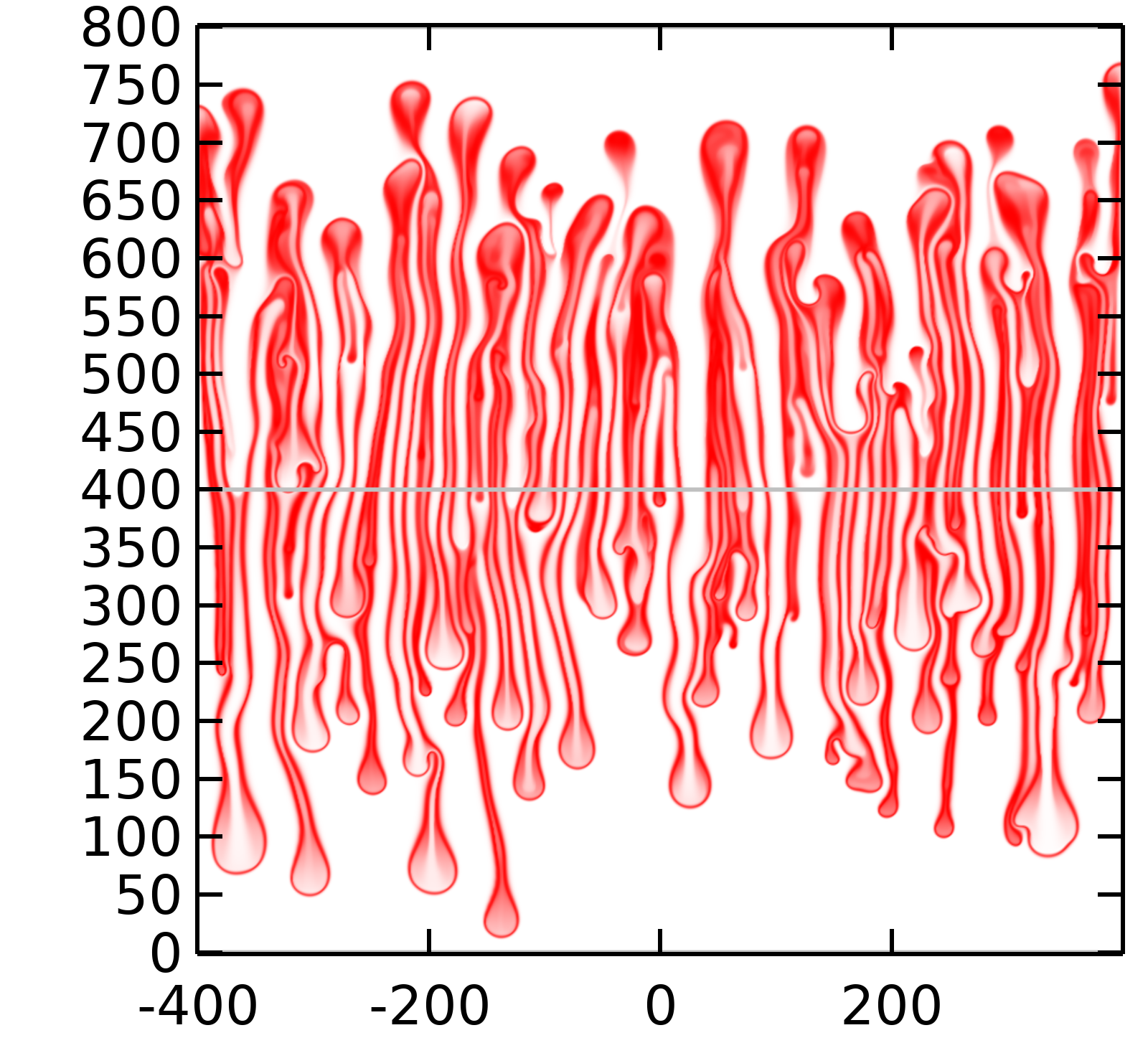}
    \includegraphics[trim=31 14 0 0,clip,width=0.242\textwidth]{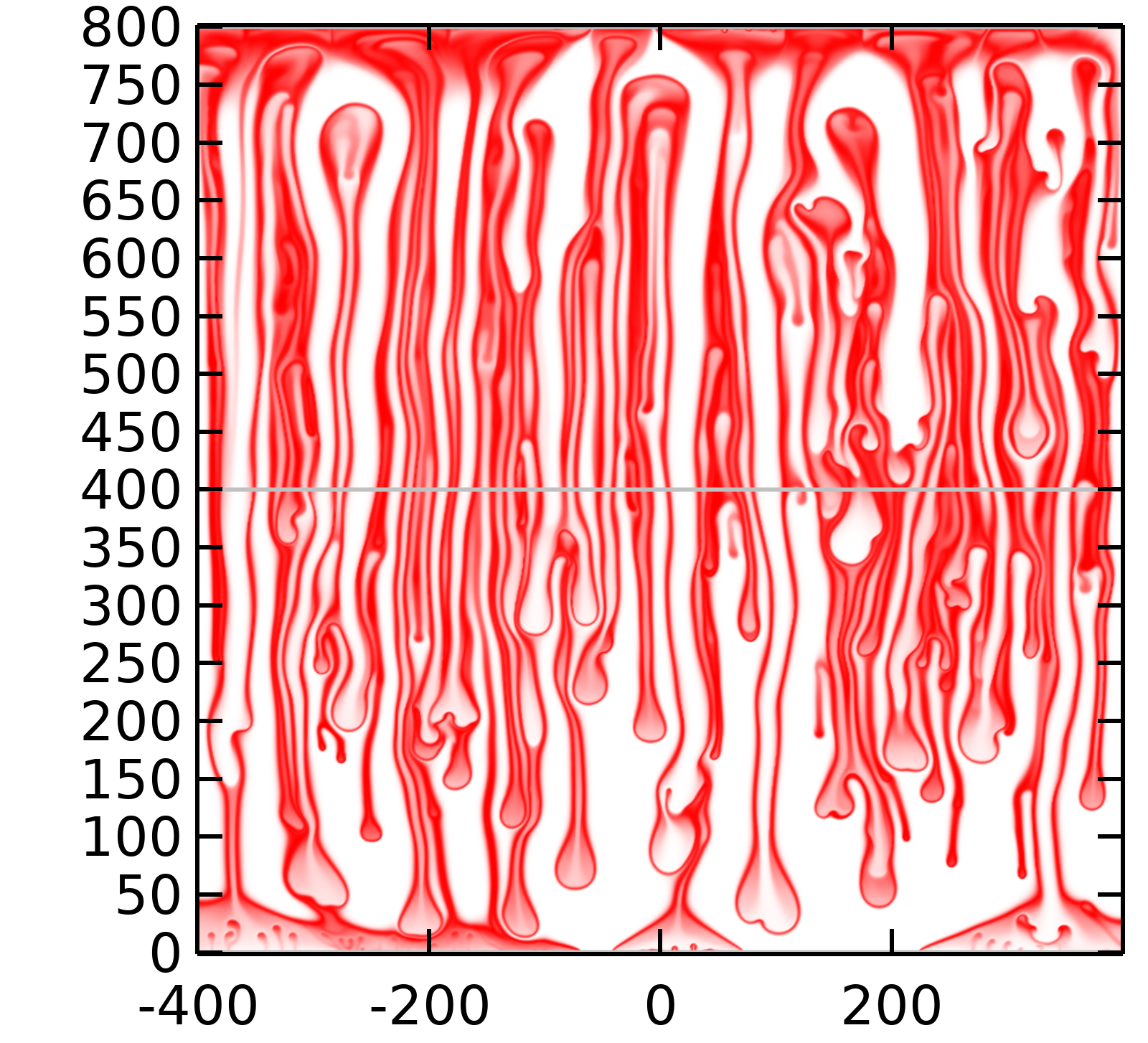}
    \begin{picture}(0,0)
    \put(-450,20){\makebox(0,0)[]{case-III}}
    \put(-450,10){\makebox(0,0)[]{\scriptsize ($R_{A}=2,R_{C}=0,\alpha=2$)}}
   \end{picture}

    \includegraphics[trim=31 14 0 0,clip,width=0.242\textwidth]{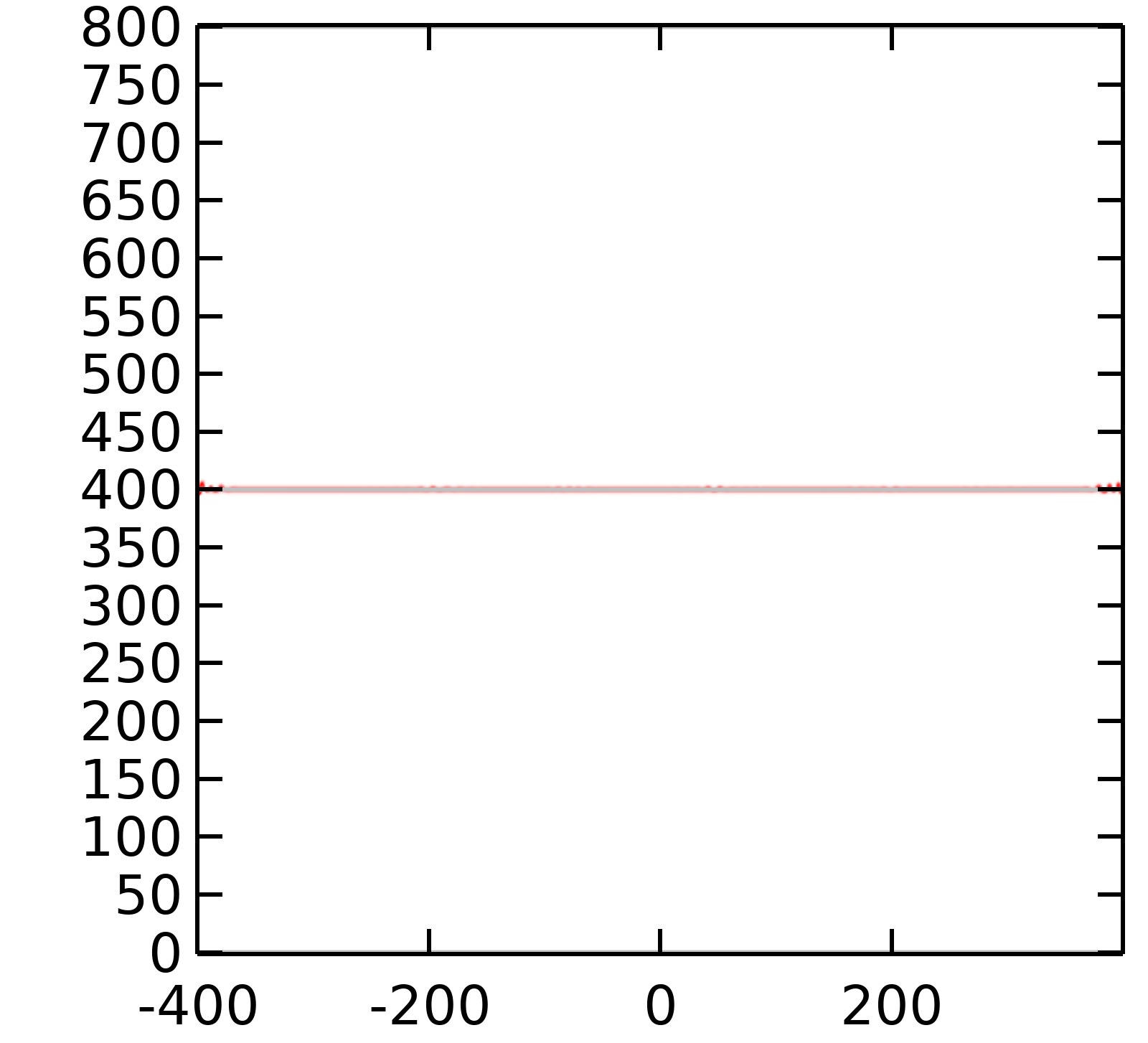}
    \includegraphics[trim=31 14 0 0,clip,width=0.242\textwidth]{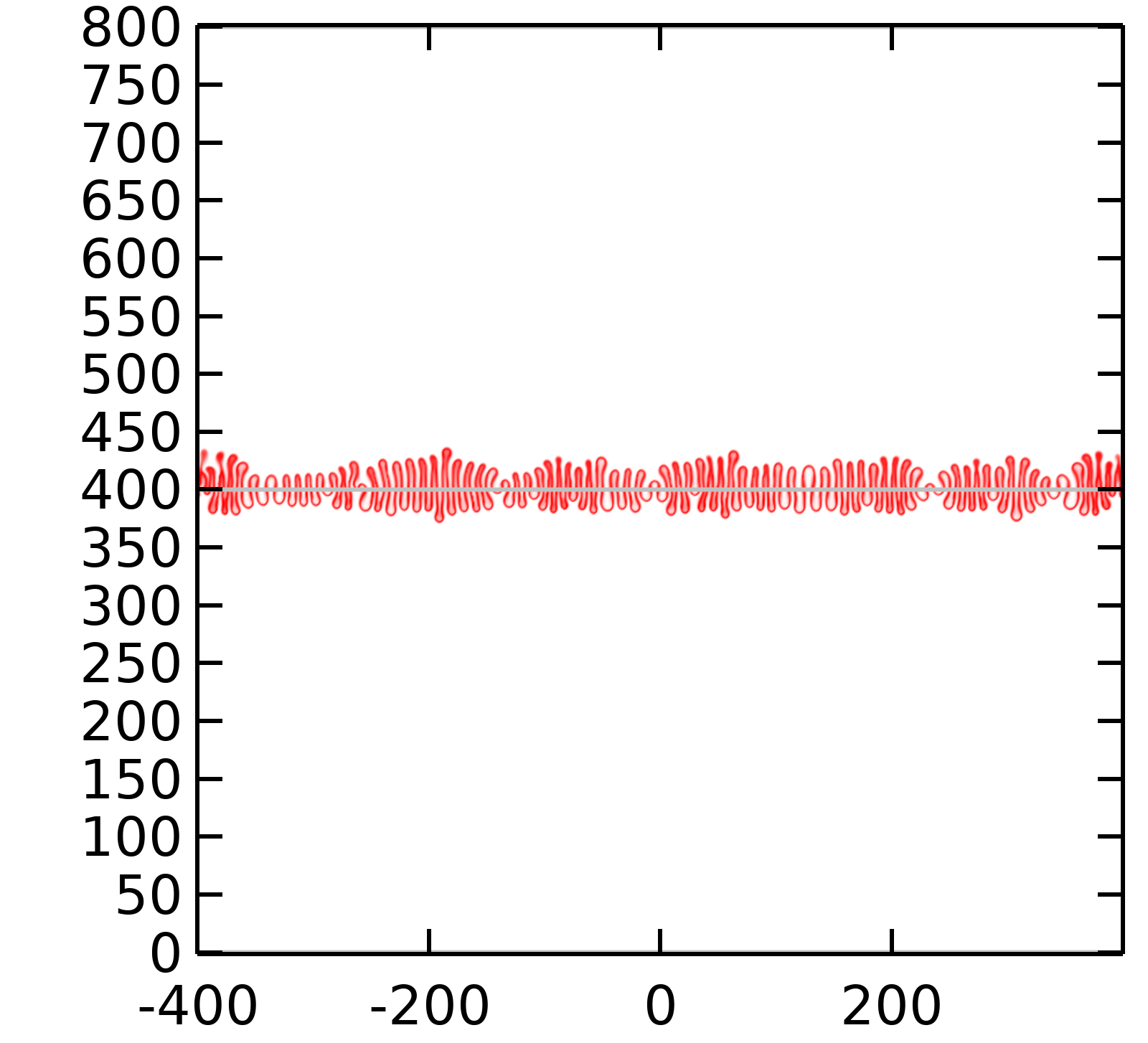}
    \includegraphics[trim=31 14 0 0,clip,width=0.242\textwidth]{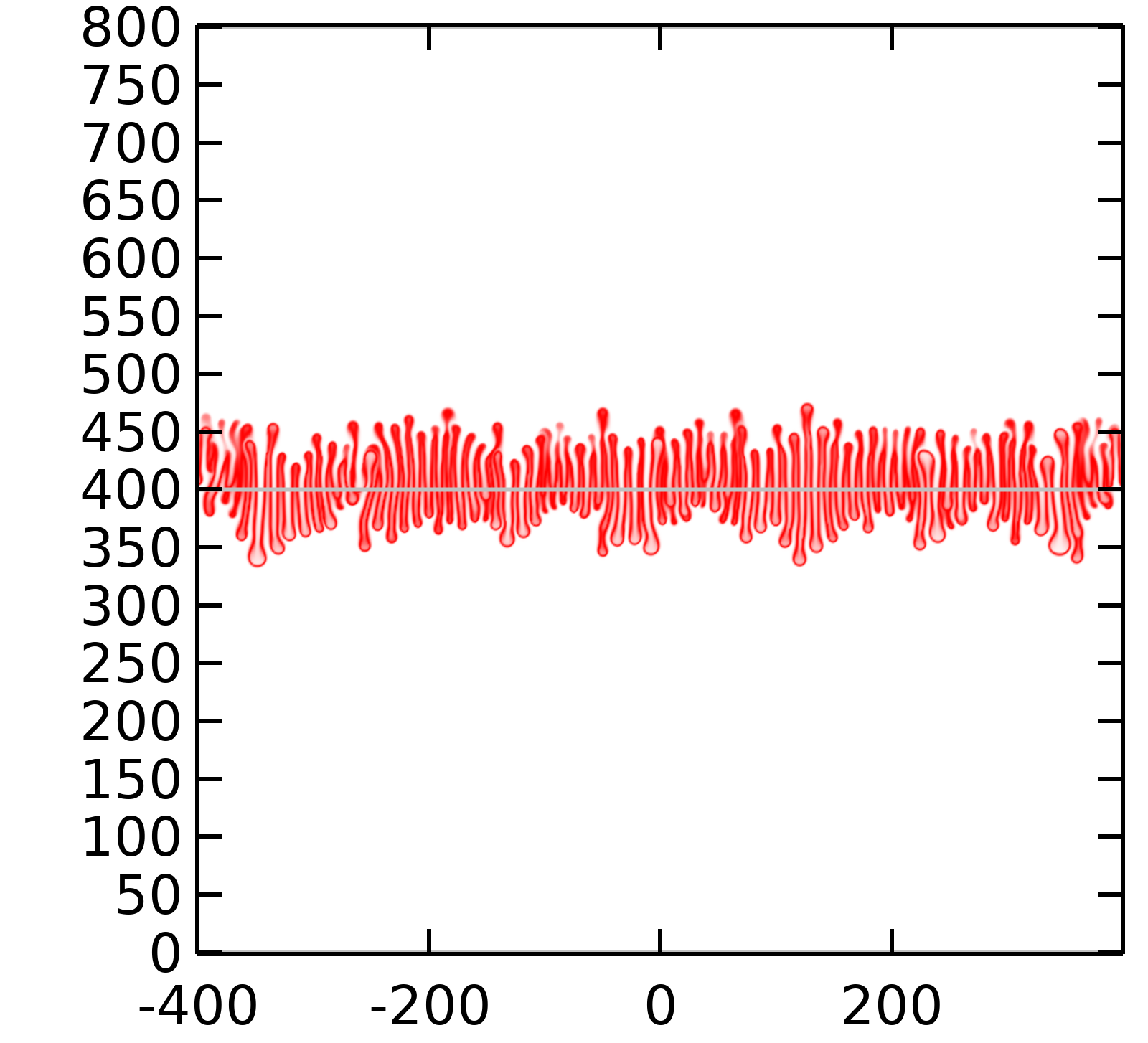}
    \includegraphics[trim=31 14 0 0,clip,width=0.242\textwidth]{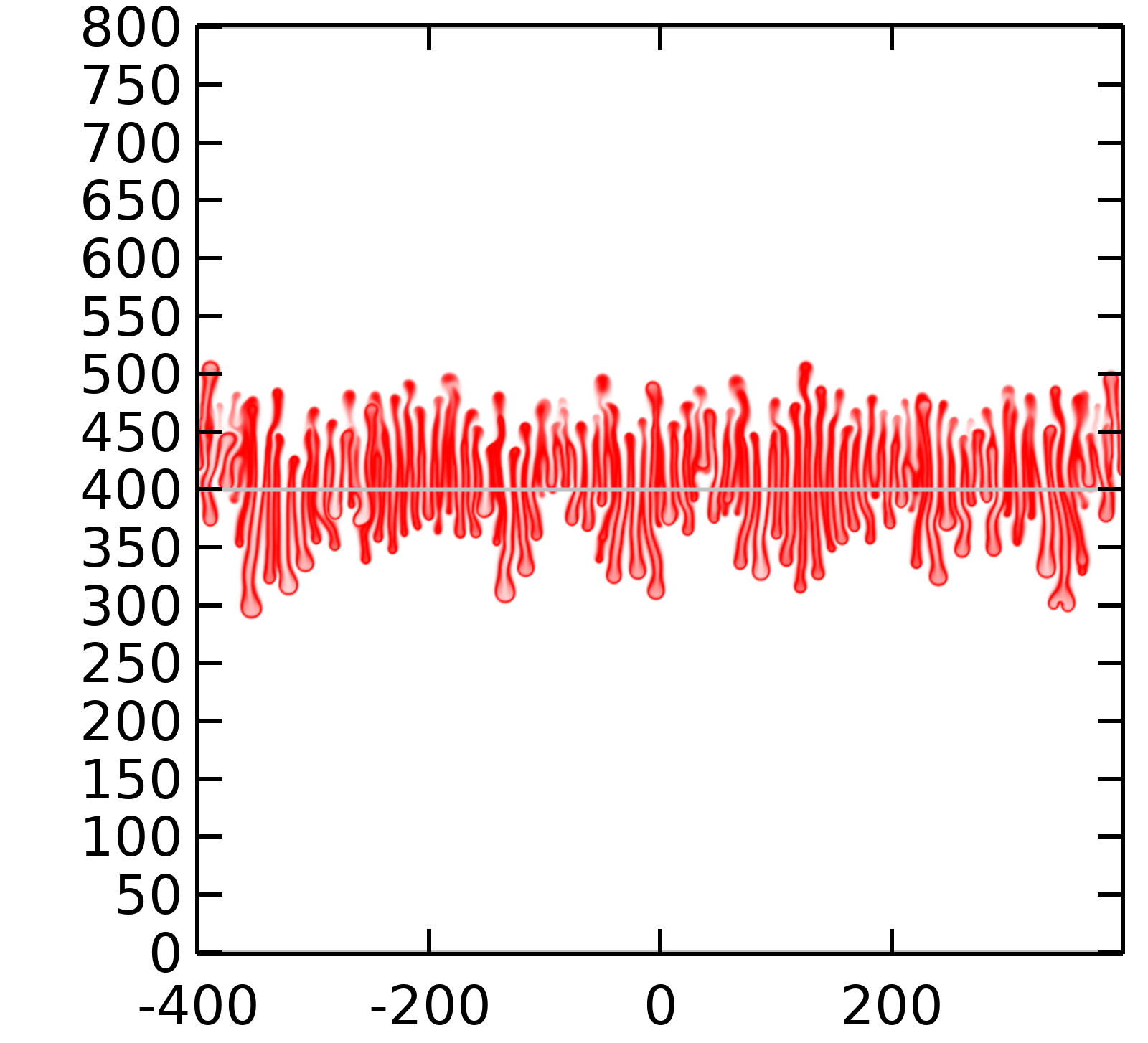}
    \begin{picture}(0,0)
    \put(-450,20){\makebox(0,0)[]{case-IV}}
    \put(-450,10){\makebox(0,0)[]{\scriptsize ($R_{A}=1,R_{C}=0,\alpha=4$)}}
   \end{picture}

    \includegraphics[trim=31 14 0 0,clip,width=0.242\textwidth]{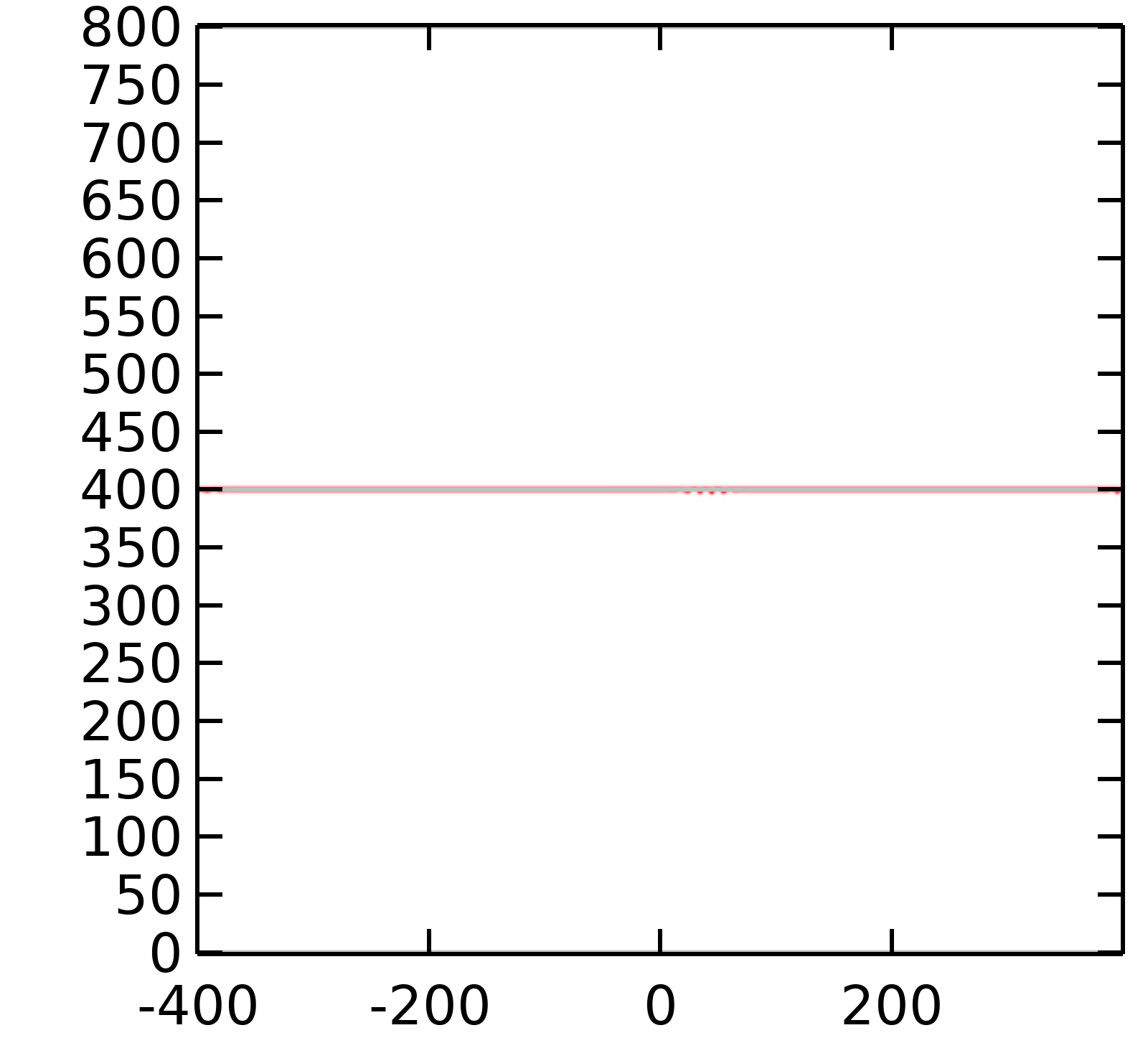}
    \includegraphics[trim=31 14 0 0,clip,width=0.242\textwidth]{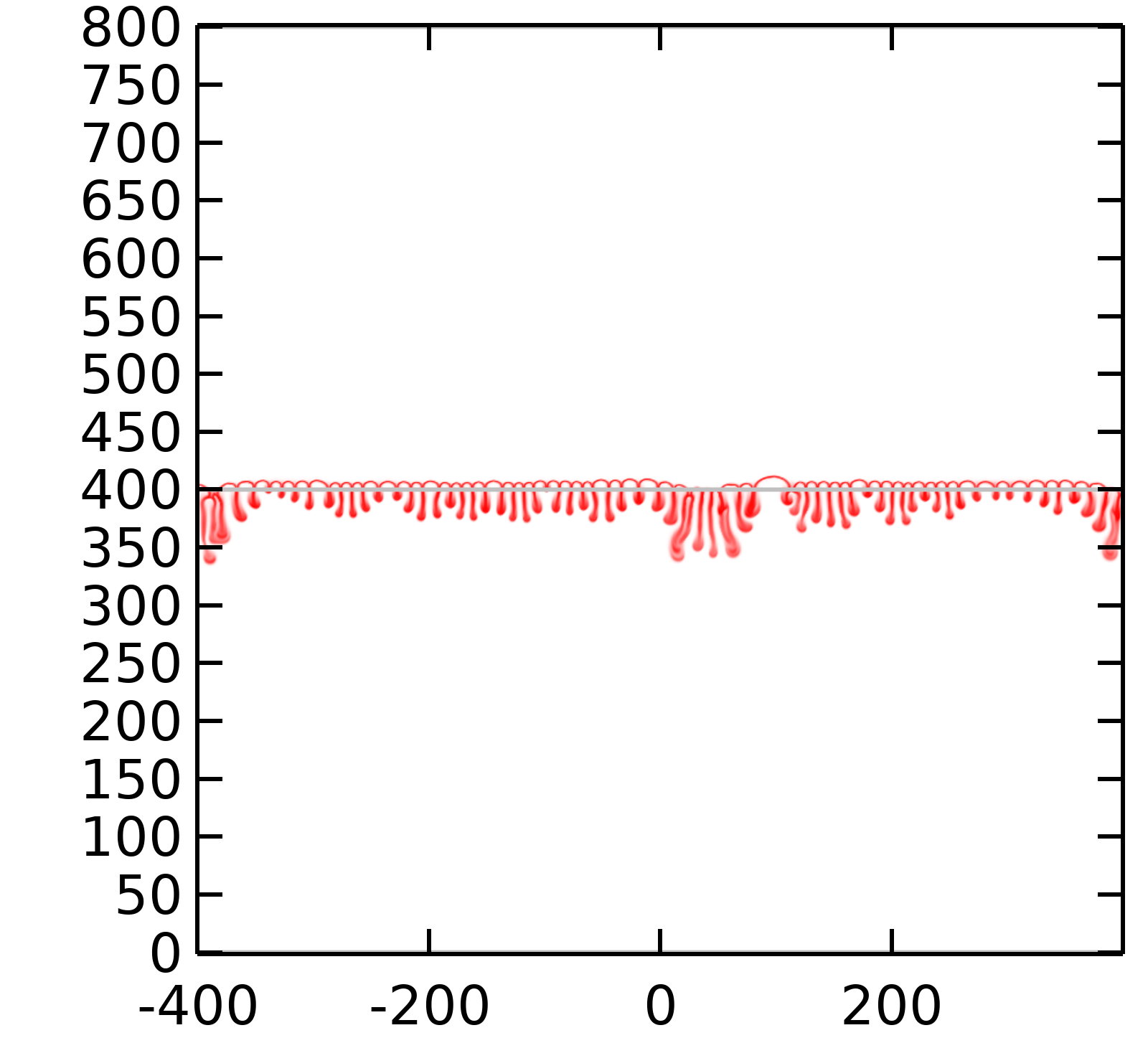}
    \includegraphics[trim=31 14 0 0,clip,width=0.242\textwidth]{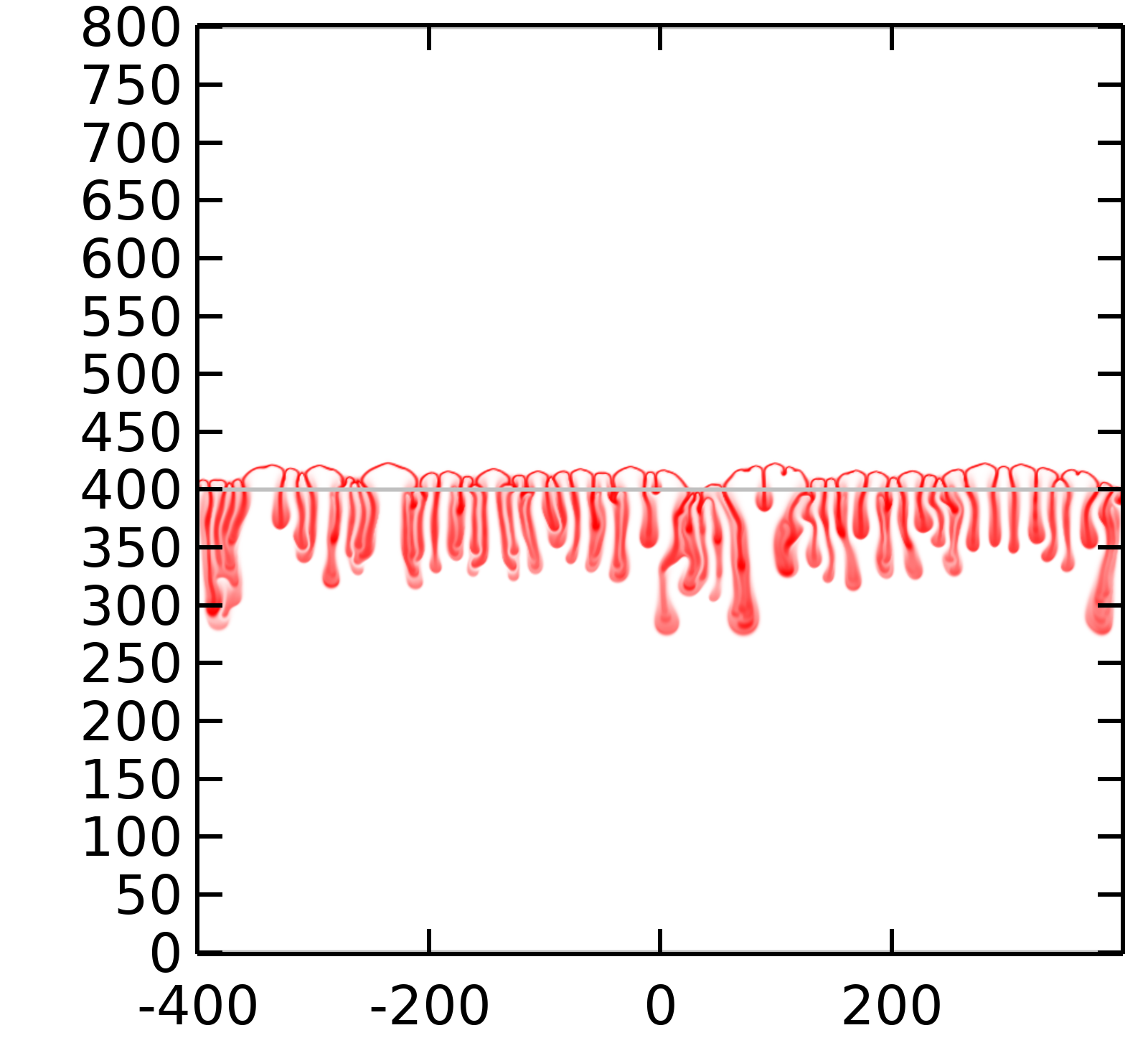}
    \includegraphics[trim=31 14 0 0,clip,width=0.242\textwidth]{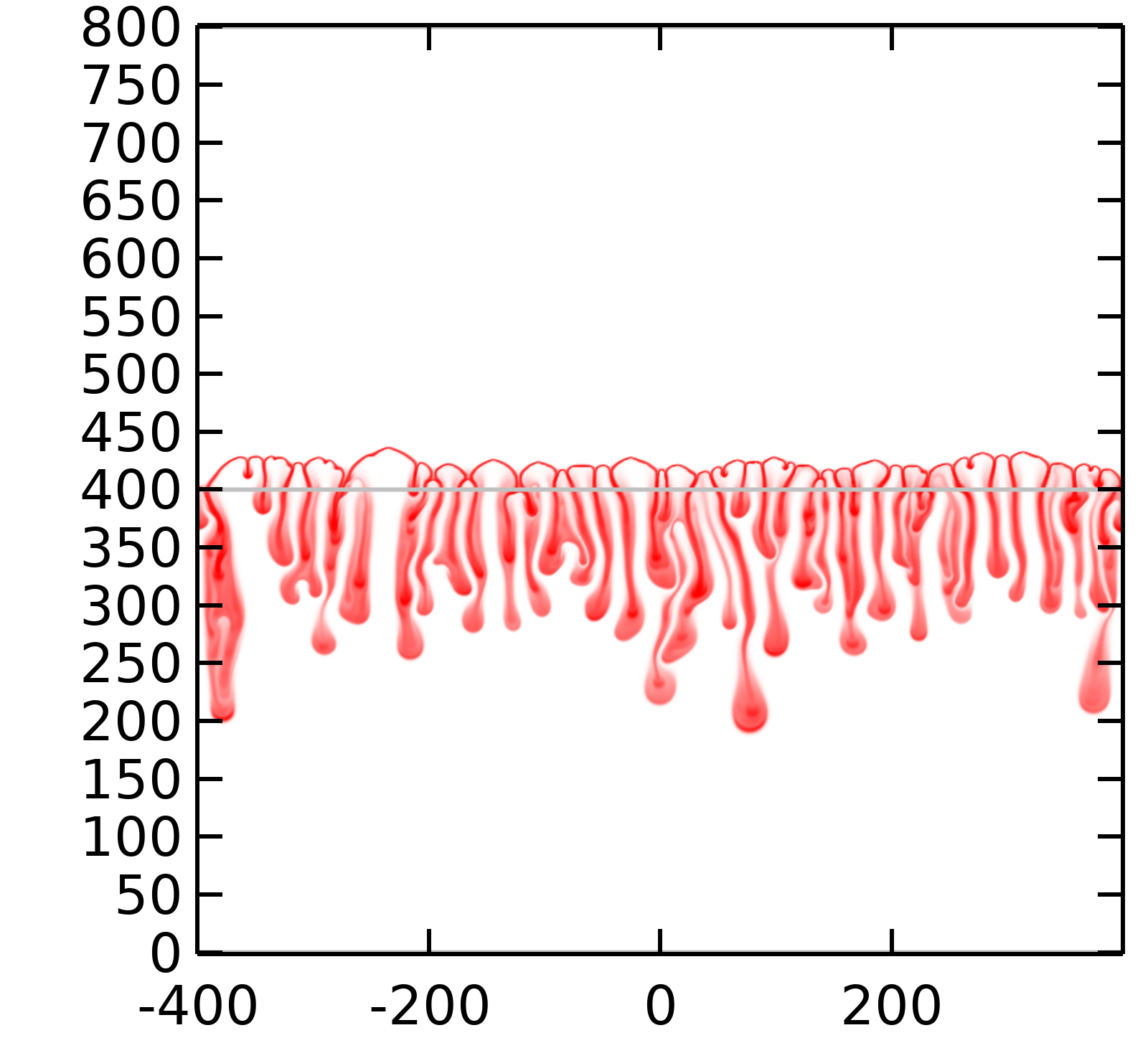}
    \begin{picture}(0,0)
    \put(-450,20){\makebox(0,0)[]{case-V}}
    \put(-450,10){\makebox(0,0)[]{\scriptsize ($R_{A}=0,R_{C}=2,\alpha=0$)}}
   \end{picture}

    \includegraphics[trim=31 14 0 0,clip,width=0.242\textwidth]{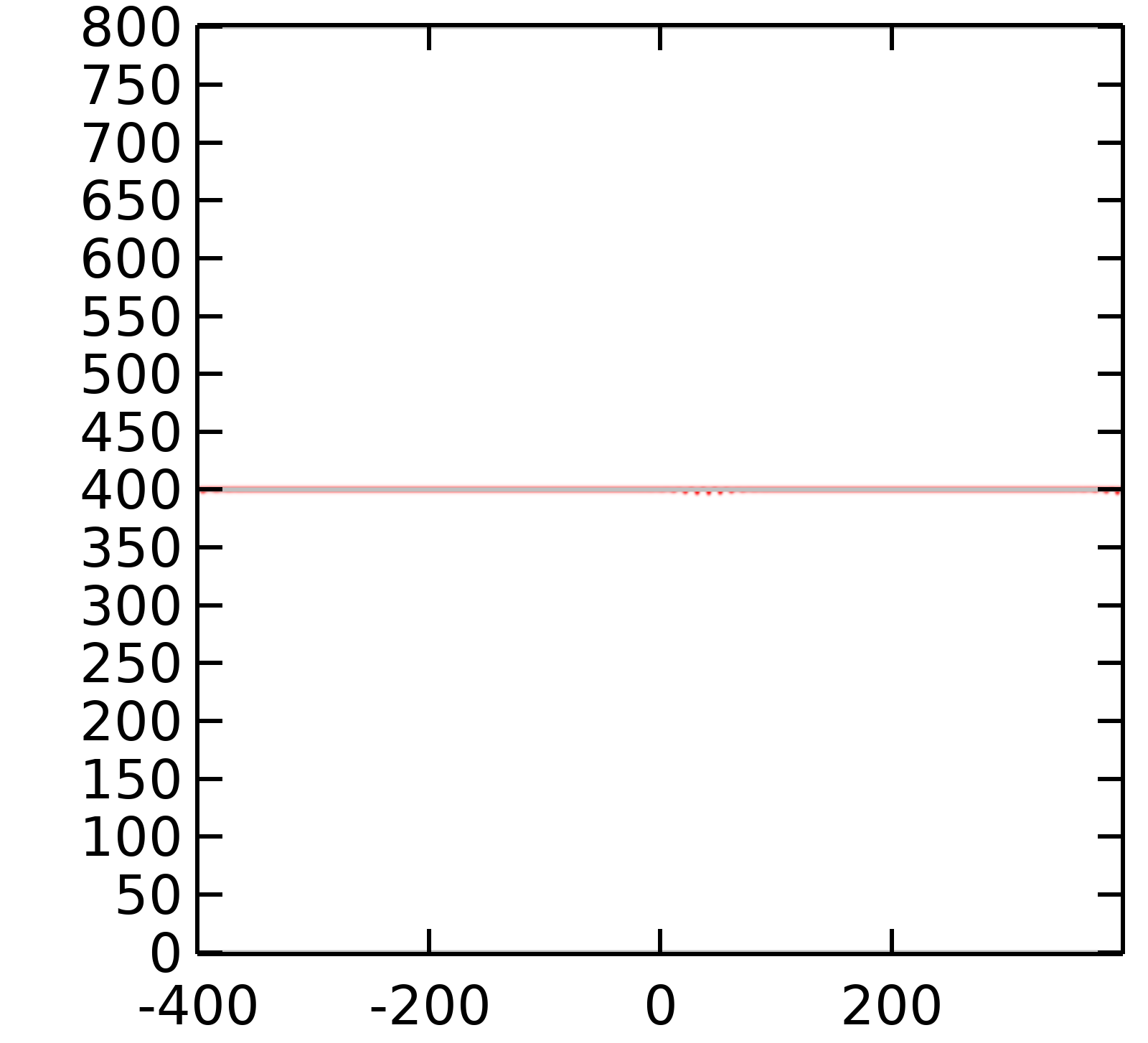}
    \includegraphics[trim=31 14 0 0,clip,width=0.242\textwidth]{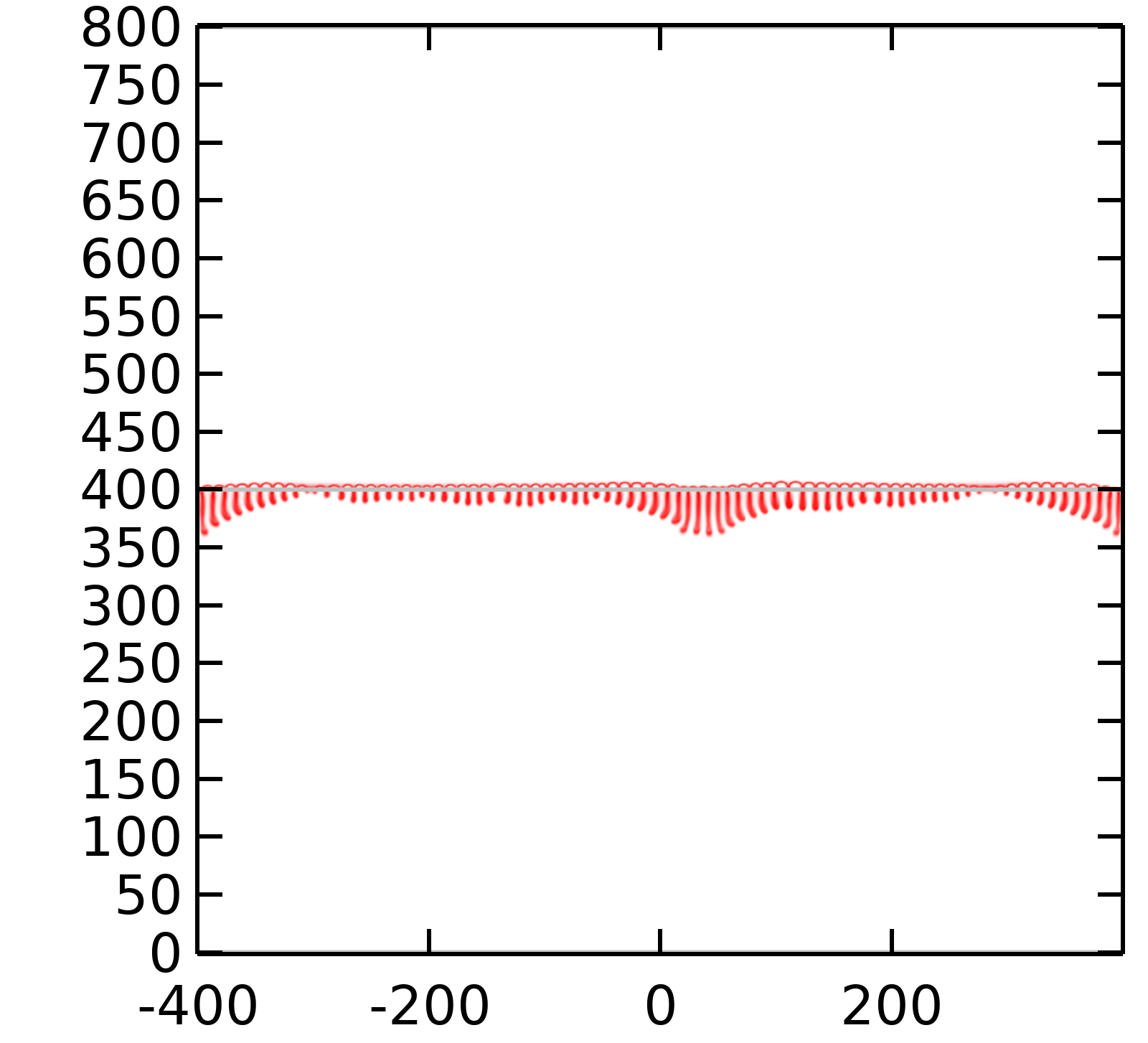}
    \includegraphics[trim=31 14 0 0,clip,width=0.242\textwidth]{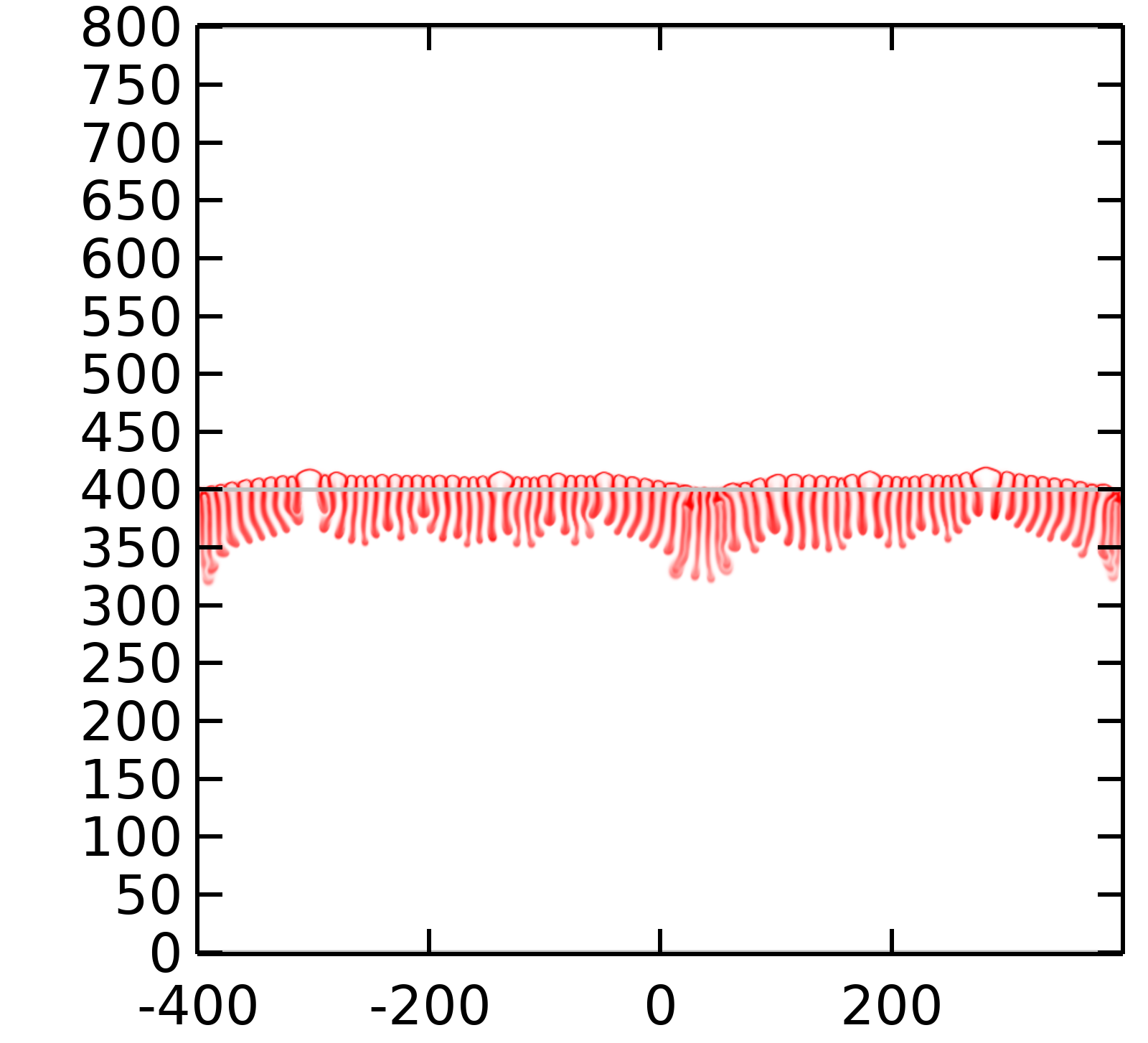}
    \includegraphics[trim=31 14 0 0,clip,width=0.242\textwidth]{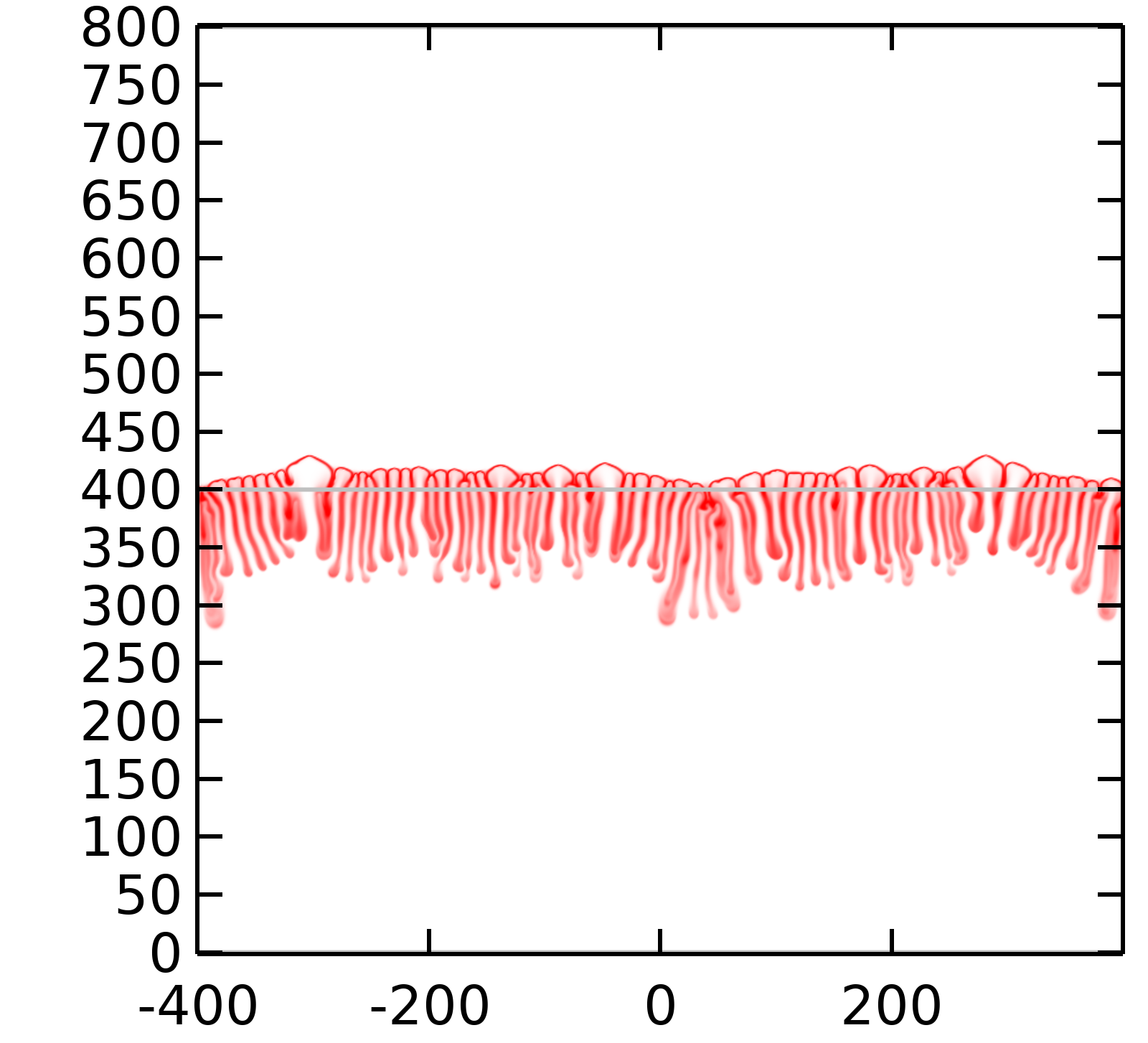}
    \begin{picture}(0,0)
    \put(-450,20){\makebox(0,0)[]{case-VI}}
    \put(-450,10){\makebox(0,0)[]{\scriptsize ($R_{A}=0,R_{C}=2,\alpha=2$)}}
   \end{picture}

    \begin{picture}(0,0)
    \put(-95,650){\makebox(0,0)[]{$t=250$}}
    \put(30,650){\makebox(0,0)[]{$t=500$}}
    \put(160,650){\makebox(0,0)[]{$t=750$}}
    \put(290,650){\makebox(0,0)[]{$t=1000$}}
   \end{picture}
    \hspace{0.38 in}\includegraphics[trim=31 10 0 140,clip,width=0.3\textwidth]{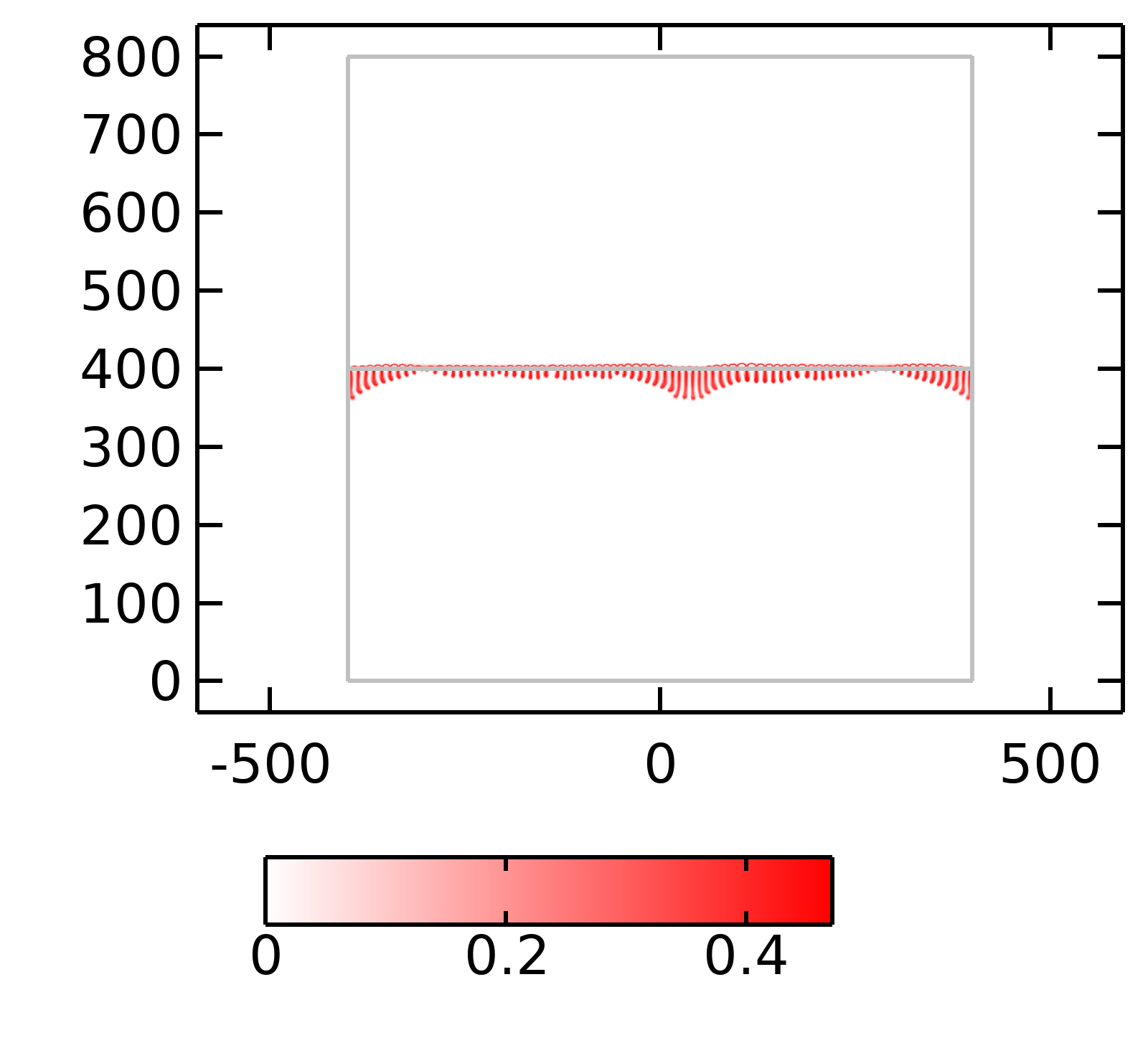}
   \vspace{-0.2 cm}\caption{Spatio-temporal evolution of the reaction product concentration profile \( c \) for cases II to VI, shown from the top to bottom panels, respectively. The axis labels are same as Fig. \ref{fig:abc}.}
    \label{fig:flat}
\end{figure}
Fig. \ref{fig:abc} shows the spatio-temporal evolution of the reactants \( A \) (top panel) and \( B \) (middle panel), along with the product concentration profiles \( C \) (bottom panel). Case I is considered here, where \( R_A = 2 \) and \(  R_C = 0 \), with the top fluid being denser than the lower reactant fluid \( B \). The system undergoes the reaction \( A + B \to C \) with a kinetic constant \( k = 1 \). There is no change in permeability due to the reaction (\( \alpha = 0 \)), and the gravity vector is a unit downward-pointing vector.  The initial conditions for the concentrations are given by \eqref{e71}, while the initial velocity vector is set to zero. The boundary conditions are specified as in \eqref{ic's}. This setup represents a case where classical Rayleigh-Taylor fingering is inherently present in the system. However, at the interface of the fingers, the reaction modulates the pattern formations, as clearly observed in the product concentration profiles shown in the bottom panel of Fig. \ref{fig:abc}.

Fig. \ref{fig:flat} summarizes the pattern formations for the remaining five cases through the spatio-temporal evolution of the product concentration profiles. In case-II, where there is no density stratification (\( R_A = R_C = 0 \)) and the reaction only alters permeability (\( \alpha = 4 \)), the system exhibits a flat front of reaction-diffusion without density-driven fingering formations (see the top panel of Fig. \ref{fig:flat}). However, permeability changes significantly affect reactive density fingers. For example, in case-III (Fig. \ref{fig:flat}), we use the same parameters as in case-I, which included density fingering, but modify the permeability parameter to \( \alpha = 2 \) instead of \( \alpha = 0 \). It is evident that this permeability change alters both the mixing length and the wavelength of the fingers, resulting in narrower fingers. To observe this, compare the panel for \( t = 500 \) in the product concentration profiles \( C \) of Fig. \ref{fig:abc} with the panel for case-III in Fig. \ref{fig:flat} at the same time (\( t = 500 \)). This shows that the local velocity change in Brinkman convection, induced by permeability alteration, significantly influences the fingers' properties and flow mixing dynamics. Similarly, we present the finger growth for case-IV, where the strength of density stratification is reduced while the effects of permeability change by the reaction product are increased compared to case-III. Specifically, we set \( R_A = 1 \) and \( \alpha = 2 \). This modification further alters the mixing length and narrows the fingers by reducing the wavelengths even more (see the panel for case-IV and compare it with that for case-III in Fig. \ref{fig:flat}).

Next, we consider \( R_A = 0 \), meaning there is no inherent density stratification in the flow before the reaction, as both reactants have the same density. However, setting \( R_C = 2 \) implies that the reaction product is denser than both reactants. This defines case-V, where there is no permeability change caused by the reaction product (\( \alpha = 0 \)). In this scenario, fingering instability is induced solely by local density changes at the reactive interface. The panel for case-V in Fig. \ref{fig:flat} shows that the \( A-C \) interface remains nearly flat, with fingers growing downward. This highlights a major change in the type of pattern formation and finger structure, which also agrees qualitatively with those observed in pure Darcian flow \cite{Almarcha_2010}. Similarly, the wavelengths of the reaction-induced density fingers and the mixing length are further reduced when the reaction also alters the permeability (\( \alpha = 2 \)), as shown in the panel for case-VI in Fig. \ref{fig:flat}.

\begin{figure}
    \centering
    \includegraphics[trim=31 14 0 0,clip,width=0.192\textwidth]{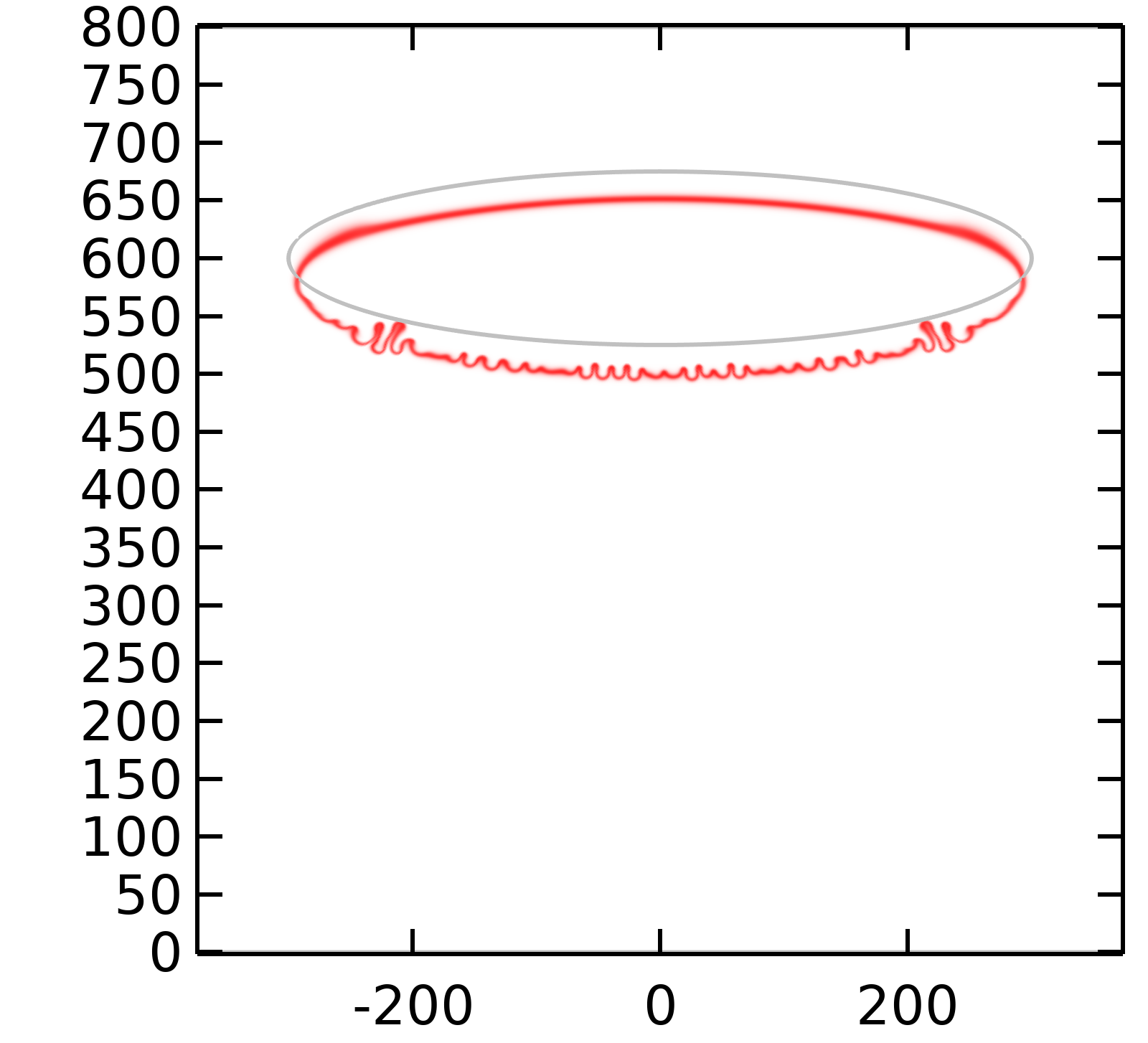}
    \includegraphics[trim=31 14 0 0,clip,width=0.192\textwidth]{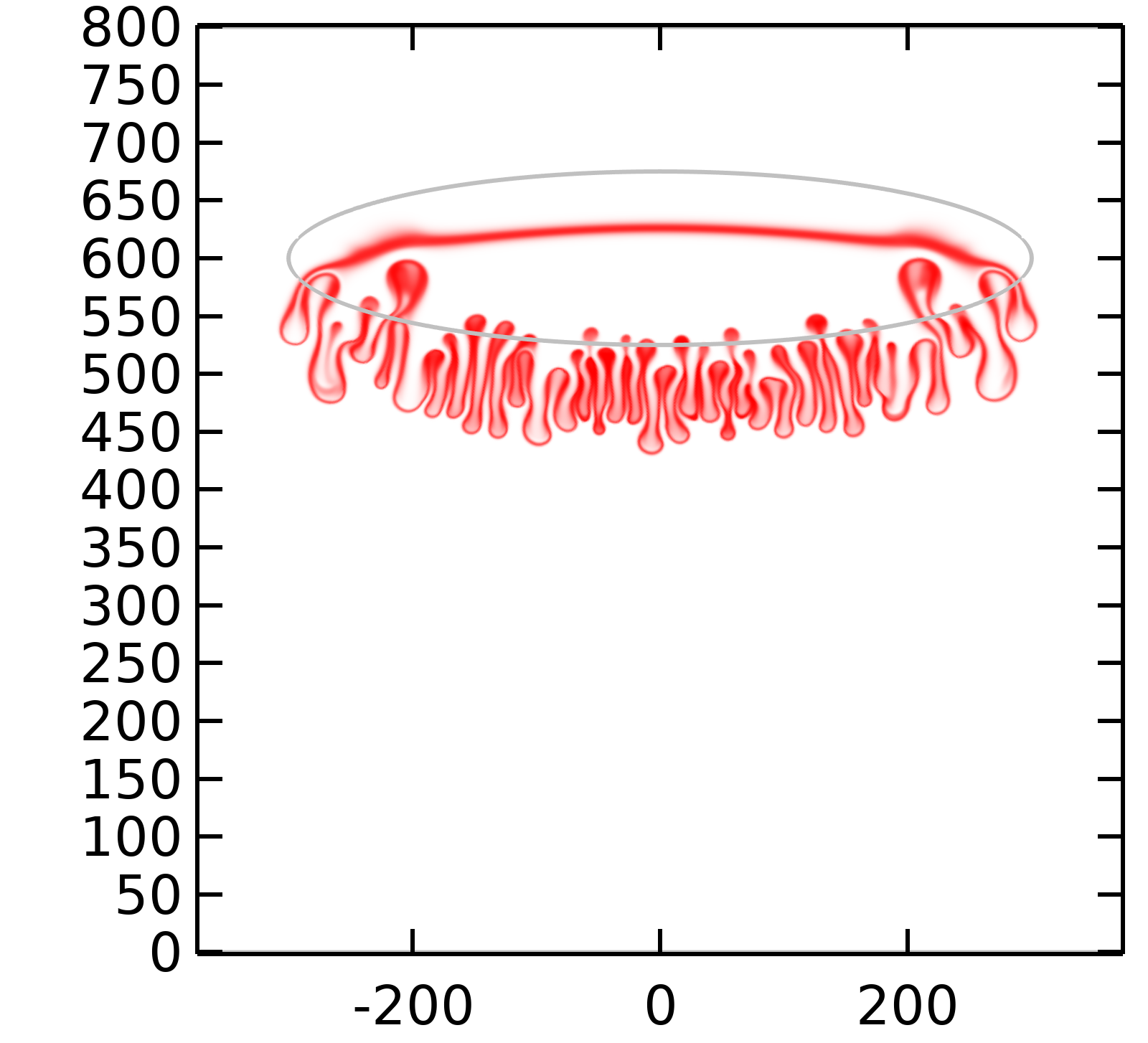}
    \includegraphics[trim=31 14 0 0,clip,width=0.192\textwidth]{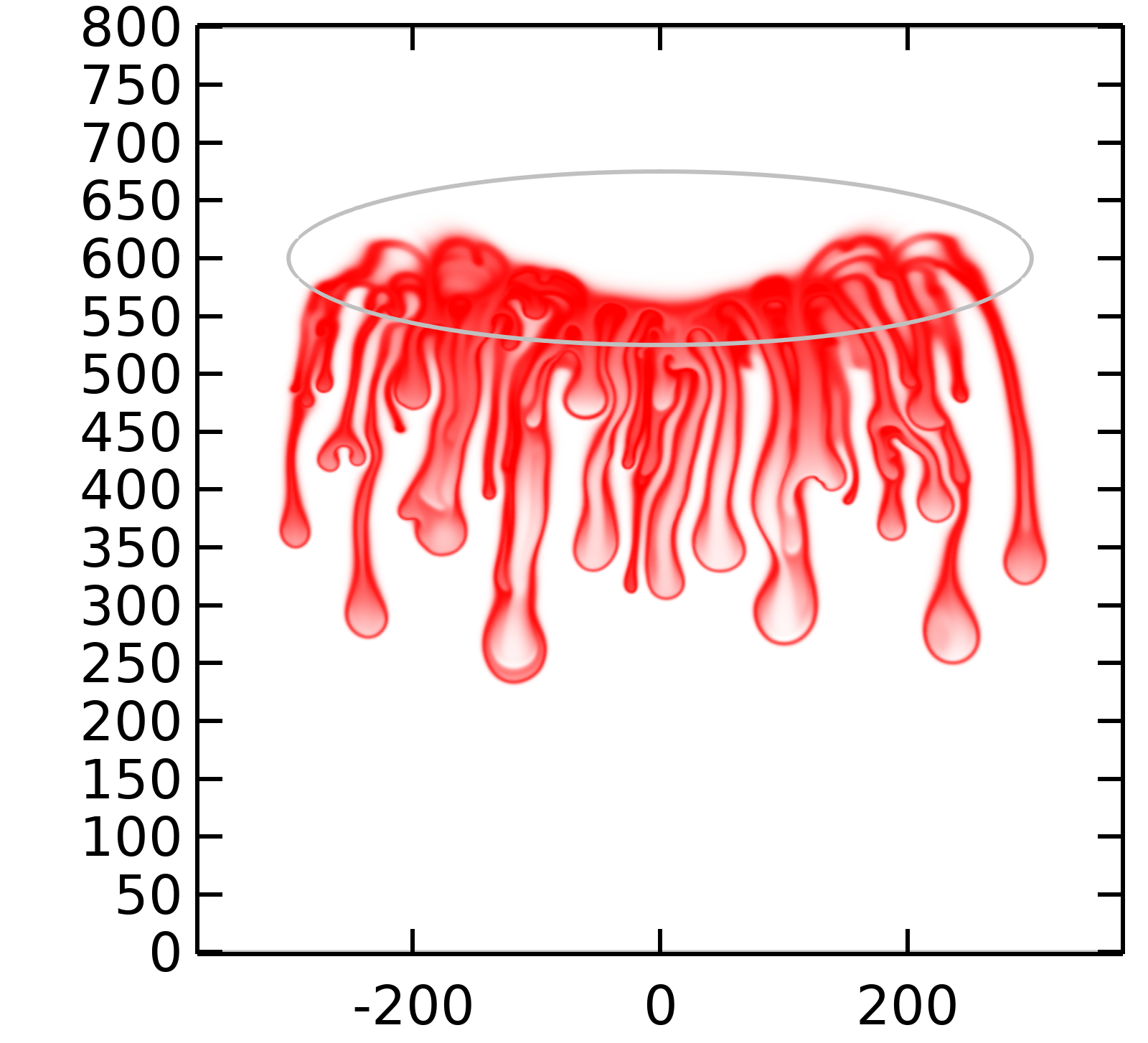}
    \includegraphics[trim=31 14 0 0,clip,width=0.192\textwidth]{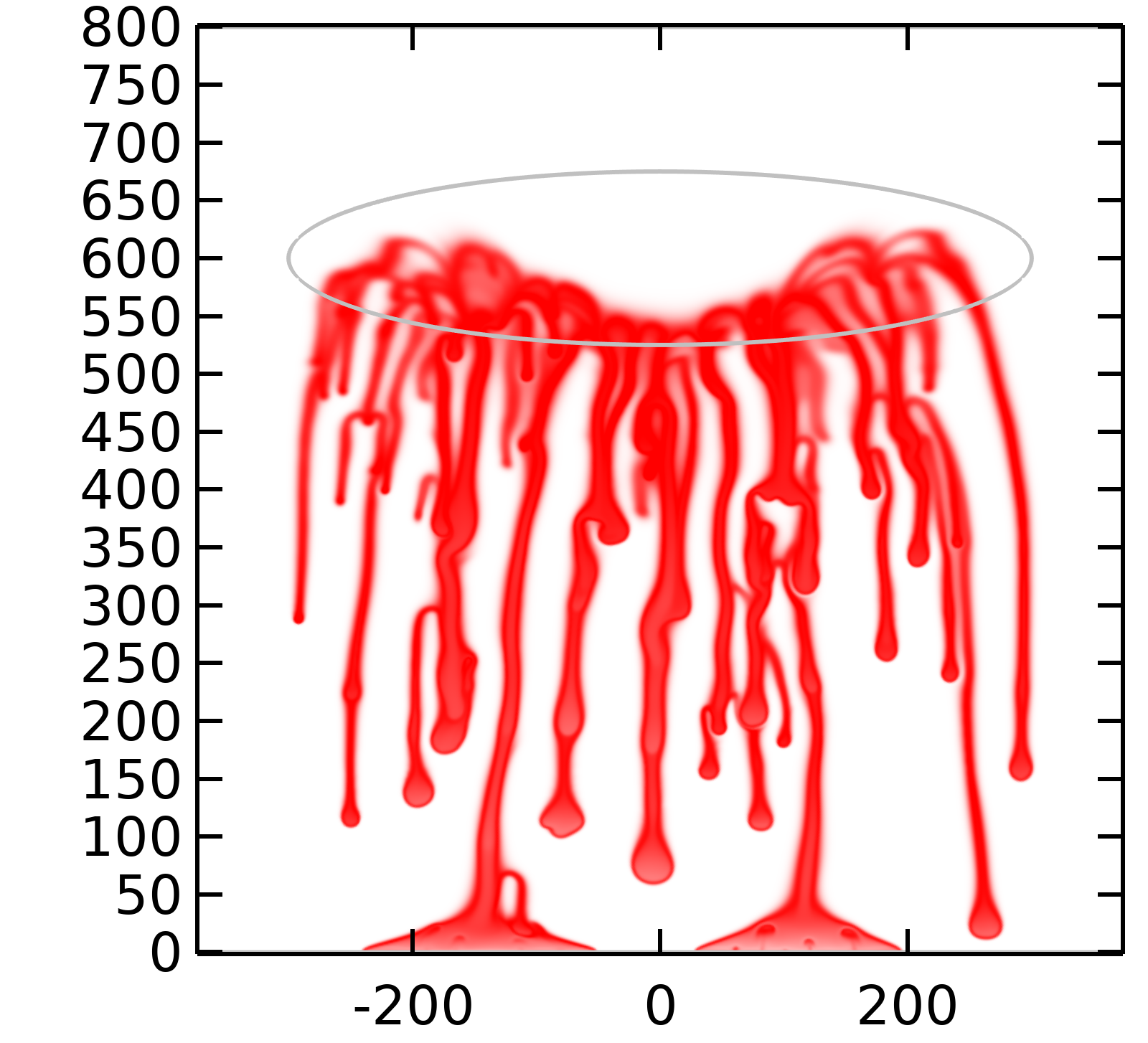}
     \includegraphics[trim=31 14 0 0,clip,width=0.192\textwidth]{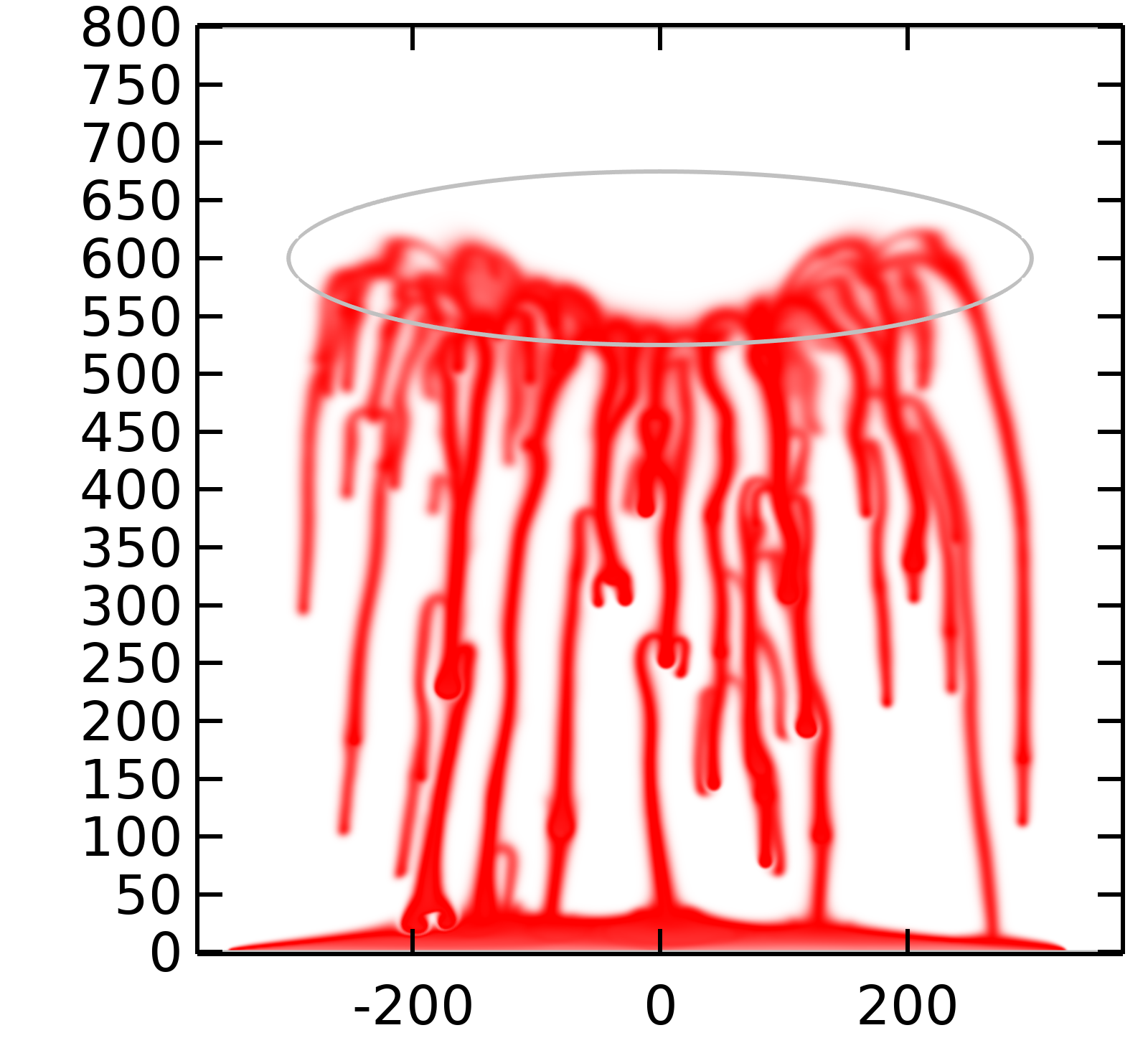}
    \begin{picture}(0,0)
    \put(-460,26){\makebox(0,0)[]{case-I}}
    \put(-460,16){\makebox(0,0)[]{\scriptsize $(R_{A}=2,R_{C}=0,\alpha=0)$}}
   \end{picture}
   
    \includegraphics[trim=31 14 0 0,clip,width=0.192\textwidth]{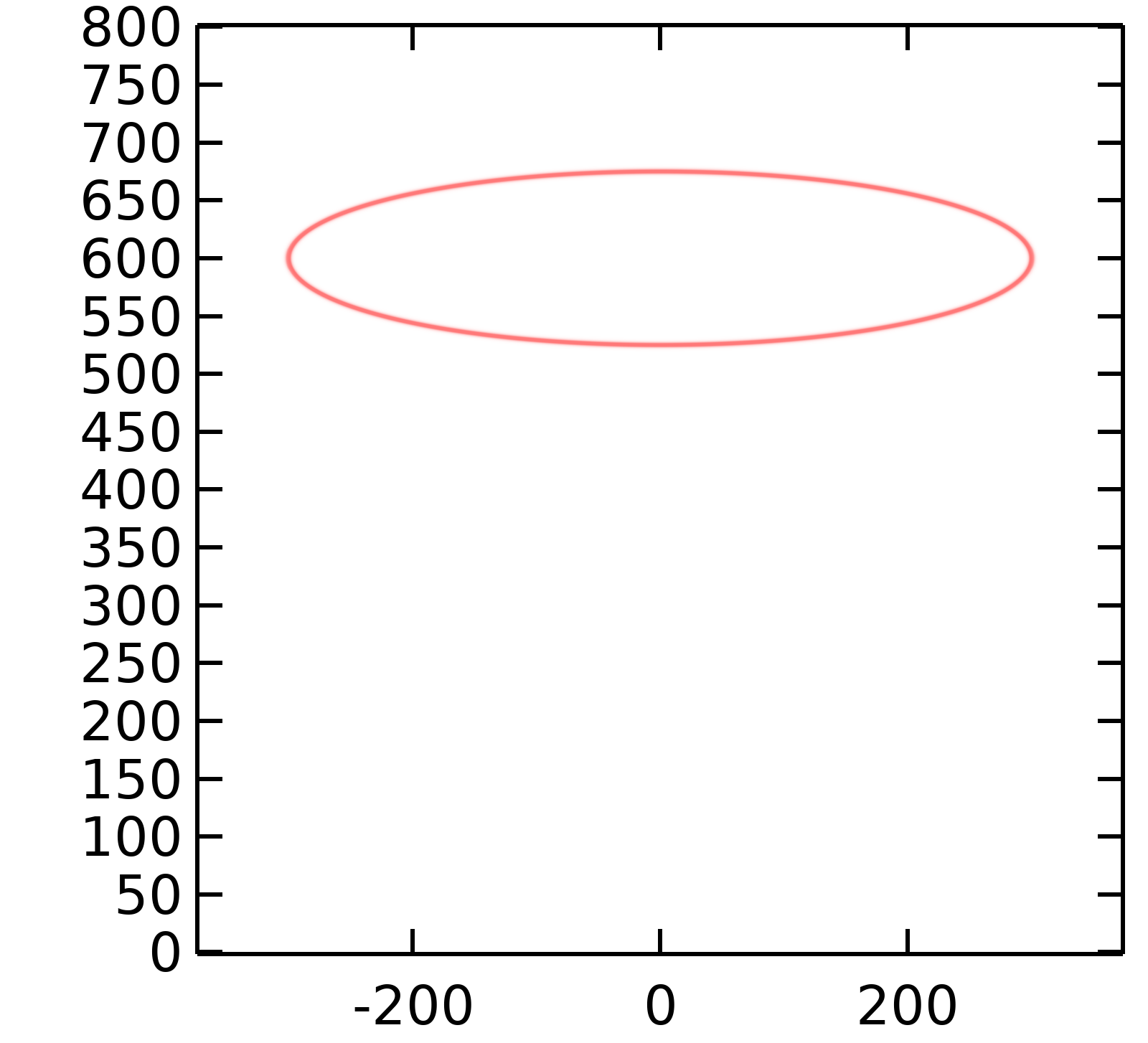}
    \includegraphics[trim=31 14 0 0,clip,width=0.192\textwidth]{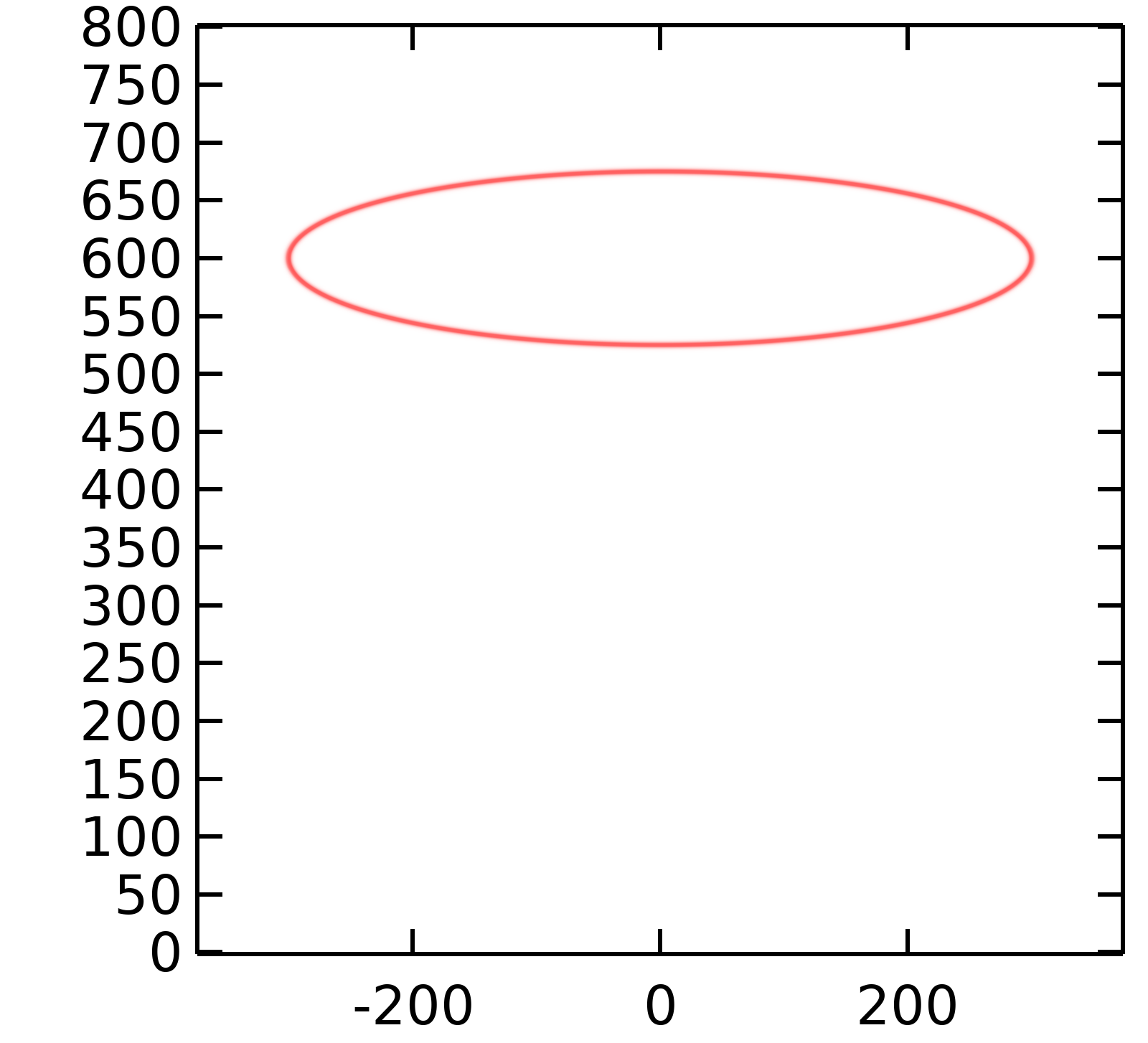}
    \includegraphics[trim=31 14 0 0,clip,width=0.192\textwidth]{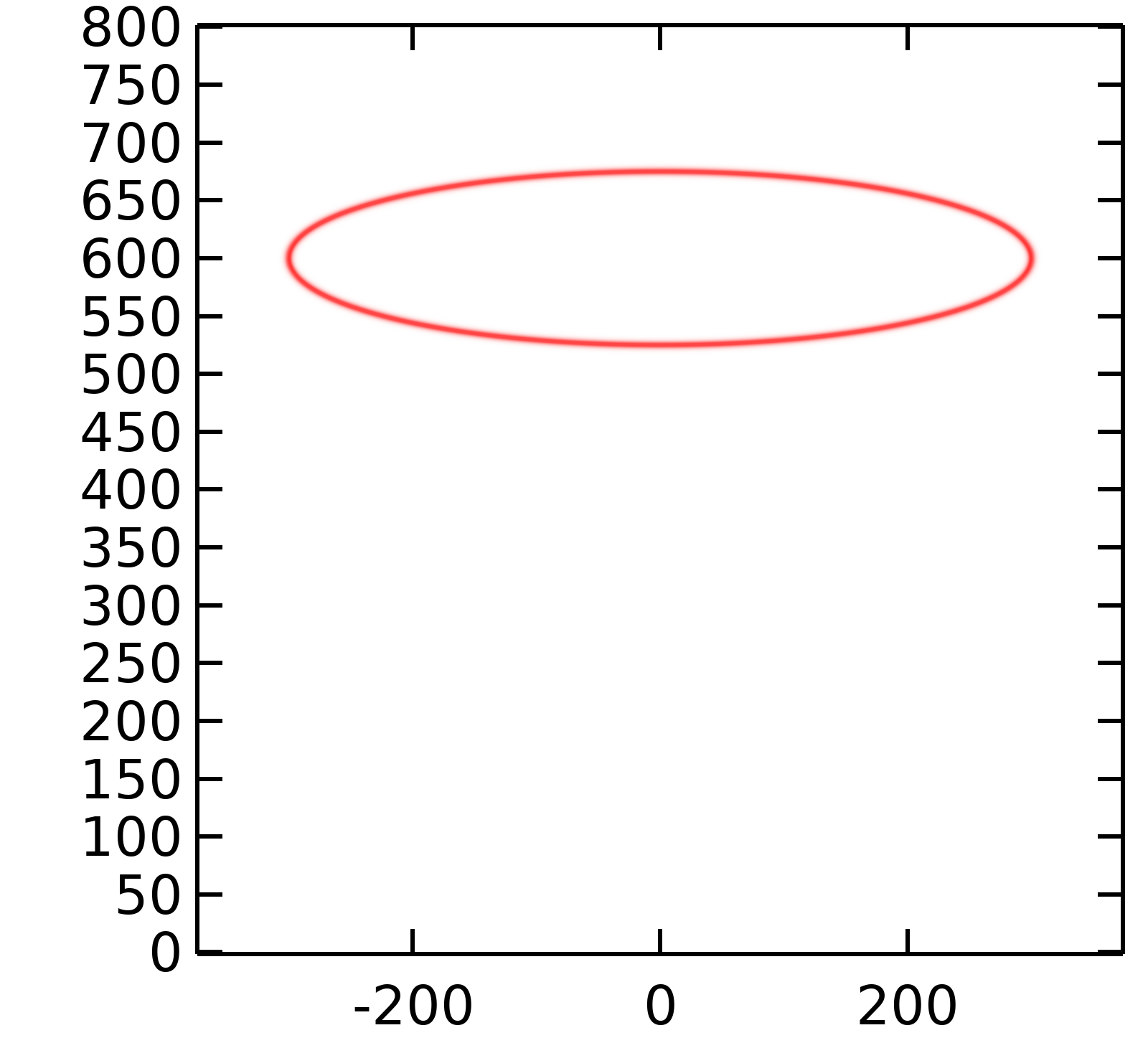}
    \includegraphics[trim=31 14 0 0,clip,width=0.192\textwidth]{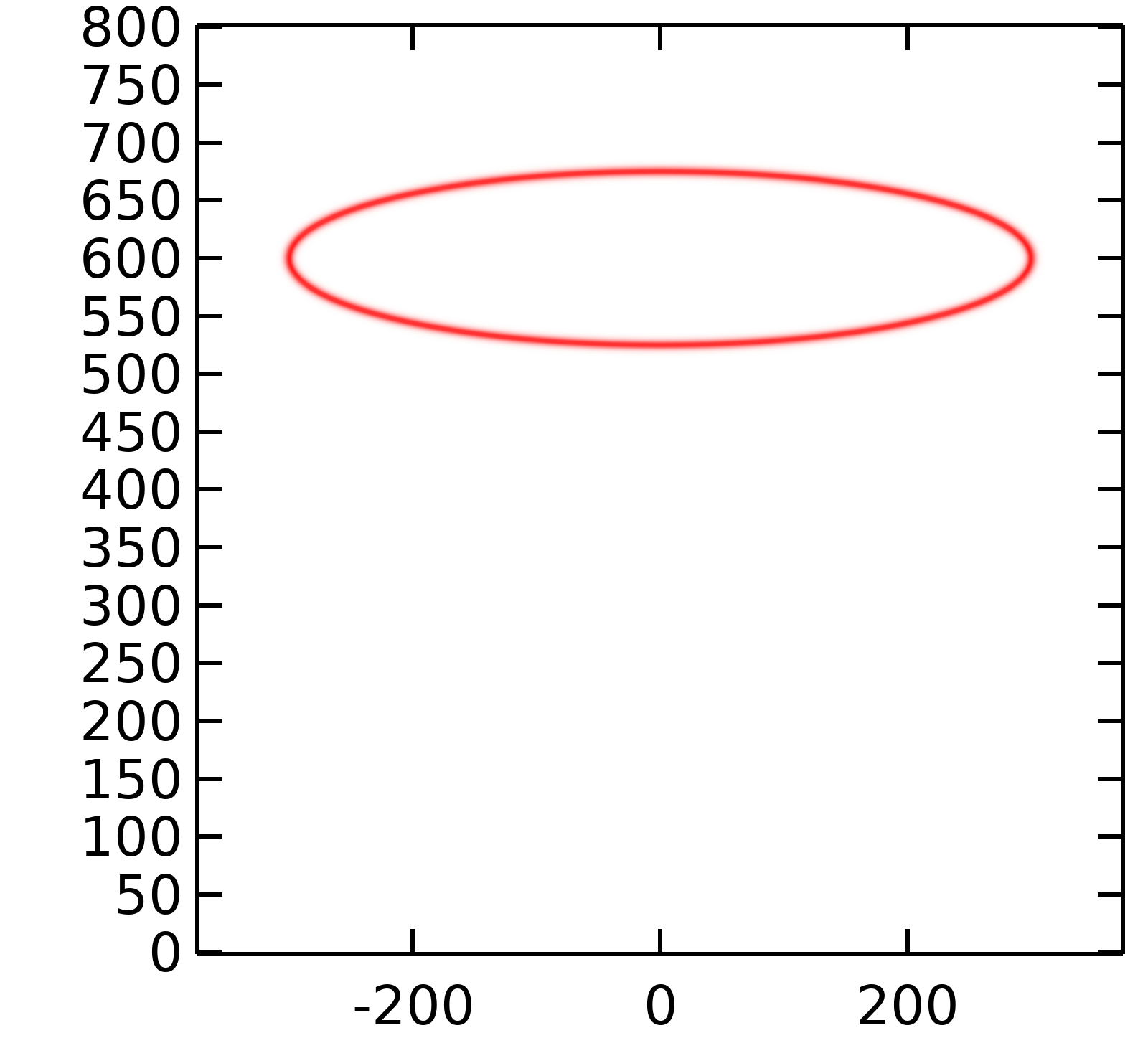}
     \includegraphics[trim=31 14 0 0,clip,width=0.192\textwidth]{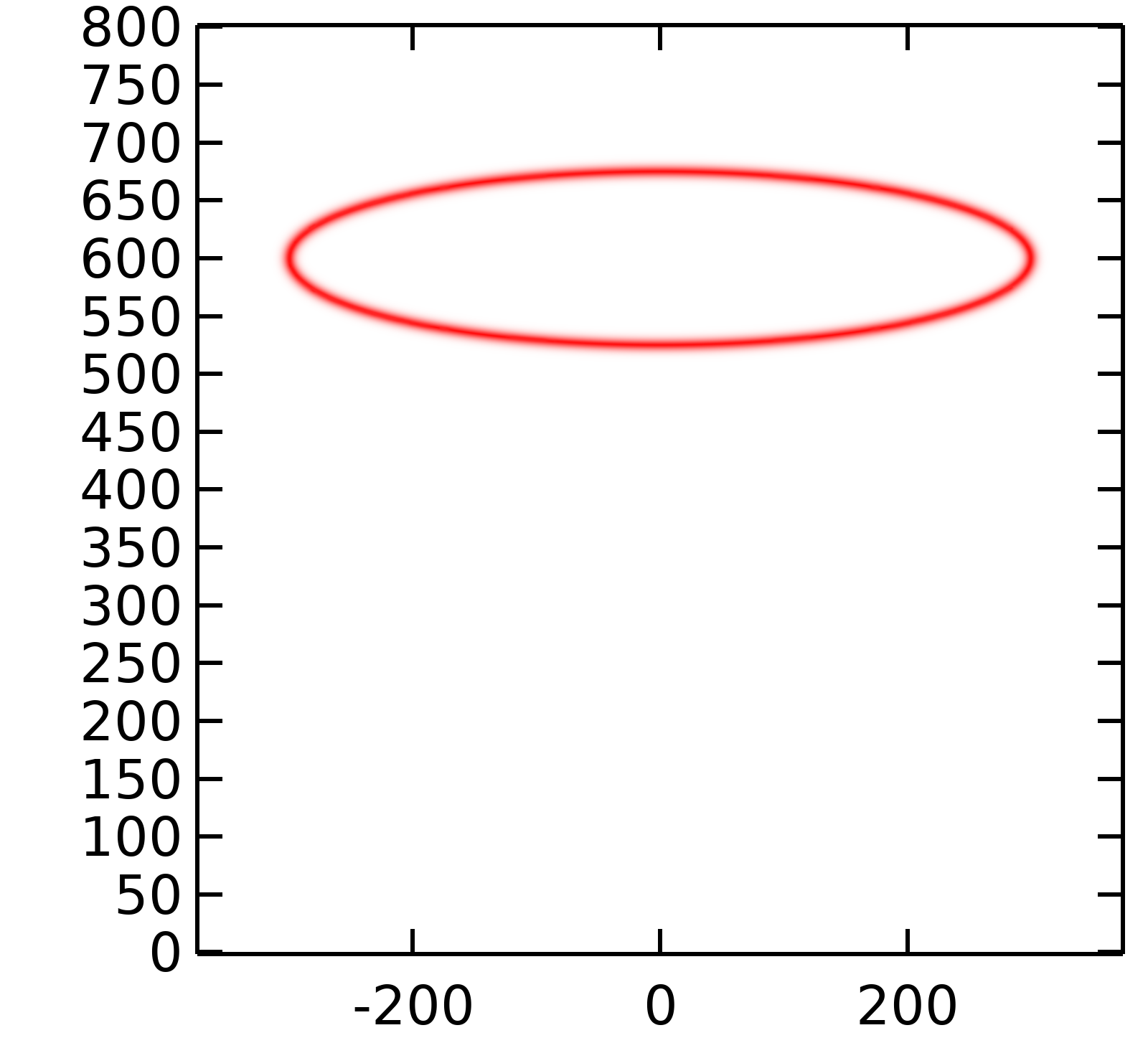}
     \begin{picture}(0,0)
    \put(-460,26){\makebox(0,0)[]{case-II}}
    \put(-460,16){\makebox(0,0)[]{\scriptsize $(R_{A}=0,R_{C}=0,\alpha=4)$}}
   \end{picture}

    \includegraphics[trim=31 14 0 0,clip,width=0.192\textwidth]{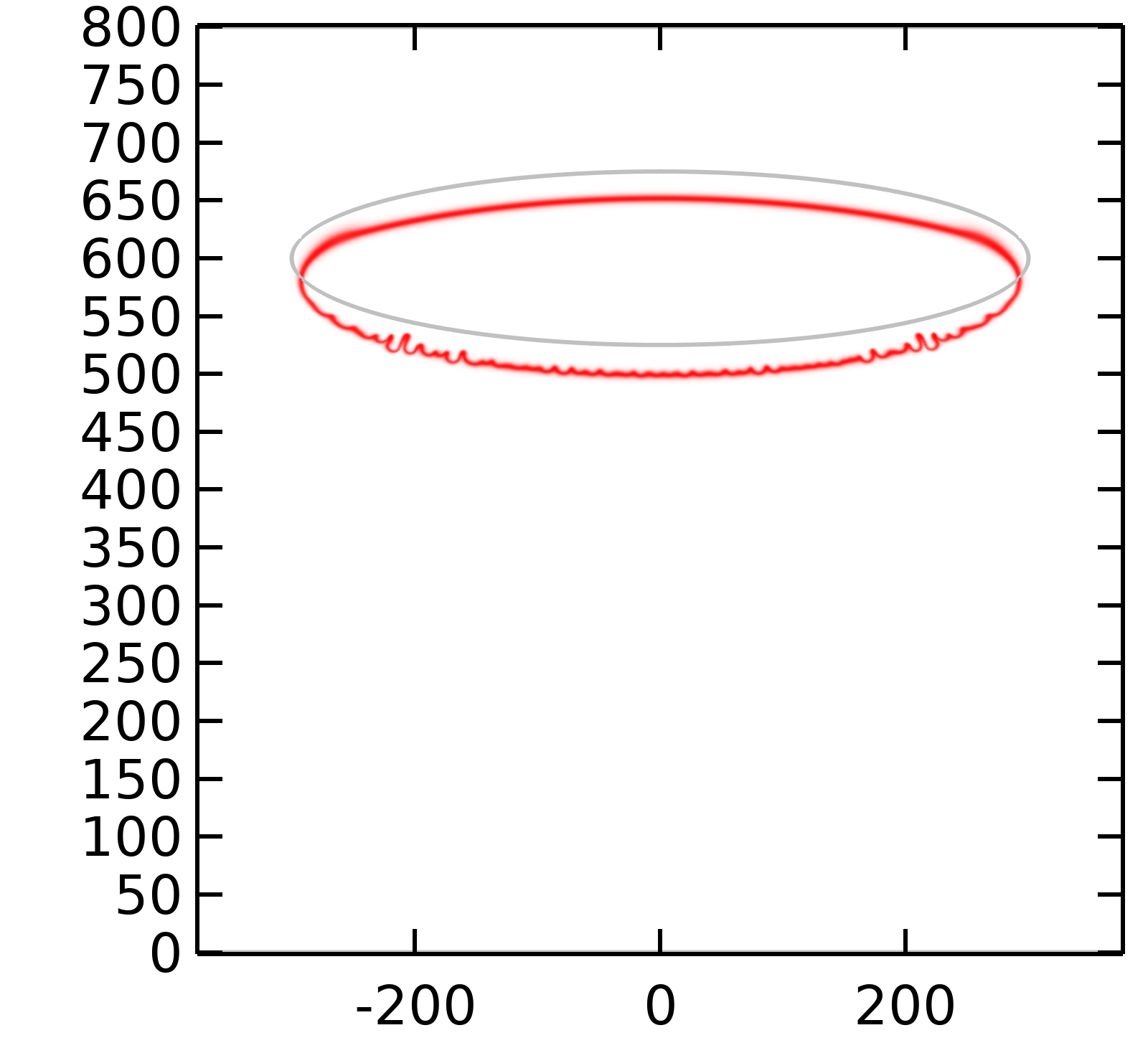}
    \includegraphics[trim=31 14 0 0,clip,width=0.192\textwidth]{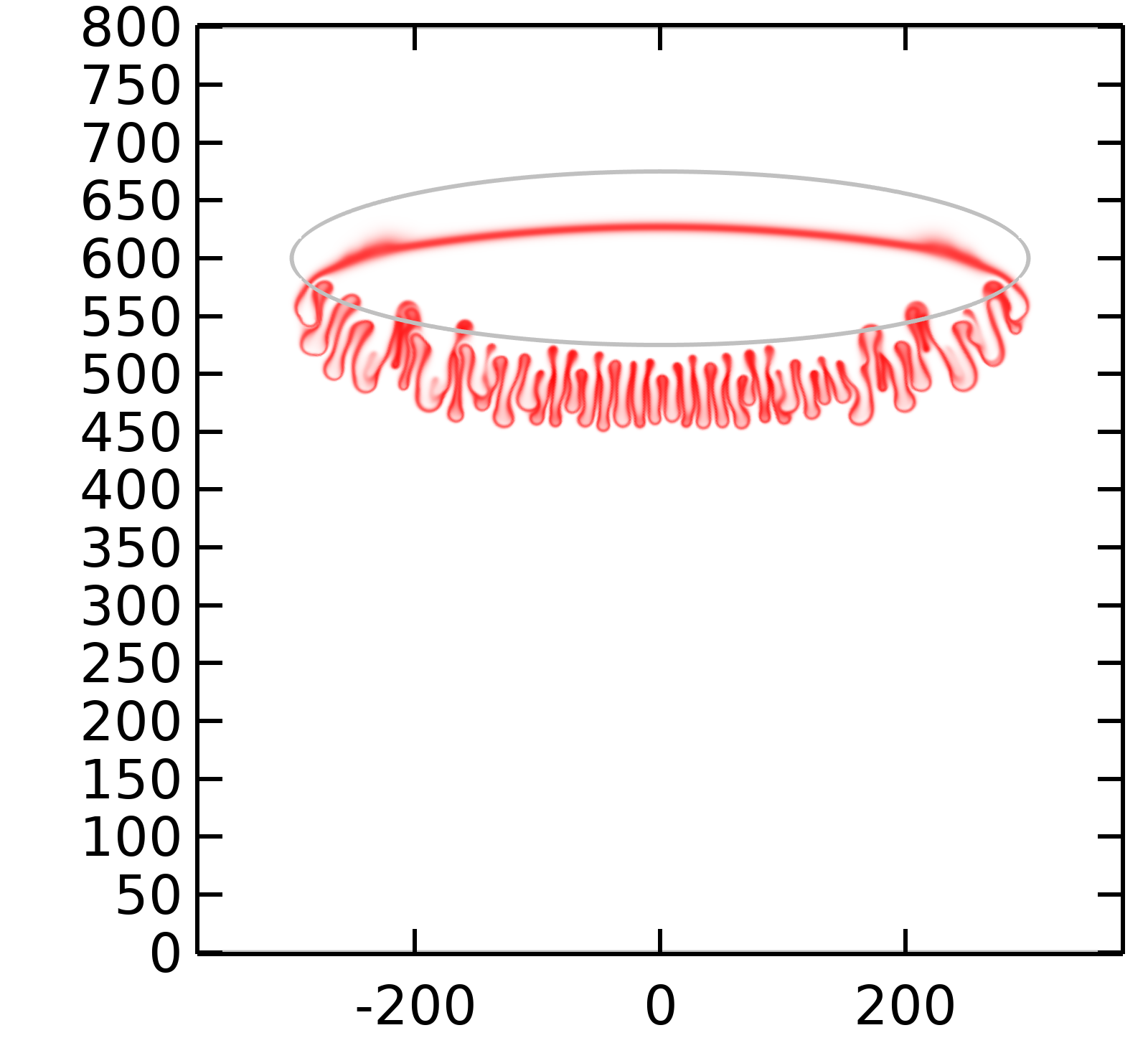}
    \includegraphics[trim=31 14 0 0,clip,width=0.192\textwidth]{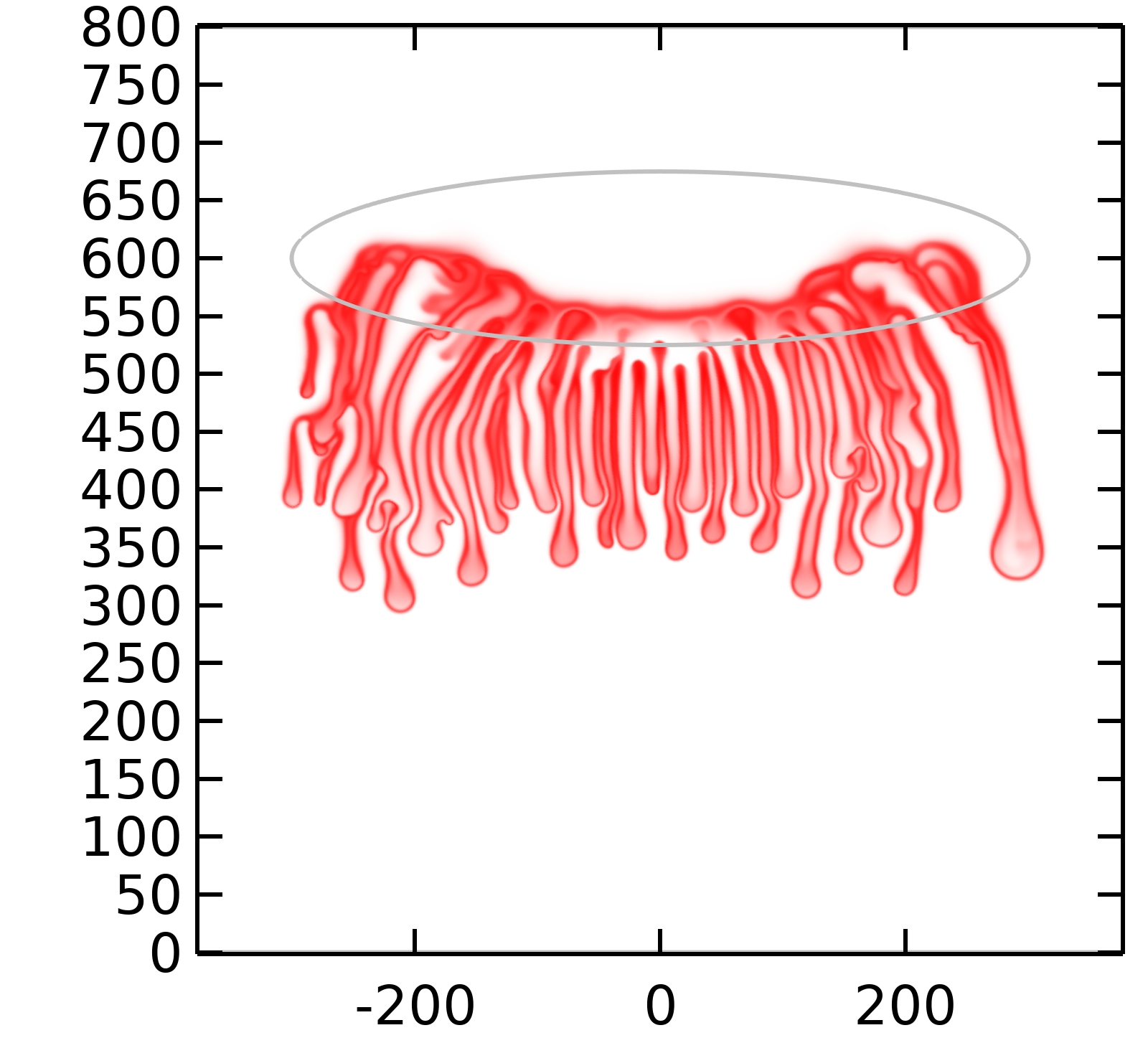}
    \includegraphics[trim=31 14 0 0,clip,width=0.192\textwidth]{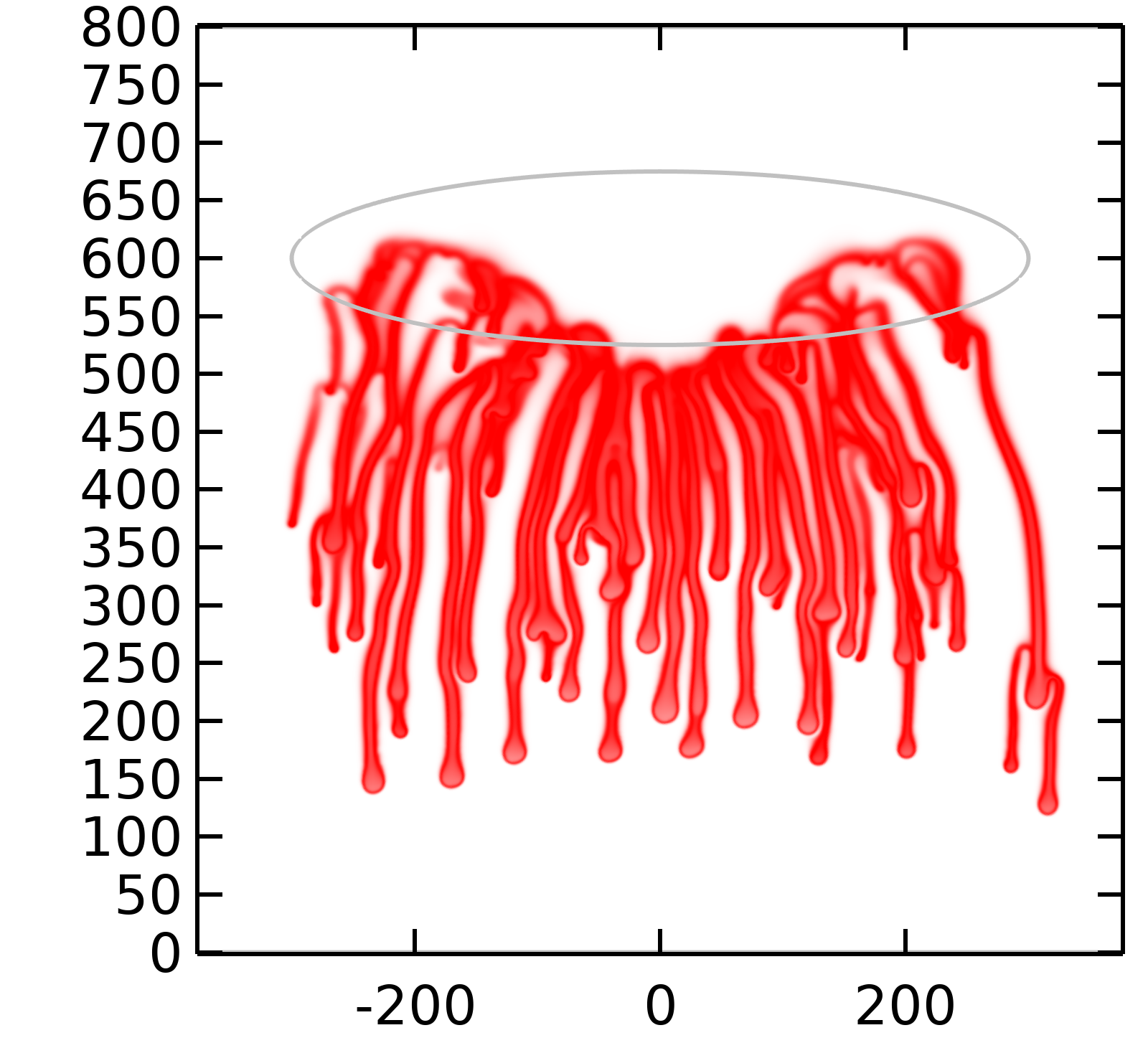}
     \includegraphics[trim=31 14 0 0,clip,width=0.192\textwidth]{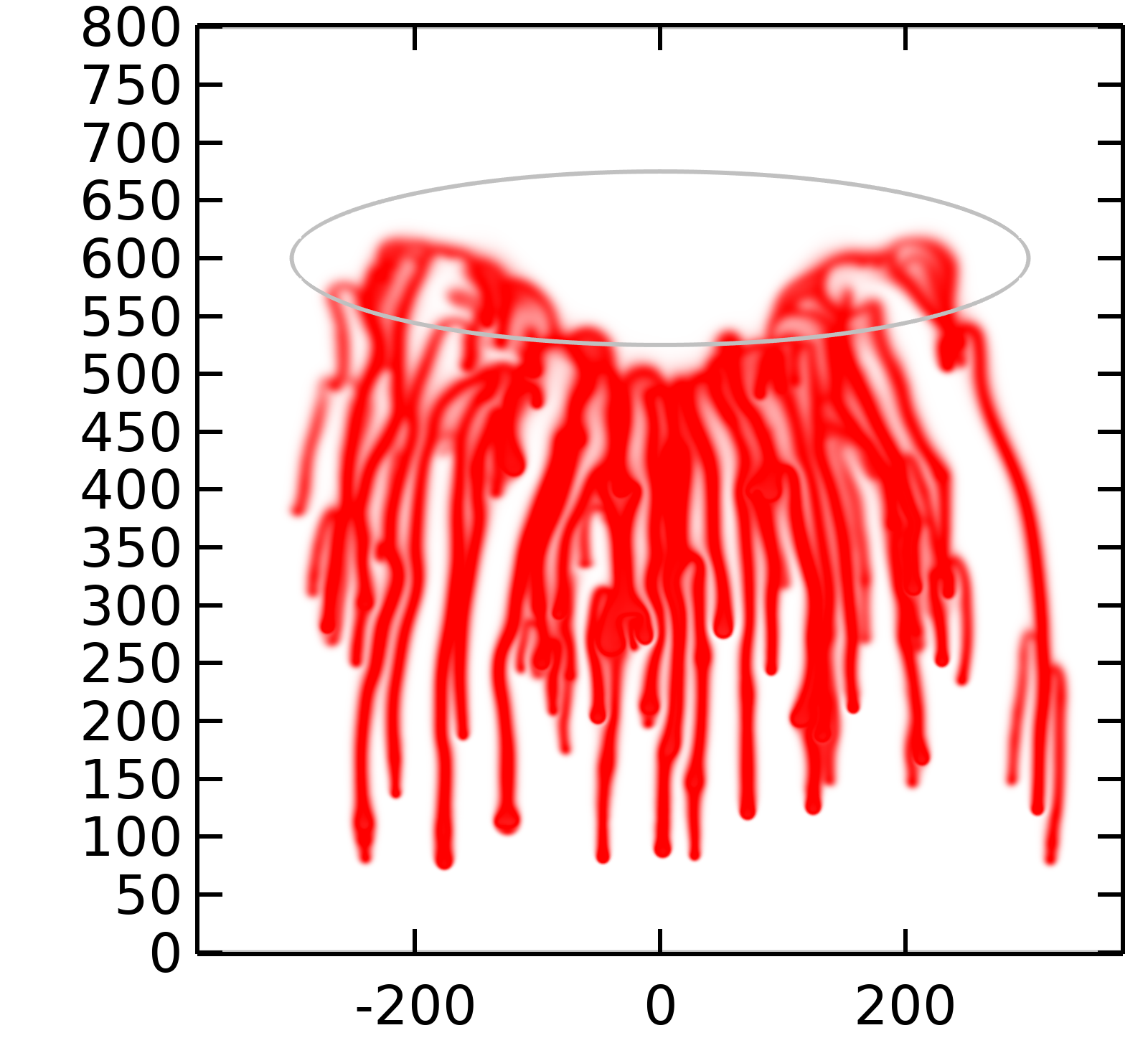}
    \begin{picture}(0,0)
    \put(-460,26){\makebox(0,0)[]{case-III}}
    \put(-460,16){\makebox(0,0)[]{\scriptsize $(R_{A}=2,R_{C}=0,\alpha=2)$}}
   \end{picture}

    \includegraphics[trim=31 14 0 0,clip,width=0.192\textwidth]{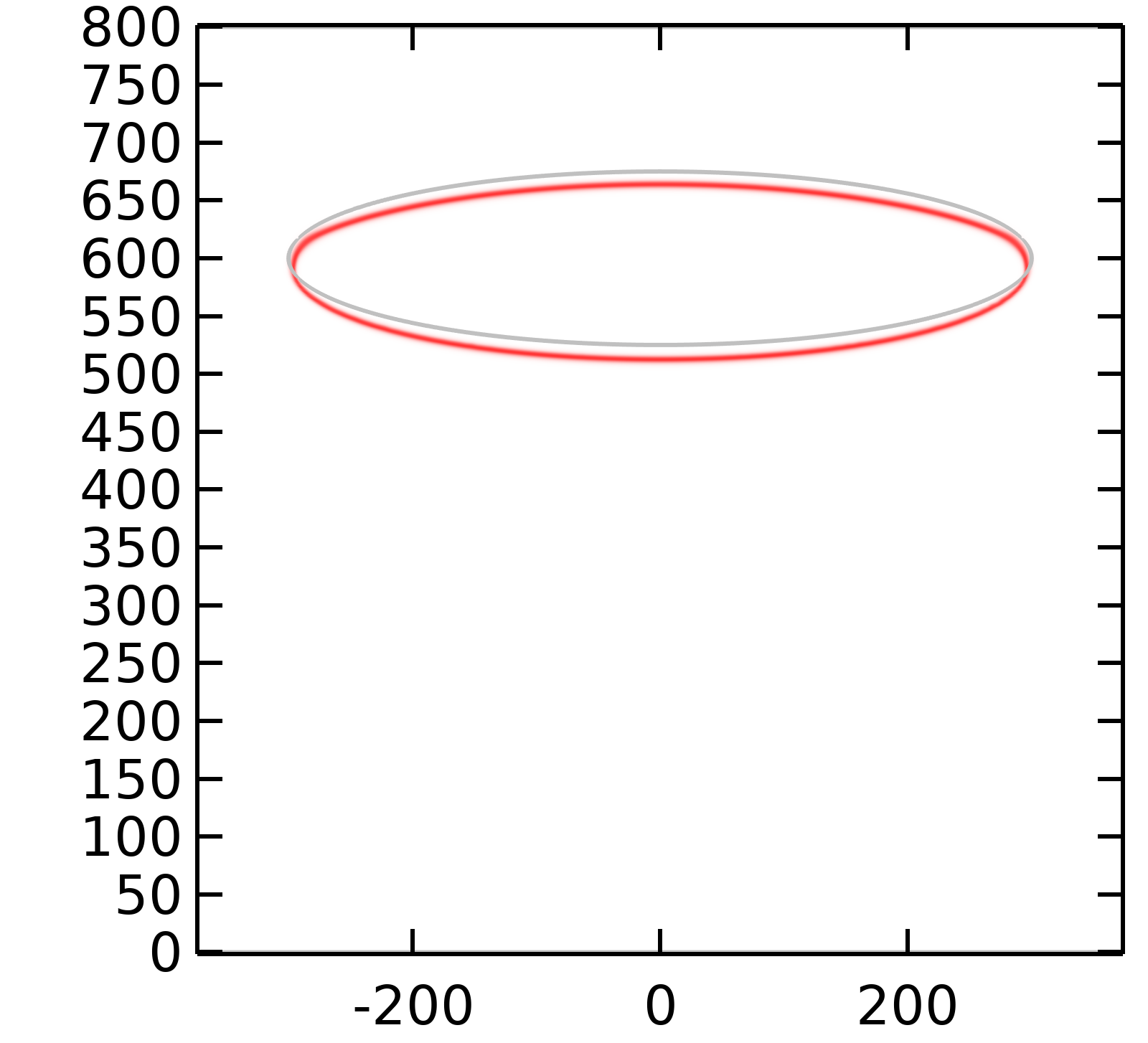}
    \includegraphics[trim=31 14 0 0,clip,width=0.192\textwidth]{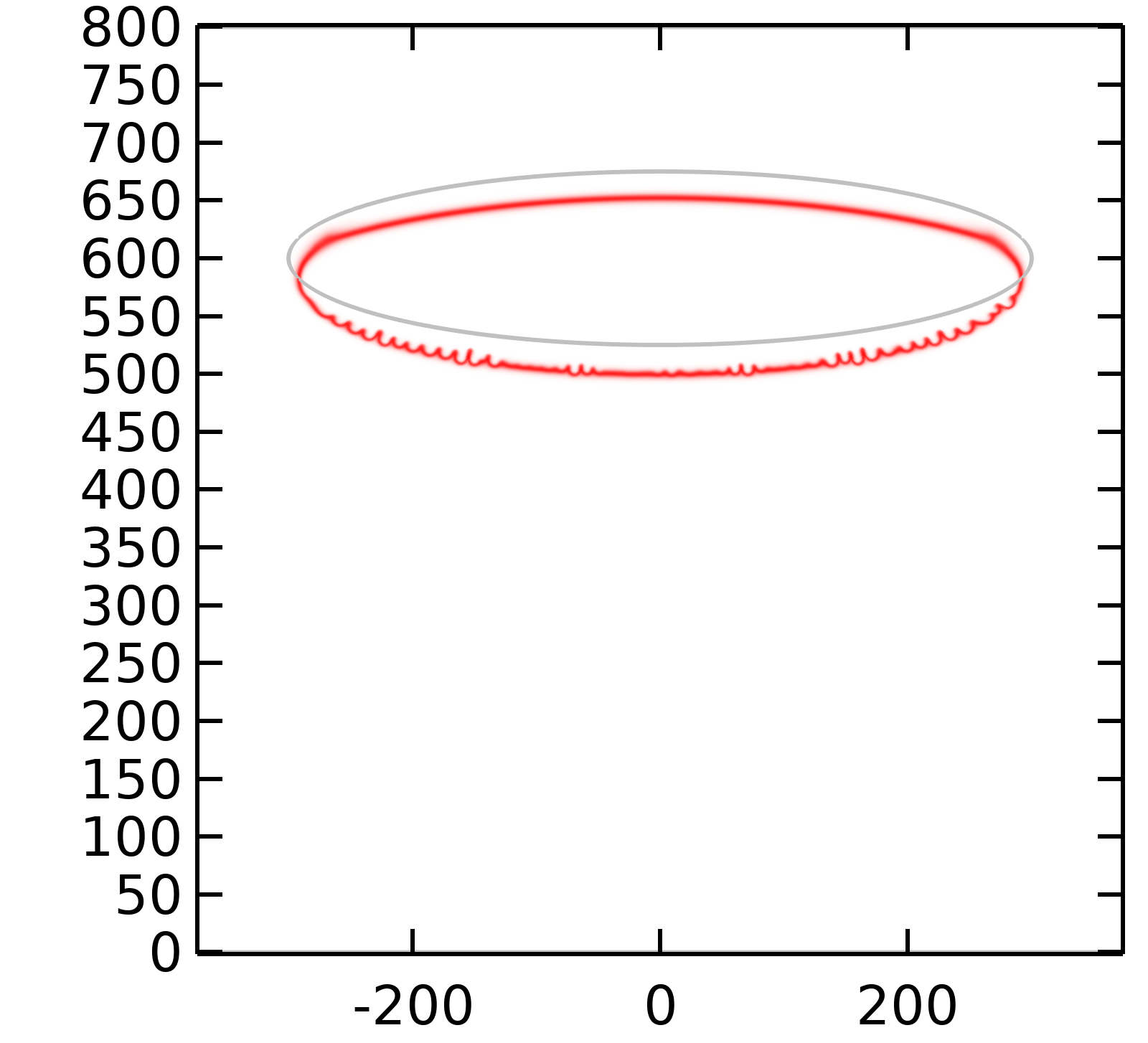}
    \includegraphics[trim=31 14 0 0,clip,width=0.192\textwidth]{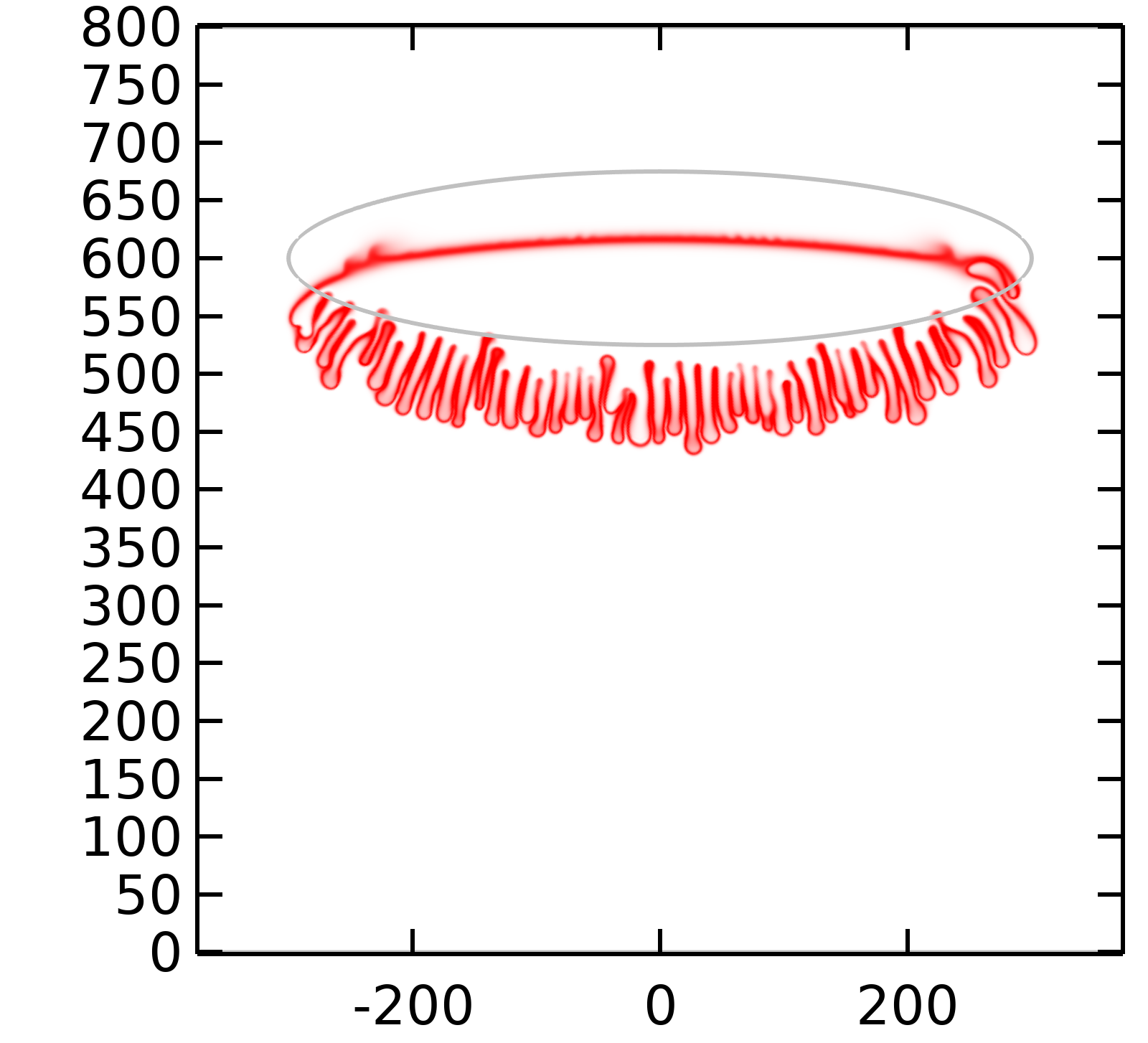}
    \includegraphics[trim=31 14 0 0,clip,width=0.192\textwidth]{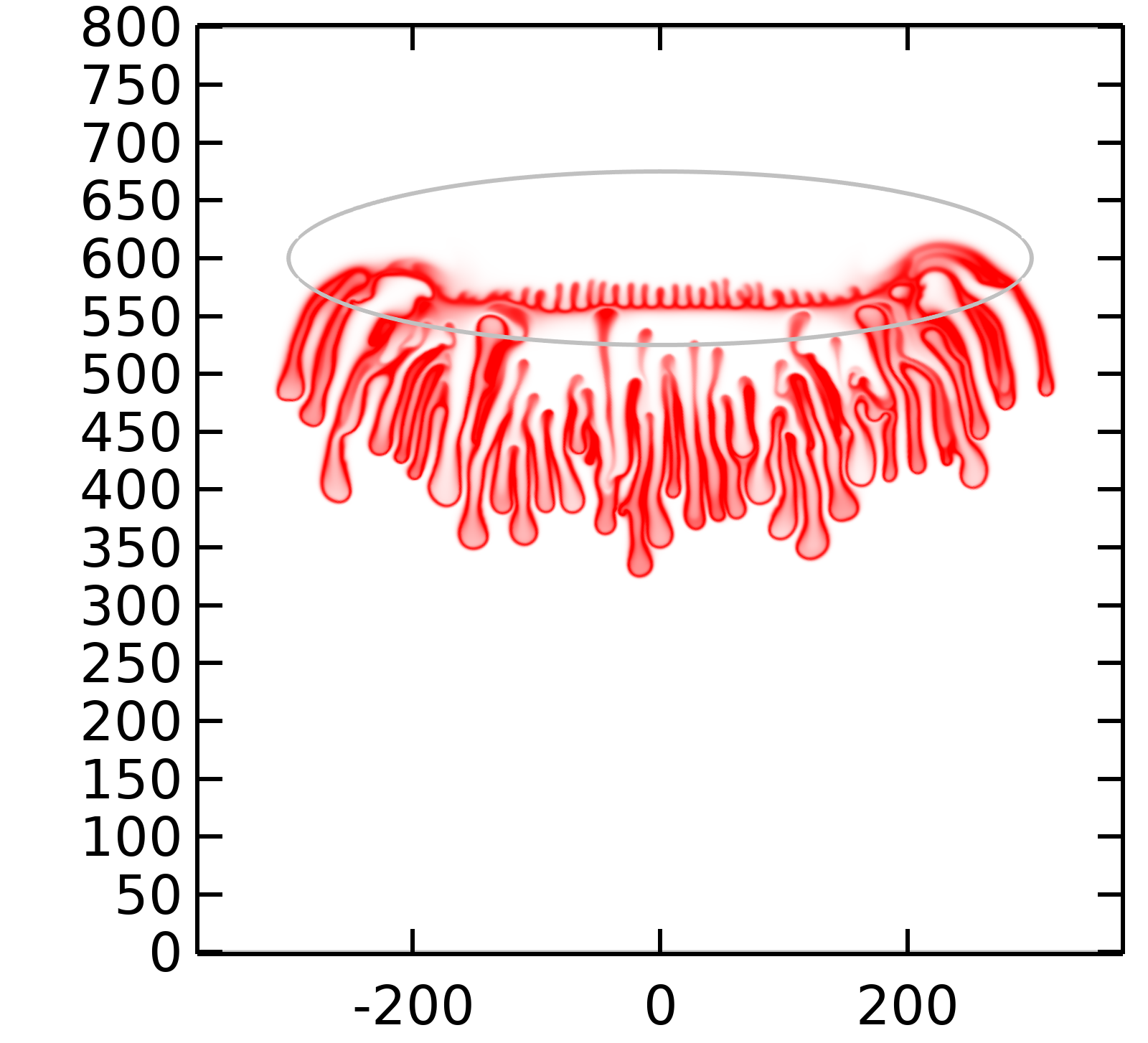}
     \includegraphics[trim=31 14 0 0,clip,width=0.192\textwidth]{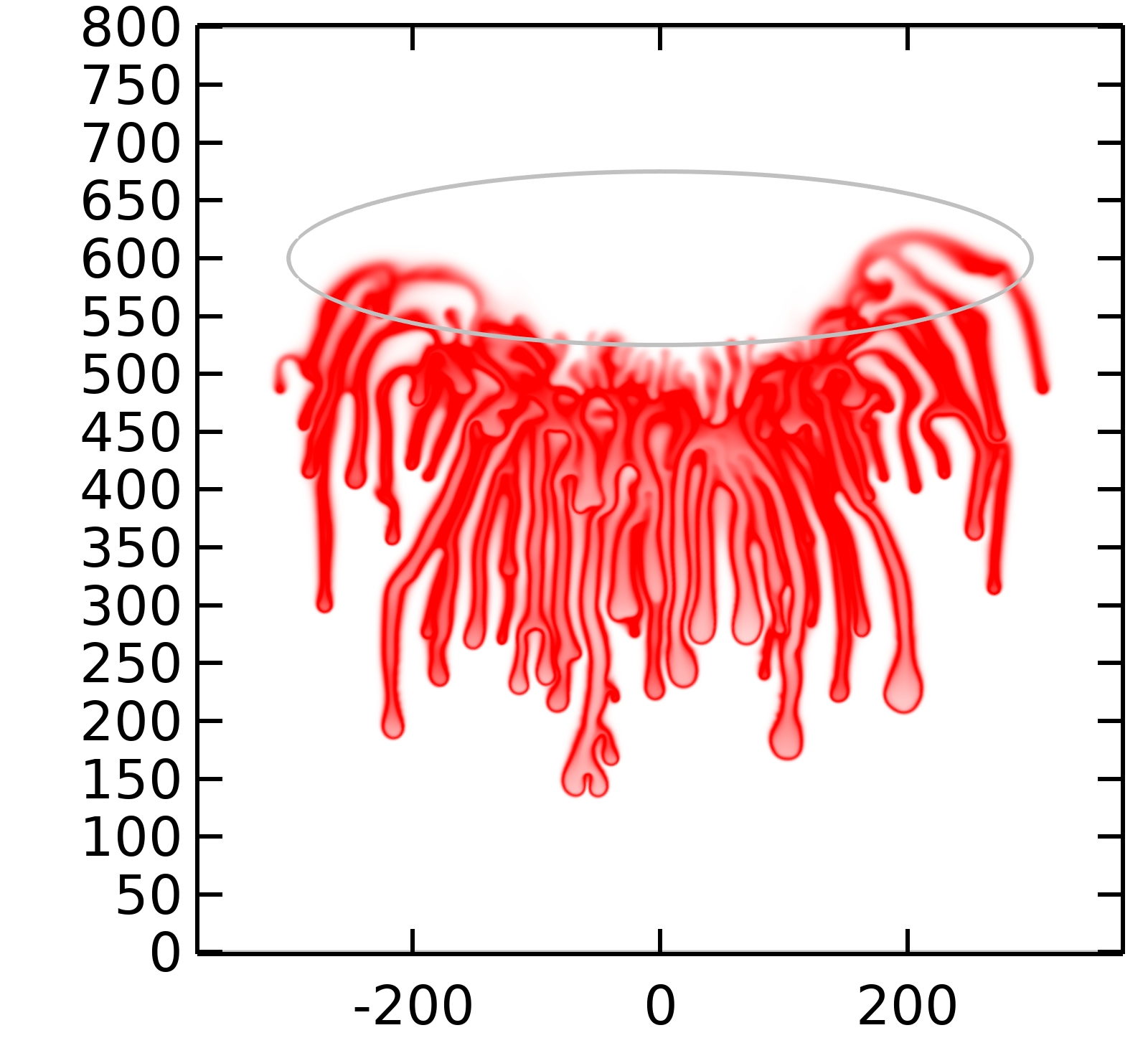}
     \begin{picture}(0,0)
    \put(-460,26){\makebox(0,0)[]{case-IV}}
    \put(-460,16){\makebox(0,0)[]{\scriptsize $(R_{A}=1,R_{C}=0,\alpha=4)$}}
   \end{picture}

    \includegraphics[trim=31 14 0 0,clip,width=0.192\textwidth]{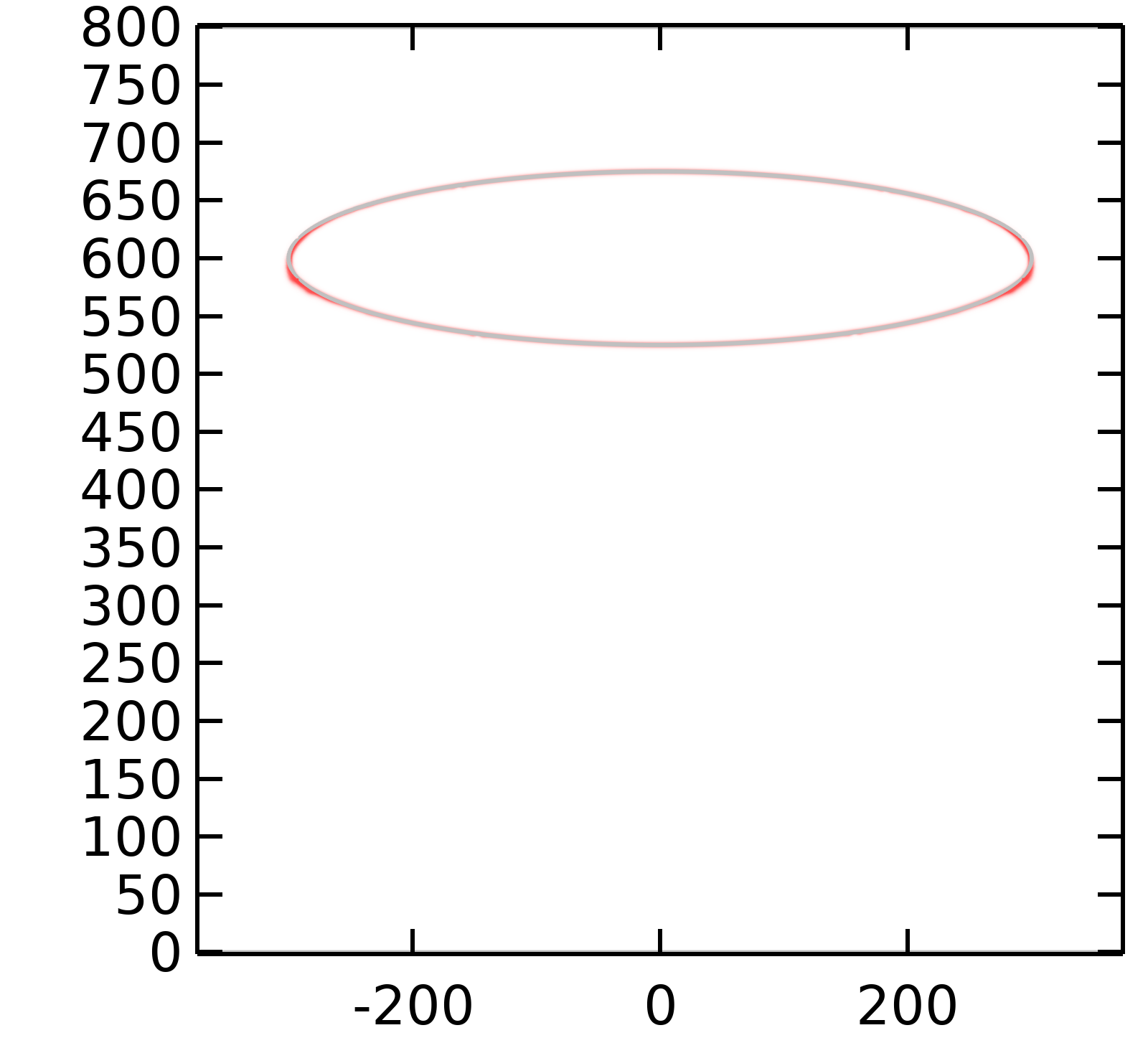}
    \includegraphics[trim=31 14 0 0,clip,width=0.192\textwidth]{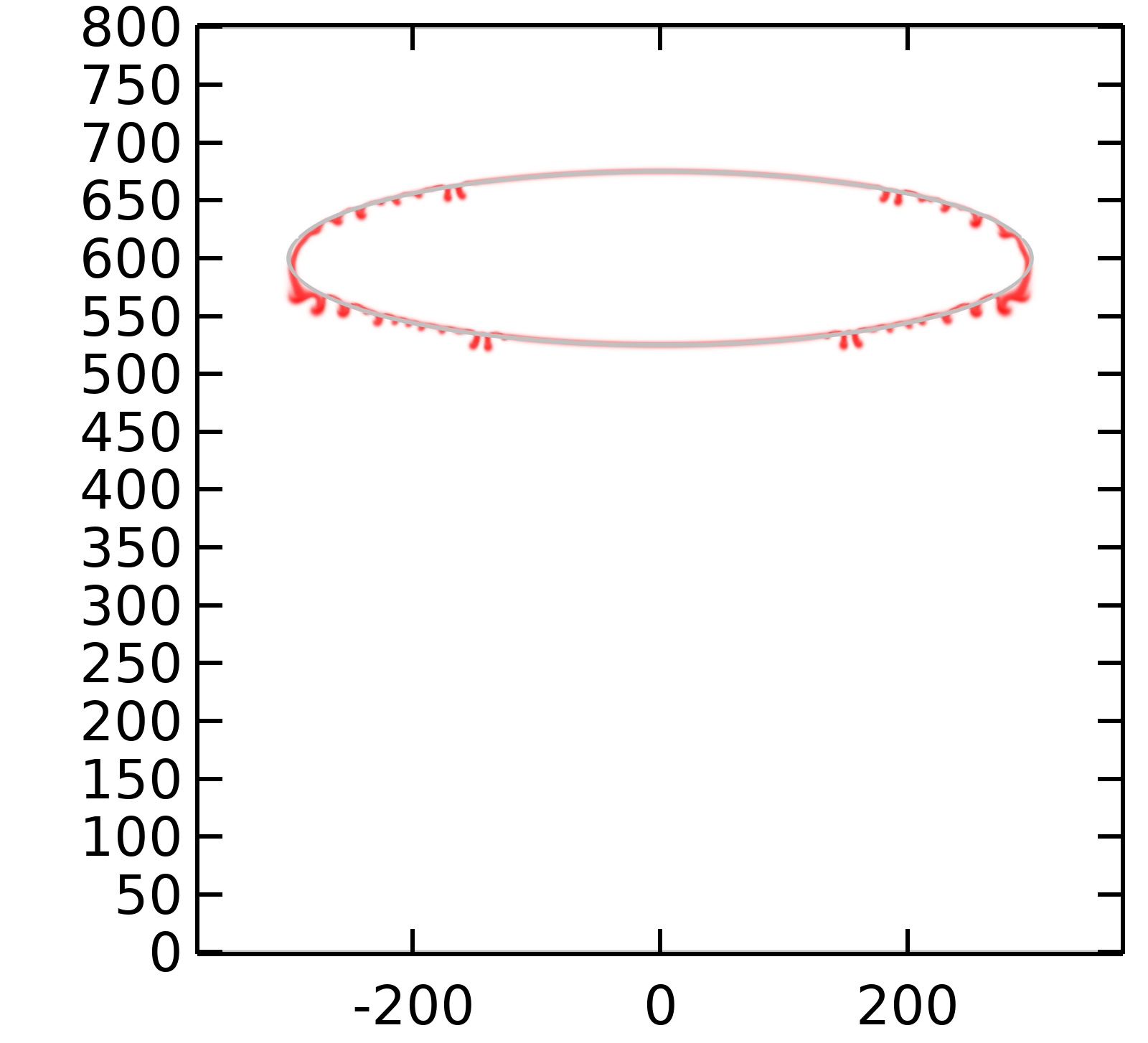}
    \includegraphics[trim=31 14 0 0,clip,width=0.192\textwidth]{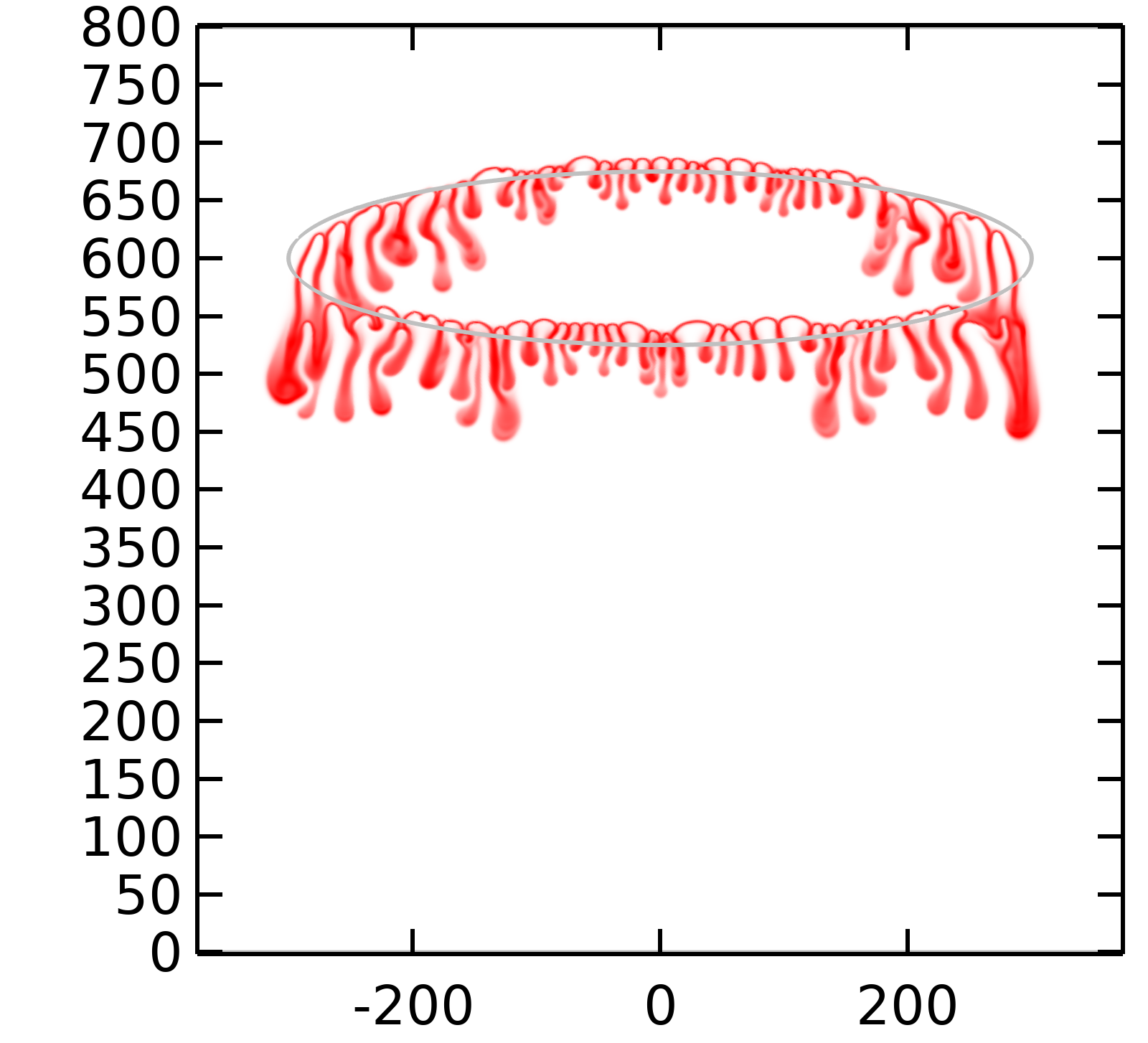}
    \includegraphics[trim=31 14 0 0,clip,width=0.192\textwidth]{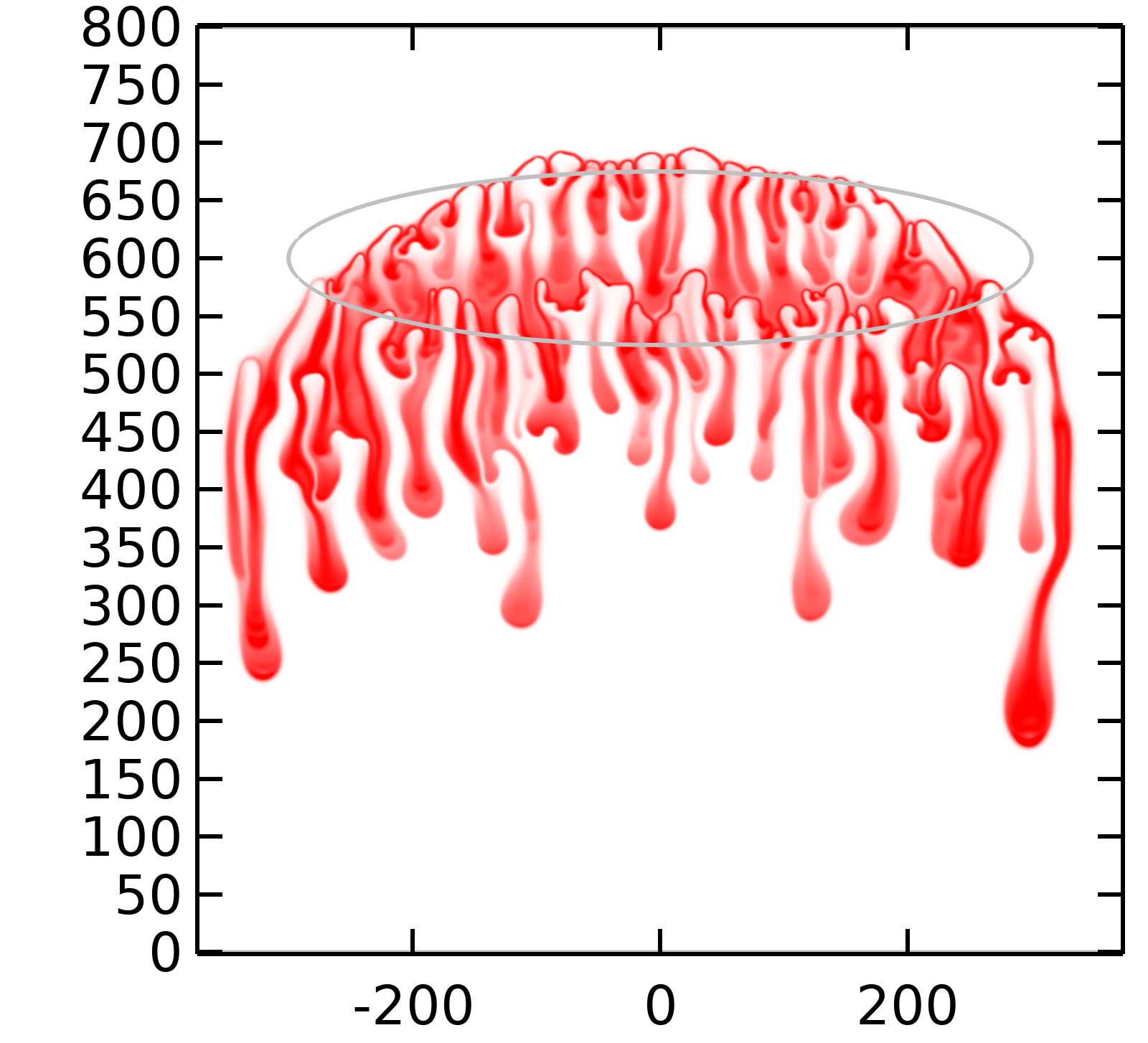}
     \includegraphics[trim=31 14 0 0,clip,width=0.192\textwidth]{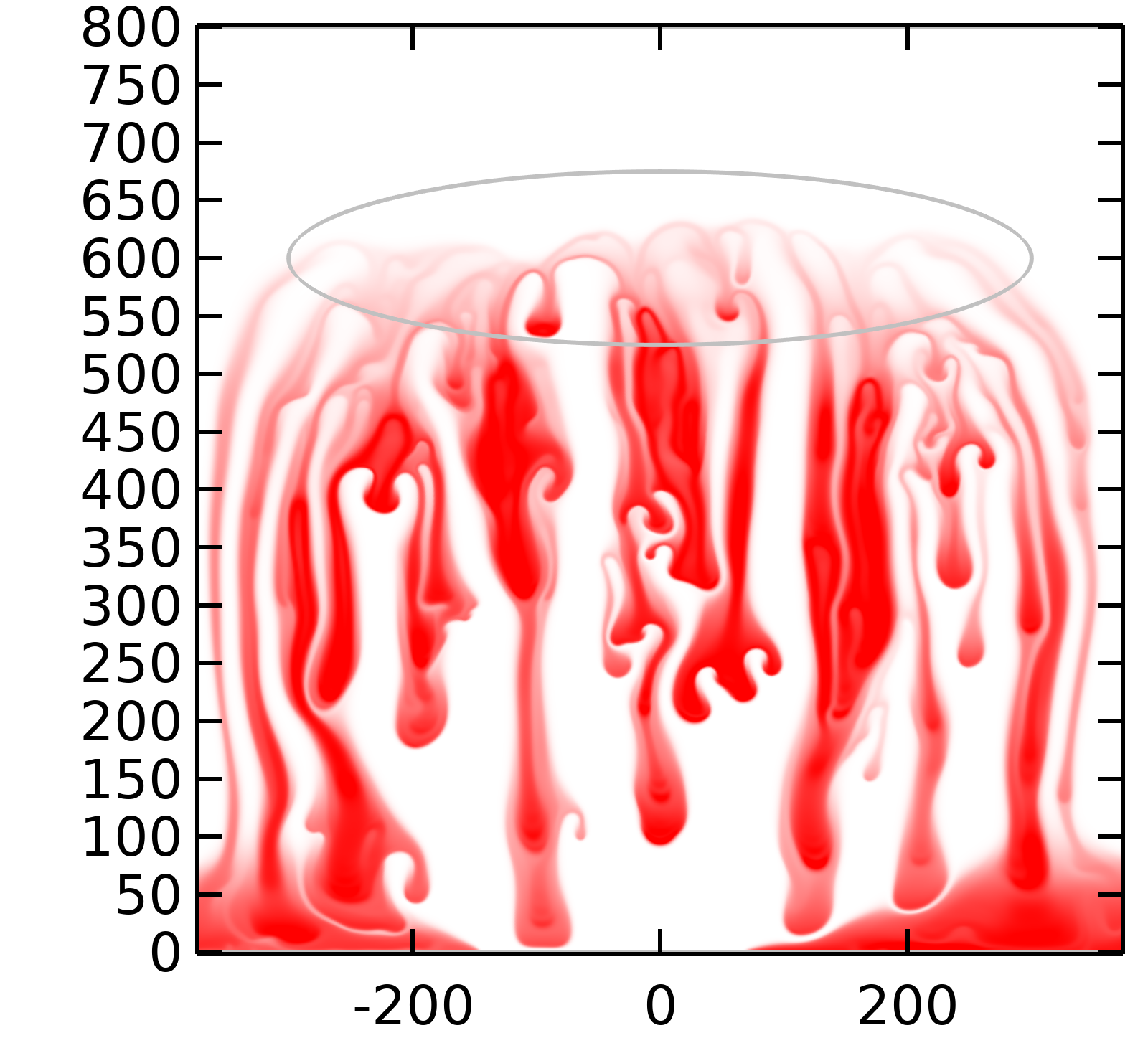}
    \begin{picture}(0,0)
    \put(-460,26){\makebox(0,0)[]{case-V}}
    \put(-460,16){\makebox(0,0)[]{\scriptsize $(R_{A}=0,R_{C}=2,\alpha=0)$}}
   \end{picture}

    \includegraphics[trim=31 14 0 0,clip,width=0.192\textwidth]{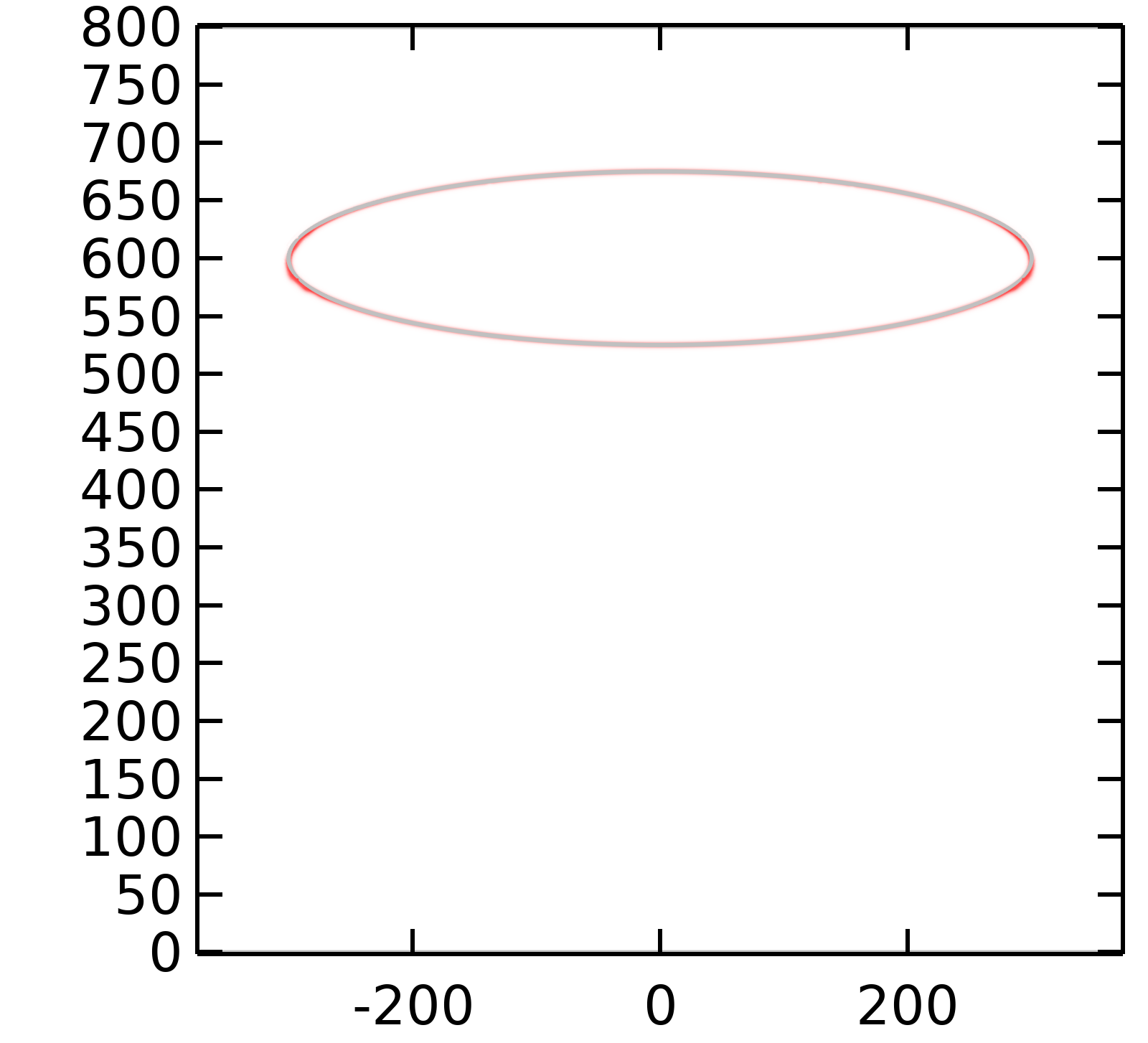}
    \includegraphics[trim=31 14 0 0,clip,width=0.192\textwidth]{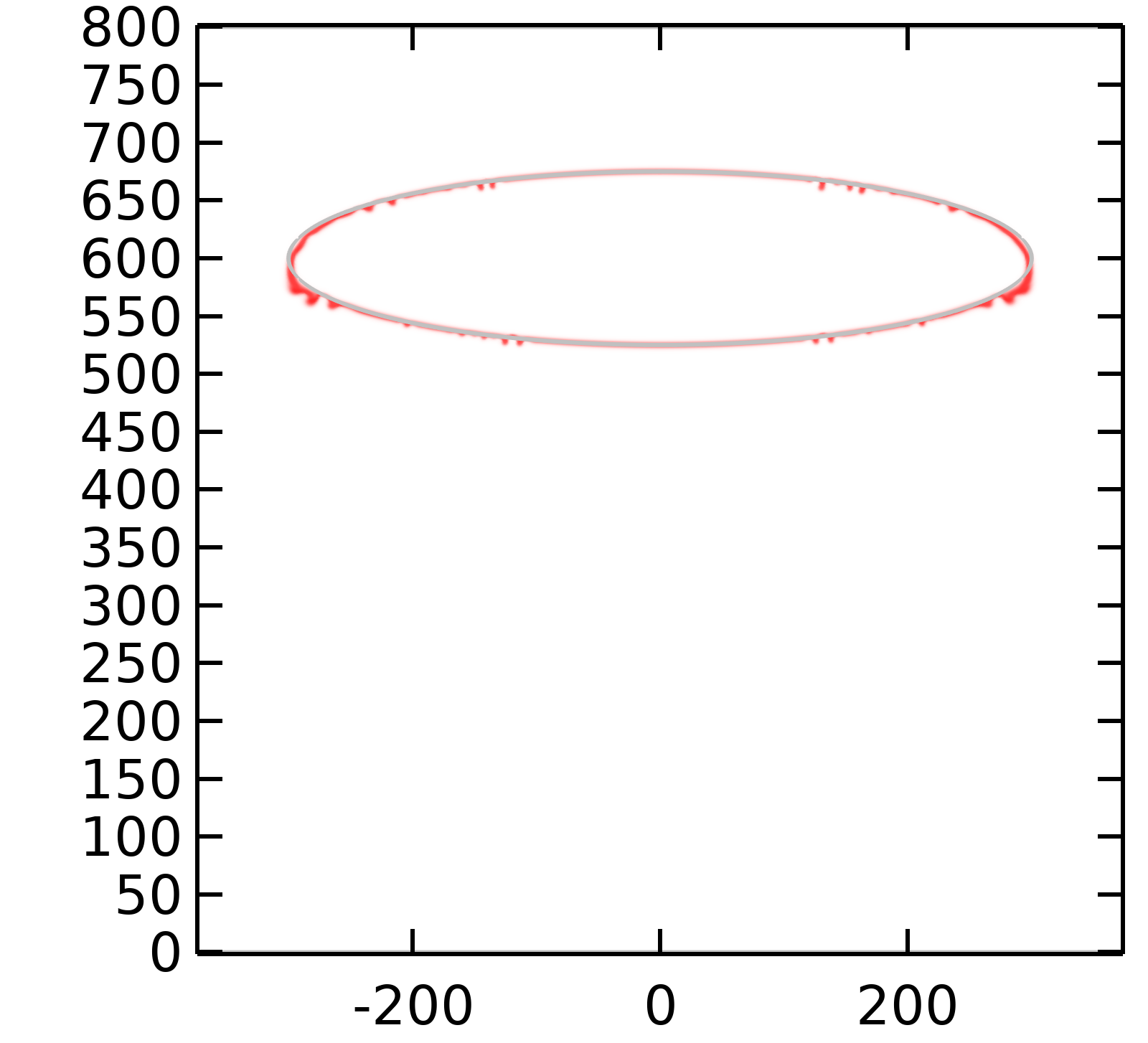}
    \includegraphics[trim=31 14 0 0,clip,width=0.192\textwidth]{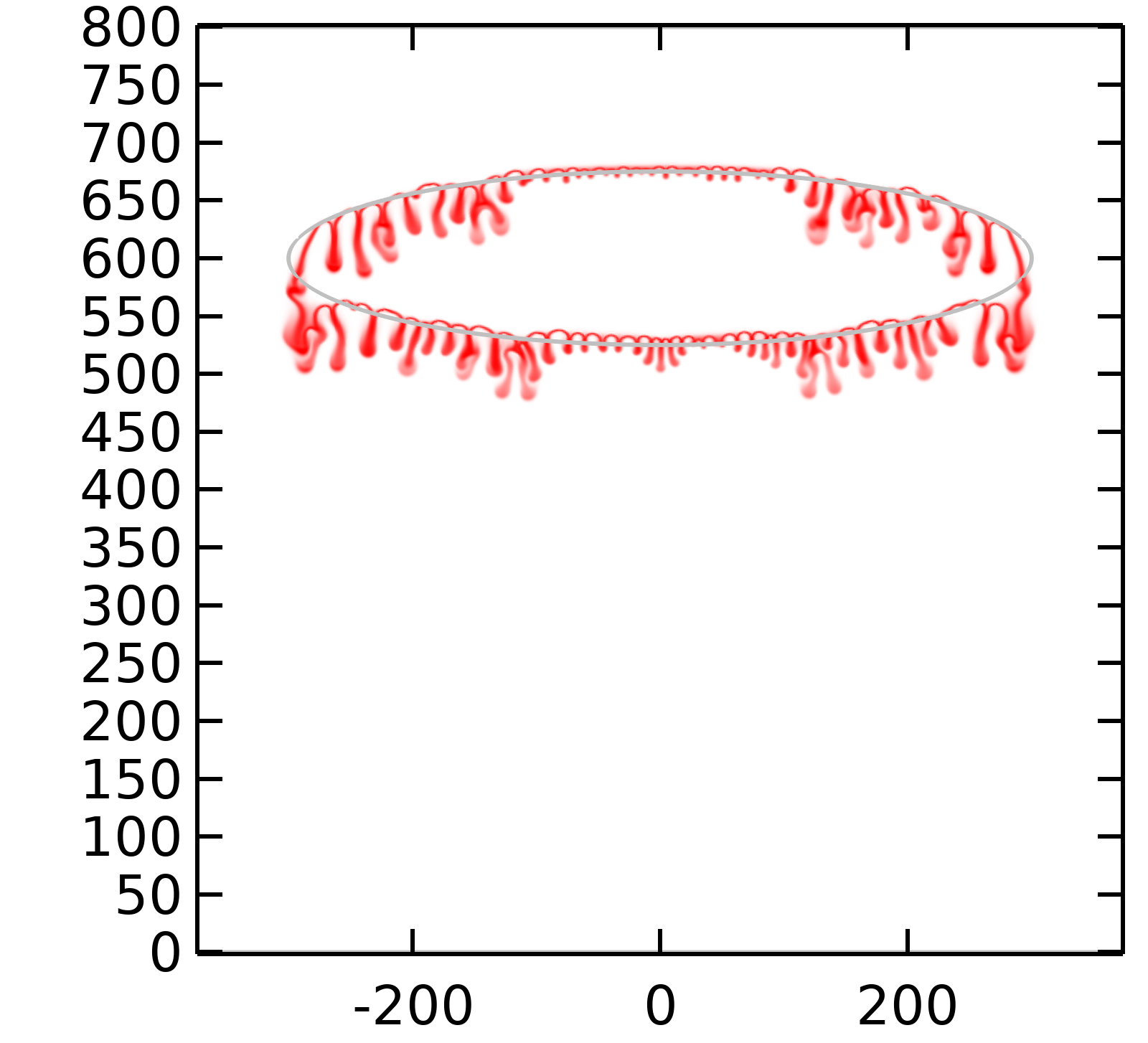}
    \includegraphics[trim=31 14 0 0,clip,width=0.192\textwidth]{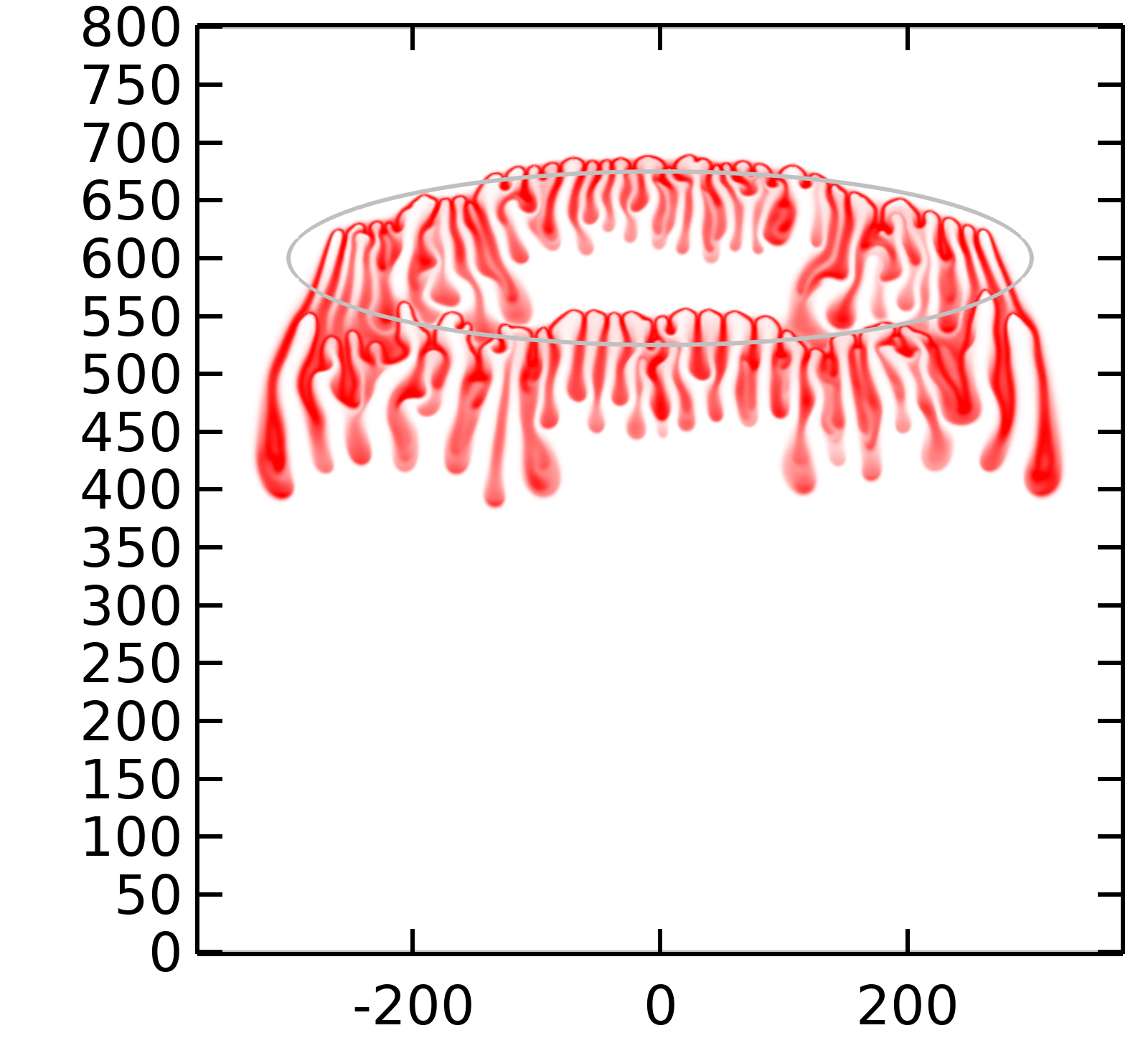}
     \includegraphics[trim=31 14 0 0,clip,width=0.192\textwidth]{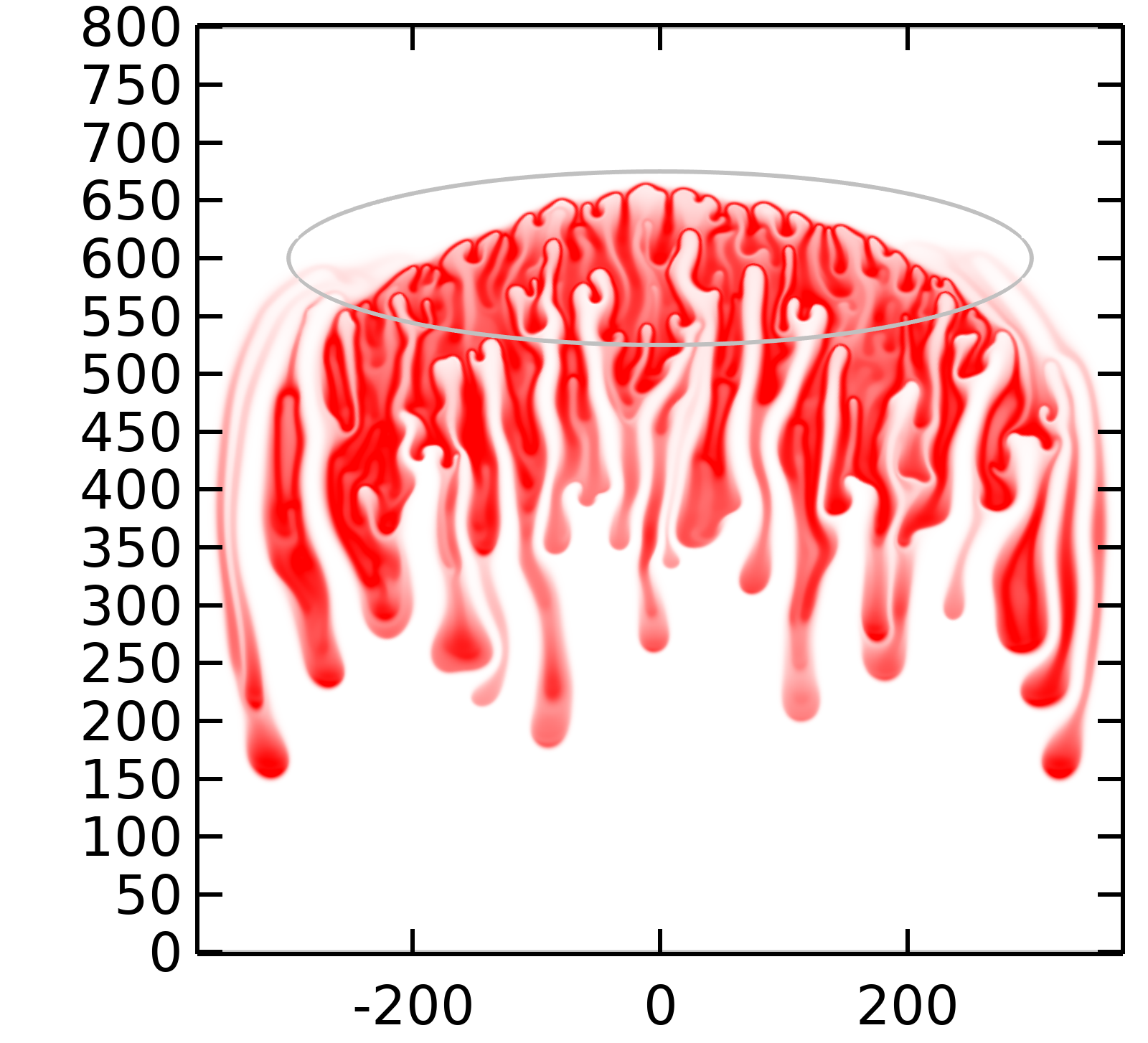}
    \begin{picture}(0,0)
    \put(-460,26){\makebox(0,0)[]{case-VI}}
    \put(-460,16){\makebox(0,0)[]{\scriptsize $(R_{A}=0,R_{C}=2,\alpha=2)$}}
   \end{picture}
  
    \begin{picture}(0,0)
    \put(-106,622){\makebox(0,0)[]{$t=100$}}
    \put(-5,622){\makebox(0,0)[]{$t=200$}}
    \put(96,622){\makebox(0,0)[]{$t=500$}}
    \put(195,622){\makebox(0,0)[]{$t=1000$}}
    \put(295,622){\makebox(0,0)[]{$t=2000$}}
   \end{picture}
    \hspace{0.38 in}\includegraphics[trim=31 10 0 140,clip,width=0.3\textwidth]{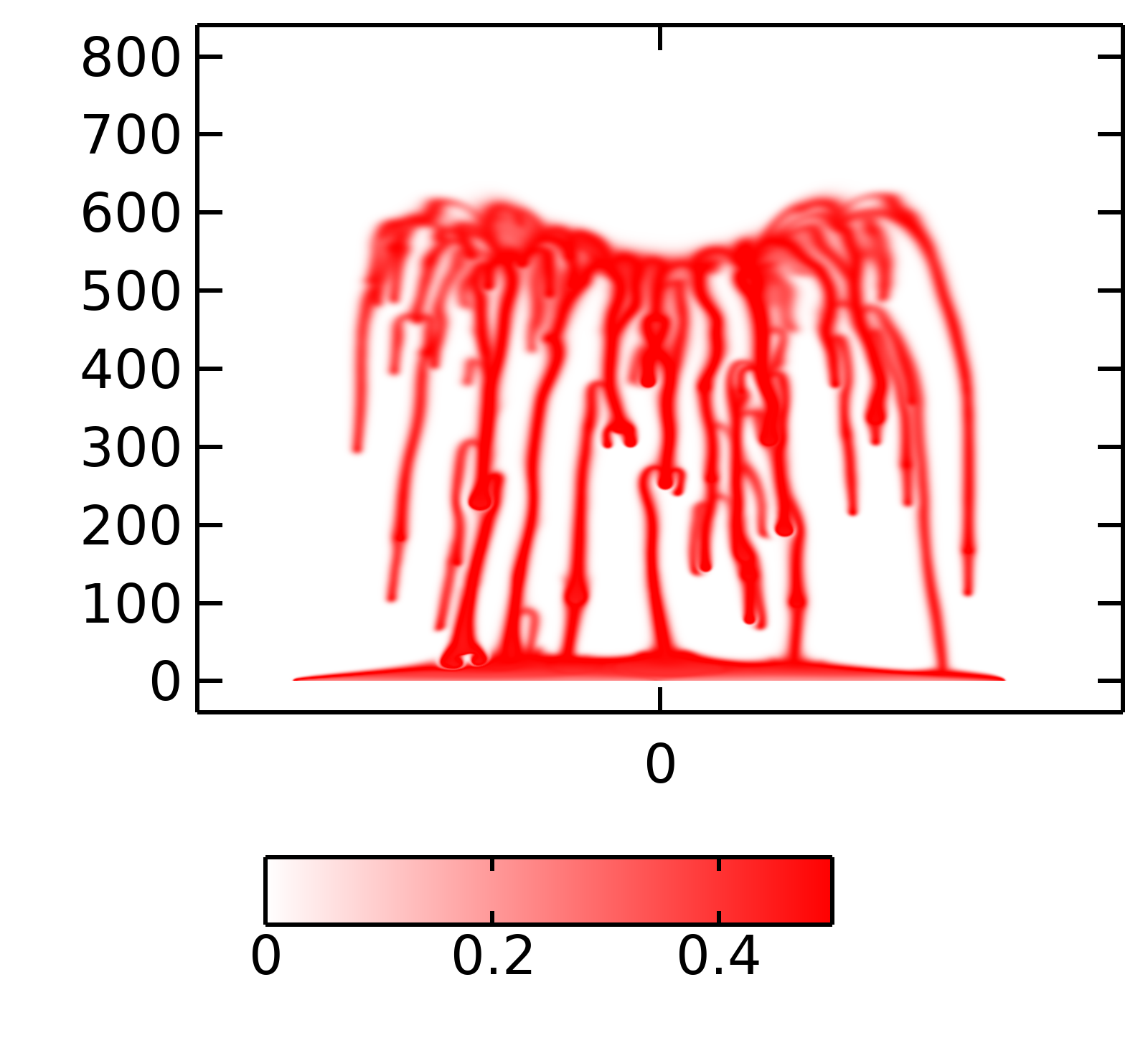}
   \vspace{-0.2 cm}\caption{Spatio-temporal evolution of the reaction product concentration profile \( c \) for cases I to VI, shown from the top to bottom panels, respectively. The axis labels are same as Fig. \ref{fig:abc}.}
    \label{fig:el}
\end{figure}

\subsection{2D initial elliptic interface}
Next, we examine the dependence of fingering instability on the initial elliptical interface, demonstrating that, for the problem of interest, the solution depends continuously on the initial data, thereby supporting our theoretical findings. The domain is initialized with an elliptic blob filled with reactant \( A \), centered at \((x, y) = (0, 600)\), with a semi-major axis of \( x_0 = 300 \) and a semi-minor axis of \( y_0 = 150 \). Outside the blob, the reactant \( B \) fills the rectangular porous medium. 
We chose this blob problem because it has significant applications in various scientific and engineering fields, such as CO\(_2\) sequestration and aquifer remediation. Similar problems have been studied previously in the presence of pure Darcy flow without reaction by \cite{Jha_2023, Pramanik_2015, SaloSalgado2024} and the references therein. 
The initial conditions are given as:
\begin{align}\label{e72}
a_0(\boldsymbol{x}) &= 
\begin{cases} 
1, & \text{if } \frac{x^2}{x_0^2} + \frac{(y-600)^2}{y_0^2} \leq 1 \\
0, & \text{otherwise}
\end{cases},~
b_0(\boldsymbol{x}) &= 
\begin{cases} 
0, & \text{if } \frac{x^2}{x_0^2} + \frac{(y-600)^2}{y_0^2} \leq 1 \\ 
1, & \text{otherwise} 
\end{cases},~
c_0(\boldsymbol{x}) &= 0.
\end{align}   
The reaction modifies the flow dynamics based on the parameters discussed in the six previously categorized cases. The initial gravitational and velocity fields are the same as those used in the flat interface case with the same boundary conditions as in \eqref{bc's}. Due to curvature at the interface, additional corner mesh refinement is added in COMSOL near the initial interface with \( h = 0.25 \), while the rest of the geometry uses \( h = 1.5 \). This results in a degree of freedom (DoF) of \( 7615615 \) for all simulations involving the 2D initial elliptic interface case. Since the reaction occurs around the elliptic interface, two distinct interfaces are observed: \(B\)-\(C\)-\(A\) for \(y > 600\), and \(A\)-\(C\)-\(B\) for \(y < 600\) (see, for example, column for \(t=100\) in Fig. \ref{fig:el}). The density fingering patterns at these interfaces exhibit unique features that depend on the governing flow parameters.

Notably, in Case I, where density stratification occurs (\(R_A > 0\)), two distinct interfaces are observed. The upper interface (\(y > 600\)) moves downward in the frame of the falling denser elliptic bulb; however, locally, no density stratification is observed at this interface as it does not interact with the less dense reactant \(B\) ($t=100$, panel for Case-I in Fig. \ref{fig:el}). In contrast, the lower interface (\(y < 600\)) is unstable, exhibiting density fingering as the denser fluid penetrates into the less dense reactant \(B\). Reaction products form at the lower interface, and the interaction between the two interfaces can be traced through the concentration of these products. 
In Case II, no density stratification occurs, and the dynamics are dominated by reaction-diffusion spreading. In case-III, the permeability change slows the descent of the bulb, preventing the fingers from reaching the lower boundary even at \(t = 2000\), unlike in Case I (compare panels of the first and third row in Fig. \ref{fig:el}). Interestingly, in Case IV, upward-pointing fingers are observed at \(t = 1000\) due to significant local permeability changes caused by the reaction products. This demonstrates that permeability changes not only resist downward motion but can also drive upward finger formation in this reactive flow system.

When there is no initial density stratification, and density stratification is induced solely by the formation of a denser product due to reaction, the pattern formation becomes even more striking. In this scenario, the product \(C\), being denser, causes both initial interfaces, \(B\)-\(C\)-\(A\) for \(y > 600\) and \(A\)-\(C\)-\(B\) for \(y < 600\), to exhibit density fingering (see panels for Case V and Case VI in Fig. \ref{fig:el}). The wavelength of these reaction-induced fingers is narrower compared to earlier cases. At \(t = 1000\), the fingers in Case V do not reach the lower boundary as the natural Rayleigh-Taylor fingers did in Case I. However, at later times, more reaction products are deposited at the bottom wall than in Case I (compare panels for \(t = 2000\) in Case I and Case VI in Fig. \ref{fig:el}). In the last case (Case VI), the falling speed of these reaction-induced fingers is reduced by the local change of permeability ($\alpha=1$), preventing them from reaching the lower boundary even at \(t = 2000\). This makes the pattern formation in Case VI distinctly different from the other cases.

\subsection{3D reaction-induced density fingering}
\begin{figure}
    \centering
    \includegraphics[width=0.33\linewidth]{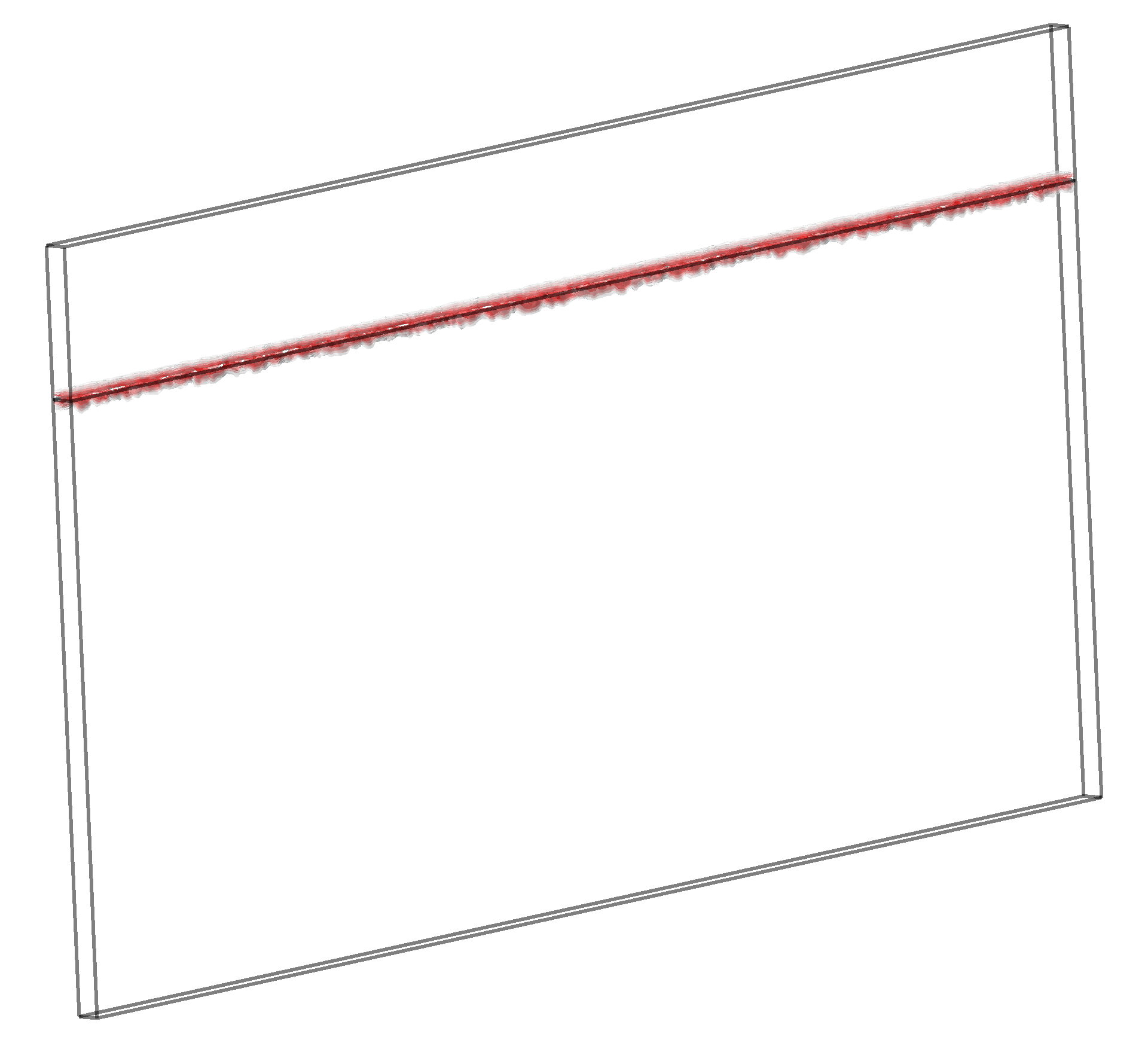}
    \includegraphics[width=0.33\linewidth]{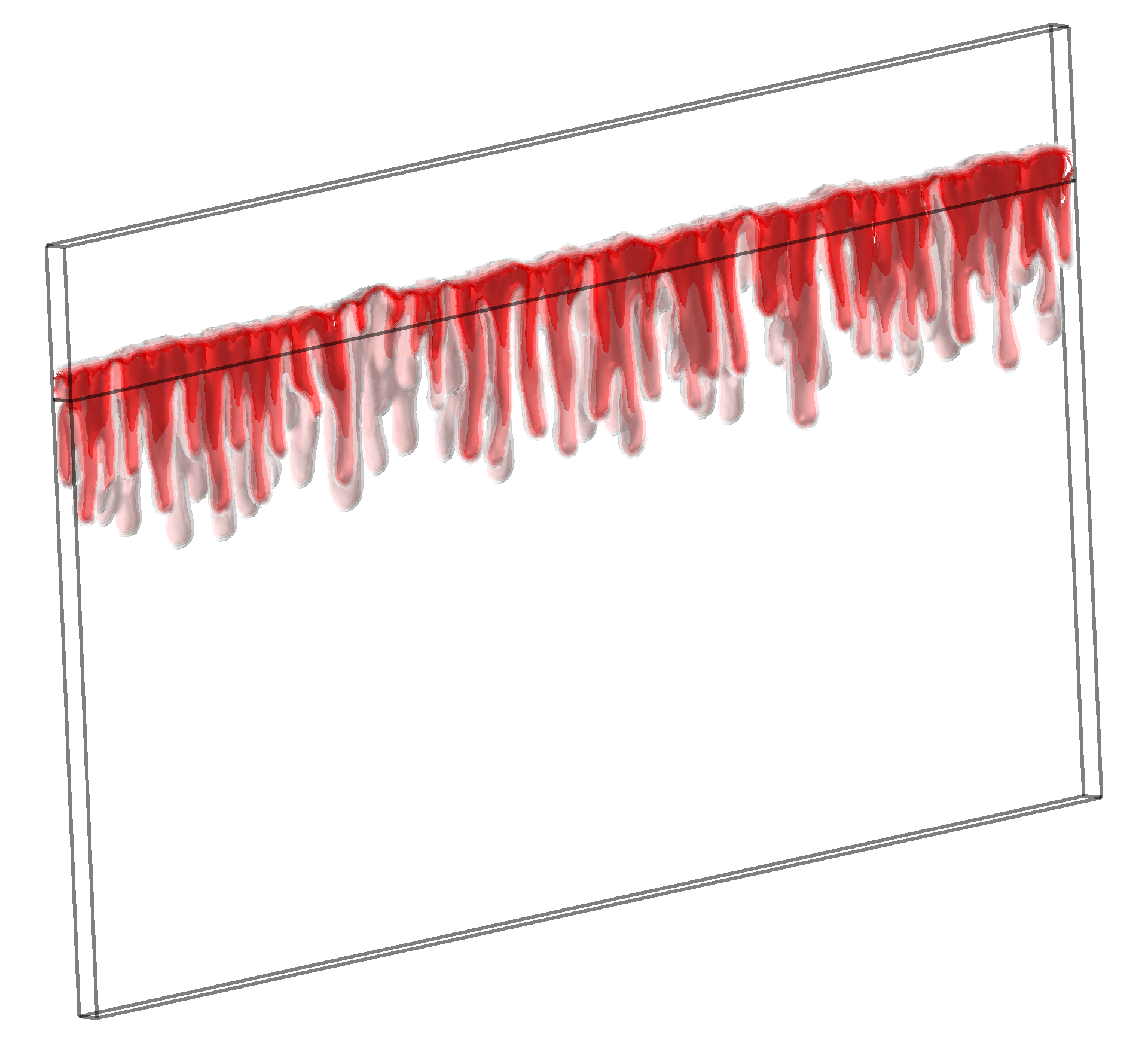}
    \includegraphics[width=0.33\linewidth]{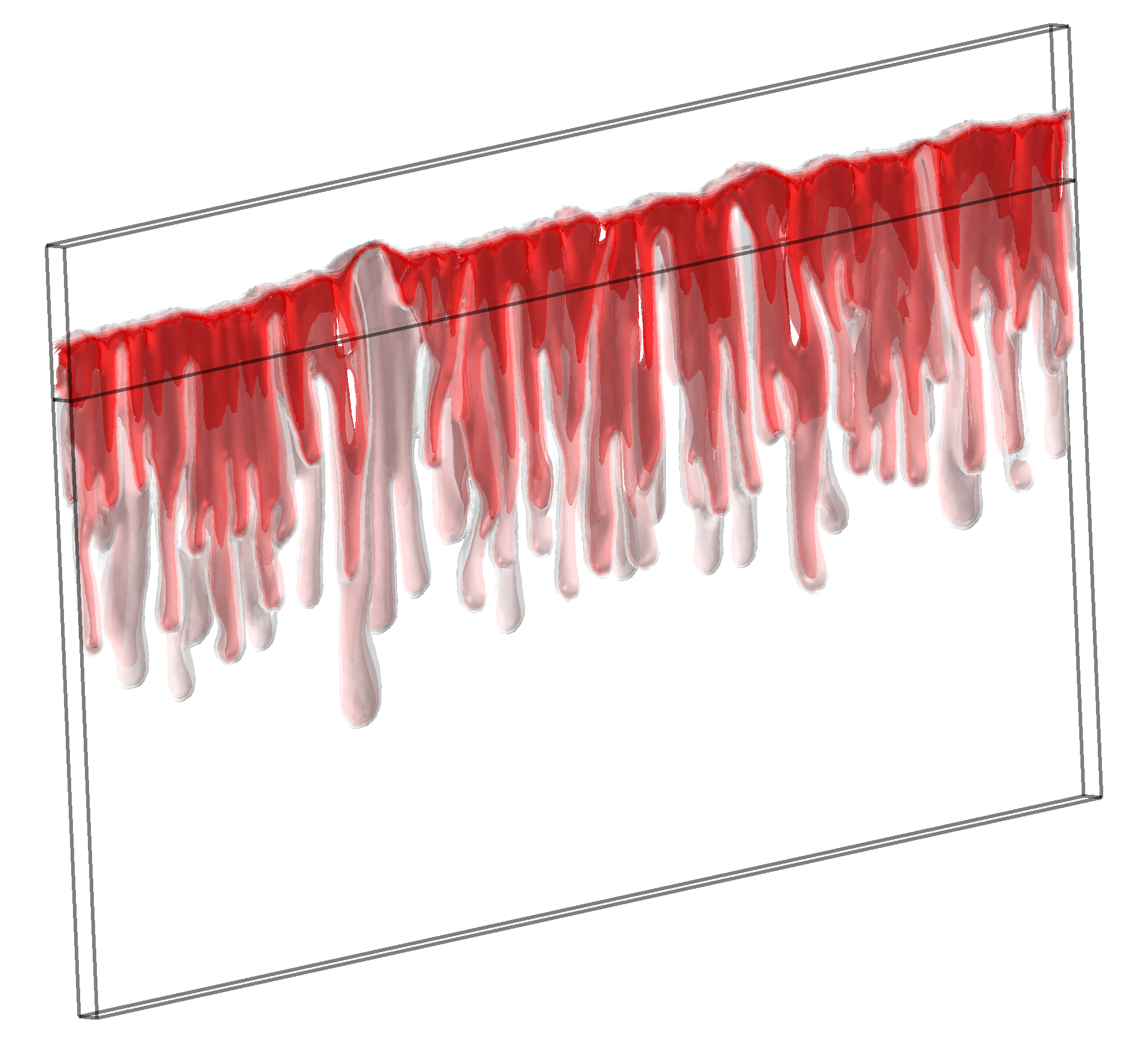}
    \includegraphics[width=0.33\linewidth]{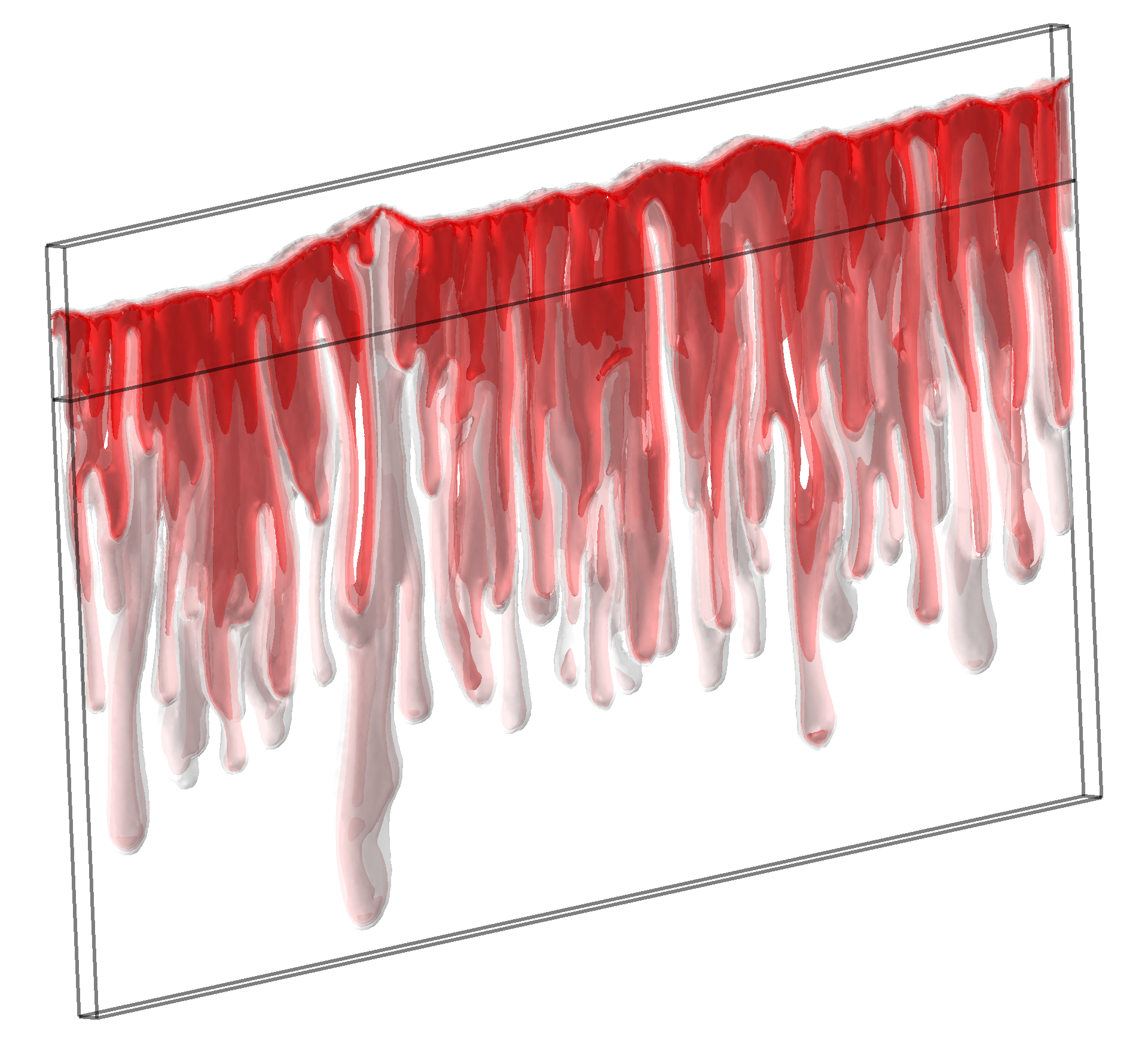}
    
    \begin{picture}(0,0)
    \put(-80,325){\makebox(0,0)[]{$t=50$}}
    \put(80,325){\makebox(0,0)[]{$t=500$}}
    \put(-80,165){\makebox(0,0)[]{$t=1000$}}
    \put(80,165){\makebox(0,0)[]{$t=1500$}}
    
    \put(170,140){\includegraphics[trim=170 22 10 38,clip, width=0.07\linewidth]{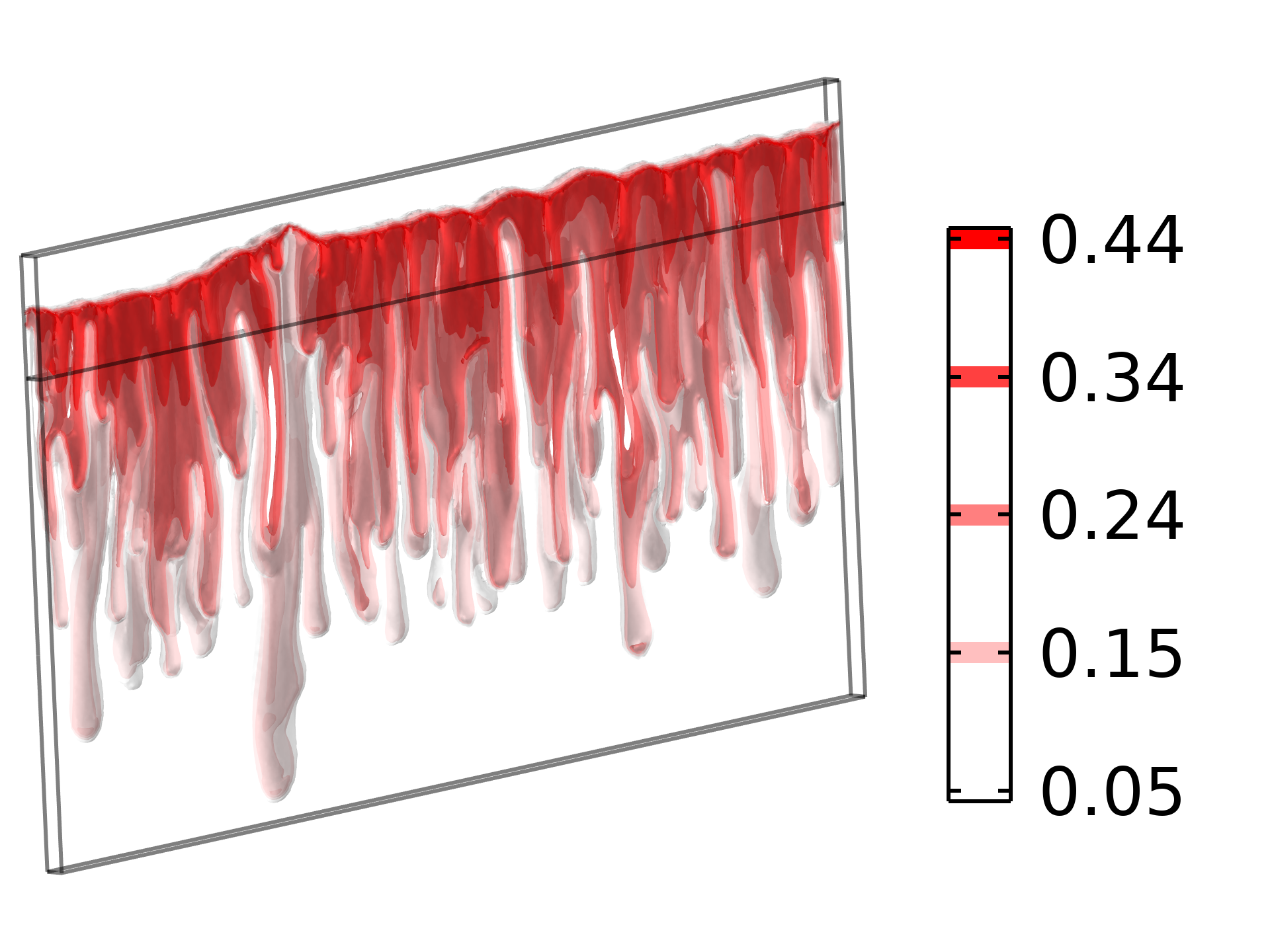}}
   \end{picture}
    \vspace{-0.7 cm}\caption{Spatio-temporal evolution of the reaction product concentration's five iso-surfaces, as indicated by the color bar, for the case where \(R_A = 0\), \(R_c = 3\), \(\alpha = 0\), and \(d = 0.01\). In this scenario, the reaction solely induces 3D fingers by changing the product density at the reactive interface.}
    \label{fig:3D_flat}
\end{figure}
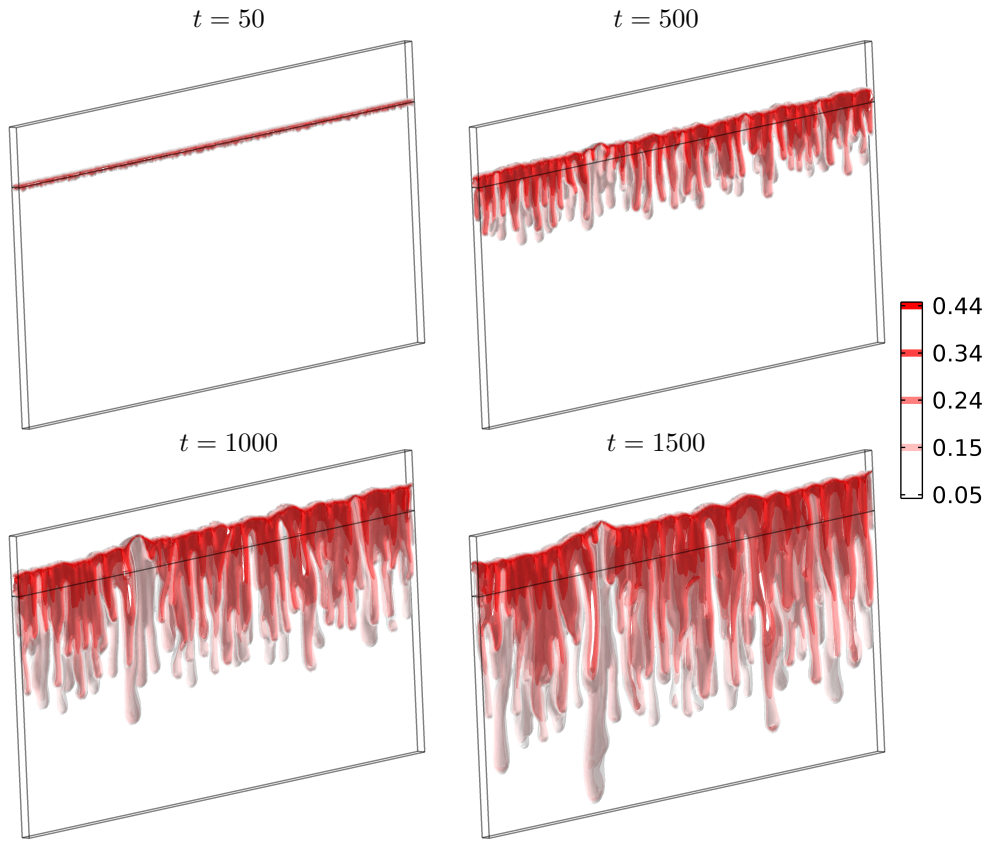
In this section, we present 3D simulation results demonstrating the continuous dependence of the solutions on the initial data to support our theoretical findings. The parameters are fixed as \(R_A = 0\), \(R_c = 3\), \(\alpha = 0\), and \(d = 0.01\). Similar effects, as observed in 2D, arising from variations in the density of reactant $A$ and permeability are also observed in 3D but are omitted here for brevity.

In Fig. \ref{fig:3D_flat}, we illustrate the spatio-temporal evolution of five iso-surfaces of the product concentration within the domain
$\Omega = \{ (x, y, z) \mid x \in [-1000, 1000], \, y \in [0, 1000], \, z \in [0, 30] \}.
$ This type of 3D domain is chosen to mimic a Hele-Shaw cell commonly used in experiments \cite{Paoli_2022}, consisting of two glass plates with a thin gap in the \( z \)-direction.
The initial conditions are specified as in \eqref{e71}, where the 3D vector is \(\boldsymbol{x} = (x, y, z)\), and the parameters are set as \(\alpha_0 = \beta_0 = 1\), \(y_0 = 800\), and \(\delta = 10^{-5}\). The initial gravitational field is again chosen to be a unit vector pointing downwards, with a zero initial velocity vector, and the boundary conditions remain the same as in \eqref{bc's}.
For the simulations, a tetrahedral mesh with a maximum element size of \(h = 7.5\) is used, resulting in a fixed degree of freedom, \(\text{DoF} = 31,276,714\). The velocity field \(\boldsymbol{u}_h\) is discretized using Lagrange \(\mathcal{P}_2\) elements, while \(\mathcal{P}_1\) elements are used for the pressure field \(p_h\) and the concentrations \(a_h\), \(b_h\), and \(c_h\), consistent with the 2D simulations. Time stepping is performed using the second order generalized-alpha method, with a time step size of \(\Delta t = 10^{-5}\).

Fig. \ref{fig:3D_flat} shows that, as time advances, the reaction induces 3D fingers at the interface, propagating downward. Since the reaction product is denser than the iso-dense reactants, the dynamics of the falling fingers resemble those observed in Case V (Fig. \ref{fig:flat}). However, in the 3D geometry, the reaction occurs along the surface of the fingers, creating distinct features. 
\begin{figure}
    \centering
    \includegraphics[width=0.33\linewidth]{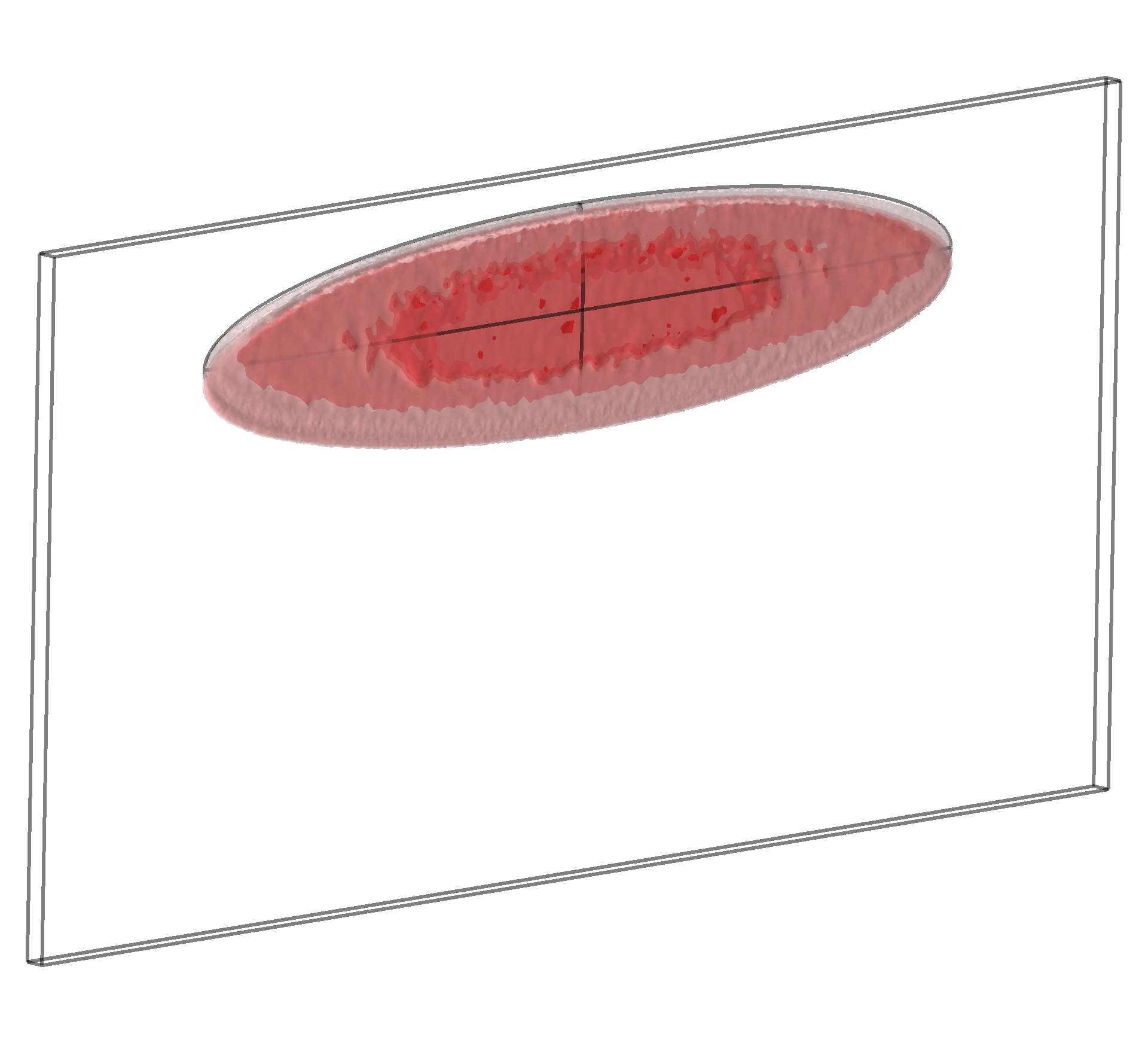}
    \includegraphics[width=0.33\linewidth]{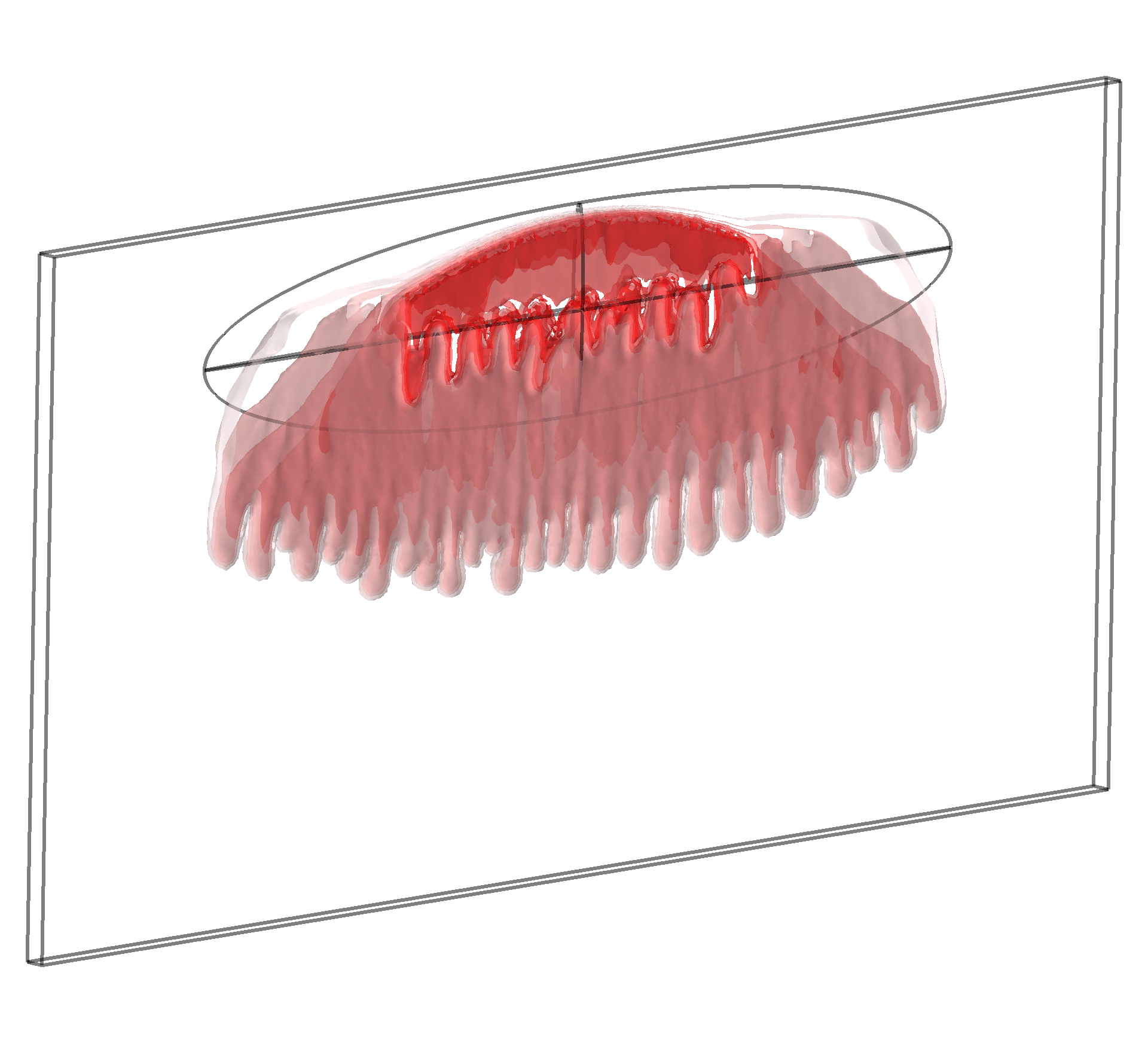}
    \includegraphics[width=0.33\linewidth]{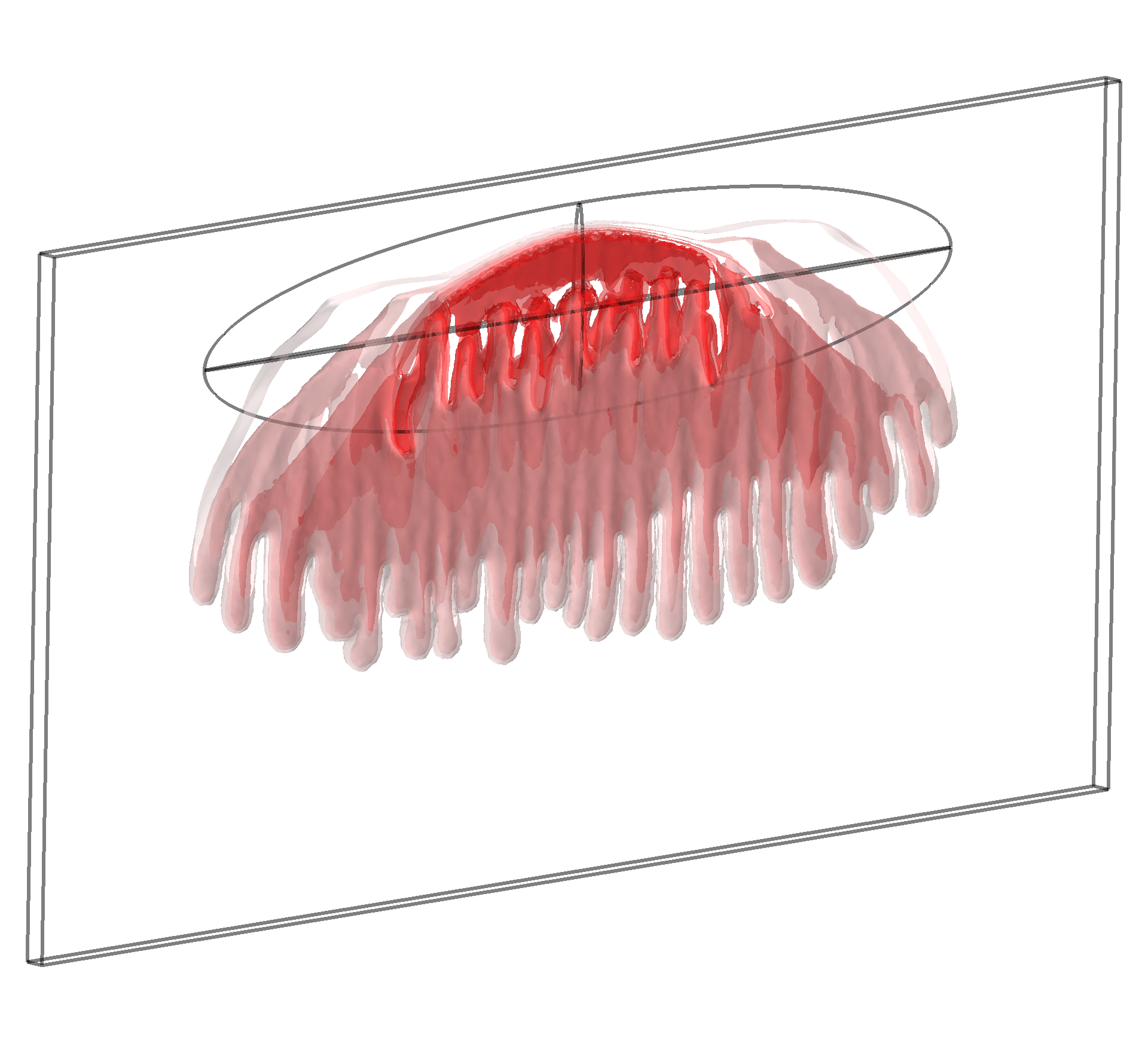}
    \includegraphics[width=0.33\linewidth]{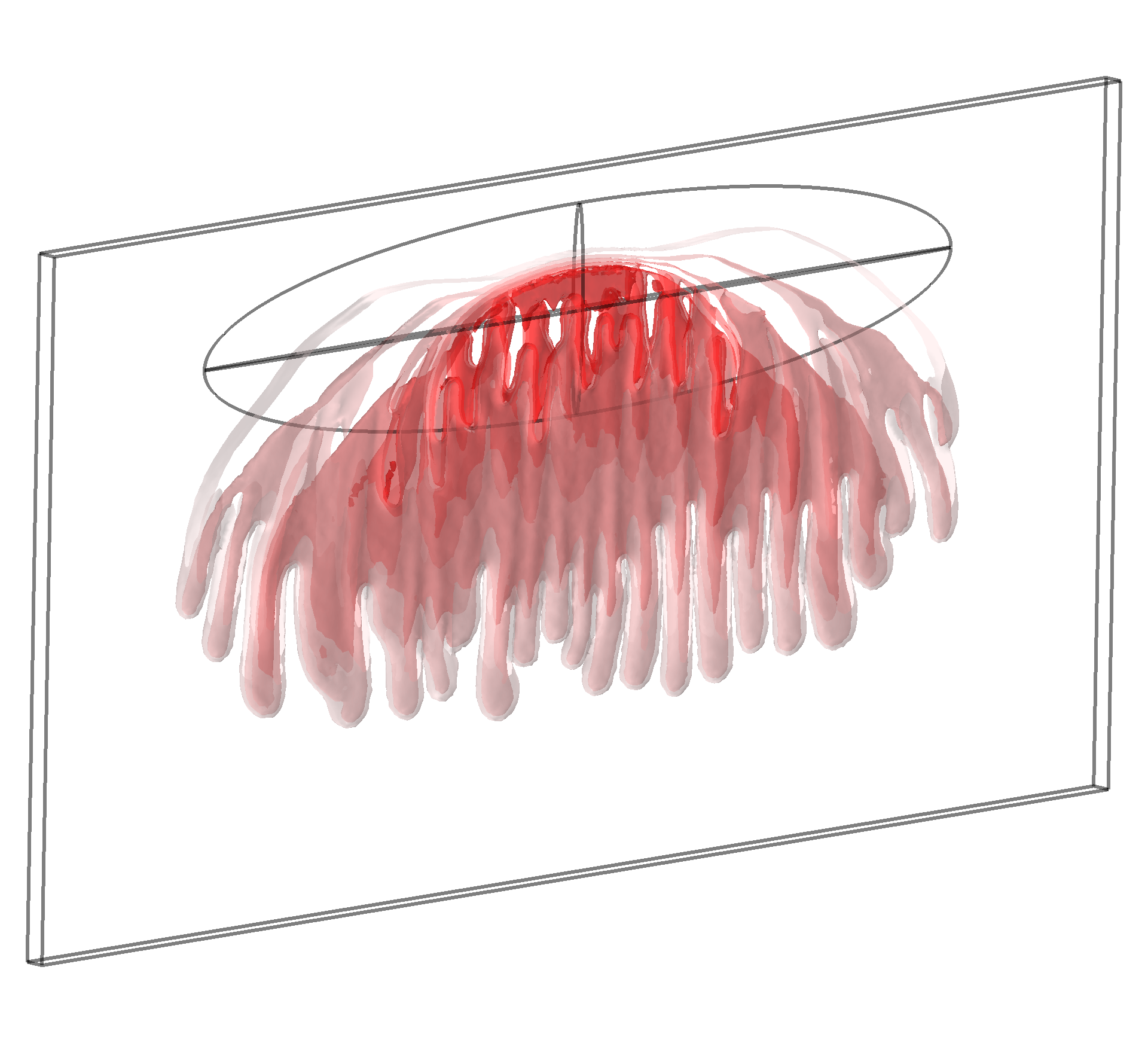}

    \begin{picture}(0,0)
    \put(-80,325){\makebox(0,0)[]{$t=50$}}
    \put(80,325){\makebox(0,0)[]{$t=500$}}
    \put(-80,165){\makebox(0,0)[]{$t=1000$}}
    \put(80,165){\makebox(0,0)[]{$t=1500$}}
    
    \put(170,140){\includegraphics[trim=170 25 10 40,clip, width=0.07\linewidth]{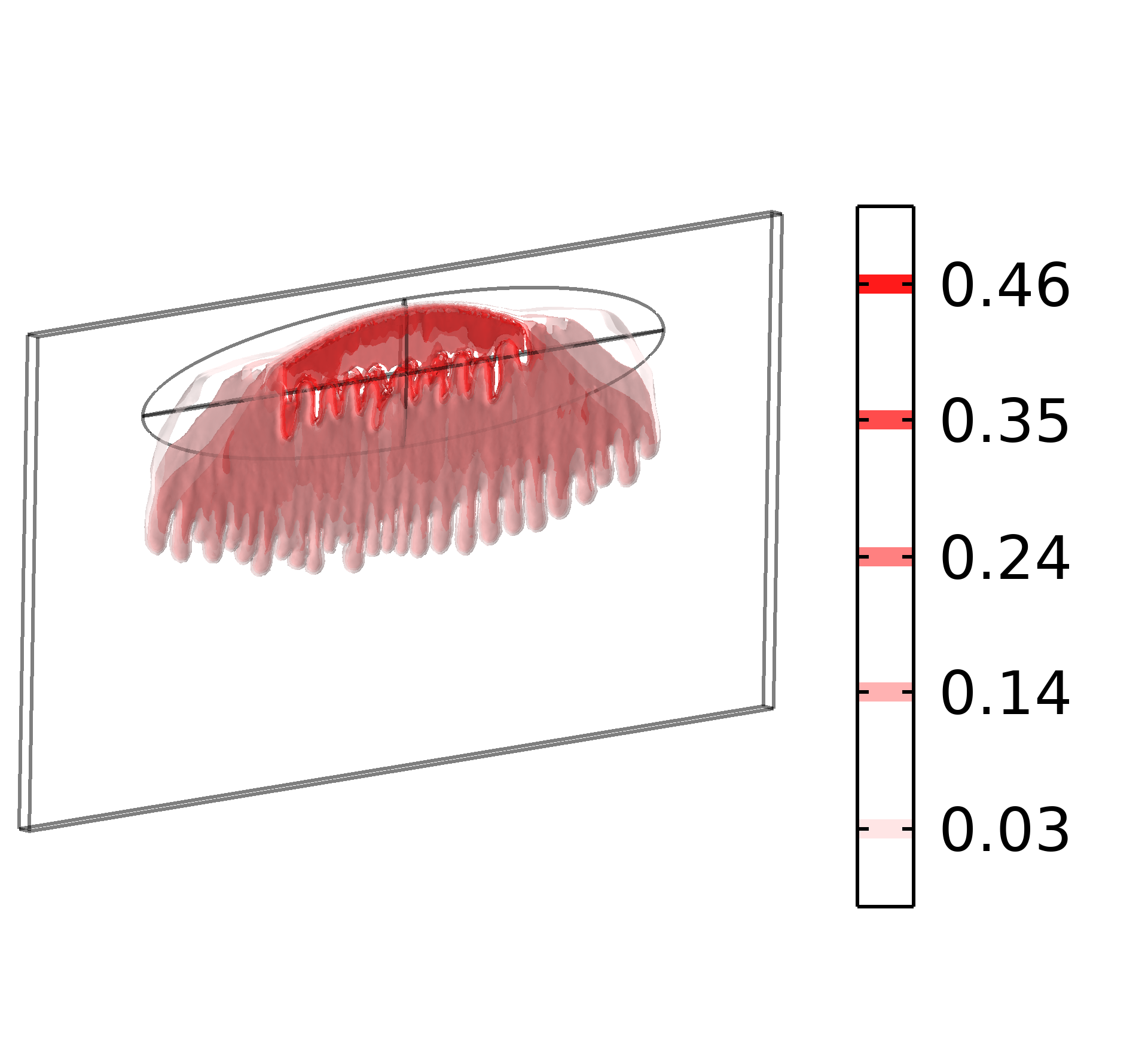}}
   \end{picture}
    \vspace{-0.7 cm} \caption{Spatio-temporal evolution of the reaction product concentration's five iso-surfaces, as indicated by the color bar, for the case where \(R_A = 0\), \(R_c = 3\), \(\alpha = 0\), and \(d = 0.01\). In this scenario, the reaction solely induces 3D fingers by changing the product density at the reactive interface.}
    \label{fig:3D_el}
\end{figure}
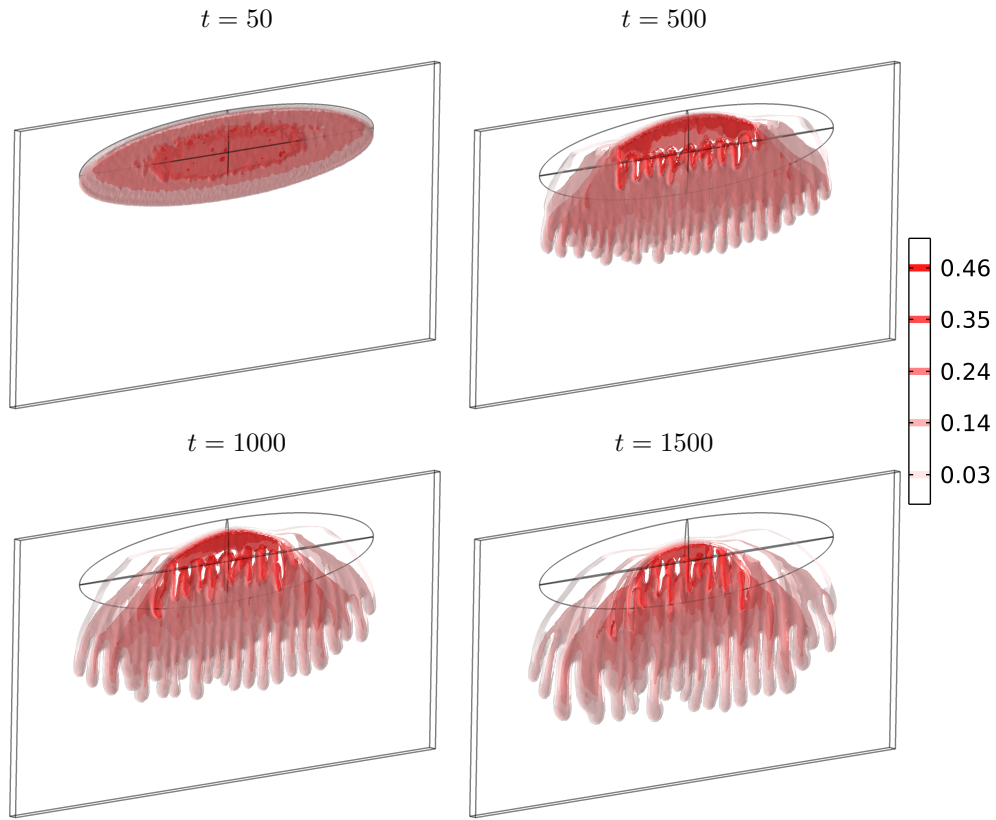
For the same parameters, we use a different initial configuration for the concentrations to demonstrate the continuous dependence of the solution in 3D geometry, while the rest of the setup remains the same as in the 3D flat interface case (Fig. \ref{fig:3D_flat}). Unlike the initial 2D elliptical bulb of reactant \(A\), in 3D, an initial ellipsoidal bulb of reactant \(A\) reacts within a 3D Hele-Shaw shell filled with reactant \(B\). Fig. \ref{fig:3D_el} shows the temporal evolution of the iso-surfaces of the product concentration, where the following initial conditions for the concentrations are employed:
\begin{align}\label{e82}
\begin{aligned}
    a_0(\boldsymbol{x}) &= 
    \begin{cases} 
        1, & \text{if } \frac{x^2}{x_0^2} + \frac{(y-600)^2}{y_0^2} + \frac{(z-15)^2}{z_0^2} \leq 1, \\ 
        0, & \text{otherwise}.
    \end{cases} \\
    b_0(\boldsymbol{x}) &= 
    \begin{cases} 
        0, & \text{if } \frac{x^2}{x_0^2} + \frac{(y-800)^2}{y_0^2} + \frac{(z-15)^2}{z_0^2} \leq 1, \\ 
        1, & \text{otherwise}.
    \end{cases} \\
    c_0(\boldsymbol{x}) &= 0.
\end{aligned}
\end{align}       
With \(x\)-major, \(y\)-major, and \(z\)-major axis lengths \(x_0 = 700\), \(y_0 = 150\), and \(z_0 = 10\), respectively, the reaction progresses over time, producing \(C\) across the surface of the ellipsoid. Fingering subsequently occurs at both interfaces, penetrating the reactant \(B\) present below in the domain (Fig. \ref{fig:3D_el}). 
The dynamics differ significantly from the 3D flat interface case due to the presence of two reactive interfaces, which eventually merge through the interaction of falling fingers, as discussed previously. This highlights the strong dependence of the solution on the initial conditions.
\section{Concluding remarks}
In this paper, we proved the well-posedness of a model that uses the time-dependent Darcy-Brinkman equation to describe the flow of reactive solutions \(A\), \(B\), and \(C\) undergoing the reaction \(A + B \to C\). We also showed that the solute concentrations remain bounded. To validate the continuous dependence of the solution on the initial data, we ran numerical simulations. These simulations showed that the fingering patterns are very different when starting with a flat reactive interface compared to an elliptical one, revealing the distinct physical behavior of the system. We examined the effects of heterogeneous permeability, finding results consistent with experiments that show precipitation effects \cite{Binda_2017}, which are modeled as an effect of mobility ratio in \cite{nagatsu2014hydrodynamic}. Additionally, we extended the model to three-dimensional space and simulated under different initial conditions. The 3D simulations showed different fingering patterns compared to 2D, with the reaction occurring on the surfaces of fingers rather than along a line interface. This work provides a framework for future studies, including extending the model to cases where viscosity and density vary simultaneously, as suggested in \cite{Jha_2023}. This would make the model more realistic and applicable to complex real-world scenarios.\\

\noindent \textbf{Acknowledgments:} M.M. acknowledges partial support from the FIST program, DST, Government of India (Ref: SR/FST/MS-I/2018/22(C)). S.K. acknowledges UGC, Government of India, for a research fellowship (Ref: 1145/CSIR-UGC NET June 2019).

\bibliography{references}
\end{document}